\documentclass[12pt,a4paper]{article}

\pdfoutput=1

\usepackage{amsfonts,amssymb,amsthm,amsmath,latexsym}
\usepackage{subeqnarray,eqsection,indent,eepicemu,url,cite,bm}
\usepackage{multirow}  
\usepackage{tensor}  
\usepackage{graphicx} 
\usepackage{tikz}
\usepackage[low-sup]{subdepth}  
\usepackage{float}  

\usepackage{color}
\def\blue{\textcolor{blue}}
\def\red{\textcolor{red}}
\def\green{\textcolor{green}}

\date{March 11, 2020 \\[1.5mm] revised February 8, 2022}

\begin{document}

\title{\vspace*{-2cm}Some multivariate master polynomials for \\
       permutations, set partitions, and perfect matchings, \\
       and their continued fractions}

\author{
     \\
     {\small Alan D.~Sokal}                  \\[-2mm]
     {\small\it Department of Mathematics}   \\[-2mm]
     {\small\it University College London}   \\[-2mm]
     {\small\it Gower Street}                \\[-2mm]
     {\small\it London WC1E 6BT}             \\[-2mm]
     {\small\it UNITED KINGDOM}              \\[-2mm]
     {\small\tt sokal@math.ucl.ac.uk}        \\[-2mm]
     {\protect\makebox[5in]{\quad}}  
     \\[-2mm]
     {\small\it Department of Physics}       \\[-2mm]
     {\small\it New York University}         \\[-2mm]
     {\small\it 726 Broadway}          \\[-2mm]
     {\small\it New York, NY 10003}      \\[-2mm]
     {\small\it USA}      \\[-2mm]
     {\small\tt sokal@nyu.edu}               \\
     \\
     {\small Jiang Zeng}                   \\[-2mm]
     {\small\it Univ Lyon, Universit\'e Claude Bernard Lyon 1}   \\[-2mm]
     {\small\it CNRS UMR 5208, Institut Camille Jordan}   \\[-2mm]
     {\small\it 43 blvd.~du 11 novembre 1918} \\[-2mm]
     {\small\it F-69622 Villeurbanne Cedex}  \\[-2mm]
     {\small\it FRANCE}                    \\[-2mm]
     {\small\tt zeng@math.univ-lyon1.fr}   \\[-2mm]
     \\
}

\maketitle
\thispagestyle{empty}   

\begin{abstract}
We find Stieltjes-type and Jacobi-type continued fractions
for some ``master polynomials'' that enumerate permutations,
set partitions or perfect matchings with a large (sometimes infinite) number
of simultaneous statistics.
Our results contain many previously obtained identities
as special cases, providing a common refinement of all of them.
\end{abstract}

\bigskip
\noindent
{\bf Key Words:}  Permutation, set partition, perfect matching,
generating polynomial, continued fraction, S-fraction, J-fraction,
Dyck path, Motzkin path.

\bigskip
\noindent
{\bf Mathematics Subject Classification (MSC 2010) codes:}
05A19 (Primary);
05A05, 05A10, 05A15, 05A18, 05A30, 30B70, 33C45, 33D45 (Secondary).

\clearpage

\newtheorem{theorem}{Theorem}[section]
\newtheorem{proposition}[theorem]{Proposition}
\newtheorem{lemma}[theorem]{Lemma}
\newtheorem{corollary}[theorem]{Corollary}
\newtheorem{definition}[theorem]{Definition}
\newtheorem{conjecture}[theorem]{Conjecture}
\newtheorem{question}[theorem]{Question}
\newtheorem{problem}[theorem]{Problem}
\newtheorem{example}[theorem]{Example}

\renewcommand{\theenumi}{\alph{enumi}}
\renewcommand{\labelenumi}{(\theenumi)}
\def\eop{\hbox{\kern1pt\vrule height6pt width4pt
depth1pt\kern1pt}\medskip}
\def\prf{\par\noindent{\bf Proof.\enspace}\rm}
\def\rmk{\par\medskip\noindent{\bf Remark\enspace}\rm}

\newcommand{\bigdash}{%
\smallskip\begin{center} \rule{5cm}{0.1mm} \end{center}\smallskip}

\newcommand{\be}{\begin{equation}}
\newcommand{\ee}{\end{equation}}
\newcommand{\<}{\langle}
\renewcommand{\>}{\rangle}
\newcommand{\widebar}{\overline}
\def\reff#1{(\protect\ref{#1})}
\def\spose#1{\hbox to 0pt{#1\hss}}
\def\ltapprox{\mathrel{\spose{\lower 3pt\hbox{$\mathchar"218$}}
    \raise 2.0pt\hbox{$\mathchar"13C$}}}
\def\gtapprox{\mathrel{\spose{\lower 3pt\hbox{$\mathchar"218$}}
    \raise 2.0pt\hbox{$\mathchar"13E$}}}
\def\textprime{${}^\prime$}
\def\proof{\par\medskip\noindent{\sc Proof.\ }}
\def\firstproof{\par\medskip\noindent{\sc First Proof.\ }}
\def\secondproof{\par\medskip\noindent{\sc Second Proof.\ }}
\def\qed{ $\square$ \bigskip}
\newcommand{\myendremark}{ $\blacksquare$ \bigskip}
\def\proofof#1{\bigskip\noindent{\sc Proof of #1.\ }}
\def\firstproofof#1{\bigskip\noindent{\sc First Proof of #1.\ }}
\def\secondproofof#1{\bigskip\noindent{\sc Second Proof of #1.\ }}
\def\thirdproofof#1{\bigskip\noindent{\sc Third Proof of #1.\ }}
\def\half{ {1 \over 2} }
\def\third{ {1 \over 3} }
\def\twothird{ {2 \over 3} }
\def\smfrac#1#2{{\textstyle{#1\over #2}}}
\def\smhalf{ {\smfrac{1}{2}} }
\newcommand{\real}{\mathop{\rm Re}\nolimits}
\renewcommand{\Re}{\mathop{\rm Re}\nolimits}
\newcommand{\imag}{\mathop{\rm Im}\nolimits}
\renewcommand{\Im}{\mathop{\rm Im}\nolimits}
\newcommand{\sgn}{\mathop{\rm sgn}\nolimits}
\newcommand{\tr}{\mathop{\rm tr}\nolimits}
\newcommand{\supp}{\mathop{\rm supp}\nolimits}
\newcommand{\disc}{\mathop{\rm disc}\nolimits}
\newcommand{\diag}{\mathop{\rm diag}\nolimits}
\newcommand{\tridiag}{\mathop{\rm tridiag}\nolimits}
\newcommand{\AZ}{\mathop{\rm AZ}\nolimits}
\newcommand{\perm}{\mathop{\rm perm}\nolimits}
\def\hboxscript#1{ {\hbox{\scriptsize\em #1}} }
\renewcommand{\emptyset}{\varnothing}
\newcommand{\eqdef}{\stackrel{\rm def}{=}}

\newcommand{\restrict}{\!\upharpoonright\!}

\newcommand{\doublehat}[1]{{\widehat{\widehat{#1\, }}}}
\newcommand{\doubletilde}[1]{{\widetilde{\widetilde{#1\, }}}}

\newcommand{\compinv}{{\langle -1 \rangle}}   

\newcommand{\scra}{{\mathcal{A}}}
\newcommand{\scrb}{{\mathcal{B}}}
\newcommand{\scrc}{{\mathcal{C}}}
\newcommand{\scrd}{{\mathcal{D}}}
\newcommand{\scre}{{\mathcal{E}}}
\newcommand{\scrf}{{\mathcal{F}}}
\newcommand{\scrg}{{\mathcal{G}}}
\newcommand{\scrh}{{\mathcal{H}}}
\newcommand{\scri}{{\mathcal{I}}}
\newcommand{\scrk}{{\mathcal{K}}}
\newcommand{\scrl}{{\mathcal{L}}}
\newcommand{\scrm}{{\mathcal{M}}}
\newcommand{\scrn}{{\mathcal{N}}}
\newcommand{\scro}{{\mathcal{O}}}
\newcommand{\scrp}{{\mathcal{P}}}
\newcommand{\scrq}{{\mathcal{Q}}}
\newcommand{\scrr}{{\mathcal{R}}}
\newcommand{\scrs}{{\mathcal{S}}}
\newcommand{\scrt}{{\mathcal{T}}}
\newcommand{\scrv}{{\mathcal{V}}}
\newcommand{\scrw}{{\mathcal{W}}}
\newcommand{\scrz}{{\mathcal{Z}}}

\newcommand{\ahat}{{\widehat{a}}}
\newcommand{\Zhat}{{\widehat{Z}}}
\renewcommand{\k}{{\mathbf{k}}}
\newcommand{\n}{{\mathbf{n}}}
\newcommand{\vv}{{\mathbf{v}}}
\newcommand{\bv}{{\mathbf{v}}}
\newcommand{\w}{{\mathbf{w}}}
\newcommand{\x}{{\mathbf{x}}}
\newcommand{\bz}{{\mathbf{z}}}
\newcommand{\bw}{{\mathbf{w}}}
\newcommand{\cc}{{\mathbf{c}}}
\newcommand{\zero}{{\mathbf{0}}}
\newcommand{\one}{{\mathbf{1}}}
\newcommand{\bmm}{ {\bf m} }

\newcommand{\C}{{\mathbb C}}
\newcommand{\D}{{\mathbb D}}
\newcommand{\Z}{{\mathbb Z}}
\newcommand{\N}{{\mathbb N}}
\newcommand{\Q}{{\mathbb Q}}
\newcommand{\PP}{{\mathbb P}}
\newcommand{\R}{{\mathbb R}}
\newcommand{\RR}{{\mathbb R}}
\newcommand{\E}{{\mathbb E}}

\newcommand{\Sym}{{\mathfrak{S}}}
\newcommand{\SymB}{{\mathfrak{B}}}

\newcommand{\myle}{\preceq}
\newcommand{\myge}{\succeq}
\newcommand{\mygt}{\succ}

\newcommand{\B}{{\sf B}}
\newcommand{\OB}{{\sf OB}}
\newcommand{\OS}{{\sf OS}}
\newcommand{\OO}{{\sf O}}
\newcommand{\SP}{{\sf SP}}
\newcommand{\OSP}{{\sf OSP}}
\newcommand{\Eu}{{\sf Eu}}
\newcommand{\ERR}{{\sf ERR}}
\newcommand{\sfB}{{\sf B}}
\newcommand{\sfD}{{\sf D}}
\newcommand{\sfE}{{\sf E}}
\newcommand{\sfG}{{\sf G}}
\newcommand{\sfJ}{{\sf J}}
\newcommand{\sfP}{{\sf P}}
\newcommand{\sfQ}{{\sf Q}}
\newcommand{\sfS}{{\sf S}}
\newcommand{\sfT}{{\sf T}}
\newcommand{\sfW}{{\sf W}}
\newcommand{\sfMV}{{\sf MV}}
\newcommand{\AMV}{{\sf AMV}}
\newcommand{\BM}{{\sf BM}}

\newcommand{\emIB}{{\hbox{\em IB}}}
\newcommand{\emIP}{{\hbox{\em IP}}}
\newcommand{\emOB}{{\hbox{\em OB}}}
\newcommand{\emSC}{{\hbox{\em SC}}}

\newcommand{\stat}{{\rm stat}}
\newcommand{\cyc}{{\rm cyc}}
\newcommand{\Asc}{{\rm Asc}}
\newcommand{\asc}{{\rm asc}}
\newcommand{\Des}{{\rm Des}}
\newcommand{\des}{{\rm des}}
\newcommand{\Exc}{{\rm Exc}}
\newcommand{\exc}{{\rm exc}}
\newcommand{\aexc}{{\rm aexc}}
\newcommand{\Wex}{{\rm Wex}}
\newcommand{\wex}{{\rm wex}}
\newcommand{\Fix}{{\rm Fix}}
\newcommand{\fix}{{\rm fix}}
\newcommand{\bfix}{{\mathbf{fix}}}
\newcommand{\lev}{{\rm lev}}
\newcommand{\lrmax}{{\rm lrmax}}
\newcommand{\rlmax}{{\rm rlmax}}
\newcommand{\Rec}{{\rm Rec}}
\newcommand{\rec}{{\rm rec}}
\newcommand{\Arec}{{\rm Arec}}
\newcommand{\arec}{{\rm arec}}
\newcommand{\ERec}{{\rm ERec}}
\newcommand{\erec}{{\rm erec}}
\newcommand{\erecop}{{\rm erecop}}
\newcommand{\erecin}{{\rm erecin}}
\newcommand{\nerecop}{{\rm nerecop}}
\newcommand{\nerecin}{{\rm nerecin}}
\newcommand{\brec}{{\rm brec}}
\newcommand{\brecop}{{\rm brecop}}
\newcommand{\brecin}{{\rm brecin}}
\newcommand{\nbrecop}{{\rm nbrecop}}
\newcommand{\nbrecin}{{\rm nbrecin}}
\newcommand{\EArec}{{\rm EArec}}
\newcommand{\earec}{{\rm earec}}
\newcommand{\recarec}{{\rm recarec}}
\newcommand{\nonrec}{{\rm nonrec}}
\newcommand{\nrar}{{\rm nrar}}
\newcommand{\ereccval}{{\rm ereccval}}
\newcommand{\ereccdrise}{{\rm ereccdrise}}
\newcommand{\eareccpeak}{{\rm eareccpeak}}
\newcommand{\eareccdfall}{{\rm eareccdfall}}
\newcommand{\rar}{{\rm rar}}
\newcommand{\nrcpeak}{{\rm nrcpeak}}
\newcommand{\nrcval}{{\rm nrcval}}
\newcommand{\nrcdrise}{{\rm nrcdrise}}
\newcommand{\nrcdfall}{{\rm nrcdfall}}
\newcommand{\nrfix}{{\rm nrfix}}
\newcommand{\Cpeak}{{\rm Cpeak}}
\newcommand{\cpeak}{{\rm cpeak}}
\newcommand{\Cval}{{\rm Cval}}
\newcommand{\cval}{{\rm cval}}
\newcommand{\Cdasc}{{\rm Cdasc}}
\newcommand{\cdasc}{{\rm cdasc}}
\newcommand{\Cddes}{{\rm Cddes}}
\newcommand{\cddes}{{\rm cddes}}
\newcommand{\cdrise}{{\rm cdrise}}
\newcommand{\cdfall}{{\rm cdfall}}
\newcommand{\cros}{{\rm cros}}
\newcommand{\cross}{{\rm cross}}
\newcommand{\nest}{{\rm nest}}
\newcommand{\ucross}{{\rm ucross}}
\newcommand{\ucrosscval}{{\rm ucrosscval}}
\newcommand{\ucrosscdrise}{{\rm ucrosscdrise}}
\newcommand{\lcross}{{\rm lcross}}
\newcommand{\lcrosscpeak}{{\rm lcrosscpeak}}
\newcommand{\lcrosscdfall}{{\rm lcrosscdfall}}
\newcommand{\unest}{{\rm unest}}
\newcommand{\unestcval}{{\rm unestcval}}
\newcommand{\unestcdrise}{{\rm unestcdrise}}
\newcommand{\lnest}{{\rm lnest}}
\newcommand{\lnestcpeak}{{\rm lnestcpeak}}
\newcommand{\lnestcdfall}{{\rm lnestcdfall}}
\newcommand{\ujoin}{{\rm ujoin}}
\newcommand{\ljoin}{{\rm ljoin}}
\newcommand{\upsnest}{{\rm upsnest}}
\newcommand{\lpsnest}{{\rm lpsnest}}
\newcommand{\psnest}{{\rm psnest}}
\newcommand{\ccc}{{\rm cc}}
\newcommand{\ecp}{{\rm ecp}}
\newcommand{\ecpar}{{\rm ecpar}}
\newcommand{\ecpnar}{{\rm ecpnar}}
\newcommand{\ocp}{{\rm ocp}}
\newcommand{\ocpar}{{\rm ocpar}}
\newcommand{\ocpnar}{{\rm ocpnar}}
\newcommand{\ecvr}{{\rm ecvr}}
\newcommand{\ecvnr}{{\rm ecvnr}}
\newcommand{\ocvr}{{\rm ocvr}}
\newcommand{\ocvnr}{{\rm ocvnr}}
\newcommand{\ecr}{{\rm ecr}}
\newcommand{\ocr}{{\rm ocr}}
\newcommand{\ene}{{\rm ene}}
\newcommand{\oone}{{\rm one}}
\newcommand{\Peak}{{\rm Peak}}
\newcommand{\peak}{{\rm peak}}
\newcommand{\Val}{{\rm Val}}
\newcommand{\val}{{\rm val}}
\newcommand{\Dasc}{{\rm Dasc}}
\newcommand{\dasc}{{\rm dasc}}
\newcommand{\Ddes}{{\rm Ddes}}
\newcommand{\ddes}{{\rm ddes}}
\newcommand{\interval}{{\rm int}}
\newcommand{\inv}{{\rm inv}}
\newcommand{\invtilde}{{\widetilde{\rm inv}}}
\newcommand{\maj}{{\rm maj}}
\newcommand{\majhat}{{\widehat{\rm maj}}}
\newcommand{\ls}{{\rm ls}}
\newcommand{\lb}{{\rm lb}}
\newcommand{\rs}{{\rm rs}}
\newcommand{\rb}{{\rm rb}}
\newcommand{\crr}{{\rm cr}}
\newcommand{\crrtilde}{{\widetilde{{\rm cr}}}}
\newcommand{\crosshat}{{\widehat{\rm cr}}}
\newcommand{\crop}{{\rm crop}}
\newcommand{\crin}{{\rm crin}}
\newcommand{\nee}{{\rm ne}}
\newcommand{\neetilde}{{\widetilde{{\rm ne}}}}
\newcommand{\neop}{{\rm neop}}
\newcommand{\nein}{{\rm nein}}
\newcommand{\psne}{{\rm psne}}
\newcommand{\psnetilde}{{\widetilde{{\rm psne}}}}
\newcommand{\crne}{{\rm crne}}
\newcommand{\qne}{{\rm qne}}
\newcommand{\qnetilde}{{\widetilde{{\rm qne}}}}
\newcommand{\itilde}{{\widetilde{i}}}
\newcommand{\ktilde}{{\widetilde{k}}}
\newcommand{\pitilde}{{\widetilde{\pi}}}
\newcommand{\ov}{{\rm ov}}
\newcommand{\ovin}{{\rm ovin}}
\newcommand{\ovinrev}{{\rm ovinrev}}
\newcommand{\ovtilde}{{\widetilde{{\rm ov}}}}
\newcommand{\ovbar}{{\overline{{\rm ov}}}}
\newcommand{\cov}{{\rm cov}}
\newcommand{\covin}{{\rm covin}}
\newcommand{\covinrev}{{\rm covinrev}}
\newcommand{\covtilde}{{\widetilde{{\rm cov}}}}
\newcommand{\covbar}{{\overline{{\rm cov}}}}
\newcommand{\qcov}{{\rm qcov}}
\newcommand{\qcovtilde}{{\widetilde{{\rm qcov}}}}
\newcommand{\pscov}{{\rm pscov}}
\newcommand{\rodd}{{\rm rodd}}
\newcommand{\reven}{{\rm reven}}
\newcommand{\lodd}{{\rm lodd}}
\newcommand{\leven}{{\rm leven}}
\newcommand{\sg}{{\rm sg}}
\newcommand{\bl}{{\rm bl}}
\newcommand{\tran}{{\rm tr}}
\newcommand{\area}{{\rm area}}
\newcommand{\ret}{{\rm ret}}
\newcommand{\peaks}{{\rm peaks}}
\newcommand{\hl}{{\rm hl}}
\newcommand{\sll}{{\rm sll}}
\newcommand{\negg}{{\rm neg}}
\newcommand{\imp}{{\rm imp}}
\newcommand{\den}{{\rm den}}
\newcommand{\denbis}{{\rm denbis}}

\newcommand{\sfa}{{{\sf a}}}
\newcommand{\sfb}{{{\sf b}}}
\newcommand{\sfc}{{{\sf c}}}
\newcommand{\sfd}{{{\sf d}}}
\newcommand{\sfe}{{{\sf e}}}
\newcommand{\sff}{{{\sf f}}}
\newcommand{\sfg}{{{\sf g}}}
\newcommand{\sfh}{{{\sf h}}}
\newcommand{\sfi}{{{\sf i}}}
\newcommand{\bsfa}{{\mbox{\textsf{\textbf{a}}}}}
\newcommand{\bsfb}{{\mbox{\textsf{\textbf{b}}}}}
\newcommand{\bsfc}{{\mbox{\textsf{\textbf{c}}}}}
\newcommand{\bsfd}{{\mbox{\textsf{\textbf{d}}}}}
\newcommand{\bsfe}{{\mbox{\textsf{\textbf{e}}}}}
\newcommand{\bsff}{{\mbox{\textsf{\textbf{f}}}}}
\newcommand{\bsfg}{{\mbox{\textsf{\textbf{g}}}}}
\newcommand{\bsfh}{{\mbox{\textsf{\textbf{h}}}}}
\newcommand{\bsfi}{{\mbox{\textsf{\textbf{i}}}}}

\newcommand{\ba}{{\bm{a}}}
\newcommand{\bahat}{{\widehat{\bm{a}}}}
\newcommand{\bb}{{\bm{b}}}
\newcommand{\bc}{{\bm{c}}}
\newcommand{\bff}{{\bm{f}}}
\newcommand{\bg}{{\bm{g}}}
\newcommand{\br}{{\bm{r}}}
\newcommand{\bu}{{\bm{u}}}
\newcommand{\bA}{{\bm{A}}}
\newcommand{\bB}{{\bm{B}}}
\newcommand{\bC}{{\bm{C}}}
\newcommand{\bE}{{\bm{E}}}
\newcommand{\bF}{{\bm{F}}}
\newcommand{\bI}{{\bm{I}}}
\newcommand{\bJ}{{\bm{J}}}
\newcommand{\bM}{{\bm{M}}}
\newcommand{\bN}{{\bm{N}}}
\newcommand{\bP}{{\bm{P}}}
\newcommand{\bQ}{{\bm{Q}}}
\newcommand{\bS}{{\bm{S}}}
\newcommand{\bT}{{\bm{T}}}
\newcommand{\bW}{{\bm{W}}}
\newcommand{\bIB}{{\bm{IB}}}
\newcommand{\bOB}{{\bm{OB}}}
\newcommand{\bOS}{{\bm{OS}}}
\newcommand{\bERR}{{\bm{ERR}}}
\newcommand{\bSP}{{\bm{SP}}}
\newcommand{\bMV}{{\bm{MV}}}
\newcommand{\bBM}{{\bm{BM}}}
\newcommand{\balpha}{{\bm{\alpha}}}
\newcommand{\bbeta}{{\bm{\beta}}}
\newcommand{\bgamma}{{\bm{\gamma}}}
\newcommand{\bdelta}{{\bm{\delta}}}
\newcommand{\bepsilon}{{\bm{\epsilon}}}
\newcommand{\bzeta}{{\bm{\zeta}}}
\newcommand{\bomega}{{\bm{\omega}}}
\newcommand{\bone}{{\bm{1}}}
\newcommand{\bzero}{{\bm{0}}}

\newcommand{\Cbar}{{\overline{C}}}
\newcommand{\Dbar}{{\overline{D}}}
\newcommand{\dbar}{{\overline{d}}}
\def\Ctilde{{\widetilde{C}}}
\def\Ftilde{{\widetilde{F}}}
\def\Gtilde{{\widetilde{G}}}
\def\Htilde{{\widetilde{H}}}
\def\Ptilde{{\widetilde{P}}}
\def\Chat{{\widehat{C}}}
\def\ctilde{{\widetilde{c}}}
\def\zbar{{\overline{Z}}}
\def\pitilde{{\widetilde{\pi}}}

\newcommand{\sech}{{\rm sech}}

%
%
\newcommand{\sn}{{\rm sn}}
\newcommand{\cn}{{\rm cn}}
\newcommand{\dn}{{\rm dn}}
\newcommand{\sm}{{\rm sm}}
\newcommand{\cm}{{\rm cm}}

%
%
\newcommand{\zfz}{ {{}_0 \! F_0} }
\newcommand{\zfo}{ {{}_0 \! F_1} }
\newcommand{\ofz}{ {{}_1 \! F_0} }
\newcommand{\oft}{ {{}_1 \! F_2} }

%
%
\newcommand{\FHyper}[2]{ {\tensor[_{#1 \!}]{F}{_{#2}}\!} }
\newcommand{\phiHyper}[2]{ {\tensor[_{#1}]{\phi}{_{#2}}} }
\newcommand{\FHYPER}[5]{ {\FHyper{#1}{#2} \!\biggl(
   \!\!\begin{array}{c} #3 \\[1mm] #4 \end{array}\! \bigg|\, #5 \! \biggr)} }
\newcommand{\phiHYPER}[6]{ {\phiHyper{#1}{#2} \!\left(
   \!\!\begin{array}{c} #3 \\ #4 \end{array}\! ;\, #5, \, #6 \! \right)\!} }
\newcommand{\ofo}{ {\FHyper{1}{1}} }
\newcommand{\tfo}{ {\FHyper{2}{1}} }
\newcommand{\FHYPERbottomzero}[3]{ {\FHyper{#1}{0} \!\biggl(
   \!\!\begin{array}{c} #2 \\[1mm] \hbox{---} \end{array}\! \bigg|\, #3 \! \biggr)} }
\newcommand{\FHYPERtopzero}[3]{ {\FHyper{0}{#1} \!\biggl(
   \!\!\begin{array}{c} \hbox{---} \\[1mm] #2 \end{array}\! \bigg|\, #3 \! \biggr)} }

%
%
\newcommand{\stirlingsubset}[2]{\genfrac{\{}{\}}{0pt}{}{#1}{#2}}
\newcommand{\stirlingcycleold}[2]{\genfrac{[}{]}{0pt}{}{#1}{#2}}
\newcommand{\stirlingcycle}[2]{\genfrac{[}{]}{0pt}{}{#1}{#2}}
\newcommand{\assocstirlingsubset}[3]{{\genfrac{\{}{\}}{0pt}{}{#1}{#2}}_{\! \ge #3}}
\newcommand{\genstirlingsubset}[4]{{\genfrac{\{}{\}}{0pt}{}{#1}{#2}}_{\! #3,#4}}
\newcommand{\euler}[2]{\genfrac{\langle}{\rangle}{0pt}{}{#1}{#2}}
\newcommand{\eulergen}[3]{{\genfrac{\langle}{\rangle}{0pt}{}{#1}{#2}}_{\! #3}}
\newcommand{\eulersecond}[2]{\left\langle\!\! \euler{#1}{#2} \!\!\right\rangle}
\newcommand{\eulersecondgen}[3]{{\left\langle\!\! \euler{#1}{#2} \!\!\right\rangle}_{\! #3}}
\newcommand{\binomvert}[2]{\genfrac{\vert}{\vert}{0pt}{}{#1}{#2}}
\newcommand{\binomsquare}[2]{\genfrac{[}{]}{0pt}{}{#1}{#2}}


\newenvironment{sarray}{
             \textfont0=\scriptfont0
             \scriptfont0=\scriptscriptfont0
             \textfont1=\scriptfont1
             \scriptfont1=\scriptscriptfont1
             \textfont2=\scriptfont2
             \scriptfont2=\scriptscriptfont2
             \textfont3=\scriptfont3
             \scriptfont3=\scriptscriptfont3
           \renewcommand{\arraystretch}{0.7}
           \begin{array}{l}}{\end{array}}

\newenvironment{scarray}{
             \textfont0=\scriptfont0
             \scriptfont0=\scriptscriptfont0
             \textfont1=\scriptfont1
             \scriptfont1=\scriptscriptfont1
             \textfont2=\scriptfont2
             \scriptfont2=\scriptscriptfont2
             \textfont3=\scriptfont3
             \scriptfont3=\scriptscriptfont3
           \renewcommand{\arraystretch}{0.7}
           \begin{array}{c}}{\end{array}}

\newcommand*{\Scale}[2][4]{\scalebox{#1}{$#2$}}

\clearpage

\tableofcontents

\clearpage

\section{Introduction}

If $(a_n)_{n \ge 0}$ is a sequence of combinatorial numbers or polynomials
with $a_0 = 1$, it is often fruitful to seek to express its
ordinary generating function as a continued fraction of either
Stieltjes (S) type,
\be
   \sum_{n=0}^\infty a_n t^n
   \;=\;
   \cfrac{1}{1 - \cfrac{\alpha_1 t}{1 - \cfrac{\alpha_2 t}{1 - \cdots}}}
   \label{def.Stype}
   \;\;,
\ee
or Jacobi (J) type,
\be
   \sum_{n=0}^\infty a_n t^n
   \;=\;
   \cfrac{1}{1 - \gamma_0 t - \cfrac{\beta_1 t^2}{1 - \gamma_1 t - \cfrac{\beta_2 t^2}{1 - \cdots}}}
   \label{def.Jtype}
   \;\;.
\ee
(Both sides of these expressions are to be interpreted as
formal power series in the indeterminate $t$.)
This line of investigation goes back at least to
Euler \cite{Euler_1760,Euler_1788},
but it gained impetus following Flajolet's \cite{Flajolet_80}
seminal discovery that any S-type (resp.\ J-type) continued fraction
can be interpreted combinatorially as a generating function
for Dyck (resp.\ Motzkin) paths with suitable weights for each rise and fall
(resp.\ each rise, fall and level step).
There are now literally dozens of sequences $(a_n)_{n \ge 0}$
of combinatorial numbers or polynomials for which
a continued-fraction expansion of the type \reff{def.Stype} or \reff{def.Jtype}
is explicitly known.

Our approach in this paper will be (in part) to run this program in reverse:
we start from a continued fraction in which the coefficients $\balpha$
(or $\bbeta$ and $\bgamma$) contain indeterminates in a nice pattern,
and we attempt to find a combinatorial interpretation for the
resulting polynomials $a_n$ --- namely,
as enumerating permutations, set partitions or perfect matchings
according to some natural multivariate statistics.
We call our $a_n$ ``master polynomials'' because our continued fractions
will contain the maximum number of independent indeterminates
consistent with the given pattern.
As a consequence, our results will contain many previously obtained identities
as special cases, providing a common refinement of all of them.

For future reference, let us recall the formula
\cite[p.~21]{Wall_48} \cite[p.~V-31]{Viennot_83}
for the contraction of an S-fraction to a J-fraction:
\reff{def.Stype} and \reff{def.Jtype} are equal if
\begin{subeqnarray}
   \gamma_0  & = &  \alpha_1
       \slabel{eq.contraction_even.coeffs.a}   \\
   \gamma_n  & = &  \alpha_{2n} + \alpha_{2n+1}  \qquad\hbox{for $n \ge 1$}
       \slabel{eq.contraction_even.coeffs.b}   \\
   \beta_n  & = &  \alpha_{2n-1} \alpha_{2n}
       \slabel{eq.contraction_even.coeffs.c}
 \label{eq.contraction_even.coeffs}
\end{subeqnarray}

The plan of this paper is as follows:
In Sections~\ref{sec.intro.permutations}--\ref{sec.intro.perfectmatchings}
we present our results for permutations, set partitions and perfect matchings,
respectively.
In Section~\ref{sec.prelim} we review briefly
the two main ingredients of our proofs:
namely, the combinatorial interpretation of continued fractions
in terms of Dyck and Motzkin paths \cite{Flajolet_80},
and the notion
(due to Flajolet \cite{Flajolet_80} and Viennot \cite{Viennot_83})
of labeled Dyck or Motzkin paths.
Finally, in Sections~\ref{sec.proofs.permutations} and
\ref{sec.proofs.setpartitions}
we supply the proofs for permutations and set partitions, respectively:
they employ bijections onto labeled Motzkin paths.
The proofs for perfect matchings will have already been presented
in Section~\ref{sec.intro.perfectmatchings},
as corollaries of the results for set partitions and permutations.


\section{Permutations: Statement of results}  \label{sec.intro.permutations}

\subsection{S-fraction}  \label{subsec.perms.S}

Euler \cite[section~21]{Euler_1760}\footnote{
   The paper \cite{Euler_1760},
   which is E247 in Enestr\"om's \cite{Enestrom_13} catalogue,
   was probably written circa 1746;
   it~was presented to the St.~Petersburg Academy in 1753,
   and published in 1760.
}
showed that the generating function of the factorials
can be represented as an S-type continued fraction
\be
   \sum_{n=0}^\infty n! \: t^n
   \;=\;
   \cfrac{1}{1 - \cfrac{1t}{1 - \cfrac{1t}{1 - \cfrac{2t}{1- \cfrac{2t}{1-\cdots}}}}}
 \label{eq.nfact.contfrac}
\ee
with coefficients $\alpha_{2k-1} = k$, $\alpha_{2k} = k$.
Inspired by \reff{eq.nfact.contfrac},
let us introduce the polynomials $P_n(x,y,u,v)$ defined
by the continued fraction
\be
   \sum_{n=0}^\infty P_n(x,y,u,v) \: t^n
   \;=\;
   \cfrac{1}{1 - \cfrac{xt}{1 - \cfrac{yt}{1 - \cfrac{(x+u)t}{1- \cfrac{(y+v)t}{1 - \cfrac{(x+2u)t}{1 - \cfrac{(y+2v)t}{1-\cdots}}}}}}}
 \label{eq.eulerian.fourvar.contfrac}
\ee
with coefficients
\begin{subeqnarray}
   \alpha_{2k-1}  & = &  x + (k-1) u \\
   \alpha_{2k}    & = &  y + (k-1) v
 \label{def.weights.eulerian.fourvar}
\end{subeqnarray}
Clearly $P_n(x,y,u,v)$ is a homogeneous polynomial of degree $n$;
it therefore has three ``truly independent'' variables.
Since $P_n(1,1,1,1) = n!$, which enumerates permutations of an $n$-element set,
it is plausible to expect that $P_n(x,y,u,v)$ enumerates permutations of $[n]$
according to some natural trivariate statistic.
Our first result gives two alternative versions of this trivariate statistic:

\begin{theorem}[S-fraction for permutations]
   \label{thm.perms.S}
The polynomials $P_n(x,y,u,v)$ defined by
\reff{eq.eulerian.fourvar.contfrac}/\reff{def.weights.eulerian.fourvar}
have the combinatorial interpretations
\be
 \hbox{(a)}\quad
   P_n(x,y,u,v)
   \;=\;
   \sum_{\sigma \in \Sym_n}
      x^{\arec(\sigma)} y^{\erec(\sigma)}
         u^{n - \exc(\sigma) - \arec(\sigma)} v^{\exc(\sigma) - \erec(\sigma)}
 \label{eq.eulerian.fourvar.arec}
\ee
and
\be
 \hbox{(b)}\quad
   P_n(x,y,u,v)
   \;=\;
   \sum_{\sigma \in \Sym_n}
      x^{\cyc(\sigma)} y^{\erec(\sigma)}
         u^{n - \exc(\sigma) - \cyc(\sigma)} v^{\exc(\sigma) - \erec(\sigma)}
 \label{eq.eulerian.fourvar.cyc}
\ee
\end{theorem}

\noindent
We will prove parts (a) and (b)
in Sections~\ref{subsec.permutations.J}
and \ref{subsec.permutations.J.v2}, respectively,
as special cases of more general results.

To explain the symbols used here
--- and others to be used subsequently ---
let~us define the relevant permutation statistics.
Given a permutation $\sigma \in \Sym_n$,
an index $i \in [n]$ (or a value $\sigma(i) \in [n]$) is called a
\begin{itemize}
   \item {\em record}\/ (rec) (or {\em left-to-right maximum}\/)
         if $\sigma(j) < \sigma(i)$ for all $j < i$
      [note in particular that the index 1 is always a record
       and that the value $n$ is always a record];
   \item {\em antirecord}\/ (arec) (or {\em right-to-left minimum}\/)
         if $\sigma(j) > \sigma(i)$ for all $j > i$
      [note in particular that the index $n$ is always an antirecord
       and that the value 1 is always an antirecord];
   \item {\em exclusive record}\/ (erec) if it is a record and not also
         an antirecord;
   \item {\em exclusive antirecord}\/ (earec) if it is an antirecord and not also
      a record;
   \item {\em record-antirecord}\/ (rar) (or {\em pivot}\/)
      if it is both a record and an antirecord;
   \item {\em neither-record-antirecord}\/ (nrar) if it is neither a record
      nor an antirecord.
\end{itemize}
Every index $i$ thus belongs to exactly one of the latter four types;
we refer to this classification as the {\em record classification}\/.
We denote the number of cycles, records, antirecords, \ldots\  in $\sigma$
by $\cyc(\sigma)$, $\rec(\sigma)$, $\arec(\sigma)$, \ldots, respectively.

Next we say that an index $i \in [n]$ is a
\begin{itemize}
   \item {\em cycle peak}\/ (cpeak) if $\sigma^{-1}(i) < i > \sigma(i)$;
   \item {\em cycle valley}\/ (cval) if $\sigma^{-1}(i) > i < \sigma(i)$;
   \item {\em cycle double rise}\/ (cdrise) if $\sigma^{-1}(i) < i < \sigma(i)$;
   \item {\em cycle double fall}\/ (cdfall) if $\sigma^{-1}(i) > i > \sigma(i)$;
   \item {\em fixed point}\/ (fix) if $\sigma^{-1}(i) = i = \sigma(i)$.
\end{itemize}
Clearly every index $i$ belongs to exactly one of these five types;
we refer to this classification as the {\em cycle classification}\/.
A rougher classification is that an index $i \in [n]$
(or a value $\sigma(i) \in [n]$) is an
\begin{itemize}
   \item {\em excedance}\/ (exc) if $\sigma(i) > i$
     [i.e.\ $i$ is either a cycle valley or a cycle double rise];
   \item {\em anti-excedance}\/ (aexc) if $\sigma(i) < i$
     [i.e.\ $i$ is either a cycle peak or a cycle double fall];
   \item {\em fixed point}\/ (fix) if $\sigma(i) = i$.
\end{itemize}
Clearly every index $i$ belongs to exactly one of these three types.
We also say that $i$ is a {\em weak excedance}\/ if $\sigma(i) \ge i$,
and a {\em weak anti-excedance}\/ if $\sigma(i) \le i$.

The record and cycle classifications of indices are related as follows:
\begin{quote}
\begin{itemize}
   \item[(a)]  Every record is a weak excedance,
      and every exclusive record is an excedance.
   \item[(b)]  Every antirecord is a weak anti-excedance,
      and every exclusive antirecord is an anti-excedance.
   \item[(c)]  Every record-antirecord is a fixed point.
\end{itemize}
\end{quote}
Furthermore,
\begin{quote}
\begin{itemize}
   \item[(d)]  The largest (resp.\ smallest) element of a cycle
         of length $\ge 2$ is always a cycle peak (resp.\ cycle valley).
\end{itemize}
\end{quote}
and hence in particular
\begin{quote}
\begin{itemize}
   \item[(d${}'$)]  Every cycle contains at least one non-excedance,
         and at least one non-anti-excedance.
\end{itemize}
\end{quote}
It follows that $\exc - \erec$,
$n - \exc - \arec$ and
$n - \exc - \cyc$ are all nonnegative,
so that the right-hand sides of \reff{eq.eulerian.fourvar.arec}
and \reff{eq.eulerian.fourvar.cyc} are indeed polynomials.

\subsection{Examples}

Many special cases of Theorem~\ref{thm.perms.S} were previously known.
For instance:

\begin{itemize}
   \item The Stirling cycle polynomials
\be
    P_n(x,1,1,1)
    \;=\;  S_n(x)
    \;=\;  \sum_{k=0}^n \stirlingcycle{n}{k} \: x^k
    \;=\;  x (x+1) \cdots (x+n-1)
    \;,
\ee
where $\stirlingcycle{n}{k}$ denotes the number of permutations of $[n]$
with $k$ cycles (or $k$ antirecords\footnote{
   Foata's fundamental transformation
   \cite[section~I.3]{Foata_70} \cite[pp.~17--18]{Stanley_86}
   \cite[section~3.3.1]{Bona_12}
   shows that cyc and rec (or equivalently arec)
   are equidistributed on $\Sym_n$.
});
or their homogenized version
\be
    P_n(x,y,y,y)
    \;=\;  S_n(x,y)
    \;=\;  \sum_{k=0}^n \stirlingcycle{n}{k} \: x^k \, y^{n-k}
    \;=\;  x (x+y) \cdots (x+(n-1)y)
    \;.
 \label{eq.stirlingcycle.homo}
\ee
The continued fraction \reff{eq.eulerian.fourvar.contfrac}
for this case was found already by Euler
\cite[section~26]{Euler_1760} \cite{Euler_1788}.\footnote{
   The paper \cite{Euler_1788},
   which is E616 in Enestr\"om's \cite{Enestrom_13} catalogue,
   was apparently presented to the St.~Petersburg Academy in 1776,
   and published posthumously in 1788.
}
   \item The Eulerian polynomials
\be
   P_n(1,y,1,y)
   \;=\;  A_n(y)
   \;=\;  \sum_{k=0}^n \euler{n}{k} \, y^k
\ee
where $\euler{n}{k}$ denotes the number of permutations of $[n]$
with $k$ excedances (or $k$ descents);
or their homogenized version
\be
   P_n(x,y,x,y)
   \;=\;  A_n(x,y)
   \;=\;  \sum_{k=0}^n \euler{n}{k} \, x^{n-k} \, y^k
   \;.
 \label{eq.eulerian.homo}
\ee
The continued fraction \reff{eq.eulerian.fourvar.contfrac}
for this case was found by Stieltjes
\cite[section~79]{Stieltjes_1894}.\footnote{
   Stieltjes does not specifically mention the Eulerian polynomials,
   but he does state that the continued fraction
   is the formal Laplace transform of
   $(1-y) / (e^{t(y-1)} - y)$,
   which is well known to be the exponential generating function
   of the Eulerian polynomials.
   Stieltjes also refrains from showing the proof:
   ``Pour abr\'eger, je supprime toujours les artifices qu'il faut employer
     pour obtenir la transformation de l'int\'egrale d\'efinie
     en fraction continue'' (!).
   But a proof is sketched, albeit also without much explanation,
   in the book of Wall \cite[pp.~207--208]{Wall_48}.
   The J-fraction corresponding to
   the contraction \reff{eq.contraction_even.coeffs} of this S-fraction
   was proven, by combinatorial methods, by Flajolet
   \cite[Theorem~3B(ii) with a slight typographical error]{Flajolet_80}.
   Also, Dumont \cite[Propositions~2 and 7]{Dumont_86}
   gave a direct combinatorial proof of the S-fraction,
   based on an interpretation of the Eulerian polynomials
   in terms of ``bipartite involutions of $[2n]$''
   and a bijection of these onto Dyck paths.
}
   \item A two-variable combination of the Stirling cycle and Eulerian
     polynomials
\cite{Carlitz_60,Dillon_68,Foata_70,Carlitz_77,Brenti_00,Ksavrelof_03,Savage_12,Ma_12}
\be
   P_n(x,y,1,y)
   \;=\;  F_n(x,y)
   \;=\; \sum_{\sigma \in \Sym_n} x^{\cyc(\sigma)} y^{\exc(\sigma)}
\ee
or its homogenized version
\be
   P_n(x,y,u,y)
   \;=\;  F_n(x,y,u)
   \;=\; \sum_{\sigma \in \Sym_n} x^{\cyc(\sigma)} y^{\exc(\sigma)}
                 u^{n - \exc(\sigma) - \cyc(\sigma)}
   \;.
\ee
The continued fraction \reff{eq.eulerian.fourvar.contfrac}
for this case was again found by Stieltjes
\cite[section~81]{Stieltjes_1894}.\footnote{
   Once again, Stieltjes does not specifically mention the polynomials,
   but he states that the continued fraction
   is the formal Laplace transform of
   $[ (1-y) / (e^{t(y-1)} - y) ]^{x}$,
   which is the exponential generating function of the polynomials
   $P_n(x,y,1,y)$.
}
These polynomials have a nice explicit formula \cite{Dillon_68,Salas_private}
that can be derived by simple algebra from
\cite[eqns.~(1.3)/(3.3)/(3.5)]{Dillon_68} \cite[Corollary~2.3]{Ma_12}:
\be
   P_n(x,y,u,y)
   \;=\;
   \sum_{k=0}^n \stirlingsubset{n}{k} \, (y-u)^{n-k} \,
      \prod\limits_{j=0}^{k-1} (x+ju)
 \label{eq.twovar.Ma.Salas}
\ee
where $\stirlingsubset{n}{k}$ denotes the number of partitions
of an $n$-element set into $k$ nonempty blocks.
When $u=y$ this reduces to \reff{eq.stirlingcycle.homo};
when $u=x$ it becomes [using \reff{eq.eulerian.homo}]
the well-known \cite[eqns.~(6.39)/(6.40)]{Graham_94} identity
\be
   \sum_{k=0}^n \euler{n}{k} \, x^{n-k} \, y^k
   \;=\;
   \sum_{k=0}^n k! \stirlingsubset{n}{k} (y-x)^{n-k}  x^k
 \label{eq.orderedbell.identity}
\ee
that relates the Eulerian polynomials to the ordered Bell polynomials.\footnote{
   The ordered Bell polynomials appear already
   (albeit without the combinatorial interpretation)
   in Euler's book {\em Foundations of Differential Calculus,
    with Applications to Finite Analysis and Series}\/,
   first published in 1755 \cite[paragraph~172]{Euler_1755}.
   This book is E212 in Enestr\"om's \cite{Enestrom_13} catalogue.
   Furthermore, the identity \reff{eq.orderedbell.identity}
   appears already there \cite[paragraphs~172 and 173]{Euler_1755};
   it was rediscovered a century-and-a-half later by Frobenius \cite{Frobenius_10}.
   See also \cite[pp.~150--151]{Flajolet_82} for a simple bijective proof.
}
   \item The record-antirecord permutation polynomials
\be
   P_n(a,b,1,1)
   \;=\;
   \sum_{\sigma \in \Sym_n} a^{\arec(\sigma)} \, b^{\erec(\sigma)}
 \label{def.recarec}
\ee
or their homogenized version
\begin{subeqnarray}
   P_n(a,b,c,c)
   & = &
   \sum_{\sigma \in \Sym_n}
      a^{\arec(\sigma)} \, b^{\erec(\sigma)} \,
                          c^{n - \arec(\sigma) - \erec(\sigma)}
      \\[2mm]
   & = &
   \sum_{\sigma \in \Sym_n}
      a^{\arec(\sigma)} \, b^{\erec(\sigma)} \, c^{\nrar(\sigma)}
   \;.
 \label{def.recarec.homo}
\end{subeqnarray}
Dumont and Kreweras \cite{Dumont_88} proved the continued fraction
\begin{subeqnarray}
   \sum_{n=0}^\infty P_n(a,b,1,1) \, t^n
   & = &
   \cfrac{1}{1 - \cfrac{at}{1 - \cfrac{bt}{1 - \cfrac{(a+1)t}{1- \cfrac{(b+1)t}{1 - \cdots}}}}}
   \qquad
         \slabel{eq.dumont.a}  \\[3mm]
   & = &
   {\FHYPERbottomzero{2}{a,b}{t}
    \:\bigg/\:
    \FHYPERbottomzero{2}{a,b-1}{t}
   }
   \;\,,
   \qquad
         \slabel{eq.dumont.b}
         \label{eq.dumont}
\end{subeqnarray}
where the second equality is a classic
\cite[section~92]{Wall_48} \cite[Theorem~6.5]{Jones_80}
corollary of Gauss' continued fraction
for ratios of contiguous $\tfo$.\footnote{
   Dumont and Kreweras \cite{Dumont_88} stated their result
   in terms of records and exclusive antirecords,
   which is of course equivalent to \reff{def.recarec} via the bijection
   $\sigma \mapsto R \circ \sigma \circ R$ with $R(i) = n+1-i$
   (i.e.~reversal combined with complementation).
}
These polynomials are also essentially identical
to the Martin--Kearney \cite{Martin_10} polynomials:
see \cite{Elvey-Sokal_selfcon} for details.
See also Section~\ref{subsec.intro.permutations.pq.S} below
for a $q$-generalization of \reff{def.recarec}/\reff{eq.dumont}.
   \item  As a special case of \reff{def.recarec.homo},
the polynomials \cite[A145879/A202992]{OEIS}
\be
   P_n(x,x,u,u)
   \;=\;
   \sum_{k=0}^n T(n,k) \, x^{n-k} u^k
\ee
where $T(n,k)$ is the number of permutations $\sigma \in \Sym_n$
having exactly $k$ indices that are the middle point of a pattern 321
(clearly $0 \le k \le n-2$ when $n \ge 2$).
In particular, $T(n,0)$ is the number of 321-avoiding permutations,
which equals the Catalan number $C_n$;
so these polynomials interpolate between $C_n$ and $n!$.\footnote{
   We thank Andrew Elvey Price for drawing our attention
   to these polynomials.
}
   \item  As another special case of \reff{def.recarec.homo},
the homogenized Narayana polynomials \cite[A001263/A090181]{OEIS}
\begin{subeqnarray}
   P_n(x,y,0,0)
   & = &
   \sum_{\sigma \in \Sym_n(321)}  x^{\arec(\sigma)} y^{\erec(\sigma)}
         \slabel{eq.narayana.a} \\[2mm]
   & = &
   \sum_{\sigma \in \Sym_n(321)}  x^{\arec(\sigma)} y^{\exc(\sigma)}
         \slabel{eq.narayana.b} \\[2mm]
   & = &
   \sum_{k=0}^n {1 \over n} \binom{n}{k} \binom{n}{k-1} \, x^k y^{n-k}
      \;,
 \label{eq.narayana}
\end{subeqnarray}
which count 321-avoiding permutations according to
the number of antirecords or exclusive records or excedances
(among many other combinatorial interpretations
 \cite{Sulanke_98,Sulanke_99}).\footnote{
   For a 321-avoiding permutation ---
   that is, one in which there are no neither-record-antirecords ---
   the containments given in (a)--(c) above are equalities:
   that is, an index $i$ is a record
   (resp.\ exclusive record, antirecord, exclusive antirecord)
   if and only if it is a weak excedance
   (resp.\ excedance, weak anti-excedance, anti-excedance).
   See Section~\ref{subsec.permutations.321} below
   for further enumerative results on 321-avoiding permutations.
}
These interpretations of the Narayana numbers
were found by Vella \cite[Proposition~2.12]{Vella_03}
and Elizalde \cite[Proposition~2.7(1)]{Elizalde_04}.
%
%
\end{itemize}


\subsection{First J-fraction}  \label{subsec.perms.J.1}

In Section~\ref{subsec.perms.S} we classified indices in a permutation
according to their record status
(exclusive record, exclusive antirecord, record-antirecord
or neither-record-antirecord) 
and also according to their cycle status
(cycle peak, cycle valley, cycle double rise, cycle double fall or fixed point).
Applying now both classifications simultaneously,
we obtain 10 disjoint categories:
\begin{itemize}
   \item ereccval:  exclusive records that are also cycle valleys;
   \item ereccdrise:  exclusive records that are also cycle double rises;
   \item eareccpeak:  exclusive antirecords that are also cycle peaks;
   \item eareccdfall:  exclusive antirecords that are also cycle double falls;
   \item rar:  record-antirecords (these are always fixed points);
   \item nrcpeak:  neither-record-antirecords that are also cycle peaks;
   \item nrcval:  neither-record-antirecords that are also cycle valleys;
   \item nrcdrise:  neither-record-antirecords that are also cycle double rises;
   \item nrcdfall:  neither-record-antirecords that are also cycle double falls;
   \item nrfix:  neither-record-antirecords that are also fixed points.
\end{itemize}
Clearly every index $i$ belongs to exactly one of these 10 types;
we call this the {\em record-and-cycle classification}\/.
The master polynomial encoding all these statistics is
\begin{eqnarray}
   & &
   Q_n(x_1,x_2,y_1,y_2,z,u_1,u_2,v_1,v_2,w)
   \;=\;
       \nonumber \\[4mm]
   & & \qquad\qquad
   \sum_{\sigma \in \Sym_n}
   x_1^{\eareccpeak(\sigma)} x_2^{\eareccdfall(\sigma)} 
   y_1^{\ereccval(\sigma)} y_2^{\ereccdrise(\sigma)} 
   z^{\rar(\sigma)} \:\times
       \qquad\qquad
       \nonumber \\[-1mm]
   & & \qquad\qquad\qquad\:
   u_1^{\nrcpeak(\sigma)} u_2^{\nrcdfall(\sigma)}
   v_1^{\nrcval(\sigma)} v_2^{\nrcdrise(\sigma)}
   w^{\nrfix(\sigma)}
   \;.
 \label{def.Qn}
\end{eqnarray}
(Thus, the variables $x_1,y_1,u_1,v_1$ are associated to cycle peaks and valleys,
 $x_2,y_2,u_2,v_2$ to cycle double rises and falls,
 and $z,w$ to fixed points.)
It turns out that these 10-variable homogeneous polynomials
have a beautiful J-fraction.

But we can go farther, by further refining the classification of fixed points.
If $i$ is a fixed point of $\sigma$, we define its {\em level}\/ by
\be
   \lev(i,\sigma)
   \;\eqdef\;
   \# \{j < i \colon\:  \sigma(j) > i \}
   \;=\;
   \# \{j > i \colon\:  \sigma(j) < i \}
   \;.
 \label{def.level}
\ee
[The two expressions are equal because $\sigma$ is a bijection
 from $[1,i) \cup (i,n]$ to itself.]
Note that for $\sigma \in \Sym_n$,
we have $0 \le \lev(i,\sigma) \le \min(i-1,n-i) \le \lfloor (n-1)/2 \rfloor$.
Clearly, a fixed point $i$ is a record-antirecord if and only if its level is 0,
and a neither-record-antirecord if and only if its level is $\ge 1$.
Let us now count the number of fixed points of each level:
for $\sigma \in \Sym_n$ and $\ell \ge 0$ we define
\be
   \fix(\sigma,\ell)
   \;\eqdef\;
   \# \{i \in [n] \colon\:  \sigma(i) = i \hbox{ and } \lev(i,\sigma) = \ell \}
   \;.
\ee
We then introduce indeterminates $\bw = (w_\ell)_{\ell \ge 0}$ and write
\be
   \bw^{\bfix(\sigma)}
   \;\eqdef\;
   \prod_{\ell=0}^\infty w_\ell^{\fix(\sigma,\ell)}
   \;=\;
   \prod_{i \in {\rm fix}} w_{\lev(i,\sigma)}
   \;.
\ee
The master polynomial encoding all these (now infinitely many) statistics is
\begin{eqnarray}
   & &
   Q_n(x_1,x_2,y_1,y_2,u_1,u_2,v_1,v_2,\bw)
   \;=\;
       \nonumber \\[4mm]
   & & \qquad\qquad
   \sum_{\sigma \in \Sym_n}
   x_1^{\eareccpeak(\sigma)} x_2^{\eareccdfall(\sigma)} 
   y_1^{\ereccval(\sigma)} y_2^{\ereccdrise(\sigma)} 
   \:\times
       \qquad\qquad
       \nonumber \\[-1mm]
   & & \qquad\qquad\qquad\:
   u_1^{\nrcpeak(\sigma)} u_2^{\nrcdfall(\sigma)}
   v_1^{\nrcval(\sigma)} v_2^{\nrcdrise(\sigma)}
   \bw^{\bfix(\sigma)}
   \;,
 \label{def.Qnbis}
\end{eqnarray}
which reduces to \reff{def.Qn} if we set $w_0 = z$
and $w_\ell = w$ for $\ell \ge 1$.
We then have:

\begin{theorem}[First J-fraction for permutations]
   \label{thm.perm.Jtype}
The ordinary generating function of the polynomials \reff{def.Qnbis} has the
J-type continued fraction
\begin{eqnarray}
   & & \hspace*{-7mm}
   \sum_{n=0}^\infty Q_n(x_1,x_2,y_1,y_2,u_1,u_2,v_1,v_2,\bw) \: t^n
   \;=\;
       \nonumber \\
   & & 
   \cfrac{1}{1 - w_0 t - \cfrac{x_1 y_1 t^2}{1 -  (x_2\!+\!y_2\!+\!w_1) t - \cfrac{(x_1\!+\!u_1)(y_1\!+\!v_1) t^2}{1 - (x_2\!+\!u_2\!+\!y_2\!+\!v_2\!+\!w_2)t - \cfrac{(x_1\!+\!2u_1)(y_1\!+\!2v_1) t^2}{1 - \cdots}}}}
       \nonumber \\[1mm]
   \label{eq.thm.perm.Jtype}
\end{eqnarray}
with coefficients
\begin{subeqnarray}
   \gamma_0  & = &   w_0   \\[1mm]
   \gamma_n  & = &   [x_2 + (n-1)u_2]  \,+\,  [y_2 + (n-1)v_2]  \,+\,  w_n
        \qquad\hbox{for $n \ge 1$}   \\[1mm]
   \beta_n   & = &   [x_1 + (n-1)u_1] \: [y_1 + (n-1)v_1]
 \label{def.weights.perm.Jtype}
\end{subeqnarray}
\end{theorem}

\noindent
We will prove this theorem in Section~\ref{subsec.permutations.J},
as a special case of a more general result.

%

\bigskip

{\bf Remark.}
The continued fraction \reff{eq.thm.perm.Jtype}
shows that $Q_n$ depends on its arguments
only via the combinations (\ref{def.weights.perm.Jtype}a,b,c).
In particular, it depends on $x_2,y_2$ only via the combination $x_2 + y_2$,
and on $u_2,v_2$ only via the combination $u_2 + v_2$;
consequently (but more weakly),
it is symmetric under $x_2 \leftrightarrow y_2$
and indepdendently under $u_2 \leftrightarrow v_2$.
Furthermore, it is symmetric under $(x_1,u_1) \leftrightarrow (y_1,v_1)$.
It would be interesting to try to understand combinatorially
(directly at the level of permutations) why these properties hold.

One very special case of these properties is easy to understand combinatorially:
$Q_n$ is invariant under the simultaneous interchange
$(x_1,u_1,x_2,u_2) \leftrightarrow (y_1,v_1,y_2,v_2)$.
This is because the bijection
$\sigma \mapsto R \circ \sigma \circ R$ with $R(i) = n+1-i$
(i.e.\ reversal combined with complementation)
interchanges cycle peaks with cycle valleys,
cycle double rises with cycle double falls,
and records with antirecords
(while preserving the number of fixed points at each level).
But the more specfic properties encoded in \reff{def.weights.perm.Jtype}
remain mysterious.
%
\myendremark

\medskip

Comparing \reff{def.Qnbis} with \reff{eq.eulerian.fourvar.arec},
we see that $Q_n(x_1,x_2,y_1,y_2,u_1,u_2,v_1,v_2,\bw)$
reduces to $P_n(x,y,u,v)$ if we set
\be
   x_1 = x_2 = w_0 = x, \quad
   y_1 = y_2 = y, \quad
   u_1 = u_2 = w_1 = w_2 = \ldots = u, \quad
   v_1 = v_2 = v  \;.
 \label{def.specialization.xyuv}
\ee
With this specialization
the J-fraction coefficients \reff{def.weights.perm.Jtype} reduce to
\begin{subeqnarray}
   \gamma_0  & = &   x   \\[1mm]
   \gamma_n  & = &   (x + nu)  \,+\,  [y + (n-1)v]
        \qquad\hbox{for $n \ge 1$}   \\[1mm]
   \beta_n   & = &   [x + (n-1)u] \: [y + (n-1)v]
 \label{def.weights.perm.Jtype.specialization.xyuv}
\end{subeqnarray}
which are precisely those that
arise as the contraction \reff{eq.contraction_even.coeffs}
of an S-type continued fraction with coefficients
\reff{def.weights.eulerian.fourvar}.
So Theorem~\ref{thm.perms.S}(a) is an immediate consequence of
a very special case of Theorem~\ref{thm.perm.Jtype}.

%

\subsection{Second J-fraction}  \label{subsec.perms.J.2}

The generalization of Theorem~\ref{thm.perms.S}(b) is less satisfying,
because cyc does not seem to mesh well with the record classification:
even the three-variable polynomials
\be
   \widehat{P}_n(x,y,\lambda)
   \;=\;
   \sum_{\sigma \in \Sym_n}
   x^{\arec(\sigma)} y^{\erec(\sigma)} \lambda^{\cyc(\sigma)}
\ee
do not have a J-fraction with polynomial coefficients
(starting at $\gamma_2$ we get rational functions).\footnote{
   The first coefficients are
   $$
      \gamma_0 \,=\, \lambda x ,\quad
      \beta_1  \,=\, \lambda xy ,\quad
      \gamma_1 \,=\, \lambda + x + y ,\quad
      \beta_2  \,=\, \lambda + x + y + \lambda xy
      \;,
   $$
   followed by
   $$ \gamma_2  \;=\;
      {(x+y)(3+xy) \:+\: (2 + x + x^2 + y + 4xy + y^2) \lambda
                   \:+\: (1 + xy) \lambda^2
       \over
       \lambda + x + y + \lambda xy
      }
     \;.
   $$
   It can then be shown that
   \begin{itemize}
      \item[(a)] $\gamma_2$ is not a polynomial in $x$
          (when $y$ and $\lambda$ are given fixed real values)
          unless $\lambda \in \{-1,0,+1\}$ or $y \in \{-1,+1\}$
          or $\lambda y = -1$.
      \item[(b)] $\gamma_2$ is not a polynomial in $y$
          unless $\lambda \in \{-1,0,+1\}$ or $x \in \{-1,+1\}$
          or $\lambda x = -1$.
      \item[(c)] $\gamma_2$ is not a polynomial in $\lambda$
          unless $x \in \{-1,+1\}$ or $y \in \{-1,+1\}$
          or $x+y = 0$ or $xy = -1$.
   \end{itemize}
}
However, cyc does mesh well with the complete cycle classification
(cpeak, cdfall, cval, cdrise, fix),
as was shown by one of us \cite[Th\'eor\`eme~3]{Zeng_93}
more than two decades ago.\footnote{
   The paper \cite{Zeng_93} used traditional combinatorial methods
   to establish an exponential generating function
   for permutations with these weights,
   and then used algebraic methods to transform this
   exponential generating function into a continued fraction
   for the ordinary generating function.
   Here, by contrast, we employ bijections onto labeled Motzkin paths
   to establish the continued fraction directly.
   We think it is instructive to compare these two quite different
   methods of proof.

   A special case of \cite[Th\'eor\`eme~3]{Zeng_93}
   with one fewer ``truly independent'' variable
   --- namely, including cval, cdrise, fix and cyc,
   but not distinguishing cpeak from cdfall ---
   was obtained recently by Elizalde \cite[eqn.~(3)]{Elizalde_18},
   using a bijection onto labeled Motzkin paths
   that is essentially the same as ours.
   We will discuss the connection with Elizalde's work
   in Sections~\ref{subsec.permutations.inv} and \ref{subsec.permutations.J.v2}.
}
In fact, it seems that cyc {\em almost}\/ meshes with the
complete record-and-cycle classification;
and we can also include the refined classification of fixed points.
Let us define the polynomial
\begin{eqnarray}
   & &
   \widehat{Q}_n(x_1,x_2,y_1,y_2,u_1,u_2,v_1,v_2,\bw,\lambda)
   \;=\;
       \nonumber \\[4mm]
   & & \qquad\qquad
   \sum_{\sigma \in \Sym_n}
   x_1^{\eareccpeak(\sigma)} x_2^{\eareccdfall(\sigma)} 
   y_1^{\ereccval(\sigma)} y_2^{\ereccdrise(\sigma)} 
   \:\times
       \qquad\qquad
       \nonumber \\[-1mm]
   & & \qquad\qquad\qquad\:
   u_1^{\nrcpeak(\sigma)} u_2^{\nrcdfall(\sigma)}
   v_1^{\nrcval(\sigma)} v_2^{\nrcdrise(\sigma)}
   \bw^{\bfix(\sigma)} \lambda^{\cyc(\sigma)}
   \;,
 \label{def.Qnhat}
\end{eqnarray}
which extends \reff{def.Qnbis}
by including the factor $\lambda^{\cyc(\sigma)}$.
We find {\em empirically}\/ that we need make only one specialization
--- either $u_1 = x_1$ or $v_1 = y_1$ ---
to obtain a good J-fraction.
Let us state the latter:

\begin{conjecture}
   \label{thm.perm.Jtype.v2}
The ordinary generating function of the polynomials $\widehat{Q}_n$
specialized to $v_1 = y_1$ has the
J-type continued fraction
\begin{eqnarray}
   & & \hspace*{-7mm}
   \sum_{n=0}^\infty \widehat{Q}_n(x_1,x_2,y_1,y_2,u_1,u_2,y_1,v_2,\bw,\lambda)
       \: t^n
   \;=\;
       \nonumber \\
   & & \hspace*{-4mm}
   \cfrac{1}{1 - \lambda w_0 t - \cfrac{\lambda x_1 y_1 t^2}{1 -  (x_2\!+\!y_2\!+\!\lambda w_1) t - \cfrac{(\lambda\!+\!1) (x_1 \!+\! u_1) y_1 t^2}{1 - (x_2\!+\!y_2\!+\!u_2\!+\!v_2 \!+\! \lambda w_2)t - \cfrac{(\lambda\!+\!2) (x_1 \!+\! 2u_1) y_1 t^2}{1 - \cdots}}}}
       \nonumber \\[1mm]
   \label{eq.thm.perm.Jtype.v2}
\end{eqnarray}
with coefficients
\begin{subeqnarray}
   \gamma_0  & = &   \lambda w_0   \\[1mm]
   \gamma_n  & = &   [x_2 + (n-1)u_2]  \,+\,  [y_2 + (n-1)v_2]  \,+\,  \lambda w_n
        \qquad\hbox{for $n \ge 1$}   \\[1mm]
   \beta_n   & = &   (\lambda+n-1) \, [x_1 + (n-1)u_1] \, y_1
 \label{def.weights.perm.Jtype.v2}
\end{subeqnarray}
\end{conjecture}

\noindent
We have tested this conjecture through $n=12$.\footnote{
   The CPU time and memory required for this computation were as follows:
    \begin{center}
    \begin{tabular}{r@{\qquad}r@{\qquad}r}
       $n$ & \multicolumn{1}{c}{CPU time} & \multicolumn{1}{c}{Memory} \\
           & \multicolumn{1}{c}{(seconds)} & \multicolumn{1}{c}{(gigabytes)} \\
       \hline \\[-2mm]
        7   &      4 &   0.04\  \\
        8   &     32 &   0.06\  \\
        9   &    301 &   0.2\hphantom{4}\   \\
       10   &   3125 &   2\hphantom{.04}\   \\
       11   &  35947 &  24\hphantom{.04}\   \\ 
      12   & 456045 & 299\hphantom{.04}\    \\
    \end{tabular}
    \end{center}
   This computation was carried out in {\sc Mathematica} 11.1.0
   on a Dell PowerEdge R930 computer with 2 TB shared memory
   and four Intel Xeon E7-8891v4 CPUs running at 2.8 GHz.
   Our attempt at computing $n=13$ crashed for lack of memory.
}

%
%

Alas, we are unable at present to prove Conjecture~\ref{thm.perm.Jtype.v2};
we are only able to prove the weaker version
in which we make the {\em two}\/ specializations
$v_1 = y_1$ and $v_2 = y_2$:

\begin{theorem}[Second J-fraction for permutations]
   \label{thm.perm.Jtype.v2.weaker0}
The ordinary generating function of the polynomials $\widehat{Q}_n$
specialized to $v_1 = y_1$ and $v_2 = y_2$ has the
J-type continued fraction
\begin{eqnarray}
   & & \hspace*{-7mm}
   \sum_{n=0}^\infty \widehat{Q}_n(x_1,x_2,y_1,y_2,u_1,u_2,y_1,y_2,\bw,\lambda)
       \: t^n
   \;=\;
       \nonumber \\
   & & \hspace*{-4mm}
   \cfrac{1}{1 - \lambda w_0 t - \cfrac{\lambda x_1 y_1 t^2}{1 -  (x_2\!+\!y_2\!+\!\lambda w_1) t - \cfrac{(\lambda\!+\!1) (x_1 \!+\! u_1) y_1 t^2}{1 - (x_2\!+\!2y_2\!+\!u_2 \!+\! \lambda w_2)t - \cfrac{(\lambda\!+\!2) (x_1 \!+\! 2u_1) y_1  t^2}{1 - \cdots}}}}
       \nonumber \\[1mm]
   \label{eq.thm.perm.Jtype.v2.weaker0}
\end{eqnarray}
with coefficients
\begin{subeqnarray}
   \gamma_0  & = &   \lambda w_0   \\[1mm]
   \gamma_n  & = &   [x_2 + (n-1)u_2]  \,+\,  ny_2  \,+\,  \lambda w_n
        \qquad\hbox{for $n \ge 1$}   \\[1mm]
   \beta_n   & = &   (\lambda+n-1) \, [x_1 + (n-1)u_1] \, y_1
 \label{def.weights.perm.Jtype.v2.weaker0}
\end{subeqnarray}
\end{theorem}

\noindent
We will prove this theorem in Section~\ref{subsec.permutations.J.v2},
as a special case of a more general result.

Comparing \reff{def.Qnhat} with \reff{eq.eulerian.fourvar.cyc},
we see that if we set
\be
   x_1 = x_2 = y,   \quad
   u_1 = u_2 = v,   \quad
   y_1 = y_2 = v_1 = v_2 = w_\ell = u,   \quad
   \lambda = x/u  \;,
 \label{def.specialization.xyuv.cyc}
\ee
then $\widehat{Q}_n(x_1,x_2,y_1,y_2,u_1,u_2,v_1,v_2,\bw,\lambda)$
reduces to
\be
   \sum_{\sigma \in \Sym_n}
      x^{\cyc(\sigma)} y^{\earec(\sigma)}
        u^{n - \aexc(\sigma) - \cyc(\sigma)} v^{\aexc(\sigma) - \earec(\sigma)}
   \;,
\ee
which is of course equivalent to \reff{eq.eulerian.fourvar.cyc}
via the bijection $\sigma \mapsto R \circ \sigma \circ R$ with $R(i) = n+1-i$,
which interchanges $\earec \leftrightarrow \erec$
and $\aexc \leftrightarrow \exc$ while preserving $\cyc$.
With this specialization
the J-fraction coefficients \reff{def.weights.perm.Jtype.v2.weaker0} reduce to
\reff{def.weights.perm.Jtype.specialization.xyuv},
which in turn arise by contraction
from the S-fraction with coefficients \reff{def.weights.eulerian.fourvar}.
So Theorem~\ref{thm.perms.S}(b) is an immediate consequence of
a very special case of Theorem~\ref{thm.perm.Jtype.v2.weaker0}.

Here is another interesting specialization of
Theorem~\ref{thm.perm.Jtype.v2.weaker0}:
let us take
\be
   x_1 = zu_1, \quad
   x_2 = zw, \quad
   v_1 = y_1, \quad
   v_2 = y_2, \quad
   w_0 = z, \quad
   w_\ell = w \hbox{ for } \ell \ge 1
  \;,
\ee
so that $\widehat{Q}_n(x_1,x_2,y_1,y_2,u_1,u_2,v_1,v_2,\bw,\lambda)$
reduces to
\be
   \sum_{\sigma \in \Sym_n}
   z^{\arec(\sigma)}
   y_1^{\cval(\sigma)} y_2^{\cdrise(\sigma)} 
   u_1^{\cpeak(\sigma)} u_2^{\nrcdfall(\sigma)}
   w^{\nrfix(\sigma) + \eareccdfall(\sigma)}
   \lambda^{\cyc(\sigma)}
   \;.
   \quad
 \label{def.Qnhat.arec-special}
\ee
Then these polynomials have a J-fraction with coefficients
\begin{subeqnarray}
   \gamma_0  & = &   \lambda z   \\[1mm]
   \gamma_n  & = &   (\lambda + z) w \,+\, (n-1)u_2 \,+\, ny_2
        \qquad\hbox{for $n \ge 1$}   \\[1mm]
   \beta_n   & = &   (\lambda+n-1) \, (z+n-1) \, u_1 y_1
\end{subeqnarray}
which are symmetric under $z \leftrightarrow \lambda$.
This symmetry of the joint distribution of
\linebreak
${(\arec,\cval,\cdrise,\cpeak,\nrcdfall,\nrfix+\eareccdfall,\cyc)}$
under $\arec \leftrightarrow \cyc$
generalizes the symmetry of $(\arec,\cyc)$
that was found by Cori \cite[Theorem~2]{Cori_09b}.
It would be interesting to find a bijective explanation of this symmetry.

%
%
%

\subsection[$p,q$-generalizations of the first J-fraction]{$\bm{p,q}$-generalizations of the first J-fraction}
   \label{subsec.intro.permutations.pq}

We can further extend Theorems~\ref{thm.perms.S}(a) and \ref{thm.perm.Jtype}
by introducing a $p,q$-generalization.
Recall that for integer $n \ge 0$ we define
\be
   [n]_{p,q}
   \;=\;
   {p^n - q^n \over p-q}
   \;=\;
   \sum\limits_{j=0}^{n-1} p^j q^{n-1-j}
\ee
where $p$ and $q$ are indeterminates;
it is a homogeneous polynomial of degree $n-1$ in $p$ and $q$,
which is symmetric in $p$ and $q$.
In particular, $[0]_{p,q} = 0$ and $[1]_{p,q} = 1$;
and for $n \ge 1$ we have the recurrence
\be
   [n]_{p,q}
   \;=\;
   p \, [n-1]_{p,q} \,+\, q^{n-1}
   \;=\;
   q \, [n-1]_{p,q} \,+\, p^{n-1}
   \;.
 \label{eq.recurrence.npq}
\ee
If $p=1$, then $[n]_{1,q}$ is the well-known $q$-integer
\be
   [n]_q
   \;=\;  [n]_{1,q}
   \;=\; {1 - q^n \over 1-q}
   \;=\;  \begin{cases}
               0  & \textrm{if $n=0$}  \\
               1+q+q^2+\ldots+q^{n-1}  & \textrm{if $n \ge 1$}
          \end{cases}
\ee
If $p=0$, then
\be
   [n]_{0,q}
   \;=\;  \begin{cases}
               0  & \textrm{if $n=0$}  \\
               q^{n-1}  & \textrm{if $n \ge 1$}
          \end{cases}
\ee

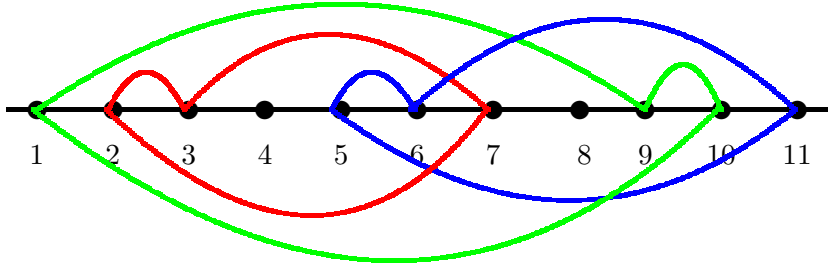
\begin{figure}[t]
\centering
\vspace*{4cm}
\begin{picture}(60,20)(120, -65)
\setlength{\unitlength}{2mm}
\linethickness{.5mm}
\put(-2,0){\line(1,0){54}}
\put(0,0){\circle*{1,3}}\put(0,0){\makebox(0,-6)[c]{\small 1}}
\put(5,0){\circle*{1,3}}\put(5,0){\makebox(0,-6)[c]{\small 2}}
\put(10,0){\circle*{1,3}}\put(10,0){\makebox(0,-6)[c]{\small 3}}
\put(15,0){\circle*{1,3}}\put(15,0){\makebox(0,-6)[c]{\small 4}}
\put(20,0){\circle*{1,3}}\put(20,0){\makebox(0,-6)[c]{\small 5}}
\put(25,0){\circle*{1,3}}\put(25,0){\makebox(0,-6)[c]{\small 6}}
\put(30,0){\circle*{1,3}}\put(30,0){\makebox(0,-6)[c]{\small 7}}
\put(35,0){ \circle*{1,3}}\put(36,0){\makebox(0,-6)[c]{\small 8}}
\put(40,0){\circle*{1,3}}\put(40,0){\makebox(0,-6)[c]{\small 9}}
\put(45,0){\circle*{1,3}}\put(45,0){\makebox(0,-6)[c]{\small 10}}
\put(50,0){\circle*{1,3}}\put(50,0){\makebox(0,-6)[c]{\small 11}}
\green{\qbezier(0,0)(20,14)(40,0)
\qbezier(40,0)(42.5,6)(45,0)}
\red{\qbezier(4,0)(6.5,5)(9,0)
\qbezier(9,0)(18,10)(29,0)}
\blue{\qbezier(18,0)(20.5,5)(23.5,0)
\qbezier(23.2,0)(36,12)(48.5,0)
\qbezier(18,0)(34,-12)(48.5,0)}
\red{\qbezier(2.5,0)(17,-14)(27.5,0)}
\green{\qbezier(-3,0)(22,-20)(42,0)}
\end{picture}
\caption{
   Pictorial representation of the permutation
   $\sigma = 9\,3\,7\,4\,6\,11\,2\,8\,10\,1\,5
           = (1,9,10)\,(2,3,7)\,(4)\,(5,6,11)\,(8) \in \Sym_{11}$.
 \label{fig.pictorial}
}
\end{figure}

The statistics on permutations corresponding to the variables $p$ and $q$
will be {\em crossings}\/ and {\em nestings}\/, defined as follows:
First we associate to each permutation $\sigma \in \Sym_n$
a pictorial representation (Figure~\ref{fig.pictorial})
by placing vertices $1,2,\ldots,n$ along a horizontal axis
and then drawing an arc from $i$ to $\sigma(i)$
above (resp.\ below) the horizontal axis
in case $\sigma(i) > i$ [resp.\ $\sigma(i) < i$];
if $\sigma(i) = i$ we do not draw any arc.
Each vertex thus has either
out-degree = in-degree = 1 (if it is not a fixed point) or
out-degree = in-degree = 0 (if it is a fixed point).
Of course, the arrows on the arcs are redundant,
because the arrow on an arc above (resp.\ below) the axis
always points to the right (resp.\ left).

We then say that a quadruplet $i < j < k < l$ forms an
\begin{itemize}
   \item {\em upper crossing}\/ (ucross) if $k = \sigma(i)$ and $l = \sigma(j)$;
   \item {\em lower crossing}\/ (lcross) if $i = \sigma(k)$ and $j = \sigma(l)$;
   \item {\em upper nesting}\/  (unest)  if $l = \sigma(i)$ and $k = \sigma(j)$;
   \item {\em lower nesting}\/  (lnest)  if $i = \sigma(l)$ and $j = \sigma(k)$.
\end{itemize}
We also consider some ``degenerate'' cases with $j=k$,
by saying that a triplet $i < j < l$ forms an
\begin{itemize}
   \item {\em upper joining}\/ (ujoin) if $j = \sigma(i)$ and $l = \sigma(j)$
      [i.e.\ the index $j$ is a cycle double rise];
   \item {\em lower joining}\/ (ljoin) if $i = \sigma(j)$ and $j = \sigma(l)$
      [i.e.\ the index $j$ is a cycle double fall];
   \item {\em upper pseudo-nesting}\/ (upsnest)
      if $l = \sigma(i)$ and $j = \sigma(j)$;
   \item {\em lower pseudo-nesting}\/ (lpsnest)
      if $i = \sigma(l)$ and $j = \sigma(j)$.
\end{itemize}
These are clearly degenerate cases of crossings and nestings, respectively.
See Figure~\ref{fig.crossnest}.
Note that $\upsnest(\sigma) = \lpsnest(\sigma)$ for all $\sigma$,
since for each fixed point~$j$,
the number of pairs $(i,l)$ with $i < j < l$ such that $l = \sigma(i)$
has to equal the number of such pairs with $i = \sigma(l)$
[it is just $\lev(j,\sigma)$ as defined in \reff{def.level}];
we therefore write these two statistics simply as
\be
   \psnest(\sigma) \;\eqdef\; \upsnest(\sigma) \;=\;  \lpsnest(\sigma)
   \;.
\ee
And of course $\ujoin = \cdrise$ and $\ljoin = \cdfall$.

\begin{figure}[p]
\centering
\begin{picture}(30,15)(145, 10)
\setlength{\unitlength}{1.5mm}
\linethickness{.5mm}
\put(2,0){\line(1,0){28}}
\put(5,0){\circle*{1,3}}\put(5,0){\makebox(0,-6)[c]{\small $i$}}
\put(12,0){\circle*{1,3}}\put(12,0){\makebox(0,-6)[c]{\small $j$}}
\put(19,0){\circle*{1,3}}\put(19,0){\makebox(0,-6)[c]{\small $k$}}
\put(26,0){\circle*{1,3}}\put(26,0){\makebox(0,-6)[c]{\small $l$}}
\red{\qbezier(5,0)(12,10)(19,0)}
\blue{\qbezier(11,0)(18,10)(25,0)}
\put(15,-6){\makebox(0,-6)[c]{\small Upper crossing}}
\put(43,0){\line(1,0){28}}
\put(47,0){\circle*{1,3}}\put(47,0){\makebox(0,-6)[c]{\small $i$}}
\put(54,0){\circle*{1,3}}\put(54,0){\makebox(0,-6)[c]{\small $j$}}
\put(61,0){\circle*{1,3}}\put(61,0){\makebox(0,-6)[c]{\small $k$}}
\put(68,0){\circle*{1,3}}\put(68,0){\makebox(0,-6)[c]{\small $l$}}
\red{\qbezier(47,0)(54,-10)(61,0)}
\blue{\qbezier(53,0)(60,-10)(67,0)}
\put(57,-6){\makebox(0,-6)[c]{\small Lower crossing}}
\end{picture}
\\[3.5cm]
\begin{picture}(30,15)(145, 10)
\setlength{\unitlength}{1.5mm}
\linethickness{.5mm}
\put(2,0){\line(1,0){28}}
\put(5,0){\circle*{1,3}}\put(5,0){\makebox(0,-6)[c]{\small $i$}}
\put(12,0){\circle*{1,3}}\put(12,0){\makebox(0,-6)[c]{\small $j$}}
\put(19,0){\circle*{1,3}}\put(19,0){\makebox(0,-6)[c]{\small $k$}}
\put(26,0){\circle*{1,3}}\put(26,0){\makebox(0,-6)[c]{\small $l$}}
\red{\qbezier(5,0)(15.5,10)(26,0)}
\blue{\qbezier(11,0)(14.5,5)(18,0)}
\put(15,-6){\makebox(0,-6)[c]{\small Upper nesting}}
\put(43,0){\line(1,0){28}}
\put(47,0){\circle*{1,3}}\put(47,0){\makebox(0,-6)[c]{\small $i$}}
\put(54,0){\circle*{1,3}}\put(54,0){\makebox(0,-6)[c]{\small $j$}}
\put(61,0){\circle*{1,3}}\put(61,0){\makebox(0,-6)[c]{\small $k$}}
\put(68,0){\circle*{1,3}}\put(68,0){\makebox(0,-6)[c]{\small $l$}}
\red{\qbezier(47,0)(57,-10)(68,0)}
\blue{\qbezier(53,0)(56.5,-5)(60.5,0)}
\put(57,-6){\makebox(0,-6)[c]{\small Lower nesting}}
\end{picture}
\\[3.5cm]
\begin{picture}(30,15)(145, 10)
\setlength{\unitlength}{1.5mm}
\linethickness{.5mm}
\put(2,0){\line(1,0){28}}
\put(5,0){\circle*{1,3}}\put(5,0){\makebox(0,-6)[c]{\small $i$}}
\put(15.5,0){\circle*{1,3}}\put(15.5,0){\makebox(0,-6)[c]{\small $j$}}
\put(26,0){\circle*{1,3}}\put(26,0){\makebox(0,-6)[c]{\small $l$}}
\red{\qbezier(5,0)(10.5,10)(15.5,0)}
\blue{\qbezier(14.75,0)(19.5,10)(25.25,0)}
\put(15,-6){\makebox(0,-6)[c]{\small Upper joining}}
\put(43,0){\line(1,0){28}}
\put(47,0){\circle*{1,3}}\put(47,0){\makebox(0,-6)[c]{\small $i$}}
\put(57.5,0){\circle*{1,3}}\put(57.5,0){\makebox(0,-6)[c]{\small $j$}}
\put(68,0){\circle*{1,3}}\put(68,0){\makebox(0,-6)[c]{\small $l$}}
\red{\qbezier(47,0)(52,-10)(57.25,0)}
\blue{\qbezier(56.75,0)(62,-10)(67,0)}
\put(57,-6){\makebox(0,-6)[c]{\small Lower joining}}
\end{picture}
\\[3.5cm]
\begin{picture}(30,15)(145, 10)
\setlength{\unitlength}{1.5mm}
\linethickness{.5mm}
\put(2,0){\line(1,0){28}}
\put(5,0){\circle*{1,3}}\put(5,0){\makebox(0,-6)[c]{\small $i$}}
\put(15.5,0){\circle*{1,3}}\put(15.5,0){\makebox(0,-6)[c]{\small $j$}}
\put(26,0){\circle*{1,3}}\put(26,0){\makebox(0,-6)[c]{\small $l$}}
\red{\qbezier(5,0)(15.5,10)(26,0)}
\put(15,-6){\makebox(0,-6)[c]{\small Upper pseudo-nesting}}
\put(43,0){\line(1,0){28}}
\put(47,0){\circle*{1,3}}\put(47,0){\makebox(0,-6)[c]{\small $i$}}
\put(57.5,0){\circle*{1,3}}\put(57.5,0){\makebox(0,-6)[c]{\small $j$}}
\put(68,0){\circle*{1,3}}\put(68,0){\makebox(0,-6)[c]{\small $l$}}
\red{\qbezier(47,0)(57.5,-10)(68,0)}
\put(57,-6){\makebox(0,-6)[c]{\small Lower pseudo-nesting}}
\end{picture}
\vspace*{3cm}
%
\caption{
   Crossing, nesting, joining and pseudo-nesting.
 \label{fig.crossnest}
}
\end{figure}

\begin{figure}[p]
\centering
\begin{picture}(30,15)(145, 10)
\setlength{\unitlength}{1.5mm}
\linethickness{.5mm}
\put(2,0){\line(1,0){28}}
\put(5,0){\circle*{1,3}}\put(5,0){\makebox(0,-6)[c]{\small $i$}}
\put(12,0){\circle*{1,3}}\put(12,0){\makebox(0,-6)[c]{\small $j$}}
\put(19,0){\circle*{1,3}}\put(19,0){\makebox(0,-6)[c]{\small $k$}}
\put(26,0){\circle*{1,3}}\put(26,0){\makebox(0,-6)[c]{\small $l$}}
\red{\qbezier(5,0)(12,10)(19,0)}
\blue{\qbezier(11,0)(18,10)(25,0)}
\blue{\qbezier(10.5,0)(11.75,-1)(13,-2)}
\put(15,-6){\makebox(0,-6)[c]{\small Upper crossing of type cval}}
\put(43,0){\line(1,0){28}}
\put(47,0){\circle*{1,3}}\put(47,0){\makebox(0,-6)[c]{\small $i$}}
\put(54,0){\circle*{1,3}}\put(54,0){\makebox(0,-6)[c]{\small $j$}}
\put(61,0){\circle*{1,3}}\put(61,0){\makebox(0,-6)[c]{\small $k$}}
\put(68,0){\circle*{1,3}}\put(68,0){\makebox(0,-6)[c]{\small $l$}}
\red{\qbezier(47,0)(54,10)(61,0)}
\blue{\qbezier(53,0)(60,10)(67,0)}
\blue{\qbezier(50,2)(51.25,1)(52.5,0)}
\put(57,-6){\makebox(0,-6)[c]{\small Upper crossing of type cdrise}}
\end{picture}
\\[3.5cm]
\begin{picture}(30,15)(145, 10)
\setlength{\unitlength}{1.5mm}
\linethickness{.5mm}
\put(2,0){\line(1,0){28}}
\put(5,0){\circle*{1,3}}\put(5,0){\makebox(0,-6)[c]{\small $i$}}
\put(12,0){\circle*{1,3}}\put(12,0){\makebox(0,-6)[c]{\small $j$}}
\put(19,0){\circle*{1,3}}\put(19,0){\makebox(0,-6)[c]{\small $k$}}
\put(26,0){\circle*{1,3}}\put(26,0){\makebox(0,-6)[c]{\small $l$}}
\red{\qbezier(5,0)(12,-10)(19,0)}
\red{\qbezier(16,2)(17.25, 1)(18.5,0)}
\blue{\qbezier(10.5,0)(18,-10)(25,0)}
\put(15,-6){\makebox(0,-6)[c]{\small Lower crossing of type cpeak}}
\put(43,0){\line(1,0){28}}
\put(47,0){\circle*{1,3}}\put(47,0){\makebox(0,-6)[c]{\small $i$}}
\put(54,0){\circle*{1,3}}\put(54,0){\makebox(0,-6)[c]{\small $j$}}
\put(61,0){\circle*{1,3}}\put(61,0){\makebox(0,-6)[c]{\small $k$}}
\put(68,0){\circle*{1,3}}\put(68,0){\makebox(0,-6)[c]{\small $l$}}
\red{\qbezier(47,0)(54,-10)(61,0)}
\blue{\qbezier(53,0)(60,-10)(67,0)}
\red{\qbezier(59.5,0)(60.75,-1)(62,-2)}
\put(57,-6){\makebox(0,-6)[c]{\small Lower crossing of type cdfall}}
\end{picture}
\\[3.5cm]
\begin{picture}(30,15)(145, 10)
\setlength{\unitlength}{1.5mm}
\linethickness{.5mm}
\put(2,0){\line(1,0){28}}
\put(5,0){\circle*{1,3}}\put(5,0){\makebox(0,-5)[c]{\small $i$}}
\put(12,0){\circle*{1,3}}\put(12,0){\makebox(0,-5)[c]{\small $j$}}
\put(19,0){\circle*{1,3}}\put(19,0){\makebox(0,-5)[c]{\small $k$}}
\put(26,0){\circle*{1,3}}\put(26,0){\makebox(0,-5)[c]{\small $l$}}
\red{\qbezier(5,0)(15.5,10)(26,0)}
\blue{\qbezier(11,0)(14.5,5)(18,0)}
\blue{\qbezier(10.25,0)(11.5,-1)(12.75,-2)}
\put(15,-6){\makebox(0,-6)[c]{\small Upper nesting of type cval}}
\put(43,0){\line(1,0){28}}
\put(47,0){\circle*{1,3}}\put(47,0){\makebox(0,-6)[c]{\small $i$}}
\put(54,0){\circle*{1,3}}\put(54,0){\makebox(0,-6)[c]{\small $j$}}
\put(61,0){\circle*{1,3}}\put(61,0){\makebox(0,-6)[c]{\small $k$}}
\put(68,0){\circle*{1,3}}\put(68,0){\makebox(0,-6)[c]{\small $l$}}
\red{\qbezier(47,0)(57,10)(68,0)}
\blue{\qbezier(53,0)(56.5,5)(60.5,0)}
\blue{\qbezier(50,2)(51.25,1)(52.5,0)}
\put(57,-6){\makebox(0,-6)[c]{\small Upper nesting of type cdrise}}
\end{picture}
\\[3.5cm]
\begin{picture}(30,15)(145, 10)
\setlength{\unitlength}{1.5mm}
\linethickness{.5mm}
\put(2,0){\line(1,0){28}}
\put(5,0){\circle*{1,3}}\put(5,0){\makebox(0,-5)[c]{\small $i$}}
\put(12,0){\circle*{1,3}}\put(12,0){\makebox(0,-5)[c]{\small $j$}}
\put(19,0){\circle*{1,3}}\put(19,0){\makebox(0,-5)[c]{\small $k$}}
\put(26,0){\circle*{1,3}}\put(26,0){\makebox(0,-5)[c]{\small $l$}}
\red{\qbezier(5,0)(15.5,-10)(26,0)}
\blue{\qbezier(11,0)(14.5,-5)(18,0)}
\blue{\qbezier(14.75,2)(16,1)(17.25,0)}
\put(15,-6){\makebox(0,-6)[c]{\small Lower nesting of type cpeak}}
\put(43,0){\line(1,0){28}}
\put(47,0){\circle*{1,3}}\put(47,0){\makebox(0,-5)[c]{\small $i$}}
\put(54,0){\circle*{1,3}}\put(54,0){\makebox(0,-5)[c]{\small $j$}}
\put(61,0){\circle*{1,3}}\put(61,0){\makebox(0,-5)[c]{\small $k$}}
\put(68,0){\circle*{1,3}}\put(68,0){\makebox(0,-5)[c]{\small $l$}}
\red{\qbezier(47,0)(57,-10)(68,0)}
\blue{\qbezier(53,0)(56.5,-5)(60.5,0)}
\blue{\qbezier(59.5,0)(60.75, -1)(62,-2)}
\put(57,-6){\makebox(0,-6)[c]{\small Lower nesting of type cdfall}}
\end{picture}
\vspace*{3cm}
%
%
%
\caption{
   Refined categories of crossing and nesting.
 \label{fig.refined_crossnest}
}
\end{figure}

Note also that
\begin{subeqnarray}
   \ucross(\sigma)  & = &  \lcross(\sigma^{-1})    \\[1mm]   
   \unest(\sigma)   & = &  \lnest(\sigma^{-1})    \\[1mm]   
   \ujoin(\sigma)   & = &  \ljoin(\sigma^{-1})    \\[1mm]   
   \psnest(\sigma)  & = &  \psnest(\sigma^{-1})
 \label{eq.crossnest.inverses}
\end{subeqnarray}

We can further refine the four crossing/nesting categories
by examining more closely the status of the inner index ($j$ or $k$)
whose {\em outgoing}\/ arc belonged to the crossing or nesting:
we say that a quadruplet $i < j < k < l$ forms an
\begin{itemize}
   \item {\em upper crossing of type cval}\/ (ucrosscval)
       if $k = \sigma(i)$ and $l = \sigma(j)$ and ${\sigma^{-1}(j) > j}$;
   \item {\em upper crossing of type cdrise}\/ (ucrosscdrise)
       \hbox{if $k = \sigma(i)$ and $l = \sigma(j)$ and ${\sigma^{-1}(j) < j}$;}
   \item {\em lower crossing of type cpeak}\/ (lcrosscpeak)
       if $i = \sigma(k)$ and $j = \sigma(l)$ and ${\sigma^{-1}(k) < k}$;
   \item {\em lower crossing of type cdfall}\/ (lcrosscdfall)
       if $i = \sigma(k)$ and $j = \sigma(l)$ and ${\sigma^{-1}(k) > k}$;
   \item {\em upper nesting of type cval}\/  (unestcval)
       if $l = \sigma(i)$ and $k = \sigma(j)$ and ${\sigma^{-1}(j) > j}$;
   \item {\em upper nesting of type cdrise}\/  (unestcdrise)
       if $l = \sigma(i)$ and $k = \sigma(j)$ and ${\sigma^{-1}(j) < j}$;
   \item {\em lower nesting of type cpeak}\/  (lnestcpeak)
       if $i = \sigma(l)$ and $j = \sigma(k)$ and ${\sigma^{-1}(k) < k}$;
   \item {\em lower nesting of type cdfall}\/  (lnestcdfall)
       if $i = \sigma(l)$ and $j = \sigma(k)$ and ${\sigma^{-1}(k) > k}$.
\end{itemize}
See Figure~\ref{fig.refined_crossnest}.
Please note that for the ``upper'' quantities
the distinguished index
(i.e.\ the one for which we examine both $\sigma$ and $\sigma^{-1}$)
is in second position ($j$),
while for the ``lower'' quantities
the distinguished index is in third position ($k$).

The master polynomial encoding all these statistics (but no others yet) is
\begin{eqnarray}
   & &
   P_n(p_{+1},p_{+2},p_{-1},p_{-2},q_{+1},q_{+2},q_{-1},q_{-2},r_+,r_-,s)
   \;=\;
       \nonumber \\[4mm]
   & & \qquad\quad
   \sum_{\sigma \in \Sym_n}
   p_{+1}^{\ucrosscval(\sigma)}
   p_{+2}^{\ucrosscdrise(\sigma)}
   p_{-1}^{\lcrosscpeak(\sigma)}
   p_{-2}^{\lcrosscdfall(\sigma)}
          \:\times
       \qquad\qquad
       \nonumber \\[-1mm]
   & & \qquad\qquad\quad
   q_{+1}^{\unestcval(\sigma)}
   q_{+2}^{\unestcdrise(\sigma)}
   q_{-1}^{\lnestcpeak(\sigma)} 
   q_{-2}^{\lnestcdfall(\sigma)} 
          \:\times
       \qquad\qquad
       \nonumber \\[3mm]
   & & \qquad\qquad\quad
   r_+^{\ujoin(\sigma)}  r_-^{\ljoin(\sigma)}
   s^{\psnest(\sigma)}
   \;.
       \qquad\quad
 \label{def.Pn.pq}
\end{eqnarray}
It turns out that these 11-variable polynomials have a beautiful J-fraction:
   
\begin{theorem}[J-fraction for crossing and nesting statistics]
   \label{thm.perm.pq.Jtype}
The ordinary generating function of the polynomials $P_n$
defined by \reff{def.Pn.pq} has the
J-type continued fraction
\begin{eqnarray}
   & & \hspace*{-7mm}
   \sum_{n=0}^\infty
   P_n(p_{+1},p_{+2},p_{-1},p_{-2},q_{+1},q_{+2},q_{-1},q_{-2},r_+,r_-,s) \: t^n
   \;=\;
       \nonumber \\
   & & \!\!\! 
   \cfrac{1}{1 - t - \cfrac{t^2}{1 -  (r_+\!+\!r_-\!+\!s) t - \cfrac{[2]_{p_{+1},q_{+1}} [2]_{p_{-1},q_{-1}} t^2}{1 - ([2]_{p_{+2},q_{+2}} r_+\!+\! [2]_{p_{-2},q_{-2}} r_-\!+\!s^2)t - \cfrac{[3]_{p_{+1},q_{+1}} [3]_{p_{-1},q_{-1}} t^2}{1 - \cdots}}}}
       \nonumber \\[1mm]
   \label{eq.thm.perm.pq.Jtype}
\end{eqnarray}
with coefficients
\begin{subeqnarray}
   \gamma_n  & = &
      [n]_{p_{+2},q_{+2}} r_+ \,+\, [n]_{p_{-2},q_{-2}} r_- \,+\, s^n
                             \\[1mm]
   \beta_n   & = &   [n]_{p_{+1},q_{+1}} \, [n]_{p_{-1},q_{-1}}
 \label{def.weights.perm.pq.Jtype}
\end{subeqnarray}
\end{theorem}

\noindent
We will prove this theorem in Section~\ref{subsec.permutations.J},
as a special case of a more general result.

The continued fraction
\reff{eq.thm.perm.pq.Jtype}/\reff{def.weights.perm.pq.Jtype}
of course has the symmetry
$(p_{+1},p_{+2},q_{+1},q_{+2},r_+)
 \leftrightarrow
 (p_{-1},p_{-2},q_{-1},q_{-2},r_-)$,
which is obvious from the definition \reff{def.Pn.pq}
by considering the bijection $\sigma \mapsto \sigma^{-1}$.
Less trivially, it has the four separate symmetries
$p_{+1} \leftrightarrow q_{+1}$, $p_{+2} \leftrightarrow q_{+2}$,
$p_{-1} \leftrightarrow q_{-1}$ and $p_{-2} \leftrightarrow q_{-2}$;
it would be interesting to understand these combinatorially.

If one tries to expand the generating function \reff{eq.thm.perm.pq.Jtype}
into an S-fraction, one obtains coefficients $\alpha_n$ that are
rational functions rather than polynomials (starting at $n=4$).
But under the specialization
\be
   p_{+1} = p_{+2}, \quad
   q_{+1} = q_{+2}, \quad
   r_+ = 1, \quad
   r_- = p_{-1} = p_{-2}, \quad
   s = q_{-1} = q_{-2}  \;,
 \label{def.specialization.perm.pq.Stype}
\ee
the J-fraction \reff{eq.thm.perm.pq.Jtype}/\reff{def.weights.perm.pq.Jtype}
does arise as the contraction \reff{eq.contraction_even.coeffs}
of an S-fraction with polynomial coefficients:

\begin{corollary}[S-fraction for crossing and nesting statistics]
   \label{cor.perm.pq.Stype}
The ordinary generating function of the polynomials $P_n$
defined by \reff{def.Pn.pq},
under the specialization \reff{def.specialization.perm.pq.Stype},
has the S-type continued fraction
\be
   \sum_{n=0}^\infty P_n(p_+,p_+,p_-,p_-,q_+,q_+,q_-,q_-,1,p_-,q_-) \: t^n
   \;\,=\;\,
   \cfrac{1}{1 - \cfrac{t}{1 -  \cfrac{t}{1 - \cfrac{[2]_{p_-,q_-} t}{1 - \cfrac{[2]_{p_+,q_+} t}{1 - \cdots}}}}}
       \nonumber \\[1mm]
   \label{eq.cor.perm.pq.Stype}
\ee
with coefficients
\begin{subeqnarray}
   \alpha_{2k-1}  & = &  [k]_{p_-,q_-}  \\[1mm]
   \alpha_{2k}    & = &  [k]_{p_+,q_+}
 \label{def.weights.perm.pq.Stype}
\end{subeqnarray}
\end{corollary}

{\bf Remark.}
We could of course interchange $+$ and $-$ everywhere,
i.e.\ set
\be
   p_{-1} = p_{-2}, \quad
   q_{-1} = q_{-2}, \quad
   r_- = 1, \quad
   r_+ = p_{+1} = p_{+2}, \quad
   s = q_{+1} = q_{+2}  \;,
 \label{def.specialization.perm.pq.Stype.bis}
\ee
and obtain $\alpha_{2k-1} = [k]_{p_+,q_+}$, $\alpha_{2k} = [k]_{p_-,q_-}$.
\myendremark

\medskip

The next step is to try to include the variables $x,y,u,v$
from \reff{eq.eulerian.fourvar.arec}
--- or~more ambitiously,
the variables $x_1,x_2,y_1,y_2,z,u_1,u_2,v_1,v_2,w$ from \reff{def.Qn};
or yet more ambitiously,
the variables $x_1,x_2,y_1,y_2,u_1,u_2,v_1,v_2,\bw$ from \reff{def.Qnbis}.
Amazingly, this latter program works.
Let us define the
polynomial
\begin{eqnarray}
   & &
   Q_n(x_1,x_2,y_1,y_2,u_1,u_2,v_1,v_2,\bw,p_{+1},p_{+2},p_{-1},p_{-2},q_{+1},q_{+2},q_{-1},q_{-2},s)
   \;=\;
       \nonumber \\[4mm]
   & & \qquad\qquad
   \sum_{\sigma \in \Sym_n}
   x_1^{\eareccpeak(\sigma)} x_2^{\eareccdfall(\sigma)} 
   y_1^{\ereccval(\sigma)} y_2^{\ereccdrise(\sigma)} 
   \:\times
       \qquad\qquad
       \nonumber \\[-1mm]
   & & \qquad\qquad\qquad\:
   u_1^{\nrcpeak(\sigma)} u_2^{\nrcdfall(\sigma)}
   v_1^{\nrcval(\sigma)} v_2^{\nrcdrise(\sigma)}
   \bw^{\bfix(\sigma)}  \:\times
       \qquad\qquad
       \nonumber \\[3mm]
   & & \qquad\qquad\qquad\:
   p_{+1}^{\ucrosscval(\sigma)}
   p_{+2}^{\ucrosscdrise(\sigma)}
   p_{-1}^{\lcrosscpeak(\sigma)}
   p_{-2}^{\lcrosscdfall(\sigma)}
          \:\times
       \qquad\qquad
       \nonumber \\[3mm]
   & & \qquad\qquad\qquad\:
   q_{+1}^{\unestcval(\sigma)}
   q_{+2}^{\unestcdrise(\sigma)}
   q_{-1}^{\lnestcpeak(\sigma)}
   q_{-2}^{\lnestcdfall(\sigma)}
   s^{\psnest(\sigma)}
   \;.
 \label{def.Qn.BIG}
\end{eqnarray}
(We have omitted $r_+$ and $r_-$ since they are redundant
 in view of $\ujoin = \cdrise$ and $\ljoin = \cdfall$.)
We then have a beautiful J-fraction that simultaneously generalizes
Theorems~\ref{thm.perm.Jtype} and \ref{thm.perm.pq.Jtype}:

\begin{theorem}[First J-fraction for permutations, $p,q$-generalization]
   \label{thm.perm.pq.Jtype.BIG}
The ordinary generating function of the polynomials $Q_n$
defined by \reff{def.Qn.BIG} has the J-type continued fraction
\begin{eqnarray}
   & & \hspace*{-7mm}
   \sum_{n=0}^\infty
       Q_n(x_1,x_2,y_1,y_2,u_1,u_2,v_1,v_2,\bw,p_{+1},p_{+2},p_{-1},p_{-2},q_{+1},q_{+2},q_{-1},q_{-2},s) \: t^n
   \;=\;
       \nonumber \\
   & & \!\!\!\!
\Scale[0.7]{
   \cfrac{1}{1 - w_0 t - \cfrac{x_1 y_1 t^2}{1 -  (x_2\!+\!y_2\!+\!sw_1) t - \cfrac{(p_{-1} x_1\!+\! q_{-1} u_1)(p_{+1} y_1\!+\! q_{+1} v_1) t^2}{1 - (p_{-2} x_2\!+\!q_{-2} u_2\!+\!p_{+2} y_2\!+\! q_{+2} v_2\!+\! s^2 w_2)t - \cfrac{(p_{-1}^2 x_1\!+\! q_{-1} [2]_{p_{-1},q_{-1}} u_1)(p_{+1}^2 y_1\!+\! q_{+1} [2]_{p_{+1},q_{+1}} v_1) t^2}{1 - \cdots}}}}
}
       \nonumber \\[1mm]
   \label{eq.thm.perm.pq.Jtype.BIG}
\end{eqnarray}
with coefficients
\begin{subeqnarray}
   \gamma_0  & = &   w_0
      \slabel{def.weights.perm.pq.Jtype.BIG.a}  \\[1mm]
   \gamma_n  & = &   (p_{-2}^{n-1} x_2 + q_{-2} \, [n-1]_{p_{-2},q_{-2}} u_2)
                 \:+\: (p_{+2}^{n-1} y_2 + q_{+2} \, [n-1]_{p_{+2},q_{+2}} v_2) \:+\: s^n w_n
       \qquad \nonumber  \\
    & & \hspace*{3.7in}
        \quad\hbox{for $n \ge 1$}
      \slabel{def.weights.perm.pq.Jtype.BIG.b}  \\
   \beta_n   & = &   (p_{-1}^{n-1} x_1 + q_{-1} \, [n-1]_{p_{-1},q_{-1}} u_1)
                  \: (p_{+1}^{n-1} y_1 + q_{+1} \, [n-1]_{p_{+1},q_{+1}} v_1)
      \slabel{def.weights.perm.pq.Jtype.BIG.c}
 \label{def.weights.perm.pq.Jtype.BIG}
\end{subeqnarray}
\end{theorem}

\noindent
We will prove this theorem in Section~\ref{subsec.permutations.J},
as a special case of a more general result.

\bigskip

{\bf Remarks.}
1.  The continued fraction \reff{eq.thm.perm.pq.Jtype.BIG}
shows that $Q_n$ depends on its arguments
only via the combinations (\ref{def.weights.perm.pq.Jtype.BIG}a,b,c):
in particular, it is symmetric under
$(x_1,u_1,p_{-1},q_{-1}) \leftrightarrow (y_1,v_1,p_{+1},q_{+1})$
and independently under
$(x_2,u_2,p_{-2},q_{-2}) \leftrightarrow (y_2,v_2,p_{+2},q_{+2})$.
It would be interesting to try to understand combinatorially
(directly at the level of permutations) why these properties hold.
The symmetry of $Q_n$ under the simultaneous application
of {\em both}\/ of these interchanges
is an immediate consequence of the bijection
$\sigma \mapsto R \circ \sigma \circ R$ with $R(i) = n+1-i$,
which interchanges cycle peaks with cycle valleys,
cycle double rises with cycle double falls,
and records with antirecords
(while preserving the number of fixed points at each level).
But the more specfic properties encoded in \reff{def.weights.perm.pq.Jtype.BIG}
remain mysterious.

2.  Please note that when $x_1 = u_1$,
then $p_{-1}^{n-1} x_1 + q_{-1} \, [n-1]_{p_{-1},q_{-1}} u_1$
simplifies to $[n]_{p_{-1},q_{-1}} x_1$,
which has the symmetry $p_{-1} \leftrightarrow q_{-1}$;
and similarly when $x_2 = u_2$, ${y_1 = v_1}$ or $y_2 = v_2$.
It would be interesting to understand combinatorially
why these symmetries hold.
If we make all four of these specializations
(i.e.\ forgo taking account of record statistics),
then the coefficients \reff{def.weights.perm.pq.Jtype.BIG} simplify to
\begin{subeqnarray}
   \gamma_n  & = &   [n]_{p_{-2},q_{-2}} x_2
                 \:+\: [n]_{p_{+2},q_{+2}} y_2 \:+\: s^n w_n
       \\[2mm]
   \beta_n   & = &   [n]_{p_{-1},q_{-1}}
                  \, [n]_{p_{+1},q_{+1}} \, x_1 y_1
 \label{def.weights.perm.pq.Jtype.BIG.simplified}
\end{subeqnarray}

3.  Once one includes the detailed classification of fixed points
according to their level \reff{def.level},
the variable $s$ becomes redundant:
it simply sends $w_\ell \to s^\ell w_\ell$.
This reflects the fact that
\be
   \psnest(\sigma)  \;=\; \sum_{j \in {\rm fix}} \lev(j,\sigma)
   \;.
 \label{eq.redundant.s}
\ee
\myendremark

\bigskip

{\bf Some historical remarks.}
1.  The pioneering work on crossings and nestings in permutations
is that of Corteel \cite{Corteel_07}, and our presentation is
strongly inspired by hers.
However, her definitions of crossings and nestings are less refined than ours,
and also partly asymmetrical between ``upper'' and ``lower''
(i.e.\ between $\sigma$ and $\sigma^{-1}$): she defines
\begin{subeqnarray}
   \cross(\sigma)  & = & \ucross(\sigma)  + \lcross(\sigma)  + \ujoin(\sigma)  
        \\[1mm]
   \nest(\sigma)   & = & \unest(\sigma)  + \lnest(\sigma)  + \psnest(\sigma)  
 \label{def.corteel.crossnest}
\end{subeqnarray}
and subsequent workers \cite{Josuat-Verges_10,Josuat-Verges_11,Shin_10,Shin_12}
have followed her in this definition.
Here $\nest(\sigma) = \nest(\sigma^{-1})$ by \reff{eq.crossnest.inverses};
but $\cross(\sigma) \neq \cross(\sigma^{-1})$ in general,
because of the appearance of $\ujoin$ without $\ljoin$
in (\ref{def.corteel.crossnest}a).
Corteel \cite{Corteel_07} obtained the special case
of Theorem~\ref{thm.perm.pq.Jtype.BIG} with three variables $y,p,q$,
where $x_1 = x_2 = u_1 = u_2 = 1$,
$y_1 = v_1 = w_\ell = y$,
$p_{+1} = p_{+2} = p_{-1} = p_{-2} = p$,
$y_2 = v_2 = p y$,
$q_{+1} = q_{+2} = q_{-1} = q_{-2} = s = q$.\footnote{
   Her usage of the variables $p$ and $q$ is the reverse of ours
   --- i.e.\ she writes $q^{\cross(\sigma)} p^{\nest(\sigma)}$ ---
   but this does not matter in comparing her formulae to ours,
   because in her specialization the formulae are anyway symmetric
   between $p$ and~$q$.
}
Some refinements of this result were obtained by
Shin and Zeng \cite{Shin_10,Shin_12}.

2.  Our refined categories of crossing and nesting are inspired
by the ``octabasic'' permutation polynomials
of Simion and Stanton \cite[Theorem~2.2]{Simion_96}.

3.  Shortly after completing the proof of Theorem~\ref{thm.perm.pq.Jtype.BIG},
we discovered that these ideas were anticipated two decades ago
in a remarkable but apparently little-known paper of
Randrianarivony \cite{Randrianarivony_98b}.\footnote{
   According to the Web of Science,
   the paper \cite{Randrianarivony_98b} has been cited
   only once \cite{Kasraoui_11}.
}
In this paper, which was also inspired in part by the work of
Simion and Stanton \cite{Simion_94,Simion_96},
Randrianarivony \cite[Th\'eor\`eme~1]{Randrianarivony_98b}
{\em almost}\/ obtained our full Theorem~\ref{thm.perm.pq.Jtype.BIG}:
he had all our variables except $\bw$ and $s$.
That is, he included all our statistics except
the classification of fixed points by level.

One reason for the unfortunate neglect of Randrianarivony's work
may be that it was written in French.
Another reason may be that he did not write {\em explicitly}\/
any continued fraction, contenting himself with an assertion
\cite[Th\'eor\`eme~1]{Randrianarivony_98b}
about the moment sequence associated to the orthogonal polynomials
satisfying a particular three-term recurrence.
But in the introduction to his paper he stated explicitly
\cite[p.~507]{Randrianarivony_98b} that
``Le probl\`eme est alors \'equivalent \`a la donn\'ee d'une
interpr\'etation des coefficients du d\'eveloppement de Taylor
de la J-fraction continue avec les param\`etres $b_n$ et $\lambda_n$.''

4.  A very recent paper of Blitvi\'c and Steingr\'{\i}msson \cite{Blitvic_20}
also contains a large part of Theorem~\ref{thm.perm.pq.Jtype.BIG},
namely, the specialization $x_1 = u_1$, $x_2 = u_2$, $y_1 = v_1$, $y_2 = v_2$,
$w_\ell = w \:\forall \ell$ \cite[Theorem~1]{Blitvic_20}.
That is, they have included all the cycle statistics
and the refined crossing-and-nesting statistics
(including the pseudo-nesting statistic),
but not the record statistics or the
detailed classification of fixed points by level.
As a result, they obtained the J-fraction with coefficients
\reff{def.weights.perm.pq.Jtype.BIG.simplified}.

5.  Some other $q$-versions, involving the number of inversions in $\sigma$,
were obtained two decades ago by one of us \cite{Zeng_89,Zeng_95}
and recently by Elizalde \cite[eqn.~(4)]{Elizalde_18}.
All of these are special cases of Theorem~\ref{thm.perm.pq.Jtype.BIG}
(or its S-fraction corollary,
 Theorem~\ref{thm.perm.pq.Stype.BIG.1} below),
as we explain in Section~\ref{subsec.permutations.inv}.
%
\myendremark

%

\subsection[$p,q$-generalizations of the first S-fraction]{$\bm{p,q}$-generalizations of the first S-fraction}
   \label{subsec.intro.permutations.pq.S}

It turns out that suitable specializations of the J-fraction
\reff{eq.thm.perm.pq.Jtype.BIG}/\reff{def.weights.perm.pq.Jtype.BIG}
can be obtained as the contraction \reff{eq.contraction_even.coeffs}
of an S-fraction with polynomial coefficients.
%
We make the specializations
\begin{subeqnarray}
  x_1  & = &  x  \\
  x_2  & = &  p_- x  \\
  y_1 \;=\; y_2  & = &  y  \\
  u_1  & = &  u  \\
  u_2  & = &  p_- u  \\
  v_1 \;=\; v_2  & = &  v  \\
  w_0  & = &  x  \\
  w_\ell & = & u \quad\hbox{for $\ell \ge 1$}  \\
  p_{+1} \;=\; p_{+2}  & = &  p_+  \\
  p_{-1} \;=\; p_{-2}  & = &  p_-  \\
  q_{+1} \;=\; q_{+2}  & = &  q_+  \\
  q_{-1} \;=\; q_{-2}  & = &  q_-  \\
  s  & = &  q_-
 \label{eq.permutations.S-fraction.pq.specializations}
\end{subeqnarray}
[which include \reff{def.specialization.perm.pq.Stype}
 as the special case $x=y=u=v=1$].
The polynomial $Q_n$ defined in \reff{def.Qn.BIG} then reduces to
the eight-variable polynomial
\begin{eqnarray}
   P_n(x,y,u,v,p_+,p_-,q_+,q_-)
   & = &
   \sum_{\sigma \in \Sym_n}
       x^{\arec(\sigma)} y^{\erec(\sigma)}
       u^{n - \exc(\sigma) - \arec(\sigma)} v^{\exc(\sigma) - \erec(\sigma)}
       \:\times
       \nonumber \\[-1mm]
   & & \qquad\,
      p_+^{\ucross(\sigma)} p_-^{\lcross(\sigma) + \ljoin(\sigma)}
      q_+^{\unest(\sigma)}  q_-^{\lnest(\sigma) + \psnest(\sigma)}
      \;.
            \nonumber \\
 \label{def.Pn1.BIG.0}
\end{eqnarray}
and the coefficients \reff{def.weights.perm.pq.Jtype.BIG} reduce to
\begin{subeqnarray}
   \gamma_0  & = &   x   \\[1mm]
   \gamma_n  & = &   (p_{-}^{n} x + q_{-} \, [n]_{p_{-},q_{-}} u)
                 \:+\: (p_{+}^{n-1} y + q_{+} \, [n-1]_{p_{+},q_{+}} v)
        \quad\hbox{for $n \ge 1$} \qquad   \\[1mm]
   \beta_n   & = &   (p_{-}^{n-1} x + q_{-} \, [n-1]_{p_{-},q_{-}} u)
                  \: (p_{+}^{n-1} y + q_{+} \, [n-1]_{p_{+},q_{+}} v)
   \label{def.weights.perm.pq.Jtype.BIG.specialization}
\end{subeqnarray}
Therefore, the J-fraction
\reff{eq.thm.perm.pq.Jtype.BIG}/\reff{def.weights.perm.pq.Jtype.BIG}
with the specializations \reff{eq.permutations.S-fraction.pq.specializations}
is the contraction \reff{eq.contraction_even.coeffs}
of the following S-fraction:

\begin{theorem}[First S-fraction for permutations, $p,q$-generalization]
   \label{thm.perm.pq.Stype.BIG.1}
The ordinary generating function of the polynomials $P_n$
defined by \reff{def.Pn1.BIG.0} has the S-type continued fraction
\begin{eqnarray}
   & & \hspace*{-7mm}
   \sum_{n=0}^\infty
       P_n(x,y,u,v,p_+,p_-,q_+,q_-) \: t^n
   \;=\;
       \nonumber \\[-1mm]
   & & \qquad\qquad
   \cfrac{1}{1 - \cfrac{xt}{1 - \cfrac{yt}{1 - \cfrac{(p_- x + q_- u)t}{1- \cfrac{(p_+ y + q_+ v)t}{1 - \cfrac{(p_-^2 x+ q_- \, [2]_{p_-,q_-} u)t}{1 - \cfrac{(p_+^2 y+ q_+ \, [2]_{p_+,q_+} v)t}{1-\cdots}}}}}}}
   \label{eq.thm.perm.pq.Stype.BIG.1}
\end{eqnarray}
with coefficients
\begin{subeqnarray}
   \alpha_{2k-1}  & = &   p_-^{k-1} x + q_- \, [k-1]_{p_-,q_-} u  \\[1mm]
   \alpha_{2k}    & = &   p_+^{k-1} y + q_+ \, [k-1]_{p_+,q_+} v
 \label{def.weights.thm.perm.pq.Stype.BIG.1}
\end{subeqnarray}
\end{theorem}

Note that the polynomial $P_n(x,y,u,v,p_+,p_-,q_+,q_-)$
reduces to $P_n(x,y,u,v)$
[cf.\ \reff{eq.eulerian.fourvar.contfrac}/\reff{eq.eulerian.fourvar.arec}]
if we set $p_+ = p_- = q_+ = q_- = 1$,
and to $P_n(p_+,p_+,p_-,p_-,q_+,q_+,q_-,q_-,1,p_-,q_-)$
[cf.\ \reff{def.Pn.pq}/\reff{def.specialization.perm.pq.Stype}]
if we set $x=y=u=v=1$.
So Theorem~\ref{thm.perm.pq.Stype.BIG.1}
simultaneously generalizes Theorem~\ref{thm.perms.S}(a)
and Corollary~\ref{cor.perm.pq.Stype}.

One other interesting specialization of Theorem~\ref{thm.perm.pq.Stype.BIG.1} is
\be
   P_n(x,qy,1,q,q,q,q^2,q^2)
   \;=\;
   \sum_{\sigma \in \Sym_n}
       x^{\arec(\sigma)} y^{\erec(\sigma)} q^{\inv(\sigma)}
   \;,
\ee
where $\inv(\sigma)$ is the inversion number
(see Section~\ref{subsec.permutations.inv} below);
this formula follows from \reff{def.Pn1.BIG.0} and \reff{eq.prop.inv.a}.
The corresponding continued fraction
\reff{eq.thm.perm.pq.Stype.BIG.1}/\reff{def.weights.thm.perm.pq.Stype.BIG.1}
was found in \cite{Zeng_89} and has coefficients
\begin{subeqnarray}
   \alpha_{2k-1}  & = &   q^{k-1} \, (x + q + \ldots + q^{k-1})  \\[1mm]
   \alpha_{2k}    & = &   q^k \, (y + q + \ldots + q^{k-1})
 \label{def.weights.zeng89}
\end{subeqnarray}
As one might expect by analogy with \reff{eq.dumont},
the corresponding ordinary generating function \reff{eq.thm.perm.pq.Stype.BIG.1}
can be written as a ratio of basic hypergeometric functions $\phiHyper{2}{0}$,
which are defined by \cite[pp.~4--5]{Gasper_04}
\be
   \phiHYPER{2}{0}{a,b}{\hbox{---}}{q}{z}
   \;\eqdef\;
   \sum_{n=0}^\infty {(a;q)_n \: (b;q)_n \over (q;q)_n} \:
         (-1)^n \, q^{-n(n-1)/2} \, z^n
\ee
and can be related to the Heine hypergeometric function $\phiHyper{2}{1}$,
\be
   \phiHYPER{2}{1}{a,b}{c}{q}{z}
   \;\eqdef\;
   \sum_{n=0}^\infty {(a;q)_n \: (b;q)_n \over (c;q)_n \: (q;q)_n} \:  z^n
   \;,
\ee
either as a limiting case
\be
   \phiHYPER{2}{0}{a,b}{\hbox{---}}{q}{z}
   \;=\;
   \lim_{c \to\infty} \phiHYPER{2}{1}{a,b}{c}{q}{cz}
 \label{eq.2phi0.limit}
\ee
or as a specialization to $c=0$:
\be
   \phiHYPER{2}{0}{a,b}{\hbox{---}}{q}{z}
   \;=\;
   \phiHYPER{2}{1}{a^{-1},b^{-1}}{0}{q^{-1}}{abq^{-1} z}
   \;.
 \label{eq.2phi0.specialization}
\ee
Starting from Heine's \cite{Heine_1847} continued fraction
for ratios of contiguous $\phiHyper{2}{1}$
\cite[pp.~318--322]{Lorentzen_92} \cite[p.~395]{Cuyt_08},
\be
   {\phiHYPER{2}{1}{a,bq}{cq}{q}{z}
    \over
    \phiHYPER{2}{1}{a,b}{c}{q}{z}
   }
   \;=\;
   \cfrac{1}{1 - \cfrac{\alpha_1 z}{1 - \cfrac{\alpha_2 z}{1 - \cdots}}}
\ee
with coefficients
\begin{subeqnarray}
   \alpha_{2k-1}  & = &   {(1-aq^{k-1}) \, (b-cq^{k-1}) \, q^{k-1}
                           \over
                           (1-cq^{2k-2}) \, (1-cq^{2k-1})
                          }
       \\[2mm]
   \alpha_{2k}    & = &   {(1-bq^k) \, (a-cq^k) \, q^{k-1}
                           \over
                           (1-cq^{2k-1}) \, (1-cq^{2k})
                          }
 \label{def.weights.heine}
\end{subeqnarray}
and applying either \reff{eq.2phi0.limit} or \reff{eq.2phi0.specialization},
we obtain a continued fraction for ratios of contiguous $\phiHyper{2}{0}$:
\be
   {\phiHYPER{2}{0}{a,bq}{\hbox{---}}{q}{q^{-1} z}
    \over
    \phiHYPER{2}{0}{a,b}{\hbox{---}}{q}{z}
   }
   \;=\;
   \cfrac{1}{1 - \cfrac{\alpha_1 z}{1 - \cfrac{\alpha_2 z}{1 - \cdots}}}
 \label{eq.contfrac.2phi0}
\ee
with coefficients
\begin{subeqnarray}
   \alpha_{2k-1}  & = &   -(1-aq^{k-1}) \, q^{-(2k-1)}
       \\[2mm]
   \alpha_{2k}    & = &   -(1-bq^k) \, q^{-2k}
 \label{def.weights.heine.2phi0}
\end{subeqnarray}
The continued fraction \reff{def.weights.zeng89}
can then be obtained from \reff{eq.contfrac.2phi0}/\reff{def.weights.heine.2phi0}
by setting
\begin{subeqnarray}
   a & = &  r^{-1} \, [r \,+\, (1-r) \, x]  \\[1mm]
   b & = &  r \,+\, (1-r) \, y  \\[1mm]
   q & = &  r^{-1}  \\[1mm]
   z & = &  {t \over 1-r}
\end{subeqnarray}
and then renaming $r \leftarrow q$.

\bigskip

{\bf Remark.}
There is an alternative specialization of
\reff{eq.thm.perm.pq.Jtype.BIG}/\reff{def.weights.perm.pq.Jtype.BIG},
generalizing \reff{def.specialization.perm.pq.Stype.bis},
in which $+$ and $-$ are interchanged
compared to \reff{eq.permutations.S-fraction.pq.specializations},
and in which also the roles of $(x,u)$ and $(y,v)$ are interchanged:
\begin{subeqnarray}
  x_1 \;=\; x_2  & = &  x  \\
  y_1  & = &  y  \\
  y_2  & = &  p_+ y  \\
  u_1  \;=\; u_2  & = &  u  \\
  v_1  & = &  v  \\
  v_2  & = &  p_+ v  \\
  w_0  & = &  y  \\
  w_\ell & = & v  \quad\hbox{for $\ell \ge 1$}  \\
  p_{+1} \;=\; p_{+2}  & = &  p_+  \\
  p_{-1} \;=\; p_{-2}  & = &  p_-  \\
  q_{+1} \;=\; q_{+2}  & = &  q_+  \\
  q_{-1} \;=\; q_{-2}  & = &  q_-  \\
  s  & = &  q_+
 \label{eq.permutations.S-fraction.pq.specializations.alt}
\end{subeqnarray}
The resulting S-fraction then has the coefficients
\reff{def.weights.thm.perm.pq.Stype.BIG.1}
but with $(x,u,p_-,q_-) \leftrightarrow (y,v,p_+,q_+)$.
This is also a consequence of the bijection
$\sigma \mapsto R \circ \sigma \circ R$ with $R(i) = n+1-i$,
which interchanges cycle peaks with cycle valleys,
cycle double rises with cycle double falls,
records with antirecords, and upper with lower.
\myendremark

\subsection{First master J-fraction}
   \label{subsec.intro.permutations.firstmaster}

Let us now return to the J-fraction of Theorem~\ref{thm.perm.pq.Jtype.BIG},
with its 17 indeterminates (or 16 if we exclude the redundant variable $s$)
plus the infinite collection $\bw$.
It turns out that this is by no means the end of the story.
Indeed, we can go much farther, and obtain a polynomial in
five infinite families of indeterminates
$\bsfa = (\sfa_{\ell,\ell'})_{\ell,\ell' \ge 0}$,
$\bsfb = (\sfb_{\ell,\ell'})_{\ell,\ell' \ge 0}$,
$\bsfc = (\sfc_{\ell,\ell'})_{\ell,\ell' \ge 0}$,
$\bsfd = (\sfd_{\ell,\ell'})_{\ell,\ell' \ge 0}$,
$\bsfe = (\sfe_\ell)_{\ell \ge 0}$
that will have a nice J-fraction
and that will include 
the polynomial \reff{def.Qn.BIG}
as a specialization.\footnote{
   In our original version of this master J-fraction,
   the weights $\bsfa,\bsfb,\bsfc,\bsfd$
   were factorized in the form
   $\sfa_{\ell,\ell'} = \sfa^{(1)}_\ell \sfa^{(2)}_{\ell'}$, etc.
   We thank Andrew Elvey Price for suggesting the generalization
   in which this factorization is avoided.
}
The basic idea is that,
rather than counting the {\em total}\/ numbers of quadruplets
$i < j < k < l$ that form upper (resp.\ lower) crossings or nestings,
we should instead count the number of upper (resp.\ lower) crossings or nestings
that use a particular vertex $j$ (resp.\ $k$)
in second (resp.\ third) position,
and then attribute weights to the vertex $j$ (resp.\ $k$)
depending on those values.

More precisely, we define
\begin{subeqnarray}
   \ucross(j,\sigma)
   & = &
   \#\{ i<j<k<l \colon\: k = \sigma(i) \hbox{ and } l = \sigma(j) \}
         \\[2mm]
   \unest(j,\sigma)
   & = &
   \#\{ i<j<k<l \colon\: k = \sigma(j) \hbox{ and } l = \sigma(i) \}
      \slabel{def.unestjk}     \\[2mm]
   \lcross(k,\sigma)
   & = &
   \#\{ i<j<k<l \colon\: i = \sigma(k) \hbox{ and } j = \sigma(l) \}
         \\[2mm]
   \lnest(k,\sigma)
   & = &
   \#\{ i<j<k<l \colon\: i = \sigma(l) \hbox{ and } j = \sigma(k) \}
 \label{def.ucrossnestjk}
\end{subeqnarray}
Note that $\ucross(j,\sigma)$ and $\unest(j,\sigma)$ can be nonzero
only when $j$ is a cycle valley or a cycle double rise,
while $\lcross(k,\sigma)$ and $\lnest(k,\sigma)$ can be nonzero
only when $k$ is a cycle peak or a cycle double fall.
And we obviously have
\be
   \ucrosscval(\sigma)
   \;=\;
   \sum\limits_{j \in {\rm cval}}  \ucross(j,\sigma)
 \label{eq.ucrosscval.sum}
\ee
and analogously for the other seven crossing/nesting quantities
defined in Section~\ref{subsec.intro.permutations.pq}.
Recall, finally, the definition \reff{def.level}
of the level of a fixed point $j$:
\be
   \lev(j,\sigma)
   \;=\;
   \#\{ i<j<l \colon\: l = \sigma(i) \}
   \;=\;
   \#\{ i<j<l \colon\: i = \sigma(l) \}
   \;.
 \label{def.level.bis}
\ee

We now introduce five infinite families of indeterminates
$\bsfa = (\sfa_{\ell,\ell'})_{\ell,\ell' \ge 0}$,
$\bsfb = (\sfb_{\ell,\ell'})_{\ell,\ell' \ge 0}$,
$\bsfc = (\sfc_{\ell,\ell'})_{\ell,\ell' \ge 0}$,
$\bsfd = (\sfd_{\ell,\ell'})_{\ell,\ell' \ge 0}$,
$\bsfe = (\sfe_\ell)_{\ell \ge 0}$
and define the polynomial $Q_n(\bsfa,\bsfb,\bsfc,\bsfd,\bsfe)$ by
\begin{eqnarray}
   & & \hspace*{-10mm}
   Q_n(\bsfa,\bsfb,\bsfc,\bsfd,\bsfe)
   \;=\;
       \nonumber \\[4mm]
   & &
   \sum_{\sigma \in \Sym_n}
   \;\:
   \prod\limits_{i \in {\rm cval}}  \! \sfa_{\ucross(i,\sigma),\,\unest(i,\sigma)}
   \prod\limits_{i \in {\rm cpeak}} \!\!  \sfb_{\lcross(i,\sigma),\,\lnest(i,\sigma)}
       \:\times
       \qquad\qquad
       \nonumber \\[1mm]
   & & \qquad\;
   \prod\limits_{i \in {\rm cdfall}} \!\!  \sfc_{\lcross(i,\sigma),\,\lnest(i,\sigma)}
   \;
   \prod\limits_{i \in {\rm cdrise}} \!\!  \sfd_{\ucross(i,\sigma),\,\unest(i,\sigma)}
   \, \prod\limits_{i \in {\rm fix}} \sfe_{\lev(i,\sigma)}
   \;.
   \quad
 \label{def.Qn.firstmaster}
\end{eqnarray}
These polynomials then have a beautiful J-fraction:

\begin{theorem}[First master J-fraction for permutations]
   \label{thm.permutations.Jtype.final1}
The ordinary generating function of the polynomials
$Q_n(\bsfa,\bsfb,\bsfc,\bsfd,\bsfe)$
has the J-type continued fraction
\begin{eqnarray}
   & & \hspace*{-8mm}
   \sum_{n=0}^\infty Q_n(\bsfa,\bsfb,\bsfc,\bsfd,\bsfe) \: t^n
   \;=\;
       \nonumber \\
   & & \hspace*{-4mm}
\Scale[0.8]{
   \cfrac{1}{1 - \sfe_0 t - \cfrac{\sfa_{00} \sfb_{00} t^2}{1 -  (\sfc_{00} + \sfd_{00} + \sfe_1) t - \cfrac{(\sfa_{01} + \sfa_{10})(\sfb_{01} + \sfb_{10}) t^2}{1 - (\sfc_{01} + \sfc_{10} + \sfd_{01} + \sfd_{10} + \sfe_2)t - \cfrac{(\sfa_{02} + \sfa_{11} + \sfa_{20})(\sfb_{02} + \sfb_{11} + \sfb_{20}) t^2}{1 - \cdots}}}}
}
       \nonumber \\[1mm]
   \label{eq.thm.permutations.Jtype.final1}
\end{eqnarray}
with coefficients
\begin{subeqnarray}
   \gamma_n  & = &   \sfc^\star_{n-1} \,+\, \sfd^\star_{n-1} \,+\, \sfe_n
          \\[1mm]
   \beta_n   & = &   \sfa^\star_{n-1} \, \sfb^\star_{n-1}
 \label{def.weights.permutations.Jtype.final1}
\end{subeqnarray}
where
\be
   \sfa^\star_{n-1}  \;\eqdef\;  \sum_{\ell=0}^{n-1} \sfa_{\ell,n-1-\ell}
 \label{def.astar}
\ee
and likewise for $\sfb,\sfc,\sfd$.
\end{theorem}

\noindent
We will prove this theorem in Section~\ref{subsec.permutations.J}.
It is our ``master theorem'' for permutations,
from which most of the others
(namely, all those not including the statistic cyc)
can be derived.

\bigskip

{\bf Remark.}  It seems far from obvious (at least to us)
why $Q_n(\bsfa,\bsfb,\bsfc,\bsfd,\bsfe)$ depends on
$\bsfa,\bsfb,\bsfc,\bsfd,\bsfe$ only via the
combinations (\ref{def.weights.permutations.Jtype.final1}a,b).
Even some partial cases of this --- e.g.\ the fact that the
dependence on $\bsfc$ and $\bsfd$ is only via their sum $\bsfc+\bsfd$ ---
seem nontrivial.
Of course, this behavior is a consequence of
the bijection onto labeled Motzkin paths
that we will use in Section~\ref{subsec.permutations.J}
to prove Theorem~\ref{thm.permutations.Jtype.final1}.
But it would be interesting to understand it combinatorially,
directly at the level of permutations.
\myendremark

\medskip

Let us now show how to recover
$Q_n(x_1,x_2,y_1,y_2,u_1,u_2,v_1,v_2, \bw, p_{+1},p_{+2},p_{-1},p_{-2},${} $q_{+1},q_{+2},q_{-1},q_{-2},s)$
as a specialization of
$Q_n(\bsfa,\bsfb,\bsfc,\bsfd,\bsfe)$,
and thereby obtain Theorem~\ref{thm.perm.pq.Jtype.BIG}
as a special case of Theorem~\ref{thm.permutations.Jtype.final1}.
We need a simple lemma:

\begin{lemma}[Records and antirecords in terms of nestings]
   \label{lemma.rec.nest}
Let $\sigma \in \Sym_n$ and $i \in [n]$.
\begin{itemize}
   \item[(a)]  If $i$ is a cycle valley or cycle double rise,
       then $i$ is a record if and only if $\unest(i,\sigma) = 0$;
       and in this case it is an exclusive record.
   \item[(b)]  If $i$ is a cycle peak or cycle double fall,
       then $i$ is an antirecord if and only if $\lnest(i,\sigma) = 0$;
       and in this case it is an exclusive antirecord.
\end{itemize}
\end{lemma}

\proof
(a)  By hypothesis we have $\sigma(i) > i$.
Then $i$ fails to be a record if and only if there exists $j < i$
such that $\sigma(j) > \sigma(i)$;
and by \reff{def.unestjk} this is exactly the assertion that
$\unest(i,\sigma) > 0$.
The final statement follows from the fact that every record-antirecord
is a fixed point.

(b) is similar.
\qed

In view of \reff{eq.redundant.s}, \reff{eq.ucrosscval.sum}
and Lemma~\ref{lemma.rec.nest},
the specialization needed for obtaining \reff{def.Qn.BIG}
from \reff{def.Qn.firstmaster} is
\begin{subeqnarray}
   \sfa_{\ell,\ell'}
   & = &
   p_{+1}^\ell q_{+1}^{\ell'} \,\times\,
   \begin{cases}
       y_1    &  \textrm{if $\ell' = 0$}  \\
       v_1    &  \textrm{if $\ell' \ge 1$}
   \end{cases}
       \\[2mm]
   \sfb_{\ell,\ell'}
   & = &
   p_{-1}^\ell q_{-1}^{\ell'} \,\times\,
   \begin{cases}
       x_1    &  \textrm{if $\ell' = 0$}  \\
       u_1    &  \textrm{if $\ell' \ge 1$}
   \end{cases}
       \\[2mm]
   \sfc_{\ell,\ell'}
   & = &
   p_{-2}^\ell q_{-2}^{\ell'} \,\times\,
   \begin{cases}
       x_2    &  \textrm{if $\ell' = 0$}  \\
       u_2    &  \textrm{if $\ell' \ge 1$}
   \end{cases}
       \\[2mm]
   \sfd_{\ell,\ell'}
   & = &
   p_{+2}^\ell q_{+2}^{\ell'} \,\times\,
   \begin{cases}
       y_2    &  \textrm{if $\ell' = 0$}  \\
       v_2    &  \textrm{if $\ell' \ge 1$}
   \end{cases}
       \\[2mm]
   \sfe_\ell  & = &  s^\ell w_\ell
 \label{eq.Qn.BIG.specializations}
\end{subeqnarray}
We then have
\be
   \sfa^\star_{n-1}
   \;=\;
   p_{+1}^{n-1} y_1 \:+\: q_{+1} \, [n-1]_{p_{+1},q_{+1}} v_1
\ee
and similarly for $\sfb^\star_{n-1}$, $\sfc^\star_{n-1}$, $\sfd^\star_{n-1}$,
so that we obtain the weights \reff{def.weights.perm.pq.Jtype.BIG}
as a specialization of \reff{def.weights.permutations.Jtype.final1}.
This shows that Theorem~\ref{thm.perm.pq.Jtype.BIG}
is a special case of Theorem~\ref{thm.permutations.Jtype.final1}.

\bigskip

{\bf Remark.}
After discovering and proving Theorem~\ref{thm.permutations.Jtype.final1},
we realized that Flajolet had, in a certain sense, anticipated these ideas
already in his 1980 paper \cite{Flajolet_80}!
In \cite[Proposition~7A]{Flajolet_80} Flajolet gives the J-fraction
for labeled Motzkin paths in which distinct weights are assigned to steps
according to their type (rise, fall or level step),
starting height $h_{i-1}$ and label $\xi_i$
(Theorem~\ref{thm.flajolet_master_labeled_Motzkin} below).
So this is a general ``master J-fraction'' for labeled Motzkin paths,
from which one need only pull back via a bijection
to obtain master J-fractions for specific combinatorial objects
(such as permutations).
As will be seen in Section~\ref{subsec.permutations.J},
our method for proving Theorem~\ref{thm.permutations.Jtype.final1}
is precisely this.
\myendremark

\subsection{First master S-fraction}
   \label{subsec.intro.permutations.firstmaster.S}

We can also obtain a master S-fraction by specializing the parameters
in Theorem~\ref{thm.permutations.Jtype.final1}
and then applying the contraction formula \reff{eq.contraction_even.coeffs}.
There are two possibilities, depending on whether we want
the S-fraction to be ``$\bsfa$ first'' or ``$\bsfb$ first''.
Let us begin by showing the latter, as it meshes better
with the specializations \reff{eq.Qn.BIG.specializations}
and \reff{eq.permutations.S-fraction.pq.specializations}.

\bigskip

{\bf $\bsfb$ first.}
The J-fraction
\reff{eq.thm.permutations.Jtype.final1}/\reff{def.weights.permutations.Jtype.final1}
is the contraction \reff{eq.contraction_even.coeffs} of the S-fraction
\be
   \cfrac{1}{1 - \cfrac{\sfb_{00} t}{1 - \cfrac{\sfa_{00} t}{1 - \cfrac{(\sfb_{01} + \sfb_{10}) t}{1 - \cfrac{(\sfa_{01} + \sfa_{10}) t}{1 - \cdots}}}}}
\ee
with coefficients $\alpha_{2k-1} = \sfb^\star_{k-1}$
and $\alpha_{2k} = \sfa^\star_{k-1}$
if we choose $\bsfa,\bsfb,\bsfc,\bsfd,\bsfe$ so that
\be
   \sfc^\star_{n-1} \,+\, \sfd^\star_{n-1} \,+\, \sfe_n
   \;=\;
   \sfa^\star_{n-1} \,+\, \sfb^\star_n
   \quad\hbox{for all $n \ge 0$}
 \label{eq.permutations.firstmaster.S.specializations}
\ee
[where of course $\sfc^\star_{-1} = \sfd^\star_{-1} = 0$
 by the definition \reff{def.astar}].
Therefore:

\begin{theorem}[First master S-fraction for permutations]
   \label{thm.permutations.Stype.final1}
In the ring $\Z[\bsfa,\bsfb,\bsfc,\bsfd,\bsfe]$,
let $\scri$ be the ideal generated by the
relations \reff{eq.permutations.firstmaster.S.specializations}
for all $n \ge 0$.
Then
\be
   \sum_{n=0}^\infty Q_n(\bsfa,\bsfb,\bsfc,\bsfd,\bsfe) \: t^n
   \;=\;
   \cfrac{1}{1 - \cfrac{\sfb_{00} t}{1 - \cfrac{\sfa_{00} t}{1 - \cfrac{(\sfb_{01} + \sfb_{10}) t}{1 - \cfrac{(\sfa_{01} + \sfa_{10}) t}{1 - \cdots}}}}}
\ee
with coefficients
\begin{subeqnarray}
   \alpha_{2k-1} & = &  \sfb^\star_{k-1} \\[1mm]
   \alpha_{2k}   & = &  \sfa^\star_{k-1}
\end{subeqnarray}
as an identity in $\Z[\bsfa,\bsfb,\bsfc,\bsfd,\bsfe]/\scri$.
\end{theorem}

In applications we will specialize $\bsfa,\bsfb,\bsfc,\bsfd,\bsfe$
in such a way that the relations
\reff{eq.permutations.firstmaster.S.specializations} hold.
There are many ways of doing this;
for instance, we could arrange to have
\begin{subeqnarray}
   \sfd^\star_{n-1}  & = &  \sfa^\star_{n-1}
                                \quad\hbox{for all $n \ge 1$}  \\[2mm]
   \sfc^\star_{n-1} \,+\, \sfe_n
                     & = &  \sfb^\star_n
                                \qquad\hbox{for all $n \ge 0$}
 \slabel{eq.permutations.firstmaster.S.specializations.bis.b}
 \label{eq.permutations.firstmaster.S.specializations.bis}
\end{subeqnarray}
(or alternatively the same thing with $\sfc$ and $\sfd$ interchanged).
And this in turn can be done in many ways;
the simplest way to have $\sfd^\star_{n-1} = \sfa^\star_{n-1}$
for all $n$ is to have $\bsfd = \bsfa$
(that is, $\sfd_{\ell,\ell'} = \sfa_{\ell,\ell'}$ for all $\ell,\ell'$).

In particular, we can obtain Theorem~\ref{thm.perm.pq.Stype.BIG.1}
in this way, by making the specializations
\reff{eq.Qn.BIG.specializations} and
\reff{eq.permutations.S-fraction.pq.specializations}:
this leads to
\begin{subeqnarray}
   \sfa_{\ell,\ell'}  \;=\; \sfd_{\ell,\ell'}
   & = &
   p_{+}^\ell q_{+}^{\ell'} \,\times\,
   \begin{cases}
       y    &  \textrm{if $\ell' = 0$}  \\
       v    &  \textrm{if $\ell' \ge 1$}
   \end{cases}
       \\[2mm]
   \sfb_{\ell,\ell'}
   & = &
   p_{-}^\ell q_{-}^{\ell'} \,\times\,
   \begin{cases}
       x    &  \textrm{if $\ell' = 0$}  \\
       u    &  \textrm{if $\ell' \ge 1$}
   \end{cases}
       \\[2mm]
   \sfc_{\ell,\ell'}
   & = &
   p_{-}^\ell q_{-}^{\ell'} \,\times\,
   \begin{cases}
       p_- x    &  \textrm{if $\ell' = 0$}  \\
       p_- u    &  \textrm{if $\ell' \ge 1$}
   \end{cases}
       \\[2mm]
   \sfe_\ell
   & = &
   q_{-}^{\ell} \,\times\,
   \begin{cases}
       x    &  \textrm{if $\ell' = 0$}  \\
       u    &  \textrm{if $\ell' \ge 1$}
   \end{cases}
 \label{eq.Qn.BIG.specializations.S}
\end{subeqnarray}
and hence
\begin{subeqnarray}
   \sfa^\star_{n-1}  \;=\;  \sfd^\star_{n-1}
   & = &
   p_+^{n-1} y + q_+ \, [n-1]_{p_+,q_+} v
       \\[2mm]
   \sfb^\star_{n-1}
   & = &
   p_-^{n-1} x + q_- \, [n-1]_{p_-,q_-} u
       \\[2mm]
   \sfc^\star_{n-1}
   & = &
   p_-^{n} x + p_- q_- \, [n-1]_{p_-,q_-} u
\end{subeqnarray}
The equation \reff{eq.permutations.firstmaster.S.specializations.bis.b}
is then satisfied by virtue of \reff{eq.recurrence.npq}.
So Theorem~\ref{thm.perm.pq.Stype.BIG.1}
is a special case of Theorem~\ref{thm.permutations.Stype.final1}.
And Theorem~\ref{thm.perms.S}(a) is the further special case
obtained by setting $p_+ = p_- = q_+ = q_- = 1$.

\bigskip

{\bf $\bsfa$ first.}
Completely analogous considerations show that
the J-fraction
\reff{eq.thm.permutations.Jtype.final1}/\reff{def.weights.permutations.Jtype.final1}
is the contraction \reff{eq.contraction_even.coeffs} of the S-fraction
\be
   \cfrac{1}{1 - \cfrac{\sfa_{00} t}{1 - \cfrac{\sfb_{00} t}{1 - \cfrac{(\sfa_{01} + \sfa_{10}) t}{1 - \cfrac{(\sfb_{01} + \sfb_{10}) t}{1 - \cdots}}}}}
\ee
with coefficients $\alpha_{2k-1} = \sfa^\star_{k-1}$
and $\alpha_{2k} = \sfb^\star_{k-1}$
if we choose $\bsfc,\bsfd,\bsfe$ so that
\be
   \sfc^\star_{n-1} \,+\, \sfd^\star_{n-1} \,+\, \sfe_n
   \;=\;
   \sfb^\star_{n-1} \,+\, \sfa^\star_n
   \quad\hbox{for all $n \ge 0$}
   \;.
\ee
This leads to an analogue of Theorem~\ref{thm.permutations.Stype.final1}
in which the roles of $\bsfa$ and $\bsfb$ are interchanged.
And one special case of this
is obtained from \reff{eq.Qn.BIG.specializations} and
\reff{eq.permutations.S-fraction.pq.specializations.alt}.

\subsection[$p,q$-generalizations of the second J-fraction]{$\bm{p,q}$-generalizations of the second J-fraction}
   \label{subsec.intro.permutations.pq.2}

We can also make a $p,q$-generalization of the second J-fraction
involving $\cyc$ (Section~\ref{subsec.perms.J.2}).
Let us define the polynomial
\begin{eqnarray}
   & &
   \widehat{Q}_n(x_1,x_2,y_1,y_2,u_1,u_2,v_1,v_2,\bw,p_{+1},p_{+2},p_{-1},p_{-2},q_{+1},q_{+2},q_{-1},q_{-2},s,\lambda)
   \;=\;
       \nonumber \\[4mm]
   & & \qquad
   \sum_{\sigma \in \Sym_n}
   x_1^{\eareccpeak(\sigma)} x_2^{\eareccdfall(\sigma)} 
   y_1^{\ereccval(\sigma)} y_2^{\ereccdrise(\sigma)} 
   \:\times
       \qquad\qquad
       \nonumber \\[-1mm]
   & & \qquad\qquad\:
   u_1^{\nrcpeak(\sigma)} u_2^{\nrcdfall(\sigma)}
   v_1^{\nrcval(\sigma)} v_2^{\nrcdrise(\sigma)}
   \bw^{\bfix(\sigma)}  \:\times
       \qquad\qquad
       \nonumber \\[3mm]
   & & \qquad\qquad\:
   p_{+1}^{\ucrosscval(\sigma)}
   p_{+2}^{\ucrosscdrise(\sigma)}
   p_{-1}^{\lcrosscpeak(\sigma)}
   p_{-2}^{\lcrosscdfall(\sigma)}
          \:\times
       \qquad\qquad
       \nonumber \\[3mm]
   & & \qquad\qquad\:
   q_{+1}^{\unestcval(\sigma)}
   q_{+2}^{\unestcdrise(\sigma)}
   q_{-1}^{\lnestcpeak(\sigma)}
   q_{-2}^{\lnestcdfall(\sigma)}
   s^{\psnest(\sigma)}
   \lambda^{\cyc(\sigma)}
   \;,
 \label{def.Qn.BIG.cyc}
\end{eqnarray}
which extends \reff{def.Qn.BIG}
by including the factor $\lambda^{\cyc(\sigma)}$.
We refrain from attempting to generalize Conjecture~\ref{thm.perm.Jtype.v2},
and simply limit ourselves to stating the $p,q$-generalization of
Theorem~\ref{thm.perm.Jtype.v2.weaker0} that we are able to prove.
It turns out that we need to make the specializations
$v_1 = y_1$, $v_2 = y_2$, $q_{+1} = p_{+1}$ and $q_{+2} = p_{+2}$.
The result is therefore the following:

\begin{theorem}[Second J-fraction for permutations, $p,q$-generalization]
   \label{thm.perm.Jtype.v2.weaker0.pq}
The ordinary generating function of the polynomials $\widehat{Q}_n$
specialized to
$v_1 = y_1$, $v_2 = y_2$, $q_{+1} = p_{+1}$, $q_{+2} = p_{+2}$
has the J-type continued fraction
\begin{eqnarray}
   & & \hspace*{-7mm}
   \sum_{n=0}^\infty
   \widehat{Q}_n(x_1,x_2,y_1,y_2,u_1,u_2,y_1,y_2,\bw,p_{+1},p_{+2},p_{-1},p_{-2},p_{+1},p_{+2},q_{-1},q_{-2},s,\lambda)
       \: t^n
   \;=\;
       \nonumber \\
   & & \hspace*{-4mm}
\Scale[0.73]{
   \cfrac{1}{1 - \lambda w_0 t - \cfrac{\lambda x_1 y_1 t^2}{1 -  (x_2\!+\!y_2\!+\!\lambda s w_1) t - \cfrac{(\lambda\!+\!1) (p_{-1} x_1 \!+\! q_{-1} u_1) p_{+1} y_1 t^2}{1 - (p_{-2} x_2\!+\!2 p_{+2} y_2\!+\! q_{-2} u_2 \!+\! \lambda s^2 w_2)t - \cfrac{(\lambda\!+\!2) (p_{-1}^2 x_1 \!+\! q_{-1} [2]_{p_{-1},q_{-1}} u_1) p_{+1}^2 y_1  t^2}{1 - \cdots}}}}
}
       \nonumber \\[1mm]
   \label{eq.thm.perm.Jtype.v2.weaker0.pq}
\end{eqnarray}
with coefficients
\begin{subeqnarray}
   \gamma_0  & = &   \lambda w_0   \\[1mm]
   \gamma_n  & = &   (p_{-2}^{n-1} x_2 + q_{-2} \, [n-1]_{p_{-2},q_{-2}} u_2)
                 \:+\: n p_{+2}^{n-1} y_2  \:+\: \lambda s^n w_n
        \quad\hbox{for $n \ge 1$}
       \qquad \\
   \beta_n   & = &   (\lambda + n-1) \:
                     (p_{-1}^{n-1} x_1 + q_{-1} \, [n-1]_{p_{-1},q_{-1}} u_1)
                  \: p_{+1}^{n-1} y_1
 \label{def.weights.perm.Jtype.v2.weaker0.pq}
\end{subeqnarray}
\end{theorem}

\noindent
We will prove this theorem in Section~\ref{subsec.permutations.J.v2},
as a special case of a more general result.

\medskip

{\bf Remarks.}
1. This J-fraction is not invariant under the reversal $R(i) = n+1-i$,
because we have made specializations affecting cval and cdrise
(namely, $v_1 = y_1$, $v_2 = y_2$, $q_{+1} = p_{+1}$, $q_{+2} = p_{+2}$)
but have not made the analogous specializations for cpeak and cdfall
(that is, $u_1 = x_1$, $u_2 = x_2$, $q_{-1} = p_{-1}$, $q_{-2} = p_{-2}$).
However, if we do make also the latter specializations,
then the coefficients \reff{def.weights.perm.Jtype.v2.weaker0.pq} simplify to
\begin{subeqnarray}
   \gamma_n  & = &   n p_{-2}^{n-1} x_2
                 \:+\: n p_{+2}^{n-1} y_2  \:+\: \lambda s^n w_n
       \\[2mm]
   \beta_n   & = &   n \, (\lambda + n-1) \:
                     p_{-1}^{n-1} p_{+1}^{n-1} \, x_1 y_1
\end{subeqnarray}
which are now invariant under the simultaneous interchange
$(x_1,x_2,p_{-1},p_{-2}) \leftrightarrow (y_1,y_2,p_{+1},p_{+2})$.

2.  It is curious to observe that the second J-fraction
\reff{eq.thm.perm.Jtype.v2.weaker0.pq}/\reff{def.weights.perm.Jtype.v2.weaker0.pq}
can be obtained as a specialization of the first J-fraction
\reff{eq.thm.perm.pq.Jtype.BIG}/\reff{def.weights.perm.pq.Jtype.BIG}:
namely,
we specialize $v_1 = y_1$, $v_2 = y_2$, $q_{+1} = p_{+1}$, $q_{+2} = p_{+2}$
and then make the substitutions
$y_1 \leftarrow \lambda y_1$, $v_1 \leftarrow y_1$,
$w_\ell \leftarrow \lambda w_\ell$.
It would be interesting to understand this identity
directly at the level of the polynomials
$Q_n$ [cf.\ \reff{def.Qn.BIG}] and $\widehat{Q}_n$ [cf.\ \reff{def.Qn.BIG.cyc}].
\myendremark

\subsection[$p,q$-generalizations of the second S-fraction]{$\bm{p,q}$-generalizations of the second S-fraction}
   \label{subsec.intro.permutations.pq.2.S}

It turns out that suitable specializations of the J-fraction
\reff{eq.thm.perm.Jtype.v2.weaker0.pq}/\reff{def.weights.perm.Jtype.v2.weaker0.pq}
can be obtained as the contraction \reff{eq.contraction_even.coeffs}
of an S-fraction with polynomial coefficients.
It suffices to make the specializations
\begin{subeqnarray}
  x_1 \;=\; x_2  & = &  x  \\
  y_1  & = &  y  \\
  y_2  & = &  p_+ y  \\
  u_1 \;=\; u_2  & = &  u  \\
  w_\ell & = & y \quad\hbox{for all $\ell \ge 0$}  \\
  p_{+1} \;=\; p_{+2}  & = &  p_+  \\
  p_{-1} \;=\; p_{-2}  & = &  p_-  \\
  q_{-1} \;=\; q_{-2}  & = &  q_-  \\
  s  & = &  p_+
 \label{eq.permutations.secondJ-fraction.pq.specializations}
\end{subeqnarray}
The polynomial $\widehat{Q}_n$ defined in \reff{def.Qn.BIG.cyc} then reduces to
the seven-variable polynomial
\begin{eqnarray}
   & &
   \widehat{P}_n(x,y,u,p_+,p_-,q_-,\lambda)
   \;=\;
   \sum_{\sigma \in \Sym_n}
       x^{\earec(\sigma)} y^{\wex(\sigma)}
       u^{n - \earec(\sigma) - \wex(\sigma)}
       \:\times
       \nonumber \\[2mm]
   & & \qquad\qquad
      p_+^{\ucross(\sigma) + \unest(\sigma) + \ereccdrise(\sigma) + \psnest(\sigma)}
      p_-^{\lcross(\sigma)}
      q_-^{\lnest(\sigma)}
      \lambda^{\cyc(\sigma)}
      \;.
            \nonumber \\
 \label{def.Pn1.BIG.0.cyc}
\end{eqnarray}
It can then be checked that the J-fraction
\reff{eq.thm.perm.Jtype.v2.weaker0.pq}/%
\reff{def.weights.perm.Jtype.v2.weaker0.pq}
with the specializations
\reff{eq.permutations.secondJ-fraction.pq.specializations}
is the contraction \reff{eq.contraction_even.coeffs}
of the following S-fraction:

\begin{theorem}[Second S-fraction for permutations, $p,q$-generalization]
   \label{thm.perm.pq.Stype.BIG.1.cyc}
The ordinary generating function of the polynomials $\widehat{P}_n$
defined by \reff{def.Pn1.BIG.0.cyc} has the S-type continued fraction
\begin{eqnarray}
   & & \hspace*{-7mm}
   \sum_{n=0}^\infty
   \widehat{P}_n(x,y,u,p_+,p_-,q_-,\lambda) \: t^n
   \;=\;
       \nonumber \\[-1mm]
   & & \qquad\qquad
   \cfrac{1}{1 - \cfrac{\lambda y t}{1 - \cfrac{xt}{1 - \cfrac{(\lambda+1) \, p_+ y t}{1- \cfrac{(p_- x + q_- u)t}{1 - \cfrac{(\lambda+2) \, p_+^2 y t}{1 - \cfrac{(p_-^2 x+ q_- \, [2]_{p_-,q_-} u)t}{1-\cdots}}}}}}}
   \label{eq.thm.perm.pq.Stype.BIG.1.cyc}
\end{eqnarray}
with coefficients
\begin{subeqnarray}
   \alpha_{2k-1}  & = &   (\lambda+k-1) \, p_+^{k-1} y  \\[1mm]
   \alpha_{2k}    & = &   p_-^{k-1} x \:+\: q_- \, [k-1]_{p_-,q_-} u
 \label{def.weights.thm.perm.pq.Stype.BIG.1.cyc}
\end{subeqnarray}
\end{theorem}

%

\subsection{Second master J-fraction}
   \label{subsec.intro.permutations.secondmaster}

As with the first J-fraction, we can go much farther,
and obtain a J-fraction in infinitely many indeterminates.
We again introduce five infinite families of indeterminates:
$\bsfa = (\sfa_{\ell})_{\ell \ge 0}$,
$\bsfb = (\sfb_{\ell,\ell'})_{\ell,\ell' \ge 0}$,
$\bsfc = (\sfc_{\ell,\ell'})_{\ell,\ell' \ge 0}$,
$\bsfd = (\sfd_{\ell,\ell'})_{\ell,\ell' \ge 0}$,
$\bsfe = (\sfe_\ell)_{\ell \ge 0}$;
please note that $\bsfa$ now has one index rather than two.
We then define the polynomial
$\widehat{Q}_n(\bsfa,\bsfb,\bsfc,\bsfd,\bsfe,\lambda)$ by
\begin{eqnarray}
   & & \hspace*{-10mm}
   \widehat{Q}_n(\bsfa,\bsfb,\bsfc,\bsfd,\bsfe,\lambda)
   \;=\;
       \nonumber \\[4mm]
   & &
   \sum_{\sigma \in \Sym_n}
   \;\:
   \lambda^{\cyc(\sigma)} \;
   \prod\limits_{i \in {\rm cval}}  \! \sfa_{\ucross(i,\sigma)+\unest(i,\sigma)}
   \prod\limits_{i \in {\rm cpeak}} \!\!  \sfb_{\lcross(i,\sigma),\,\lnest(i,\sigma)}
       \:\times
       \qquad\qquad
       \nonumber \\[1mm]
   & & \;
   \prod\limits_{i \in {\rm cdfall}} \!\!  \sfc_{\lcross(i,\sigma),\,\lnest(i,\sigma)}
   \,
   \prod\limits_{i \in {\rm cdrise}} \!\!  \sfd_{\ucross(i,\sigma)+\unest(i,\sigma),\,\unest(\sigma^{-1}(i),\sigma)}
   \, \prod\limits_{i \in {\rm fix}} \sfe_{\lev(i,\sigma)}
   \;.
   \qquad
 \label{def.Qn.secondmaster}
\end{eqnarray}
Note that here, in contrast to the first master J-fraction,
$\widehat{Q}_n$ depends on $\ucross(i,\sigma)$ and $\unest(i,\sigma)$
only via their sum
(that is the price we have to pay in order to include the statistic cyc);
and note also the somewhat bizarre appearance of
$\unest(\sigma^{-1}(i),\sigma)$ as the second index on $\sfd$.
These polynomials have a nice J-fraction:

\begin{theorem}[Second master J-fraction for permutations]
   \label{thm.permutations.Jtype.final2}
The ordinary generating function of the polynomials
$\widehat{Q}_n(\bsfa,\bsfb,\bsfc,\bsfd,\bsfe,\lambda)$
has the J-type continued fraction
\begin{eqnarray}
   & & \hspace*{-8mm}
   \sum_{n=0}^\infty \widehat{Q}_n(\bsfa,\bsfb,\bsfc,\bsfd,\bsfe,\lambda) \: t^n
   \;=\;
       \nonumber \\
   & & \hspace*{-4mm}
\Scale[0.85]{
   \cfrac{1}{1 - \lambda\sfe_0 t - \cfrac{\lambda \sfa_0 \sfb_{00} t^2}{1 -  (\sfc_{00} + \sfd_{00} + \lambda\sfe_1) t - \cfrac{(\lambda+1) \sfa_1 (\sfb_{01} + \sfb_{10}) t^2}{1 - (\sfc_{01} + \sfc_{10} + \sfd_{10} + \sfd_{11} + \lambda\sfe_2)t - \cfrac{(\lambda+2) \sfa_2 (\sfb_{02} + \sfb_{11} + \sfb_{20}) t^2}{1 - \cdots}}}}
}
       \nonumber \\[1mm]
   \label{eq.thm.permutations.Jtype.final2}
\end{eqnarray}
with coefficients
\begin{subeqnarray}
   \gamma_n  & = &   \sfc^\star_{n-1} \,+\, \sfd^\natural_{n-1}
                                      \,+\, \lambda\sfe_n
          \\[1mm]
   \beta_n   & = &   (\lambda+n-1) \, \sfa_{n-1} \, \sfb^\star_{n-1}
 \label{def.weights.permutations.Jtype.final2}
\end{subeqnarray}
where
\begin{subeqnarray}
   \sfb^\star_{n-1}  & \eqdef &   \sum_{\ell=0}^{n-1} \sfb_{\ell,n-1-\ell}
        \\[2mm]
   \sfc^\star_{n-1}  & \eqdef &   \sum_{\ell=0}^{n-1} \sfc_{\ell,n-1-\ell}
        \\[2mm]
   \sfd^\natural_{n-1}  & \eqdef &  \sum_{\ell=0}^{n-1} \sfd_{n-1,\ell}
\end{subeqnarray}
\end{theorem}

\noindent
We will prove this theorem in Section~\ref{subsec.permutations.J.v2}.
It is our ``master theorem'' for permutation polynomials
that include the statistic cyc.

Let us now show how to recover
$\widehat{Q}_n(x_1,x_2,y_1,y_2,u_1,u_2,y_1,y_2,\bw,p_{+1},p_{+2},p_{-1},p_{-2},$ $p_{+1},p_{+2},q_{-1},q_{-2},s,\lambda)$
as a specialization of
$\widehat{Q}_n(\bsfa,\bsfb,\bsfc,\bsfd,\bsfe)$,
and thereby obtain Theorem~\ref{thm.perm.Jtype.v2.weaker0.pq}
as a special case of Theorem~\ref{thm.permutations.Jtype.final2}.
In view of \reff{eq.redundant.s}, \reff{eq.ucrosscval.sum}
and Lemma~\ref{lemma.rec.nest},
it suffices to set
\begin{subeqnarray}
   \sfa_{\ell}
   & = &
   p_{+1}^\ell y_1
       \\[2mm]
   \sfb_{\ell,\ell'}
   & = &
   p_{-1}^\ell q_{-1}^{\ell'} \,\times\,
   \begin{cases}
       x_1    &  \textrm{if $\ell' = 0$}  \\
       u_1    &  \textrm{if $\ell' \ge 1$}
   \end{cases}
       \\[2mm]
   \sfc_{\ell,\ell'}
   & = &
   p_{-2}^\ell q_{-2}^{\ell'} \,\times\,
   \begin{cases}
       x_2    &  \textrm{if $\ell' = 0$}  \\
       u_2    &  \textrm{if $\ell' \ge 1$}
   \end{cases}
       \\[2mm]
   \sfd_{\ell,\ell'}
   & = &
   p_{+2}^\ell y_2  \qquad\hbox{(no dependence on $\ell'$)}
       \\[2mm]
   \sfe_\ell  & = &  s^\ell w_\ell
 \label{eq.Qn.BIG.specializations.v2}
\end{subeqnarray}
We then have
\begin{subeqnarray}
   \sfb^\star_{n-1}
   & = &
   p_{-1}^{n-1} x_1 \:+\: q_{-1} \, [n-1]_{p_{-1},q_{-1}} u_1
        \\[2mm]
   \sfc^\star_{n-1}
   & = &
   p_{-2}^{n-1} x_2 \:+\: q_{-2} \, [n-1]_{p_{-2},q_{-2}} u_2
        \\[2mm]
   \sfd^\natural_{n-1}
   & = &
   n p_{+2}^{n-1} y_2
\end{subeqnarray}
so that we obtain the weights \reff{def.weights.perm.Jtype.v2.weaker0.pq}
as a specialization of \reff{def.weights.permutations.Jtype.final2}.
This shows that Theorem~\ref{thm.perm.Jtype.v2.weaker0.pq}
is a special case of Theorem~\ref{thm.permutations.Jtype.final2}.

\subsection{Second master S-fraction}

Analogously to what was done in
Section~\ref{subsec.intro.permutations.firstmaster.S}
to obtain the first master S-fraction,
we can also obtain a second master S-fraction by specializing the parameters
in Theorem~\ref{thm.permutations.Jtype.final2}
and then applying the contraction formula \reff{eq.contraction_even.coeffs}.
There are various ways in which this can be done,
but the one that seems to us most natural goes as follows:
{}From \reff{def.weights.permutations.Jtype.final2}
we see that $\lambda$ is associated with $\bsfe$
and with either $\bsfa$ or $\bsfb$;
but since $\bsfa$ and $\bsfe$ are both singly-indexed,
while $\bsfb$ is doubly-indexed,
it seems most natural to set $\bsfa = \bsfe$
and associate $\lambda$ with $\bsfa$.
We thus take
$\alpha_{2k-1} = (\lambda+k-1) \sfa_{k-1} = (\lambda+k-1) \sfe_{k-1}$
and $\alpha_{2k} = \sfb^\star_{k-1}$;
we must then choose $\bsfa,\bsfb,\bsfc,\bsfd$ so that
\be
   \sfc^\star_{n-1} \,+\, \sfd^\natural_{n-1}
   \;=\;
   \sfb^\star_{n-1} \,+\, n \sfe_n
   \;.
\ee
This too can be done in various ways;
the simplest seems to be to choose $\bsfb = \bsfc$
and $d_{n-1,\ell} = e_n = a_n \: \forall \ell$.
This yields:

\begin{theorem}[Second master S-fraction for permutations]
   \label{thm.permutations.Stype.final2}
The ordinary generating function of the polynomials
$\widehat{Q}_n(\bsfa,\bsfb,\bsfc,\bsfd,\bsfe,\lambda)$
specialized to $\bsfa = \bsfe$, $\bsfb = \bsfc$
and $d_{n-1,\ell} = e_n = a_n \: \forall \ell$
has the S-type continued fraction
\be
   \sum_{n=0}^\infty
   \widehat{Q}_n(\bsfa,\bsfb,\bsfb,\bsfd,\bsfa,\lambda)
       \bigr|_{d_{n-1,\ell} = a_n \: \forall \ell} \: t^n
   \;=\;
   \cfrac{1}{1 - \cfrac{\lambda\sfa_0 t}{1 - \cfrac{\sfb^\star_0 t}{1 - \cfrac{(\lambda+1) \sfa_1 t}{1 - \cfrac{\sfb^\star_1 t}{1 - \cdots}}}}}
\ee
with coefficients
\begin{subeqnarray}
   \alpha_{2k-1} & = &  (\lambda+k-1) \sfa_{k-1} \\[1mm]
   \alpha_{2k}   & = &  \sfb^\star_{k-1}
\end{subeqnarray}
\end{theorem}

If in Theorem~\ref{thm.permutations.Stype.final2}
we then make the specializations \reff{eq.Qn.BIG.specializations.v2}
and \reff{eq.permutations.secondJ-fraction.pq.specializations},
we obtain Theorem~\ref{thm.perm.pq.Stype.BIG.1.cyc}.

\subsection{Counting connected components; indecomposable permutations}  
   \label{subsec.permutations.connected}

%

%
%

Let us now show how to extend our permutation polynomials
to count also the connected components of a permutation.
As a corollary, we will obtain continued fractions for
indecomposable permutations.

A {\em divider}\/ of a permutation $\sigma \in \Sym_n$
is an index $i \in [n]$ such that $\sigma$ maps the interval $[1,i]$
into (hence onto) itself;
equivalently, $\sigma$ maps the complementary interval $[i+1,n]$
into (hence onto) itself.
Clearly, when $n = 0$ (hence $\sigma = \emptyset$) there are no dividers;
when $n \ge 1$, the index $n$ is always a divider,
and there may or may not be others.
A {\em connected component}\/ of $\sigma \in \Sym_n$
\cite[p.~262]{Comtet_74} \cite[A059438]{OEIS}
is a minimal nonempty interval $[i,j] \subseteq [n]$
such that the intervals $[1,i-1]$, $[i,j]$ and $[j+1,n]$
are all mapped by $\sigma$ into (hence onto) themselves.
If $1 \le i_1 < i_2 < \ldots < i_k = n$ are the dividers of $\sigma$,
then $[1,i_1]$, $[i_1 + 1, i_2]$, \ldots, $[i_{k-1} + 1, i_k]$
are its connected components.
So the number of connected components equals the number of dividers;
we write it as $\ccc(\sigma)$.
Thus $\ccc(\emptyset) = 0$;
for $n \ge 1$ we have $1 \le \ccc(\sigma) \le n$,
with $\ccc(\sigma) = n$ if and only if $\sigma$ is the identity permutation.
A permutation $\sigma$ is called {\em indecomposable}\/
(or {\em irreducible}\/ or {\em connected}\/) if $\ccc(\sigma) = 1$
\cite[A003319]{OEIS}
(see also \cite{Comtet_72a,Dumont_88,Cori_09b,Foata_09,Cori_17}).\footnote{
   {\bf Warnings:}
   1) B\'ona \cite[p.~162]{Bona_12} defines a permutation to be
   ``indecomposable'' if there does not exist an index $k \in [n-1]$
   such that $\sigma(i) > \sigma(j)$ (greater!!)\ whenever $i \le k < j$.
   This differs from our definition, but is related to it by the involution
   $\sigma \mapsto R \circ \sigma$ where $R(i) = n+1-i$.

   2) The term ``irreducible'' has also been employed
   \cite{Atkinson_02,Baril_16}
   to denote a {\em different}\/ class of permutations,
   namely, those in which there is no index $i$
   satisfying $\sigma(i+1) - \sigma(i) = 1$.
}


In any of the permutation polynomials studied thus far,
we can insert an additional factor $\zeta^{\ccc(\sigma)}$.
This affects the continued fractions as follows:

\begin{theorem}[Counting connected components in permutations]
   \label{thm.connected}
Consider any of the polynomials
\reff{eq.eulerian.fourvar.arec},
\reff{eq.eulerian.fourvar.cyc},
\reff{def.Qnbis},
\reff{def.Qnhat},
\reff{def.Pn.pq},
\reff{def.Qn.BIG},
\reff{def.Pn1.BIG.0},
\reff{def.Qn.firstmaster},
\reff{def.Qn.BIG.cyc}
or \reff{def.Qn.secondmaster},
and insert an additional factor $\zeta^{\ccc(\sigma)}$.
Then the continued fractions associated to the ordinary generating functions
are modified as follows:
in each S-fraction, multiply $\alpha_1$ by $\zeta$;
in each J-fraction, multiply $\gamma_0$ and $\beta_1$ by $\zeta$.
\end{theorem}

This result has an easy proof in our labeled-Motzkin-paths formalism,
as we shall remark in
Sections~\ref{subsec.permutations.J} and \ref{subsec.permutations.J.v2}.
But it also has a simple ``renewal theory'' explanation, as follows:
Given permutations $\sigma = (\sigma_1,\ldots,\sigma_m) \in \Sym_m$
and $\tau = (\tau_1,\ldots,\tau_n) \in \Sym_n$,
let us define their {\em concatenation}\/ $\sigma|\tau \in \Sym_{m+n}$
as $(\sigma_1,\ldots,\sigma_m, \tau_1+m,\ldots,\tau_n+m)$.
Multiple concatenations are defined in the obvious way.
Then every permutation can be written uniquely as a concatenation
of (zero or more) indecomposable permutations
(namely, $\sigma$ restricted to its connected components,
 with indices relabeled to start at 1).
Now let $P_n$ be any permutation polynomial based on statistics that are
additive under concatenation, and include also a factor $\zeta^{\ccc(\sigma)}$;
and let $P_n^{\rm ind}$ be the corresponding polynomial with the sum
restricted to indecomposable permutations (without the factor $\zeta$).
Now define the ordinary generating functions
\begin{subeqnarray}
   f(t)  & = &  \sum_{n=0}^\infty P_n t^n  \\[2mm]
   g(t)  & = &  \sum_{n=1}^\infty P_n^{\rm ind} t^n
\end{subeqnarray}
Then it is immediate from the foregoing that
\be
   f(t)  \;=\;  {1 \over 1 \,-\, \zeta g(t)}
   \;.
 \label{eq.perms.renewal}
\ee
Moreover, all of the statistics that have been considered here
are indeed additive under concatenation:
it is easy to see that this holds for statistics
based on the cycle structure, on the record structure,
or on crossings and nestings.
Theorem~\ref{thm.connected} is an immediate consequence.

\bigskip

{\bf Remark.}
Let us observe that, by contrast,
some of the statistics based on the {\em linear}\/ structure
of the permutation (see Section~\ref{subsec.permutations.linear})
are {\em not}\/ additive under concatenation.
For instance, the ascents in $\sigma|\tau$
include those in $\sigma$ and $\tau$
{\em plus one more}\/ at the boundary between $\sigma$ and $\tau$;
so the ascent statistic is not additive.
(It does, however, behave in a simple way under concatenation.)
On the other hand, the descents in $\sigma|\tau$
include only those in $\sigma$ and $\tau$,
so the descent statistic {\em is}\/ additive.
\myendremark

\medskip

These considerations also allow us to deduce continued fractions
for the ordinary generating functions of the polynomials $P_n^{\rm ind}$
associated to indecomposable permutations.
Indeed, it follows immediately from \reff{eq.perms.renewal}
that if $f(t)$ is the ordinary generating function associated to all
permutations (without factors $\zeta^{\ccc(\sigma)}$),
then $g(t) = 1 - 1/f(t)$ is the ordinary generating function associated to
indecomposable permutations.
Thus, if the polynomials $P_n$ (without factors $\zeta^{\ccc(\sigma)}$)
are given by an S-fraction
\be
   \sum_{n=0}^\infty P_n t^n
   \;=\;
   \cfrac{1}{1 - \cfrac{\alpha_1 t}{1 - \cfrac{\alpha_2 t}{1 - \cdots}}}
   \;\;,
 \label{eq.countconn.Pn.S}
\ee
then the polynomials $P_n^{\rm ind}$  are given by an S-fraction
\be
   \sum_{n=1}^\infty P_n^{\rm ind} t^n
   \;=\;
   \cfrac{\alpha_1 t}{1 - \cfrac{\alpha_2 t}{1 - \cfrac{\alpha_3 t}{1 - \cdots}}}
   \;\;.
 \label{eq.countconn.Pnind.S}
\ee
And if the $P_n$ are given by a J-fraction
\be
   \sum_{n=0}^\infty P_n t^n
   \;=\;
   \cfrac{1}{1 - \gamma_0 t - \cfrac{\beta_1 t^2}{1 - \gamma_1 t - \cfrac{\beta_2 t^2}{1 - \cdots}}}
   \;\;,
 \label{eq.countconn.Pn.J}
\ee
then the $P_n^{\rm ind}$  are given by a J-fraction
plus an additive linear term:
\be
   \sum_{n=1}^\infty P_n^{\rm ind} t^n
   \;=\;
   \gamma_0 t \:+\:
   \cfrac{\beta_1 t^2}{1 - \gamma_1 t - \cfrac{\beta_2 t^2}{1 - \gamma_2 t - \cfrac{\beta_3 t^2}{1 - \cdots}}}
   \;\;.
 \label{eq.countconn.Pnind.J}
\ee

\subsection{321-avoiding permutations}
     \label{subsec.permutations.321}

A permutation $\sigma$ is called {\em 321-avoiding}\/
if there do not exist indices $i < j < k$ such that
$\sigma(i) > \sigma(j) > \sigma(k)$.
In other words, $\sigma$ is 321-avoiding
if there does not exist an index $j$
that is neither a record nor an antirecord.
In detail, this means that
every cycle valley or cycle double rise is an exclusive record,
every cycle peak or cycle double fall is an exclusive antirecord,
and every fixed point is a record-antirecord.\footnote{
   The converses are of course true in every permutation,
   as was observed in Section~\ref{subsec.perms.S}.
}
It is well known \cite[Section~4.2]{Bona_12} \cite[item~115]{Stanley_15}
that the number of 321-avoiding permutations of $[n]$
is the Catalan number $C_n$.
By specializing our preceding results to suppress
neither-record-antirecords,
we can deduce continued fractions for the ordinary generating functions
of polynomials enumerating 321-avoiding permutations
with a variety of record and crossing/nesting statistics.

The simplest such result is obtained by setting
$u = v = 0$ in Theorem~\ref{thm.perms.S};
as mentioned already in \reff{eq.narayana},
this leads to the homogenized Narayana polynomials
and their ordinary generating function
\be
   \sum_{n=0}^\infty P_n(x,y,0,0) \, t^n
   \;=\;
   \cfrac{1}{1 - \cfrac{xt}{1 - \cfrac{yt}{1 -  \cfrac{xt}{1- \cfrac{yt}{1- \cdots}}}}}
 \label{eq.cfrac.narayana}
\ee
with coefficients $\alpha_{2k-1} = x$, $\alpha_{2k} = y$.
Specializing further to $x=y=1$ we obtain the well-known
S-fraction for the ordinary generating function of the Catalan numbers:
\be
   \sum_{n=0}^\infty C_n t^n
   \;=\;
   \cfrac{1}{1 - \cfrac{t}{1 - \cfrac{t}{1 -  \cfrac{t}{1- \cdots}}}}
 \label{eq.cfrac.catalan}
\ee
with coefficients $\alpha_n = 1$.

More generally, we can set $u_1 = u_2 = v_1 = v_2 = w = 0$ in \reff{def.Qn};
or even more generally, we can set
$u_1 = u_2 = v_1 = v_2 = 0$ and $w_\ell = 0$ for $\ell \ge 1$
in \reff{def.Qnbis} or \reff{def.Qn.BIG}.
Taking the latter leads to the J-fraction
\begin{eqnarray}
   & & \hspace*{-14mm}
   \sum_{n=0}^\infty
       Q_n(x_1,x_2,y_1,y_2,0,0,0,0,w_0,\bzero,p_{+1},p_{+2},p_{-1},p_{-2},q_{+1},q_{+2},q_{-1},q_{-2},s) \: t^n
   \;\:=\;\:
       \nonumber \\
   & &
   \cfrac{1}{1 - w_0 t - \cfrac{x_1 y_1 t^2}{1 -  (x_2\!+\!y_2) t - \cfrac{p_{-1} p_{+1} x_1 y_1 t^2}{1 - (p_{-2} x_2\!+\!p_{+2} y_2)t - \cfrac{p_{-1}^2 p_{+1}^2 x_1 y_1 t^2}{1 - \cdots}}}}
   \qquad
   \label{eq.thm.perm.pq.Jtype.BIG.321}
\end{eqnarray}
with coefficients
\begin{subeqnarray}
   \gamma_0  & = &   w_0
      \slabel{def.weights.perm.pq.Jtype.BIG.a.321}  \\[1mm]
   \gamma_n  & = &   p_{-2}^{n-1} x_2 \:+\: p_{+2}^{n-1} y_2
        \quad\hbox{for $n \ge 1$}
      \slabel{def.weights.perm.pq.Jtype.BIG.b.321}  \\[1mm]
   \beta_n   & = &   p_{-1}^{n-1} p_{+1}^{n-1} x_1 y_1
      \slabel{def.weights.perm.pq.Jtype.BIG.c.321}
 \label{def.weights.perm.pq.Jtype.BIG.321}
\end{subeqnarray}
Note that $q_{\pm 1}, q_{\pm 2}, s$ do not enter here;
this is because:

\begin{lemma}
   \label{lemma.321-avoiding}
A 321-avoiding permutation cannot have an upper nesting, lower nesting,
upper pseudo-nesting or lower pseudo-nesting.
\end{lemma}

\proof
Suppose that $\sigma \in \Sym_n$ has an upper nesting or pseudo-nesting,
i.e.\ there exist $i,j \in [n]$ such that $i < j \le \sigma(j) < \sigma(i)$.
Now partition $[n] = [1,j] \cup [j+1,n]$.
There is at least one arc of $\sigma$ running from $[1,j]$ to $[j+1,n]$
--- namely, the arc from $i$ to $\sigma(i)$ ---
so there must be at least one arc running in the opposite direction,
i.e.\ there exists $k > j$ such that $\sigma(k) \le j$.
And in case $\sigma(j) = j$, then we must have $\sigma(k) < j$.
So $\sigma(k) < \sigma(j)$, contradicting the hypothesis
that $\sigma$ is 321-avoiding.
An analogous proof works if $\sigma$ has a lower nesting or pseudo-nesting.
\qed

We can also get an S-fraction by setting $u=v=0$ in \reff{def.Pn1.BIG.0};
this leads to
\begin{eqnarray}
   & & \hspace*{-7mm}
   \sum_{n=0}^\infty
       P_n(x,y,0,0,p_+,p_-,q_+,q_-) \: t^n
   \;=\;
       \nonumber \\[-1mm]
   & & \qquad\qquad
   \cfrac{1}{1 - \cfrac{xt}{1 - \cfrac{yt}{1 - \cfrac{p_- x t}{1- \cfrac{p_+ y t}{1 - \cfrac{p_-^2 x t}{1 - \cfrac{p_+^2 y t}{1-\cdots}}}}}}}
   \label{eq.thm.perm.pq.Stype.BIG.1.321}
\end{eqnarray}
with coefficients
\begin{subeqnarray}
   \alpha_{2k-1}  & = &   p_-^{k-1} x   \\[1mm]
   \alpha_{2k}    & = &   p_+^{k-1} y 
 \label{def.weights.thm.perm.pq.Stype.BIG.1.321}
\end{subeqnarray}

\subsection{Cycle-alternating permutations}
     \label{subsec.permutations.cycle-alternating}

A permutation $\sigma$ is called {\em cycle-alternating}\/
if it has no cycle double rises, cycle double falls, or fixed points;
thus, each cycle of $\sigma$ is of even length (call it $2k$)
and consists of $k$ cycle valleys and $k$ cycle peaks in alternation.
Deutsch and Elizalde \cite[Proposition~2.2]{Deutsch_11}
have shown that the number of cycle-alternating permutations of $[2n]$
is the secant number $E_{2n}$ (see also Biane \cite{Biane_93}).
By specializing our preceding results to suppress
cycle double rises, cycle double falls and fixed points,
we can deduce continued fractions for the ordinary generating functions
of polynomials enumerating cycle-alternating permutations
with a variety of record and crossing/nesting statistics.

The simplest such result is obtained by setting
$x_2 = y_2 = u_2 = v_2 = w_\ell = 0$ in Theorem~\ref{thm.perm.Jtype}.
Then all the coefficients $\gamma_n$ in the J-fraction vanish,
so that we obtain an S-fraction in the variable $t^2$
that enumerates cycle-alternating permutations
according to their record statistics.
Changing $t^2$ to $t$, we have:

\begin{theorem}[S-fraction for cycle-alternating permutations]
   \label{thm.perm.cycle-alternating}
The ordinary generating function of the polynomials
\be
   Q_{2n}(x_1,0,y_1,0,u_1,0,v_1,0,\bzero)
   \;=\;
   \sum_{\sigma \in \Sym^{\rm ca}_{2n}}
   x_1^{\eareccpeak(\sigma)}
   y_1^{\ereccval(\sigma)}
   u_1^{\nrcpeak(\sigma)}
   v_1^{\nrcval(\sigma)}
 \label{def.Qnbis.cycle-alternating}
\ee
enumerating cycle-alternating permutations
according to their record statistics
has the S-type continued fraction
\be
   \sum_{n=0}^\infty Q_{2n}(x_1,0,y_1,0,u_1,0,v_1,0,\bzero) \: t^n
   \;=\;
   \cfrac{1}{1 - \cfrac{x_1 y_1 t}{1 -  \cfrac{(x_1\!+\!u_1)(y_1\!+\!v_1) t}{1 - \cfrac{(x_1\!+\!2u_1)(y_1\!+\!2v_1) t}{1 - \cdots}}}}
   \label{eq.thm.perm.Stype.cycle-alternating}
\ee
with coefficients
\be
   \alpha_n   \;=\;   [x_1 + (n-1)u_1] \: [y_1 + (n-1)v_1]
   \;.
 \label{def.weights.perm.Stype.cycle-alternating}
\ee
\end{theorem}

In particular, specializing to $x_1 = y_1 = u_1 = v_1 = 1$,
we obtain the well-known
\cite{Stieltjes_1889,Rogers_07}
\cite[Theorem~3B(iii)]{Flajolet_80}
S-fraction expansion of the ordinary generating function
of the secant numbers:
\be
   \sum_{n=0}^\infty E_{2n} t^n
   \;=\;
   \cfrac{1}{1 - \cfrac{1^2 t}{1 - \cfrac{2^2 t}{1 -  \cfrac{3^2 t}{1- \cdots}}}}
 \label{eq.cfrac.secant}
\ee
with coefficients $\alpha_n = n^2$.

More generally, we can make the corresponding specialization
in Theorem~\ref{thm.perm.pq.Jtype.BIG}
and thereby include crossing and nesting statistics:

\begin{theorem}[S-fraction for cycle-alternating permutations, $p,q$-generalization]
   \label{thm.perm.cycle-alternating.2}
The ordinary generating function of the polynomials
\begin{eqnarray}
   & &
   Q_{2n}(x_1,0,y_1,0,u_1,0,v_1,0,\bzero,p_{+1},0,p_{-1},0,q_{+1},0,q_{-1},0,0)
   \;=\;
       \nonumber \\[4mm]
   & & \qquad\qquad
   \sum_{\sigma \in \Sym^{\rm ca}_{2n}}
   x_1^{\eareccpeak(\sigma)}
   y_1^{\ereccval(\sigma)}
   u_1^{\nrcpeak(\sigma)}
   v_1^{\nrcval(\sigma)}
       \:\times
       \qquad\qquad
       \nonumber \\[0mm]
   & & \qquad\qquad\qquad\;\:
   p_{+1}^{\ucrosscval(\sigma)}
   p_{-1}^{\lcrosscpeak(\sigma)}
   q_{+1}^{\unestcval(\sigma)}
   q_{-1}^{\lnestcpeak(\sigma)}
 \label{def.Qn.BIG.cycle-alternating}
\end{eqnarray}
enumerating cycle-alternating permutations
according to their record and crossing/nesting statistics
has the S-type continued fraction
\begin{eqnarray}
   & & \hspace*{-2cm}
   \sum_{n=0}^\infty
       Q_{2n}(x_1,0,y_1,0,u_1,0,v_1,0,\bzero,p_{+1},0,p_{-1},0,q_{+1},0,q_{-1},0,0) \; t^n
   \;=\;
       \nonumber \\
   & & \!\!\!\!
   \cfrac{1}{1 - \cfrac{x_1 y_1 t}{1 - \cfrac{(p_{-1} x_1\!+\! q_{-1} u_1)(p_{+1} y_1\!+\! q_{+1} v_1) t}{1 - \cfrac{(p_{-1}^2 x_1\!+\! q_{-1} [2]_{p_{-1},q_{-1}} u_1)(p_{+1}^2 y_1\!+\! q_{+1} [2]_{p_{+1},q_{+1}} v_1) t}{1 - \cdots}}}}
   \label{eq.thm.perm.pq.Jtype.BIG.cycle-alternating}
\end{eqnarray}
with coefficients
\be
   \alpha_n   \;=\;   (p_{-1}^{n-1} x_1 + q_{-1} \, [n-1]_{p_{-1},q_{-1}} u_1)
                  \: (p_{+1}^{n-1} y_1 + q_{+1} \, [n-1]_{p_{+1},q_{+1}} v_1)
   \;.
 \label{def.weights.perm.pq.Jtype.BIG.cycle-alternating}
\ee
\end{theorem}

\medskip

Most generally, we can set $\bsfc = \bsfd = \bsfe = \bzero$
in the first master J-fraction for permutations
(Theorem~\ref{thm.permutations.Jtype.final1}) to obtain:

\begin{theorem}[First master S-fraction for cycle-alternating permutations]
   \label{thm.perm.cycle-alternating.3}
The ordinary generating function of the polynomials
\be
   Q_{2n}(\bsfa,\bsfb,\bzero,\bzero,\bzero)
   \;=\;
   \sum_{\sigma \in \Sym^{\rm ca}_{2n}}
   \;\:
   \prod\limits_{i \in {\rm cval}}  \! \sfa_{\ucross(i,\sigma),\,\unest(i,\sigma)}
   \prod\limits_{i \in {\rm cpeak}} \!\! \sfb_{\lcross(i,\sigma),\,\lnest(i,\sigma)}
 \label{def.Qn.firstmaster.cycle-alternating}
\ee
enumerating cycle-alternating permutations has the S-type continued fraction
\be
   \sum_{n=0}^\infty Q_{2n}(\bsfa,\bsfb,\bzero,\bzero,\bzero) \: t^n
   \;=\;
   \cfrac{1}{1 - \cfrac{\sfa_{00} \sfb_{00} t}{1 - \cfrac{(\sfa_{01} + \sfa_{10})(\sfb_{01} + \sfb_{10}) t}{1 - \cfrac{(\sfa_{02} + \sfa_{11} + \sfa_{20})(\sfb_{02} + \sfb_{11} + \sfb_{20}) t}{1 - \cdots}}}}
   \label{eq.thm.permutations.Stype.final1.cycle-alternating}
\ee
with coefficients
\be
   \alpha_n   \;=\;   \sfa^\star_{n-1} \, \sfb^\star_{n-1}
 \label{def.weights.permutations.Stype.final1.cycle-alternating}
\ee
where
$\displaystyle
   \sfa^\star_{n-1}  \;\eqdef\;  \sum_{\ell=0}^{n-1} \sfa_{\ell,n-1-\ell}
$
and likewise for $\sfb$.
\end{theorem}

\bigskip

We can also obtain a second S-fraction for cycle-alternating permutations
by setting $x_2 = y_2 = u_2 = v_2 = w_\ell = 0$
in Theorem~\ref{thm.perm.Jtype.v2.weaker0}:
this allows us to include the counting of cycles ($\lambda$),
but at the expense of ignoring the record status of cycle valleys
($v_1 = y_1$).  We have:

\begin{theorem}[Second S-fraction for cycle-alternating permutations]
   \label{thm.perm.cycle-alternating.second}
The ordinary generating function of the polynomials
\be
   \widehat{Q}_{2n}(x_1,0,y_1,0,u_1,0,y_1,0,\bzero,\lambda)
   \;=\;
   \sum_{\sigma \in \Sym^{\rm ca}_{2n}}
   x_1^{\eareccpeak(\sigma)}
   u_1^{\nrcpeak(\sigma)}
   y_1^{\cval(\sigma)}
   \lambda^{\cyc(\sigma)}
 \label{def.Qnbis.cycle-alternating.second}
\ee
has the S-type continued fraction
\be
   \sum_{n=0}^\infty \widehat{Q}_{2n}(x_1,0,y_1,0,u_1,0,y_1,0,\bzero,\lambda)
      \: t^n
   \;=\;
   \cfrac{1}{1 - \cfrac{\lambda x_1 y_1 t}{1 -  \cfrac{(\lambda+1)(x_1\!+\!u_1) y_1 t}{1 - \cfrac{(\lambda+2)(x_1\!+\!2u_1) y_1 t}{1 - \cdots}}}}
   \label{eq.thm.perm.Stype.cycle-alternating.second}
\ee
with coefficients
\be
   \alpha_n   \;=\;   (\lambda+n-1) \, [x_1 + (n-1)u_1] \, y_1
   \;.
 \label{def.weights.perm.Stype.cycle-alternating.second}
\ee
\end{theorem}

\noindent
Note that here the variable $y_1$ is redundant,
as it can be absorbed into $x_1$ and $u_1$.
This reflects the fact that for cycle-alternating permutations
one has $\cval = \cpeak = \eareccpeak + \nrcpeak$.

More generally, we can make the corresponding specialization
in Theorem~\ref{thm.perm.Jtype.v2.weaker0.pq}
and thereby include crossing and nesting statistics
(again subject to the specialization $v_1 = y_1$):

\begin{theorem}[Second S-fraction for cycle-alternating permutations, $p,q$-generalization]
   \label{thm.perm.cycle-alternating.2.second}
The ordinary generating function of the polynomials
\begin{eqnarray}
   & &
   \widehat{Q}_{2n}(x_1,0,y_1,0,u_1,0,y_1,0,\bzero,p_{+1},0,p_{-1},0,q_{+1},0,q_{-1},0,0,\lambda)
   \;=\;
       \nonumber \\[4mm]
   & & \qquad\qquad
   \sum_{\sigma \in \Sym^{\rm ca}_{2n}}
   x_1^{\eareccpeak(\sigma)}
   u_1^{\nrcpeak(\sigma)}
   y_1^{\cval(\sigma)}
       \:\times
       \qquad\qquad
       \nonumber \\[0mm]
   & & \qquad\qquad\qquad\;\:
   p_{+1}^{\ucrosscval(\sigma)}
   p_{-1}^{\lcrosscpeak(\sigma)}
   q_{+1}^{\unestcval(\sigma)}
   q_{-1}^{\lnestcpeak(\sigma)}
   \lambda^{\cyc(\sigma)}
   \qquad\qquad
 \label{def.Qn.BIG.cycle-alternating.second}
\end{eqnarray}
has the S-type continued fraction
\begin{eqnarray}
   & & \hspace*{-2cm}
   \sum_{n=0}^\infty
   \widehat{Q}_{2n}(x_1,0,y_1,0,u_1,0,y_1,0,\bzero,p_{+1},0,p_{-1},0,q_{+1},0,q_{-1},0,0,\lambda)
\; t^n
   \;=\;
       \nonumber \\
   & & \!\!\!\!
   \cfrac{1}{1 - \cfrac{\lambda x_1 y_1 t}{1 - \cfrac{(\lambda+1)(p_{-1} x_1\!+\! q_{-1} u_1) p_{+1} y_1 t}{1 - \cfrac{(\lambda+2)(p_{-1}^2 x_1\!+\! q_{-1} [2]_{p_{-1},q_{-1}} u_1) p_{+1}^2 y_1 t}{1 - \cdots}}}}
   \label{eq.thm.perm.pq.Jtype.BIG.cycle-alternating.second}
\end{eqnarray}
with coefficients
\be
   \alpha_n   \;=\;   (\lambda+n-1) \,
       (p_{-1}^{n-1} x_1 + q_{-1} \, [n-1]_{p_{-1},q_{-1}} u_1)
                  \: p_{+1}^{n-1} y_1
   \;.
 \label{def.weights.perm.pq.Jtype.BIG.cycle-alternating.second}
\ee
\end{theorem}

\medskip

Most generally, we can set $\bsfc = \bsfd = \bsfe = \bzero$
in the second master J-fraction for permutations
(Theorem~\ref{thm.permutations.Jtype.final2}) to obtain:

\begin{theorem}[Second master S-fraction for cycle-alternating permutations]
   \label{thm.perm.cycle-alternating.3.second}
The ordinary generating function of the polynomials
\be
   \widehat{Q}_{2n}(\bsfa,\bsfb,\bzero,\bzero,\bzero,\lambda)
   \;=\;
   \sum_{\sigma \in \Sym^{\rm ca}_{2n}}
   \;\:
   \lambda^{\cyc(\sigma)} \;
   \prod\limits_{i \in {\rm cval}}  \! \sfa_{\ucross(i,\sigma)+\unest(i,\sigma)}
   \prod\limits_{i \in {\rm cpeak}} \!\! \sfb_{\lcross(i,\sigma),\,\lnest(i,\sigma)}
 \label{def.Qn.firstmaster.cycle-alternating.second}
\ee
has the S-type continued fraction
\be
   \sum_{n=0}^\infty \widehat{Q}_{2n}(\bsfa,\bsfb,\bzero,\bzero,\bzero,\lambda) \: t^n
   \;=\;
   \cfrac{1}{1 - \cfrac{\lambda \sfa_{0} \sfb_{00} t}{1 - \cfrac{(\lambda+1) \sfa_{1} (\sfb_{01} + \sfb_{10}) t}{1 - \cfrac{(\lambda+2) \sfa_{2} (\sfb_{02} + \sfb_{11} + \sfb_{20}) t}{1 - \cdots}}}}
   \label{eq.thm.permutations.Stype.final1.cycle-alternating.second}
\ee
with coefficients
\be
   \alpha_n   \;=\;   (\lambda+n-1) \, \sfa_{n-1} \, \sfb^\star_{n-1}
 \label{def.weights.permutations.Stype.final1.cycle-alternating.second}
\ee
where
$\displaystyle
   \sfb^\star_{n-1}  \;\eqdef\;  \sum_{\ell=0}^{n-1} \sfb_{\ell,n-1-\ell}
$.
\end{theorem}

\medskip

{\bf A final remark.}
In Section~\ref{sec.intro.perfectmatchings} we will enumerate perfect matchings
--- which are a subclass of cycle-alternating permutations ---
with distinct weights for even and odd cycle peaks.
On the other hand, Dumont \cite{Dumont_79,Dumont_81a}
has shown that if one enumerates {\em all}\/ permutations
with distinct weights for even and odd cycle peaks,
one obtains the Schett polynomials,
which are closely related to the Jacobian elliptic functions.
It would be interesting to know whether any of our S-fractions
for cycle-alternating permutations can likewise be refined by
giving distinct weights for even and odd cycle peaks.
\myendremark

\subsection{A remark on the inversion statistic}
   \label{subsec.permutations.inv}

An {\em inversion}\/ of a permutation $\sigma\in \Sym_n$
is a pair $(i,j)\in [n]\times [n]$ such that $i<j$ and $\sigma(i)>\sigma(j)$.  
We write
\be
   \inv(\sigma)
   \;\eqdef\;
   \#\{(i,j)\colon\: i \le j \textrm{ and } \sigma(i)>\sigma(j) \}
 \label{def.inv}
\ee
for the number of inversions in $\sigma$.
Note that $\inv(\sigma) = \inv(\sigma^{-1})$.

De M\'edicis and Viennot \cite[Lemme 3.1]{DeMedicis_94}
gave an expression for $\inv(\sigma)$
that can be translated into our language as follows
(see also \cite[eq.~(40)]{Shin_10}):

\begin{proposition}
   \label{prop.inv}
We have
\begin{subeqnarray}
   \hspace*{-1.5cm}
   \inv
   & \!=\! &
   \exc + (\ucross + 2\,\unest) + (\lcross + \ljoin + 2\,\lnest + 2\,\lpsnest)
      \slabel{eq.prop.inv.a} \\[1mm]
   & \!=\! & \cval + \cdrise + \cdfall + \ucross + \lcross
                    + 2(\unest + \lnest + \psnest)
   \qquad
  \label{eq.prop.inv}
\end{subeqnarray}
\end{proposition}

\noindent
Note that (\ref{eq.prop.inv}b) is invariant under
$\sigma \leftrightarrow \sigma^{-1}$,
since $\cval(\sigma) = \cval(\sigma^{-1})$,
$\cdrise(\sigma) = \cdfall(\sigma^{-1})$,
$\ucross(\sigma) = \lcross(\sigma^{-1})$,
$\unest(\sigma) = \lnest(\sigma^{-1})$
and $\psnest(\sigma) = \psnest(\sigma^{-1})$.

Since the proof of Proposition~\ref{prop.inv} given in \cite{DeMedicis_94}
is rather lengthy, for completeness let us give a short proof.
This proof is based on the following pair of identities
\cite[Lemma~8]{Clarke_97}:

\begin{lemma}
   \label{lemma.prop.inv}
For any permutation $\sigma$, we have
\begin{eqnarray}
   \#\{(i,j)\colon\: i \le j < \sigma(i) \textrm{ and } \sigma(j) > j \}
   & = & 
   \#\{(i,j)\colon\: \sigma(i)<\sigma(j)\leq  i \textrm{ and } \sigma(j) > j\}
       \nonumber \\ \label{clarke2} \\
   \#\{(i,j)\colon\: i \le j < \sigma(i) \textrm{ and } \sigma(j) \le j \}
   & = & 
   \#\{(i,j)\colon\: \sigma(i)<\sigma(j)\leq  i \textrm{ and } \sigma(j) \le  j\}
       \nonumber \\ \label{clarke2bis}
\end{eqnarray}
\end{lemma}

\proofof{Lemma~\ref{lemma.prop.inv}}
The equality~\reff{clarke2} was proven by Clarke~\cite[Lemma~3]{Clarke_95}.
To prove \reff{clarke2bis},
note that $\sum\limits_{i=1}^n \, [\sigma(i) - i] = 0$
and hence
$\sum\limits_{i \colon\, \sigma(i) > i} [\sigma(i) - i]
 =\! \sum\limits_{i \colon\, \sigma(i) < i} [i - \sigma(i)]$,
or in other words
\be
   \#\{(i,j)\colon\: i \le j < \sigma(i) \}
   \;=\;
   \#\{(i,j)\colon\: \sigma(i) < \sigma(j) \le i \}
   \;.
 \label{eq.proof.lemma.prop.inv}
\ee
Subtracting \reff{clarke2} from \reff{eq.proof.lemma.prop.inv}
yields \reff{clarke2bis}.
\qed

\proofof{Proposition~\ref{prop.inv}}
Let us begin by observing that \reff{clarke2} and \reff{clarke2bis}
can be rewritten, by consideration of cases, as
\begin{eqnarray}
   \#\{(i,j)\colon\: i \le j < \sigma(i) \textrm{ and } \sigma(j) > j \}
   & \!=\! &
   \exc+\ucross+\unest
       \label{eq.clarke2.rewrite} \\[2mm]
   \#\{(i,j)\colon\: \sigma(i)<\sigma(j)\leq  i \textrm{ and } \sigma(j) \le  j\}
   & \!=\! &
   \ljoin + \lcross + \lnest + \lpsnest
       \nonumber \\[-1mm] \label{eq.clarke2bis.rewrite}
\end{eqnarray}
Let us now divide the set of inversion pairs 
$\{(i,j) \colon\: i<j \hbox{ and } \sigma(i)>\sigma(j)\}$
into three classes:
\begin{itemize}
   \item[1)]  $\sigma(i) \le i$ [hence $\sigma(j) < \sigma(i) \le i < j$];
   \item[2)]  $\sigma(i) > i$ and $\sigma(j) \ge j$
                 [hence $i < j \le \sigma(j) < \sigma(i)$];
   \item[3)]  $\sigma(i) > i$ and $\sigma(j) < j$.
\end{itemize}
The first class yields $\lnest + \lpsnest$.
The second class yields $\unest + \upsnest$.
Let us divide the third class into two subclasses:
\begin{itemize}
   \item[3a)]  $\sigma(i) > i$ and $\sigma(j) < j$ and $\sigma(i) \le j$
      [hence $\sigma(j) < \sigma(i) \le j$ and $\sigma(i) > i$];
   \item[3b)]  $\sigma(i) > i$ and $\sigma(j) < j$ and $\sigma(i) > j$
      [hence $i < j < \sigma(i)$ and $\sigma(j) < j$].
\end{itemize}
Class (3a) is the right-hand side of \reff{clarke2}
[with $i \leftrightarrow j$], which by \reff{clarke2}/\reff{eq.clarke2.rewrite}
equals ${\exc + \ucross + \unest}$.
On the other hand, the left-hand side of \reff{clarke2bis}
can be rewritten as
   $\#\{(i,j)\colon\: i < j < \sigma(i) \textrm{ and } \sigma(j) \le j \}$
(since $i=j$ contradicts the other inequalities),
which in turn equals class (3b) plus $\upsnest$;
therefore, by \reff{clarke2bis}/\reff{eq.clarke2bis.rewrite},
class (3b) yields ${\ljoin + \lcross + \lnest}$
(since $\upsnest = \lpsnest$).
Combining classes (1), (2), (3a) and (3b)
and using $\upsnest = \lpsnest$ gives (\ref{eq.prop.inv}a).
\qed


See also the end of Section~\ref{subsec.permutations.J.v2}
for an alternate proof of Proposition~\ref{prop.inv}
based on the Biane bijection to labeled Motzkin paths.

Using Proposition~\ref{prop.inv}, results involving $\inv$ can be translated
to the cycle, crossing and nesting statistics,
which in our opinion are more fundamental.
For instance, we see from (\ref{eq.prop.inv}b) that
Zeng's \cite{Zeng_89} S-fraction
for the weights $a^{\arec(\sigma)} b^{\erec(\sigma)} q^{\inv(\sigma)}$
is the special case of Theorem~\ref{thm.perm.pq.Stype.BIG.1} with
\be
   x=a ,\quad y=qb ,\quad u=1 ,\quad v=q ,\quad p_+ = p_- = q
                                         ,\quad q_+ = q_- = q^2
   \;,
\ee
as already remarked in Section~\ref{subsec.intro.permutations.pq.S}.
See also \cite[eqns.~(2.3) and (3.2)]{Zeng_95} for the special case $b=1$.
Similarly, Elizalde's \cite[eqn.~(4)]{Elizalde_18} J-fraction
for the weights
$a^{\cval(\sigma)} b^{\cdrise(\sigma)} w^{\fix(\sigma)} q^{\inv(\sigma)}$
is the special case of Theorem~\ref{thm.perm.pq.Jtype.BIG} with
\begin{eqnarray}
   & &
   x_1 = u_1 = 1,\quad
   x_2 = u_2 = q,\quad
   y_1 = v_1 = qa,\quad
   y_2 = v_2 = qb,\quad
   w_\ell = w\;\forall\ell,\quad
        \nonumber \\
   & & \quad
   p_{+1} = p_{+2} = p_{-1} = p_{-2} = q,\quad
   q_{+1} = q_{+2} = q_{-1} = q_{-2} = q^2,\quad
   s = q^2
   \;.
   \qquad\qquad
\end{eqnarray}
Indeed, we can obtain a more general J-fraction with weights
$a^{\cval(\sigma)} b^{\cdrise(\sigma)} c^{\cpeak(\sigma)}$
$d^{\cdfall(\sigma)} w^{\fix(\sigma)} q^{\inv(\sigma)}$
by specializing Theorem~\ref{thm.perm.pq.Jtype.BIG} to
\begin{eqnarray}
   & &
   x_1 = u_1 = c,\quad
   x_2 = u_2 = qd,\quad
   y_1 = v_1 = qa,\quad
   y_2 = v_2 = qb,\quad
   w_\ell = w\;\forall\ell,\quad
        \nonumber \\
   & & \quad
   p_{+1} = p_{+2} = p_{-1} = p_{-2} = q,\quad
   q_{+1} = q_{+2} = q_{-1} = q_{-2} = q^2,\quad
   s = q^2
   \;,
   \qquad\qquad
 \label{eq.spec.perm.inv}
\end{eqnarray}
yielding continued-fraction coefficients
\begin{subeqnarray}
   \gamma_n  & = &  q^n \, [n]_q \, (b+d) \:+\:  \lambda q^{2n} w
        \\[2mm]
   \beta_n   & = &  q^{2n-1} \, [n]_q^2 \, ac
\end{subeqnarray}
as a specialization of \reff{def.weights.perm.pq.Jtype.BIG.simplified}.

\bigskip

{\bf Remark.}
Note that we are {\em unable}\/ to include
an additional weight $\lambda^{\cyc(\sigma)}$ in this latter J-fraction,
because Theorem~\ref{thm.perm.Jtype.v2.weaker0.pq}
requires $q_{+1} = p_{+1}$ and $q_{+2} = p_{+2}$,
which are not the case in \reff{eq.spec.perm.inv}.
This inability to include the weight $\lambda^{\cyc(\sigma)}$
is not merely a limitation of our method of proof,
but is inherent in the problem.
Indeed, even the simpler weight $q^{\inv(\sigma)} \lambda^{\cyc(\sigma)}$
gives rise to a J-fraction with coefficients that are rational functions
rather than polynomials:
the first coefficients are
\be
   \gamma_0  \,=\,  \lambda ,\;
   \beta_1  \,=\,  \lambda q ,\;
   \gamma_1  \,=\,  q(2 + \lambda q) ,\;
   \beta_2  \,=\,  q^2 (\lambda + 2q + \lambda q^2) ,
\ee
followed by
\be
   \gamma_2  \;=\;
     {q^2 (2 + 6 \lambda q + 6 q^2 + \lambda^2 q^2 + 4 \lambda q^3 + \lambda^2 q^4)
      \over
      \lambda + 2q + \lambda q^2
     } \;.
\ee
It can then be shown that
\begin{itemize}
   \item[(a)] $\gamma_2$ is a polynomial in $\lambda$
      if and only if $q \in \{-1,0,+1,-i,i\}$; and
   \item[(b)] $\gamma_2$ is a polynomial in $q$
      if and only if $\lambda \in \{-1,0,+1\}$.\
\end{itemize}
Note also that the cases $q = \pm 1$ and $\lambda = \pm 1$
are reducible to the trivial cases $q=1$ or $\lambda=1$
by exploiting the elementary identity
$\inv(\sigma) + \cyc(\sigma) = n \bmod 2$.
\myendremark

\medskip

Similarly, Elizalde's \cite[section~4.2]{Elizalde_18} J-fraction
for 321-avoiding permutations with the weights
$a^{\cval(\sigma)} b^{\cdrise(\sigma)} w_0^{\rar(\sigma)} q^{\inv(\sigma)}$
is the special case of
\reff{eq.thm.perm.pq.Jtype.BIG.321}/\reff{def.weights.perm.pq.Jtype.BIG.321}
with
\be
   x_1 = 1,\quad
   x_2 = q,\quad
   y_1 = qa,\quad
   y_2 = qb,\quad
   p_{+1} = p_{+2} = p_{-1} = p_{-2} = q
   \;.
   \quad
\ee

Finally, Biane \cite[section~6]{Biane_93} has given a $q$-analogue
of \reff{eq.cfrac.secant}, by defining the $q$-secant numbers
in terms of cycle-alternating permutations:
\be
   E_{2n}(q)
   \;\eqdef\;
   \sum_{\sigma \in \Sym^{\rm ca}_{2n}} q^{\inv(\sigma)}
   \;.
\ee
Since a cycle-alternating permutation cannot have any
cycle double rises, cycle double falls or fixed points,
the identity \reff{eq.prop.inv} specializes in this case to
\be
   \inv  \;=\; \cval + \ucross + \lcross + 2(\unest + \lnest)
   \;.
\ee
Applying Theorem~\ref{thm.perm.cycle-alternating.2} with
\be
   x_1 = u_1 = 1 ,\quad
   y_1 = v_1 = q ,\quad
   p_{+1} = p_{-1} = q ,\quad
   q_{+1} = q_{-1} = q^2 \;,
\ee
we obtain an S-fraction with coefficients
\be
   \alpha_n  \;=\; q^{2n-1} \, [n]_q^2
   \;,
\ee
exactly as given by Biane \cite{Biane_93}.

%
%

\subsection{A remark on linear statistics} 
   \label{subsec.permutations.linear}

In this paper we have studied the classification of indices $i \in [n]$
in a permutation $\sigma \in \Sym_n$ according to {\em cyclic}\/ statistics,
i.e.\ cycle peak, cycle valley, cycle double rise, cycle double fall,
and fixed point.
An alternative classification involves {\em linear}\/ statistics,
i.e.\ classifying indices $i \in [n]$ as
\begin{itemize}
   \item {\em peak}\/:  $\sigma(i-1) < \sigma(i) > \sigma(i+1)$;
   \item {\em valley}\/:  $\sigma(i-1) > \sigma(i) < \sigma(i+1)$;
   \item {\em double rise}\/:  $\sigma(i-1) < \sigma(i) < \sigma(i+1)$;
   \item {\em double fall}\/:  $\sigma(i-1) > \sigma(i) > \sigma(i+1)$.
\end{itemize}
However, in order to define the linear statistics it is necessary
to adopt boundary conditions at the two ends ($i=0$ and $i=n+1$):
namely, at each end we can set $\sigma$ to be either
0 (or equivalently $-\infty$),
$n+1$ (or equivalently $+\infty$),
or ``undefined'' --- where an inequality involving ``undefined''
is considered to be automatically false.
(Thus, for instance, if $\sigma(0) = \hbox{undefined}$,
then we count peaks, valleys, double rises and double falls
only starting at $i=2$.)
So there are {\em nine}\/ possible combinations of boundary conditions,
of which four are essentially distinct:
\begin{itemize}
   \item $\sigma(0) = \sigma(n+1) = 0$
              [or the ``dual'' $\sigma(0) = \sigma(n+1) = n+1$];
   \item $\sigma(0) = 0$ and $\sigma(n+1) = n+1$
              [or the ``dual'' $\sigma(0) = n+1$ and $\sigma(n+1) = 0$];
   \item $\sigma(0) = 0$ and $\sigma(n+1) = \hbox{undefined}$
              [or any of the three ``duals''];
   \item $\sigma(0) = \sigma(n+1) = \hbox{undefined}$.
\end{itemize}
The first two of these boundary conditions have been extensively
studied (e.g.\ \cite{Francon_79,Zeng_93}),
and some limited sets of statistics (e.g.\ only peaks) have been studied
under all of these boundary conditions
(e.g.\ \cite{Petersen_07,Ma_12b,Zhuang_17}).

It would be an interesting project to study the joint distribution
of the four linear statistics --- possibly along with other statistics ---
under each of the four possible boundary conditions,
and to obtain continued fractions for the associated
ordinary generating functions.
We can imagine at least two ways in which this could be done:
by transforming our results for cyclic statistics
using Foata's fundamental transformation
\cite[section~1.3]{Foata_70} \cite[pp.~17--18]{Stanley_86}
\cite[section~3.3.1]{Bona_12}
or some other bijection (e.g.\ \cite[section~5]{Clarke_97});
or alternatively by imitating our proofs
in Section~\ref{sec.proofs.permutations}
but using the Fran\c{c}on--Viennot \cite{Francon_79} bijection
in place of the Foata--Zeilberger \cite{Foata_90}
and Biane \cite{Biane_93} bijections.\footnote{
   See \cite[section~5]{Clarke_97} for a discussion of
   the relationship between these three bijections.
}
We must leave this study for future work;
but we should note that an impressive start
was already made two decades ago by Randrianarivony
\cite[Th\'eor\`eme~2]{Randrianarivony_98b}.

Let us observe, finally, that there is an alternate way of getting
$p,q$-generalizations by using generalized pattern avoidance instead of
crossings and nestings
\cite{Clarke_97,Claesson_02,Corteel_07,Shin_10,Shin_12,Blitvic_20}.
These pattern-avoidance statistics mesh naturally with linear statistics.

\section{Set partitions: Statement of results}   \label{sec.intro.setpartitions}

\subsection{S-fraction}  \label{subsec.intro.setpartitions.S}

The {\em Bell number}\/ $B_n$ is, by definition,
the number of partitions of an $n$-element set into nonempty blocks;
by convention we set $B_0 = 1$.
The ordinary generating function of the Bell numbers
can be represented as an S-type continued fraction
\be
   \sum_{n=0}^\infty B_n \, t^n
   \;=\;
   \cfrac{1}{1 - \cfrac{1t}{1 - \cfrac{1 t}{1 - \cfrac{1t}{1- \cfrac{2 t}{1 - \cdots}}}}}
 \label{eq.bellnumber.contfrac}
\ee
with coefficients $\alpha_{2k-1} = 1$, $\alpha_{2k} = k$.\footnote{
   We are not sure where the S-fraction \reff{eq.bellnumber.contfrac}
   first appeared.
   The J-fraction that is equivalent to \reff{eq.bellnumber.contfrac}
   by contraction \reff{eq.contraction_even.coeffs}
   was found by Touchard \cite[section~4]{Touchard_56} in 1956
   (up to a change of variables $x = 1/t$).
   Flajolet \cite[Theorem~2(ia)]{Flajolet_80}
   gave a combinatorial proof of this J-fraction;
   and he observed \cite[pp.~141--142]{Flajolet_80}
   that this J-fraction is implicit in the three-term recurrence relation
   for the Poisson--Charlier polynomials
   \cite[p.~25, Exercise~4.10]{Chihara_78}.
   The S-fraction \reff{eq.bellnumber.contfrac}
   --- as well as some $q$-generalizations ---
   can be derived directly from the functional equation satisfied by the
   ordinary generating function $\sum_{n=0}^\infty B_n t^n$:
   see \cite{Dumont_89} \cite[proof of Lemma~3]{Zeng_95}
   for this elegant method.
   The S-fraction \reff{eq.bellnumber.contfrac}
   is also a straightforward consequence of Aigner's \cite{Aigner_99b}
   evaluation of the zero-shifted and once-shifted Hankel determinants
   of the Bell numbers.
   However, there may well be earlier references
   (for either the S-fraction or the J-fraction) of which we are unaware;
   we would be grateful to readers for drawing our attention to them.
}
Inspired by \reff{eq.bellnumber.contfrac},
let us introduce the polynomials $B_n(x,y,v)$ defined by the continued fraction
\be
   \sum_{n=0}^\infty B_n(x,y,v) \: t^n
   \;=\;
   \cfrac{1}{1 - \cfrac{xt}{1 - \cfrac{yt}{1 - \cfrac{xt}{1- \cfrac{(y+v)t}{1 - \cfrac{xt}{1 - \cfrac{(y+2v)t}{1-\cdots}}}}}}}
 \label{eq.setpartitions.contfrac}
\ee
with coefficients
\begin{subeqnarray}
   \alpha_{2k-1}  & = &  x \\
   \alpha_{2k}    & = &  y + (k-1) v
 \label{def.weights.setpartitions}
\end{subeqnarray}
Clearly $B_n(x,y,v)$ is a homogeneous polynomial of degree $n$;
it therefore has two ``truly independent'' variables.
Since $B_n(1,1,1) = B_n$,
it is plausible to expect that $B_n(x,y,v)$
enumerates partitions of the set $[n]$
according to some natural bivariate statistic.
Of course $B_n(x,y,v)$ is simply $P_n(x,y,u,v)$
[cf.\ \reff{eq.eulerian.fourvar.contfrac}]
specialized to $u=0$;
but our goal here is to interpret it in terms of set partitions,
not permutations.  Our result is:

\begin{theorem}[S-fraction for set partitions]
   \label{thm.setpartitions}
The polynomials $B_n(x,y,v)$ defined by
\reff{eq.setpartitions.contfrac}/\reff{def.weights.setpartitions}
have the combinatorial interpretation
\be
   B_n(x,y,v)
   \;=\;
   \sum_{\pi \in \Pi_n}
      x^{|\pi|} y^{\erec(\pi)} v^{n - |\pi| - \erec(\pi)}
 \label{eq.thm.setpartitions}
\ee
where $|\pi|$ [resp.\ $\erec(\pi)$]
denotes the number of blocks (resp.\ exclusive records) in $\pi$.
\end{theorem}

We need to explain what we mean by an ``exclusive record''
of a set partition $\pi$.
First of all,
given a partition $\pi$ of $[n]$, we say that an element $i \in [n]$ is
\begin{itemize}
   \item an {\em opener}\/ if it is the smallest element of a
      block of size $\ge 2$;
   \item a {\em closer}\/ if it is the largest element of a
      block of size $\ge 2$;
   \item an {\em insider}\/ if it is a non-opener non-closer element of a
      block of size $\ge 3$;
   \item a {\em singleton}\/ if it is the sole element of a block of size 1.
\end{itemize}
Clearly every element $i \in [n]$ belongs to precisely one of these
four classes.

Then we define ``exclusive record'' as follows:
For $i \in [n]$, write $B_\pi(i)$ for the block of $\pi$ containing $i$,
and then define $\sigma'(i)$ to be the next-larger element of $B_\pi(i)$
after $i$, if $i$ is not the largest element of $B_\pi(i)$,
and 0 otherwise.
We then say that an index $i \in [n]$ is an {\em exclusive record}\/ of $\pi$
if it is a nonzero record of the word $\sigma'$,
i.e.\ if $\sigma'(i) \neq 0$
and $\sigma'(j) < \sigma'(i)$ for all $j < i$.
Pictorially, we can say that $i$ is an exclusive record of $\pi$
if it is not the largest element of its block
(that is, it is either an opener or an insider)
and its right neighbor (within its block) sticks out farther to the right
than any right neighbor (within its block) of a vertex $< i$.
In Section~\ref{subsec.intro.setpartitions.firstmaster}
we will reinterpret the notion of exclusive record
in terms of nestings.

Since every exclusive record is either an opener or an insider,
while $|\pi|$ equals the number of closers plus the number of singletons,
it follows that ${n - |\pi| - \erec(\pi)}$ counts the
openers and insiders that are not exclusive records.
In particular,
\linebreak
${n - |\pi| - \erec(\pi) \ge 0}$,
so that the right-hand side of \reff{eq.thm.setpartitions}
is indeed a polynomial.

We will prove Theorem~\ref{thm.setpartitions}
in Section~\ref{subsec.setpartitions.S}.

\bigskip

{\bf Remarks.}
1.  For the special case $y=v$,
i.e.\ the Bell polynomials
\be
   B_n(x) \;=\;
   B_n(x,1,1)  \;=\;
   \sum_{\pi \in \Pi_n} x^{|\pi|}
          \;=\; \sum_{k=0}^n \stirlingsubset{n}{k} \: x^k
\ee
or their homogenization,
Flajolet \cite[Theorem~2(ib)]{Flajolet_80} found a J-type continued fraction
that is equivalent by contraction \reff{eq.contraction_even.coeffs}
to the specialization of \reff{eq.setpartitions.contfrac}.
Later, Dumont \cite{Dumont_89} found the S-fraction directly
by a functional-equation method,
and one of us used this same method to find \cite[Lemma~3]{Zeng_95}
two $q$-generalizations of the S-fraction
(these will be discussed in Section~\ref{subsec.setpartitions.inv} below).

2.  The triangular array corresponding to the polynomials $B_n(x,1,x)$
can be found in \cite[A085791]{OEIS}, but without a
combinatorial interpretation.
The sub-triangular array corresponding to the polynomials $B_n(1,1,v)$,
which interpolate between the Catalan numbers ($v=0$)
and the Bell numbers ($v=1$), is apparently not in \cite{OEIS}.
\myendremark



\subsection{J-fraction}

We can refine the polynomial $B_n(x,y,v)$
by distinguishing between singletons and blocks of size $\ge 2$;
in addition, we can distinguish between exclusive records
that are openers and those that are insiders;
and finally, among indices that are not exclusive records,
we can also distinguish those that are openers from those that are insiders.
That is, we define
\be
   B_n(x_1,x_2,y_1,y_2,v_1,v_2)
   \;=\;
   \sum_{\pi \in \Pi_n}
      x_1^{m_1(\pi)} x_2^{m_{\ge 2}(\pi)}
      y_1^{\erecin(\pi)} y_2^{\erecop(\pi)}
      v_1^{\nerecin(\pi)} v_2^{\nerecop(\pi)}
   \;,
 \label{def.Bn.x1x2y1y2}
\ee
where $m_1(\pi)$ is the number of singletons in $\pi$,
$m_{\ge 2}(\pi)$ is the number of non-singleton blocks
(or equivalently the number of closers),
$\erecin(\pi)$ is the number of insiders that are exclusive records,
$\erecop(\pi)$ is the number of openers that are exclusive records,
$\nerecin(\pi)$ is the number of insiders that are not exclusive records,
and $\nerecop(\pi)$ is the number of openers that are not exclusive records.
We then have a nice J-fraction:

\begin{theorem}[J-fraction for set partitions]
   \label{thm.setpartitions.Jtype}
The ordinary generating function of the polynomials
$B_n(x_1,x_2,y_1,y_2,v_1,v_2)$
has the J-type continued fraction
\begin{eqnarray}
   & &
   \sum\limits_{n=0}^\infty B_n(x_1,x_2,y_1,y_2,v_1,v_2) \: t^n
   \;=\;
         \nonumber \\
   & & \qquad\qquad
   \cfrac{1}{1 - x_1 t - \cfrac{x_2 y_2 t^2}{1 -  (x_1\!+\!y_1) t - \cfrac{x_2(y_2\!+\!v_2) t^2}{1 - (x_1\!+\!y_1\!+\!v_1)t - \cfrac{x_2(y_2\!+\!2v_2) t^2}{1 - \cdots}}}}
   \qquad
   \label{eq.thm.setpartitions.Jtype}
\end{eqnarray}
with coefficients
\begin{subeqnarray}
   \gamma_0  & = &   x_1   \\[1mm]
   \gamma_n  & = &   x_1 + y_1 + (n-1)v_1
        \qquad\hbox{for $n \ge 1$}   \\[1mm]
   \beta_n   & = &   x_2 \, [y_2 + (n-1)v_2]
 \label{def.weights.setpartitions.Jtype}
\end{subeqnarray}
\end{theorem}

\noindent
We will prove this theorem in Section~\ref{subsec.setpartitions.J},
as a special case of a more general result.
The case $y_1 = y_2 = v_1 = v_2$ was obtained previously by
Flajolet \cite[Theorem~2]{Flajolet_80}.

When specialized to $x_1 = x_2$ and $y_1 = y_2$ and $v_1 = v_2$,
\reff{def.Bn.x1x2y1y2} reduces to \reff{eq.thm.setpartitions},
and the J-fraction \reff{eq.thm.setpartitions.Jtype}
is the contraction \reff{eq.contraction_even.coeffs}
of the S-fraction \reff{eq.setpartitions.contfrac}.
So Theorem~\ref{thm.setpartitions} is a corollary of
a special case of Theorem~\ref{thm.setpartitions.Jtype}.

\bigskip

{\bf Warning:}  Because we have chosen to define $x_1$ and $x_2$
as conjugate to $m_1(\pi)$ and $m_{\ge 2}(\pi)$, respectively,
the meaning of the subscripts 1 and 2
in the continued fraction \reff{eq.thm.setpartitions.Jtype}
is reversed vis-\`a-vis our usage for permutations
(compare Theorem~\ref{thm.setpartitions.Jtype}
 with Theorem~\ref{thm.perm.Jtype}).

\bigskip

{\bf Remark.}
By setting $y_1 = v_1 = 0$, we can suppress insiders
and thereby restrict the sum to set partitions in which every block
is of size either 1 or 2.
(These are in obvious bijection with involutions,
i.e.\ permutations in which every cycle is of length either 1 or 2.)
The resulting J-fraction
\reff{eq.thm.setpartitions.Jtype}/\reff{def.weights.setpartitions.Jtype}
has $\gamma_n = x_1$ for all $n$,
which means \cite[eq.~(6.15)]{Aigner_01a} \cite[Proposition~4]{Barry_09}
that the polynomials $B_n(x_1,x_2,0,y_2,0,v_2)$
are the $x_1$-binomial transform of the polynomials
$B_n(0,x_2,0,y_2,0,v_2) = x_2^n B_n(0,1,0,y_2,0,v_2)$
that enumerate perfect matchings:
\be
   B_n(x_1,x_2,0,y_2,0,v_2)
   \;=\;
   \sum_{k=0}^n \binom{n}{k} \, B_k(0,x_2,0,y_2,0,v_2) \: x_1^{n-k}
   \;.
\ee
This relation is also obvious combinatorially.
See Section~\ref{sec.intro.perfectmatchings}
for more information on perfect matchings.
\myendremark

\subsection[First $p,q$-generalization: Crossings and nestings]{First $\bm{p,q}$-generalization: Crossings and nestings}
   \label{subsec.intro.setpartitions.pq}

Let $\pi = \{B_1,B_2,\ldots,B_{|\pi|}\}$ be a partition of $[n]$.
We associate to the partition~$\pi$ a~graph $\scrg_\pi$
with vertex set $[n]$ such that $i,j$ are joined by an edge
if and only if they are {\em consecutive}\/ elements
within the same block.\footnote{
   Kasraoui and Zeng \cite{Kasraoui_06} call this the
   {\em partition graph}\/ associated to $\pi$;
   Mansour \cite[Definition~3.50]{Mansour_13}
   calls it the {\em standard representation}\/ of $\pi$.
}
(The graph $\scrg_\pi$ thus has $|\pi|$~connected components
and $n-|\pi|$~edges.)
We always write an edge $e$ of $\scrg_\pi$ as a pair $(i,j)$ with $i < j$.
We then say that a quadruplet $i < j < k < l$ forms a
\begin{itemize}
  \item {\em crossing}\/ (cr) if $(i,k) \in \scrg_\pi$ and $(j,l) \in \scrg_\pi$
     [note that $i,k$ and $j,l$ must then belong to different blocks];
  \item {\em nesting}\/  (ne) if $(i,l) \in \scrg_\pi$ and $(j,k) \in \scrg_\pi$
     [note that $i,l$ and $j,k$ must then belong to different blocks].
\end{itemize}
We also consider the ``degenerate'' case with $j=k$,
by saying that a triplet $i < j < l$ forms a
\begin{itemize}
   \item {\em pseudo-nesting}\/ (psne) if $j$ is a singleton
       and $(i,l) \in \scrg_\pi$.
\end{itemize}
See Figure~\ref{fig.setpartitions.crossnest}.
We define $\crr(\pi)$ [resp.\ $\nee(\pi)$, $\psne(\pi)$] to be number of
crossings (resp.\ nestings, pseudo-nestings) in $\pi$.

\begin{figure}[p]
\centering
\begin{picture}(30,15)(145, 10)
\setlength{\unitlength}{1.5mm}
\linethickness{.5mm}
\put(2,0){\line(1,0){28}}
\put(5,0){\circle*{1,3}}\put(5,0){\makebox(0,-6)[c]{\small $i$}}
\put(12,0){\circle*{1,3}}\put(12,0){\makebox(0,-6)[c]{\small $j$}}
\put(19,0){\circle*{1,3}}\put(19,0){\makebox(0,-6)[c]{\small $k$}}
\put(26,0){\circle*{1,3}}\put(26,0){\makebox(0,-6)[c]{\small $l$}}
\red{\qbezier(5,0)(12,10)(19,0)}
\blue{\qbezier(11,0)(18,10)(25,0)}
\put(15,-6){\makebox(0,-6)[c]{\small Crossing}}
\end{picture}
\begin{picture}(30,15)(145, 10)
\setlength{\unitlength}{1.5mm}
\linethickness{.5mm}
\put(43,0){\line(1,0){28}}
\put(47,0){\circle*{1,3}}\put(47,0){\makebox(0,-6)[c]{\small $i$}}
\put(54,0){\circle*{1,3}}\put(54,0){\makebox(0,-6)[c]{\small $j$}}
\put(61,0){\circle*{1,3}}\put(61,0){\makebox(0,-6)[c]{\small $k$}}
\put(68,0){\circle*{1,3}}\put(68,0){\makebox(0,-6)[c]{\small $l$}}
\red{\qbezier(47,0)(57.5,10)(68,0)}
\blue{\qbezier(53,0)(56.5,5)(60,0)}
\put(57,-6){\makebox(0,-6)[c]{\small Nesting}}
\end{picture}
\\[2.5cm]
\begin{picture}(30,15)(145, 10)
\setlength{\unitlength}{1.5mm}
\linethickness{.5mm}
\put(23,0){\line(1,0){28}}
\put(26,0){\circle*{1,3}}\put(26,0){\makebox(0,-6)[c]{\small $i$}}
\put(36.5,0){\circle*{1,3}}\put(36.5,0){\makebox(0,-6)[c]{\small $j$}}
\put(47,0){\circle*{1,3}}\put(47,0){\makebox(0,-6)[c]{\small $l$}}
\red{\qbezier(26,0)(36.5,10)(47,0)}
\put(36,-6){\makebox(0,-6)[c]{\small Pseudo-nesting}}
\end{picture}
\vspace*{2.5cm}
\caption{
   Crossing, nesting and pseudo-nesting for set partitions.
 \label{fig.setpartitions.crossnest}
}
\end{figure}
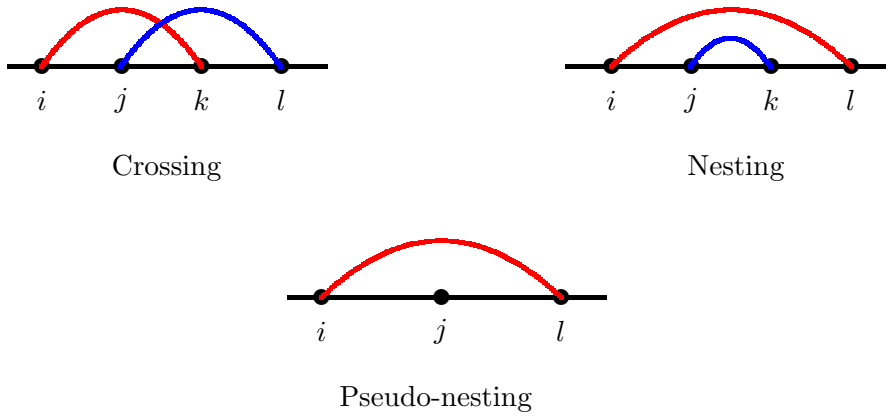

\begin{figure}[p]
\centering
\begin{picture}(30,15)(145, 10)
\setlength{\unitlength}{1.5mm}
\linethickness{.5mm}
\put(2,0){\line(1,0){28}}
\put(5,0){\circle*{1,3}}\put(5,0){\makebox(0,-6)[c]{\small $i$}}
\put(12,0){\circle*{2,4}}\put(12,0){\makebox(0,-6)[c]{\small $j$}}
\put(19,0){\circle*{1,3}}\put(19,0){\makebox(0,-6)[c]{\small $k$}}
\put(26,0){\circle*{1,3}}\put(26,0){\makebox(0,-6)[c]{\small $l$}}
\red{\qbezier(5,0)(12,10)(19,0)}
\blue{\qbezier(11,0)(18,10)(25,0)}
\put(15,-6){\makebox(0,-6)[c]{\small Crossing of opener type}}
\put(43,0){\line(1,0){28}}
\put(47,0){\circle*{1,3}}\put(47,0){\makebox(0,-6)[c]{\small $i$}}
\put(54,0){\circle*{1,3}}\put(54,0){\makebox(0,-6)[c]{\small $j$}}
\put(61,0){\circle*{1,3}}\put(61,0){\makebox(0,-6)[c]{\small $k$}}
\put(68,0){\circle*{1,3}}\put(68,0){\makebox(0,-6)[c]{\small $l$}}
\red{\qbezier(47,0)(54,10)(61,0)}
\blue{\qbezier(53,0)(60,10)(67,0)}
\blue{\qbezier(50,2)(51.25,1)(52.5,0)}
\put(57,-6){\makebox(0,-6)[c]{\small Crossing of insider type}}
\end{picture}
\\[3.5cm]
\begin{picture}(30,15)(145, 10)
\setlength{\unitlength}{1.5mm}
\linethickness{.5mm}
\put(2,0){\line(1,0){28}}
\put(5,0){\circle*{1,3}}\put(5,0){\makebox(0,-5)[c]{\small $i$}}
\put(12,0){\circle*{2,4}}\put(12,0){\makebox(0,-5)[c]{\small $j$}}
\put(19,0){\circle*{1,3}}\put(19,0){\makebox(0,-5)[c]{\small $k$}}
\put(26,0){\circle*{1,3}}\put(26,0){\makebox(0,-5)[c]{\small $l$}}
\red{\qbezier(5,0)(15.5,10)(26,0)}
\blue{\qbezier(11,0)(14.5,5)(18,0)}
\put(15,-6){\makebox(0,-6)[c]{\small Nesting of opener type}}
\put(43,0){\line(1,0){28}}
\put(47,0){\circle*{1,3}}\put(47,0){\makebox(0,-6)[c]{\small $i$}}
\put(54,0){\circle*{1,3}}\put(54,0){\makebox(0,-6)[c]{\small $j$}}
\put(61,0){\circle*{1,3}}\put(61,0){\makebox(0,-6)[c]{\small $k$}}
\put(68,0){\circle*{1,3}}\put(68,0){\makebox(0,-6)[c]{\small $l$}}
\red{\qbezier(47,0)(57,10)(68,0)}
\blue{\qbezier(53,0)(56.5,5)(60.5,0)}
\blue{\qbezier(50,2)(51.25,1)(52.5,0)}
\put(57,-6){\makebox(0,-6)[c]{\small Nesting of insider type}}
\end{picture}
\vspace*{2.5cm}
\caption{
   Refined categories of crossing and nesting for set partitions.
   Here a large circle marks an opener.
 \label{fig.setpartitions.refined_crossnest}
}
\end{figure}
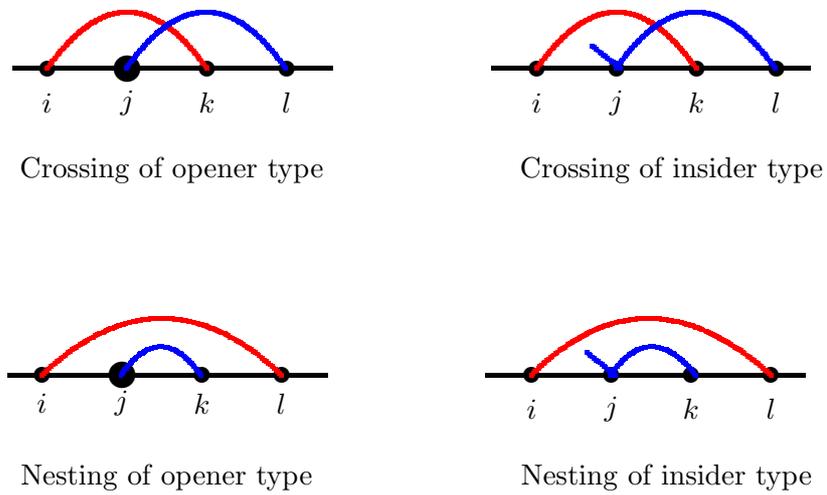

We now introduce a $p,q$-generalization of the polynomial
\reff{def.Bn.x1x2y1y2}:
\begin{eqnarray}
   & &  \hspace*{-6mm}
   B_n(x_1,x_2,y_1,y_2,v_1,v_2,p,q,r)
        \nonumber \\[2mm]
   & & 
   \,=\;
   \sum_{\pi \in \Pi_n}
      x_1^{m_1(\pi)} x_2^{m_{\ge 2}(\pi)}
      y_1^{\erecin(\pi)} y_2^{\erecop(\pi)}
      v_1^{\nerecin(\pi)} v_2^{\nerecop(\pi)}
      p^{\crr(\pi)} q^{\nee(\pi)} r^{\psne(\pi)}
   \:.
   \qquad
 \label{def.Bn.crossnest}
\end{eqnarray}
These polynomials have a nice J-fraction, as we shall see.

%

But we can go farther, and refine the categories of crossing and nesting
by analogy with what was done for permutations in
Section~\ref{subsec.intro.permutations.pq}.
Let us say that a quadruplet $i < j < k < l$ forms a
\begin{itemize}
  \item {\em crossing of opener type}\/ (crop)
     if $(i,k) \in \scrg_\pi$ and $(j,l) \in \scrg_\pi$
     and $j$ is an opener;
  \item {\em crossing of insider type}\/ (crin)
     if $(i,k) \in \scrg_\pi$ and $(j,l) \in \scrg_\pi$
     and $j$ is an insider;
  \item {\em nesting of opener type}\/  (neop)
     if $(i,l) \in \scrg_\pi$ and $(j,k) \in \scrg_\pi$
     and $j$ is an opener;
  \item {\em nesting of insider type}\/  (nein)
     if $(i,l) \in \scrg_\pi$ and $(j,k) \in \scrg_\pi$
     and $j$ is an insider.
\end{itemize}
See Figure~\ref{fig.setpartitions.refined_crossnest}.
Please note that here the distinguished index $j$
is the one in second position.
The categories crop, crin, neop, nein for set partitions
correspond, respectively,
to ucrosscval, ucrosscdrise, unestcval, unestcdrise for permutations,
under a mapping that will be discussed in
Section~\ref{subsec.intro.setpartitions.firstmaster} below.

Let us now define the refined polynomial
\begin{eqnarray}
   & &  \hspace*{-3mm}
   B_n(x_1,x_2,y_1,y_2,v_1,v_2,p_1,p_2,q_1,q_2,r)
        \nonumber \\[2mm]
   & & 
   \qquad=\;
   \sum_{\pi \in \Pi_n}
      x_1^{m_1(\pi)} x_2^{m_{\ge 2}(\pi)}
      y_1^{\erecin(\pi)} y_2^{\erecop(\pi)}
      v_1^{\nerecin(\pi)} v_2^{\nerecop(\pi)}
        \:\times
        \qquad\qquad
        \nonumber \\[1mm]
    & & \qquad\qquad\quad\;\:
      p_1^{\crin(\pi)} p_2^{\crop(\pi)}
      q_1^{\nein(\pi)} q_2^{\neop(\pi)}
      r^{\psne(\pi)}
   \:.
   \qquad
 \label{def.Bn.crossnest.refined}
\end{eqnarray}
(Thus, the variables $y_1,v_1,p_1,q_1$ are associated with insiders;
$y_2,v_2,p_2,q_2$ are associated with openers;
$x_1$ and $r$ are associated with singletons;
and $x_2$ can be interpreted as associated with closers.)
These polynomials have a nice J-fraction:

\begin{theorem}[J-fraction for set partitions, first $p,q$-generalization]
   \label{thm.setpartitions.Jtype.pq.refined}
The ordinary generating function of the polynomials
$B_n(x_1,x_2,y_1,y_2,v_1,v_2,p_1,p_2,q_1,q_2,r)$ has the J-type continued fraction
\begin{eqnarray}
   & & \hspace*{-12mm}
   \sum_{n=0}^\infty B_n(x_1,x_2,y_1,y_2,v_1,v_2,p_1,p_2,q_1,q_2,r) \: t^n
   \;=\;
       \nonumber \\
   & & 
   \cfrac{1}{1 - x_1 t - \cfrac{x_2 y_2 t^2}{1 -  (rx_1\!+\!y_1) t - \cfrac{x_2(p_2 y_2\!+\!q_2 v_2) t^2}{1 - (r^2 x_1\!+\!p_1 y_1\!+\!q_1 v_1)t - \cfrac{x_2(p_2^2 y_2\!+\! q_2 \, [2]_{p_2,q_2}v_2) t^2}{1 - \cdots}}}}
       \nonumber \\[1mm]
   \label{eq.thm.setpartitions.Jtype.pq.refined}
\end{eqnarray}
with coefficients
\begin{subeqnarray}
   \gamma_0  & = &   x_1   \\[1mm]
   \gamma_n  & = &   r^n x_1 + p_1^{n-1} y_1 + q_1 \, [n-1]_{p_1,q_1} v_1
        \qquad\hbox{for $n \ge 1$}   \\[1mm]
   \beta_n   & = &   x_2 \, (p_2^{n-1} y_2 + q_2 \, [n-1]_{p_2,q_2} v_2)
 \label{def.weights.setpartitions.Jtype.pq.refined}
\end{subeqnarray}
\end{theorem}

\noindent
We will prove this theorem in Section~\ref{subsec.setpartitions.J},
as a special case of a more general result.

Note that $\beta_n$ is homogeneous of degree $n-1$ in the pair $(p_2,q_2)$.
So if we multiply both $p_2$ and $q_2$ by $C$,
this has the effect of multiplying $\beta_n$ by $C^{n-1}$.

Note also that when $y_i=v_i$ (for $i=1$ and/or 2),
the weight $p_i^{n-1} y_i + q_i \, [n-1]_{p_i,q_i} v_i$
simplifies to $[n]_{p_i,q_i} y_i$.
In this case the weights \reff{def.weights.setpartitions.Jtype.pq.refined}
are invariant under $p_i \leftrightarrow q_i$.
(In~particular, if $y_1 = v_1$ {\em and}\/ $y_2 = v_2$,
 then the weights are invariant under the {\em independent}\/ interchanges
 $p_1 \leftrightarrow q_1$ and $p_2 \leftrightarrow q_2$.)
It~would be interesting to find a bijective proof of these symmetries.

In particular, the specialization
$y_1 = y_2 = v_1 = v_2$, $p_1 = p_2$, $q_1 = q_2$ and $r=1$
of Theorem~\ref{thm.setpartitions.Jtype.pq.refined}
was proven earlier by Kasraoui and Zeng \cite[Proposition~4.1]{Kasraoui_06}.
In this specialization, the symmetry $p \leftrightarrow q$
is also a consequence of Kasraoui--Zeng's \cite{Kasraoui_06}
bijection that interchanges crossings and nestings
while preserving various other statistics.

\bigskip

When specialized to $x_1 = x_2$, $y_1 = y_2$, $v_1 = v_2$,
$p_2 = r p_1$ and $q_2 = r q_1$,
the J-fraction \reff{eq.thm.setpartitions.Jtype.pq.refined}
arises as the contraction \reff{eq.contraction_even.coeffs}
of an S-fraction with polynomial coefficients:

\begin{corollary}[S-fraction for set partitions, first $p,q$-generalization]
   \label{cor.setpartitions.Jtype.pq}
The ordinary generating function of the polynomials
$B_n(x_1,x_2,y_1,y_2,v_1,v_2,p_1,p_2,q_1,q_2,r)$
specialized to $x_1 = x_2$, $y_1 = y_2$, $v_1 = v_2$,
$p_2 = r p_1$, $q_2 = r q_1$
has the S-type continued fraction
\be
   \sum_{n=0}^\infty B_n(x,x,y,y,v,v,p,rp,q,rq,r) \: t^n
   \;=\;
   \cfrac{1}{1 - \cfrac{xt}{1 - \cfrac{yt}{1 - \cfrac{rxt}{1- \cfrac{(py+qv)t}{1 - \cfrac{r^2 xt}{1 - \cfrac{(p^2 y+ q\, [2]_{p,q}v)t}{1-\cdots}}}}}}}
   \label{eq.cor.setpartitions.Jtype.pq}
\ee
with coefficients
\begin{subeqnarray}
   \alpha_{2k-1}  & = &  r^{k-1} x   \\[1mm]
   \alpha_{2k}    & = &  p^{k-1} y + q \, [k-1]_{p,q} v
 \label{def.weights.setpartitions.pq}
\end{subeqnarray}
\end{corollary}

\bigskip

{\bf Remark.}
Josuat-Verg\`es and Rubey \cite[Theorem~1.2]{Josuat-Verges_11}
have given an explicit formula for $[x^k] \, B_n(x,x,1,1,1,1,1,1,q,q,1)$,
which enumerates the partitions of $[n]$ with $k$ blocks
according to the number of nestings (or crossings).
This formula is analogous to, but more complicated than,
the Touchard--Riordan formula for perfect matchings
given in \reff{eq.Mn1q} below.
\myendremark

\subsection[Second $p,q$-generalization: Overlaps and coverings]{Second $\bm{p,q}$-generalization: Overlaps and coverings}
   \label{subsec.intro.setpartitions.pq.2.0}

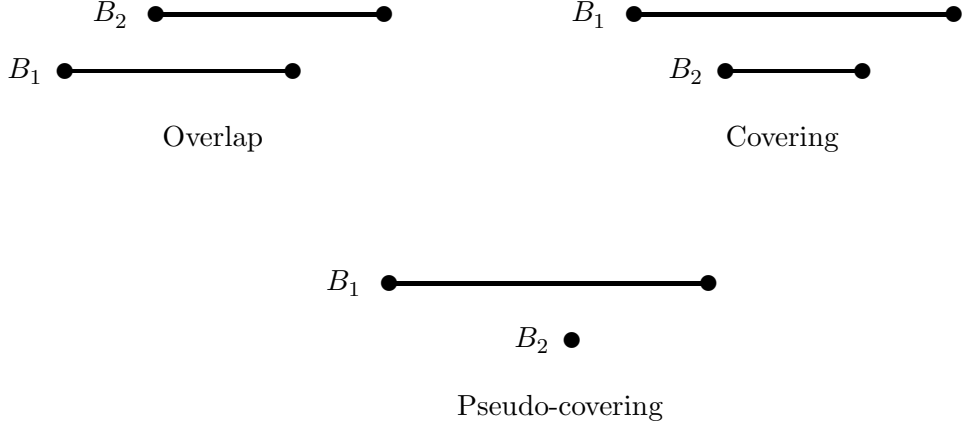
\begin{figure}[t]
\vspace*{1cm}
\centering
\begin{picture}(30,15)(145, 10)
\setlength{\unitlength}{1.5mm}
\linethickness{.5mm}
\put(4,3){\makebox(4,4)[c]{\small $B_2$}}
\put(10,5){\line(1,0){20}}
\put(10,5){\circle*{1.3}}
\put(30,5){\circle*{1.3}}
\put(-1,0){\makebox(-1,0)[c]{\small $B_1$}}
\put(2,0){\line(1,0){20}}
\put(2,0){\circle*{1.3}}
\put(22,0){\circle*{1.3}}
\put(15,-3){\makebox(0,-6)[c]{\small Overlap}}
\end{picture}
\begin{picture}(30,15)(145, 10)
\setlength{\unitlength}{1.5mm}
\linethickness{.5mm}
\put(38,3){\makebox(4,4)[c]{\small $B_1$}}
\put(44,5){\line(1,0){28}}
\put(44,5){\circle*{1.3}}
\put(72,5){\circle*{1.3}}
\put(49,0){\makebox(-1,0)[c]{\small $B_2$}}
\put(52,0){\line(1,0){12}}
\put(52,0){\circle*{1.3}}
\put(64,0){\circle*{1.3}}
\put(57,-3){\makebox(0,-6)[c]{\small Covering}}
\end{picture}
\\[3cm]
\begin{picture}(30,15)(30, 10)
\setlength{\unitlength}{1.5mm}
\linethickness{.5mm}
\put(1,3){\makebox(-10,4)[c]{\small $B_1$}}
\put(0,5){\circle*{1.3}}
\put(0,5){\line(1,0){28}}
\put(28,5){\circle*{1.3}}
\put(13,0){\makebox(-1,0)[c]{\small $B_2$}}
\put(16,0){\circle*{1.3}}
\put(15,-3){\makebox(0,-6)[c]{\small Pseudo-covering}}
\end{picture}
\vspace*{2cm}
\caption{
   Overlap, covering and pseudo-covering for set partitions.
   The line segments represent the intervals
   $[\min B_1,\max B_1]$ and $[\min B_2,\max B_2]$.
   Only the smallest and largest elements of each block are shown explicitly
   with a dot.
 \label{fig.setpartitions.ovcov}
}
\end{figure}

We can form a different type of $p,q$-generalization
by looking at the crossings and nestings of entire blocks
rather than nearest-neighbor edges.
We say that blocks $B_1$ and $B_2$ of $\pi$ form
\begin{itemize}
  \item an {\em overlap}\/ if $\min B_1 < \min B_2 < \max B_1 < \max B_2$;
  \item a {\em covering}\/ if $\min B_1 < \min B_2 < \max B_2 < \max B_1$;
  \item a {\em pseudo-covering}\/ if $\min B_1 < \min B_2 = \max B_2 < \max B_1$
     (so that here $B_2$ is a singleton).
\end{itemize}
See Figure~\ref{fig.setpartitions.ovcov}.
We write $\ov(\pi)$, $\cov(\pi)$ and $\pscov(\pi)$
for the number of overlaps, coverings and pseudo-coverings in $\pi$,
respectively.

Let us also observe that $\pscov(\pi) = \psne(\pi)$:
for if $B_2 = \{j\}$ is a singleton and $\min B_1 < j < \max B_1$,
then there is precisely one edge $(i,l) \in \scrg_\pi$ with $i,l \in B_1$
such that $i<j<l$.
So pseudo-coverings and pseudo-nestings are
two different names for the same quantity.

Let us now define two related quantities:
\begin{subeqnarray}
   \ovin(\pi)
   & = &
   \!\!
   \sum\limits_{j \in {\rm insiders}}  \!\!
   \#\{ (B_1,B_2) \colon\:  j \in B_2 \,\hbox{ and }\,
                          \min B_1 < j < \max B_1 < \max B_2  \}
         \nonumber \\[-4mm] \\[1mm]
   \covin(\pi)
   & = &
   \!\!
   \sum\limits_{j \in {\rm insiders}}  \!\!
   \#\{ (B_1,B_2) \colon\:  j \in B_2 \,\hbox{ and }\,
                          \min B_1 < j < \max B_2 < \max B_1  \}
         \nonumber \\[-4mm]
 \label{def.ovin.covin}
\end{subeqnarray}
Here $j$ is explicitly required to be an insider of the block $B_2$ (not an opener).
Let us stress, however, that {\em both}\/ possible inequalities
of $\min B_1$ and $\min B_2$ are allowed, with opposite effect:
\begin{itemize}
   \item If $\min B_1 < \min B_2$, then the pair $(B_1,B_2)$
      contributing to $\ovin$ (resp.\ $\covin$) is an overlap (resp.\ covering).
   \item If $\min B_1 > \min B_2$, then the pair $(B_1,B_2)$
      contributing to $\ovin$ (resp.\ $\covin$) is a covering (resp.\ overlap).
\end{itemize}
The motivation for these somewhat strange definitions
will become apparent in Section~\ref{subsec.intro.setpartitions.secondmaster}.

Along with these definitions, we can introduce a different notion of record
that is better adapted to overlaps and coverings.
Recall that $j$ is an exclusive record of $\pi$
if it is either an opener or an insider
and its right neighbor (within its block) sticks out farther to the right
than any previous right neighbor as we read the graph $\scrg_\pi$
from left to right.
Let us now say that $j$ is a {\em block-record}\/
if it is either an opener or an insider
and its {\em block}\/ sticks out farther to the right
than any {\em block}\/ containing a vertex $< j$.
That~is, $j$ is a block-record
if it is an opener or insider of a block $B$
and there does not exist a block $B'$ satisfying
$\min B' < j < \max B < \max B'$.
We write $\brec(\pi)$ for the number of block-records in~$\pi$;
and more specifically, we write
$\brecin(\pi)$ for the number of insiders that are block-records,
$\brecop(\pi)$ for the number of openers that are block-records,
$\nbrecin(\pi)$ for the number of insiders that are not block-records,
and $\nbrecop(\pi)$ for the number of openers that are not block-records.
In Section~\ref{subsec.intro.setpartitions.secondmaster}
we will reinterpret the notion of block-record in terms of coverings.

We now define a polynomial that is analogous to
\reff{def.Bn.crossnest.refined}
but uses block-records, overlaps, coverings and pseudo-coverings
in place of exclusive records, crossings, nestings and pseudo-nestings:
\begin{eqnarray}
   & &  \hspace*{-3mm}
   B^{(2)}_n(x_1,x_2,y_1,y_2,v_1,v_2,p_1,p_2,q_1,q_2,r)
        \nonumber \\[2mm]
   & & 
   \qquad=\;
   \sum_{\pi \in \Pi_n}
      x_1^{m_1(\pi)} x_2^{m_{\ge 2}(\pi)}
      y_1^{\brecin(\pi)} y_2^{\brecop(\pi)}
      v_1^{\nbrecin(\pi)} v_2^{\nbrecop(\pi)}
        \:\times
        \qquad\qquad
        \nonumber \\[-1mm]
    & & \qquad\qquad\quad\;\:
      p_1^{\ovin(\pi)} p_2^{\ov(\pi)}
      q_1^{\covin(\pi)} q_2^{\cov(\pi)}
      r^{\pscov(\pi)}
   \:.
   \qquad
 \label{def.Bn.ovcov.brec}
\end{eqnarray}
It turns out that these polynomials are not merely {\em analogous}\/
to \reff{def.Bn.crossnest.refined};
they are {\em identical}\/ to \reff{def.Bn.crossnest.refined}:

\begin{theorem}
   \label{thm.setpartitions.pq.ovcov.brec}
We have
\be
   B^{(2)}_n(x_1,x_2,y_1,y_2,v_1,v_2,p_1,p_2,q_1,q_2,r)
   \;=\;
   B_n(x_1,x_2,y_1,y_2,v_1,v_2,p_1,p_2,q_1,q_2,r)
   \;.
 \label{eq.thm.setpartitions.pq.ovcov.brec}
\ee
\end{theorem}

\noindent
We will prove this theorem in Sections~\ref{subsec.setpartitions.J}
and \ref{subsec.setpartitions.J.2},
by showing that both polynomials have the same J-fraction
\reff{eq.thm.setpartitions.Jtype.pq.refined}/%
\reff{def.weights.setpartitions.Jtype.pq.refined}.
It would be interesting to find a direct bijective proof
of the identity \reff{eq.thm.setpartitions.pq.ovcov.brec}.

\bigskip

{\bf Historical remark.}
Our introduction of overlaps and coverings is inspired by
the work of Flajolet and Schott \cite{Flajolet_90}
and Claesson \cite{Claesson_01}.
Flajolet and Schott considered {\em nonoverlapping}\/ partitions
(i.e.\ partitions with no overlaps)
and gave \cite[eq.~(27)]{Flajolet_90} a J-fraction for the weight $x^{|\pi|}$,
i.e.\ \reff{eq.thm.setpartitions.Jtype.pq.refined}/%
\reff{def.weights.setpartitions.Jtype.pq.refined}
reinterpreted using Theorem~\ref{thm.setpartitions.pq.ovcov.brec}
and then specialized to $x_1 = x_2 = x$, $y_1 = y_2 = v_1 = v_2 = 1$,
$p_1 = 1$, $p_2 =0$ and $q_1 = q_2 =r=1$:
namely, $\gamma_n = x+n$ and $\beta_n = x$.
Claesson considered {\em monotone}\/ partitions
(i.e.\ partitions with no coverings)
and constructed a bijection between them and nonoverlapping partitions
\cite[Proposition~6]{Claesson_01};
he also related them to permutations that avoid certain generalized patterns.
To our knowledge, no one has heretofore considered giving weights
(other than 0 or~1) to overlaps and coverings.
\myendremark


\subsection{Some useful identities}
   \label{subsec.intro.setpartitions.identities}

Before proceeding further, let us record some useful identities
that relate the statistics that have just been introduced.

To begin with, we have the trivial identities
\begin{subeqnarray}
   \erecop(\pi) + \nerecop(\pi) & = & \brecop(\pi) + \nbrecop(\pi) \\[1mm]
   \erecin(\pi) + \nerecin(\pi) & = & \brecin(\pi) + \nbrecin(\pi)
 \label{eq.identity.erec.brec}
\end{subeqnarray}
in which both sides express the total number of openers (resp.\ insiders).

Somewhat less trivially, we have:
   
\begin{lemma}[Crossings + nestings = overlaps + coverings]
   \label{lemma.cross+nest.ov+cov}
We have
\begin{subeqnarray}
   \crop(\pi) + \neop(\pi) & = & \ov(\pi) + \cov(\pi)  \\[1mm]
   \crin(\pi) + \nein(\pi) & = & \ovin(\pi) + \covin(\pi)
 \label{eq.lemma.cross+nest.ov+cov}
\end{subeqnarray}
\end{lemma}

\noindent
The proof of Lemma~\ref{lemma.cross+nest.ov+cov} is not difficult,
but we defer it to Section~\ref{subsec.intro.setpartitions.secondmaster},
where it will arise naturally as a special case of a more general identity
(see Lemma~\ref{lemma.qne}).

Finally, we have:

\begin{lemma}[Crossings and nestings modulo 2]
   \label{lemma.cross+nest.mod2}
We have
\begin{subeqnarray}
   \crr(\pi)  \;=\;
   \crin(\pi) +\, \crop(\pi)  & = &  \ov(\pi)  \quad\;\,\bmod 2   \\[1mm]
   \crin(\pi) + \neop(\pi)  & = &  \cov(\pi) \quad\bmod 2  \\[1mm]
   \crop(\pi) + \nein(\pi)  & = &  \ov(\pi) + \ovin(\pi) + \covin(\pi) \quad\bmod 2
        \\[1mm]
   \nee(\pi)  \;=\;
   \nein(\pi) + \neop(\pi)  & = &  \cov(\pi) + \ovin(\pi) + \covin(\pi) \;\;\bmod 2
        \qquad
 \label{eq.lemma.cross+nest.mod2}
\end{subeqnarray}
\end{lemma}

\begin{figure}[t]
\centering
\begin{picture}(80,20)(12, 12)
\setlength{\unitlength}{1.5mm}
\linethickness{.5mm}
\put(-15,0){\line(1,0){50}}
\put(-15,0){\circle*{1.3}}
\put(-15,0){\makebox(0,-6)[c]{\small $\min B_1$}}
\put(-5,0){\circle*{2}}
\put(5,0){\circle*{1.3}}
\put(15,0){\circle*{2}}
\put(25,0){\circle*{1.3}}
\put(35,0){\circle*{2}}
\put(25,0){\makebox(0,-6)[c]{\small $\max B_1$}}
\blue{\qbezier(-15,0)(-5.5,10)(5,0)}
\blue{\qbezier(5,0)(15,10)(24.5,0)}
\red{\qbezier(-7,0)(5,12)(13.5,0)}
\red{\qbezier(12.5,0)(24,12)(32.5,0)}
\put(-5,0){\makebox(0,-6)[c]{\small $\min B_2$}}
\put(35,0){\makebox(0,-6)[c]{\small $\max B_2$}}
\put(5,-8){(a)}
\end{picture}
\\[3cm]
\begin{picture}(80,20)(12, 12)
\setlength{\unitlength}{1.5mm}
\linethickness{.5mm}
\put(-20,0){\line(1,0){60}}
\put(-20,0){\circle*{1.3}}
\put(-20,0){\makebox(0,-6)[c]{\small $\min B_1$}}
\put(0,0){\circle*{1.3}}
\put(20,0){\circle*{1.3}}
\put(40,0){\circle*{1.3}}
\put(40,0){\makebox(0,-6)[c]{\small $\max B_1$}}
\put(-10,0){\circle*{2}}
\put(10,0){\circle*{2}}
\put(30,0){\circle*{2}}
\put(-10,0){\makebox(0,-6)[c]{\small $\min B_2$}}
\put(30,0){\makebox(0,-6)[c]{\small $\max B_2$}}
\blue{\qbezier(-20,0)(-10.5,10)(0,0)}
\blue{\qbezier(0,0)(10,10)(20,0)}
\blue{\qbezier(18,0)(29,10)(38,0)}
\red{\qbezier(-12.5,0)(-2,12)(7.2,0)}
\red{\qbezier(6.8,0)(17,12)(26.5,0)}
\put(5,-8){(b)}
\end{picture}
\vspace*{2cm}
\caption{
   Situation in the proof of Lemma~\ref{lemma.cross+nest.mod2},
   when the pair $(B_1,B_2)$ forms (a) an overlap or (b) a covering.
 \label{fig.lemma.cross+nest.mod2}
}
\end{figure}

\proof
We first prove (\ref{eq.lemma.cross+nest.mod2}a).
Each pair of crossing arcs $(i,k)$ and $(j,l)$
must belong to a pair of distinct non-singleton blocks,
call them $B_1$ and $B_2$ where $\min B_1 < \min B_2$;
and the pair $(B_1,B_2)$ must form either an overlap
(i.e.\ $\min B_1 < \min B_2 < \max B_1 < \max B_2$)
or a covering
(i.e.\ $\min B_1 < \min B_2 < \max B_2 < \max B_1$).
So we shall consider pairs of blocks $(B_1,B_2)$ of these two types,
and for each such pair we shall count modulo 2
the number of pairs of crossing arcs between $B_1$ and $B_2$.

Suppose first that the pair $(B_1,B_2)$ forms an overlap.
Then (see Figure~\ref{fig.lemma.cross+nest.mod2}a)
each arc of $B_1$ is crossed either zero or two times by an arc of $B_2$,
except for the arc of $B_1$ that is crossed by the opener arc of $B_2$,
which is crossed only once.
So the total number of crossing pairs of $B_1$ with $B_2$ is odd.

Now suppose that the pair $(B_1,B_2)$ forms a covering.
Then (see Figure~\ref{fig.lemma.cross+nest.mod2}b)
each arc of $B_1$ is crossed either zero or two times by an arc of $B_2$,
except for the two arcs of $B_1$
that are crossed by the opener and closer arcs of $B_2$,
which are each crossed only once.
So the total number of crossing pairs of $B_1$ with $B_2$ is even.

Summing over all pairs $(B_1,B_2)$ gives (\ref{eq.lemma.cross+nest.mod2}a).

Then (\ref{eq.lemma.cross+nest.mod2}b,c,d) are an immediate consequence of
(\ref{eq.lemma.cross+nest.mod2}a) and (\ref{eq.lemma.cross+nest.ov+cov}a,b).
\qed

\subsection[Third and fourth $p,q$-generalizations: Crossings, nestings, overlaps and coverings]{Third and fourth $\bm{p,q}$-generalizations: \\ Crossings, nestings, overlaps and coverings}
   \label{subsec.intro.setpartitions.pq.3}

Let us now try to go even farther,
by introducing a ridiculously general polynomial that includes
{\em both}\/ crossing-nesting and overlap-covering statistics,
and {\em both}\/ exclusive-record and block-record statistics:
\begin{eqnarray}
   & &  \hspace*{-6mm}
  \widetilde{B}_n(x_1,x_2,y_1,y_2,v_1,v_2,y'_1,y'_2,v'_1,v'_2,p_1,p_2,q_1,q_2,p'_1,p'_2,q'_1,q'_2,r)
        \nonumber \\[3mm]
   & &
   =\;
   \sum_{\pi \in \Pi_n}
      x_1^{m_1(\pi)} x_2^{m_{\ge 2}(\pi)}
      y_1^{\erecin(\pi)} y_2^{\erecop(\pi)}
      v_1^{\nerecin(\pi)} v_2^{\nerecop(\pi)}
        \:\times
        \qquad\qquad
        \nonumber \\[-3mm]
    & & \hspace*{3.8cm}
      (y'_1)^{\brecin(\pi)} (y'_2)^{\brecop(\pi)}
      (v'_1)^{\nbrecin(\pi)} (v'_2)^{\nbrecop(\pi)}
        \:\times
        \qquad\qquad
        \nonumber \\[2mm]
    & & \hspace*{3.8cm}
      p_1^{\crin(\pi)} p_2^{\crop(\pi)}
      q_1^{\nein(\pi)} q_2^{\neop(\pi)}
        \:\times
        \qquad\qquad
        \nonumber \\[2mm]
    & & \hspace*{3.8cm}
      (p'_1)^{\ovin(\pi)} (p'_2)^{\ov(\pi)}
      (q'_1)^{\covin(\pi)} (q'_2)^{\cov(\pi)}
      \: r^{\psne(\pi)}
   \:.
 \label{def.Bn.third}
\end{eqnarray}
Of course, this polynomial is much too general to admit a J-fraction
with polynomial coefficients;
indeed, even the specialization
$y'_1 = y'_2 = v'_1 = v'_2 = 1$ and $p_1 = p_2 = q_1 = q_2 = p'_1 = q'_1 = 1$
does not admit such a J-fraction.\footnote{
   The first coefficients of the J-fraction are
   $$
      \gamma_0 \,=\, x_1 ,\quad
      \beta_1  \,=\, x_2 y_2 ,\quad
      \gamma_1 \,=\, rx_1 + y_1 ,\quad
      \beta_2  \,=\, x_2 (p'_2 y_2 + q'_2 v_2)
      \;,
   $$
   but then
   $$
      \gamma_2  \;=\;
      r^2 x_1 \:+\:
    {p'_2 (y_1 v_2 + y_2 v_1) \,+\, q'_2 (y_1 y_2 + v_1 v_2)
     \over
     p'_2 y_2 + q'_2 v_2
    }
   $$
is not a polynomial.
}
Nevertheless, there are some specializations of \reff{def.Bn.third}
with a surprisingly large number of independent variables
that do have nice J-fractions.
However, it turns out that they have nice J-fractions
because they are really just the polynomial \reff{def.Bn.crossnest.refined}
in disguise.

Before explaining how all this comes about,
let us first note some simple homogeneities that arise
as a result of the identities
\reff{eq.identity.erec.brec}--\reff{eq.lemma.cross+nest.mod2}:
\begin{itemize}
   \item[(a)]  Multiplying $y_1$ and $v_1$ by $C$
      is equivalent to multiplying $y'_1$ and $v'_1$ by $C$.
   \item[(b)]  Multiplying $y_2$ and $v_2$ by $C$
      is equivalent to multiplying $y'_2$ and $v'_2$ by $C$.
   \item[(c)]  Multiplying $p_1$ and $q_1$ by $C$
      is equivalent to multiplying $p'_1$ and $q'_1$ by $C$.
   \item[(d)]  Multiplying $p_2$ and $q_2$ by $C$
      is equivalent to multiplying $p'_2$ and $q'_2$ by $C$.
   \item[(e)]  Multiplying $p'_2$ by $\epsilon = \pm 1$
      is equivalent to multiplying $p_1$ and $p_2$ by $\epsilon$.
   \item[(f)]  Multiplying $q'_2$ by $\epsilon = \pm 1$
      is equivalent to multiplying $p_1$ and $q_2$ by $\epsilon$.
\end{itemize}
In particular, in the specializations that will be considered:
\begin{itemize}
   \item[(a${}_1$)]  Setting $y_1 = v_1 = C$ is equivalent to
      setting $y_1 = v_1 = 1$ and then multiplying $y'_1$ and~$v'_1$ by $C$.
   \item[(a${}_2$)]  Setting $y'_1 = v'_1 = C$ is equivalent to
      setting $y'_1 = v'_1 = 1$ and then multiplying $y_1$ and~$v_1$ by $C$.
\end{itemize}
Similarly,
\begin{itemize}
   \item[(b${}_1$)]  Setting $y_2 = v_2 = C$ is equivalent to
      setting $y_2 = v_2 = 1$ and then multiplying $y'_2$ and~$v'_2$ by $C$.
   \item[(b${}_2$)]  Setting $y'_2 = v'_2 = C$ is equivalent to
      setting $y'_2 = v'_2 = 1$ and then multiplying $y_2$ and~$v_2$ by $C$.
\end{itemize}
Likewise,
\begin{itemize}
   \item[(c${}_1$)]  Setting $p_1 = q_1 = C$ is equivalent to
      setting $p_1 = q_1 = 1$ and then multiplying $p'_1$ and~$q'_1$ by $C$.
   \item[(c${}_2$)]  Setting $p'_1 = q'_1 = C$ is equivalent to
      setting $p'_1 = q'_1 = 1$ and then multiplying $p_1$ and~$q_1$ by $C$.
   \item[(d${}_1$)]  Setting $p_2 = q_2 = C$ is equivalent to
      setting $p_2 = q_2 = 1$ and then multiplying $p'_2$ and~$q'_2$ by $C$.
   \item[(d${}_2$)]  Setting $p'_2 = q'_2 = C$ is equivalent to
      setting $p'_2 = q'_2 = 1$ and then multiplying $p_2$ and~$q_2$ by $C$.
\end{itemize}
And more generally,
\begin{itemize}
   \item[(e/f)]  Setting $p'_2 = \epsilon C$ and $q'_2 = \widehat{\epsilon} C$
      with $\epsilon, \widehat{\epsilon} \in \{-1,1\}$
      is equivalent to setting $p'_2 = q'_2 = 1$
      and then multiplying $p_1$ by $\epsilon \widehat{\epsilon}$,
      $p_2$ by $\epsilon C$, and $q_2$ by $\widehat{\epsilon} C$.
\end{itemize}
So, in making these specializations,
we might as well take $C=1$ to simplify the formulae,
while remembering that the results actually generalize to arbitrary $C$.

We now have:

\begin{theorem}[Four equivalent specializations of $\widetilde{B}_n$]
   \label{thm.Bnhat.third.equivalences}
The following specializations
of the polynomial
$\widetilde{B}_n(x_1,x_2,y_1,y_2,v_1,v_2,y'_1,y'_2,v'_1,v'_2,p_1,p_2,q_1,q_2,p'_1,p'_2,q'_1,q'_2,r)$
are equal:
\begin{itemize}
   \item[(i)]  $y'_1 = v'_1 = 1$, $y'_2 = v'_2 = 1$,
      $p'_1 = q'_1 = 1$, $p'_2 = q'_2 = 1$.
   \item[(ii)]  $y_1 = v_1 = 1$, $y_2 = v_2 = 1$,
      $p_1 = q_1 = 1$, $p_2 = q_2 = 1$ (and then dropping primes).
   \item[(iii)]  $y'_1 = v'_1 = 1$, $y_2 = v_2 = 1$,
      $p'_1 = q'_1 = 1$, $p_2 = q_2 = 1$ (and then dropping primes).
   \item[(iv)]  $y_1 = v_1 = 1$, $y'_2 = v'_2 = 1$,
      $p_1 = q_1 = 1$, $p'_2 = q'_2 = 1$ (and then dropping primes).
\end{itemize}
\end{theorem}

Here the specialization~(i) is just the polynomial $B_n$
defined in \reff{def.Bn.crossnest.refined},
while the specialization~(ii) is the polynomial $B^{(2)}_n$
defined in \reff{def.Bn.ovcov.brec};
and the equality of these two specializations
was already stated in Theorem~\ref{thm.setpartitions.pq.ovcov.brec}.
Now Theorem~\ref{thm.Bnhat.third.equivalences}
asserts that the specializations~(iii) and (iv) are also equal to these.
Note the logic:
we can choose to count insiders
either by crossings, nestings and exclusive records
[specializations~(i) and (iii)]
or by overlaps, coverings and block-records [(ii) and (iv)];
and we can {\em independently}\/ choose to count openers
either by crossings, nestings and exclusive records [(i) and (iv)]
or by overlaps, coverings and block-records [(ii) and (iii)].
No matter which of the four choices we make, we obtain the same polynomial.
We will prove this theorem in
Sections~\ref{subsec.setpartitions.J}--\ref{subsec.setpartitions.J.3},
by showing that all four polynomials have the same J-fraction
\reff{eq.thm.setpartitions.Jtype.pq.refined}/%
\reff{def.weights.setpartitions.Jtype.pq.refined}.
It would be interesting to find a direct bijective proof of these identities.

\subsection{First master J-fraction}
   \label{subsec.intro.setpartitions.firstmaster}

But we can go much farther, and obtain a polynomial in
four infinite families of indeterminates
$\bsfa = (\sfa_{\ell,\ell'})_{\ell,\ell' \ge 0}$,
$\bsfb = (\sfb_{\ell})_{\ell \ge 0}$,
$\bsfd = (\sfd_{\ell,\ell'})_{\ell,\ell' \ge 0}$,
$\bsfe = (\sfe_\ell)_{\ell \ge 0}$
that will have a nice J-fraction
and that will include $B_n(x_1,x_2,y_1,y_2,v_1,v_2,p_1,p_2,q_1,q_2,r)$
[defined in \reff{def.Bn.crossnest.refined}]
as a specialization.\footnote{
   In our original version of this master J-fraction,
   the weights $\bsfa,\bsfd$
   were factorized in the form
   $\sfa_{\ell,\ell'} = \sfa^{(1)}_\ell \sfa^{(2)}_{\ell'}$, etc.
   We thank Andrew Elvey Price for suggesting the generalization
   in which this factorization is avoided.
   This generalization will play a key role
   in our analysis of perfect matchings
   (see the proofs of Theorems~\ref{thm.matchings} and
    \ref{thm.matching.fourvar.pq.Stype} below).
}
The basic idea is that,
rather than counting the {\em total}\/ numbers of quadruplets
$i < j < k < l$ that form crossings or nestings,
we should instead count the number of crossings or nestings
that use a particular vertex $j$ in second (or sometimes third) position,
and then attribute weights to the vertex $j$ depending on those values.

More precisely, we define
\begin{subeqnarray}
   \crr(j,\pi)
   & = &
   \#\{ i<j<k<l \colon\: (i,k) \in \scrg_\pi \hbox{ and } (j,l) \in \scrg_\pi \}
         \\[2mm]
   \nee(j,\pi)
   & = &
   \#\{ i<j<k<l \colon\: (i,l) \in \scrg_\pi \hbox{ and } (j,k) \in \scrg_\pi \}
       \label{def.cr.ne.k.pi}
       \slabel{def.ne.k.pi}
\end{subeqnarray}
Note that $\crr(j,\pi)$ and $\nee(j,\pi)$ can be nonzero
only when $j$ is either an opener or an insider;
and summing over those $j$ gives
the refined categories of crossings and nestings:
\begin{subeqnarray}
   \crop(\pi)
   & = &
   \!\!
   \sum\limits_{j \in {\rm openers}}  \!\! \crr(j,\pi)
         \\[2mm]
   \crin(\pi)
   & = &
   \!\!
   \sum\limits_{j \in {\rm insiders}}  \!\! \crr(j,\pi)
         \\[2mm]
   \neop(\pi)
   & = &
   \!\!
   \sum\limits_{j \in {\rm openers}}  \!\! \nee(j,\pi)
         \\[2mm]
   \nein(\pi)
   & = &
   \!\!
   \sum\limits_{j \in {\rm insiders}}  \!\! \nee(j,\pi)
 \label{eq.crne.total.refined}
\end{subeqnarray}

In addition, let us define
\be
   \qne(j,\pi)  \;=\;  \#\{ i<j<l \colon\: (i,l) \in \scrg_\pi  \}
   \;;
 \label{eq.qne}
\ee
we call this a {\em quasi-nesting}\/ of the vertex $j$.
Please note that here $j$ can be a vertex of any type
(but of course it must belong to a block that is different from the one
 containing $i$ and $l$).
When $j$ is a singleton, this gives the pseudo-nestings:
\be
   \psne(\pi)   \;=\;  \sum\limits_{j \in {\rm singletons}}  \qne(j,\pi)
   \;.
\ee
When $j$ is an opener or an insider, we have simply
\be
   \qne(j,\pi)  \;=\;  \crr(j,\pi)  \,+\, \nee(j,\pi)
   \;,
 \label{eq.qne.crr+nee}
\ee
so no new information is obtained.
And finally, when $j$ is a closer,
$\qne(j,\pi)$ counts the number of times that the closer $j$
occurs in {\em third}\/ position in a crossing or nesting:
when $(i,l) \in \scrg_\pi$ is a pair contributing to $\qne(j,\pi)$,
and $(h,j) \in \scrg_\pi$,
then we have either
$h < i < j < l$ (so that the quadruplet is a crossing)
or
$i < h < j < l$ (so that the quadruplet is a nesting),
but we do not keep track of which one it is.

We now introduce four infinite families of indeterminates
$\bsfa = (\sfa_{\ell,\ell'})_{\ell,\ell' \ge 0}$,
$\bsfb = (\sfb_{\ell})_{\ell \ge 0}$,
$\bsfd = (\sfd_{\ell,\ell'})_{\ell,\ell' \ge 0}$,
$\bsfe = (\sfe_\ell)_{\ell \ge 0}$
and define the polynomials
$B_n(\bsfa,\bsfb,\bsfd,\bsfe)$ by
\begin{eqnarray}
   & &  \hspace*{-7mm}
   B_n(\bsfa,\bsfb,\bsfd,\bsfe)
   \;=\;
       \nonumber \\[4mm]
   & & 
   \sum_{\pi \in \Pi_n}
   \;
   \prod\limits_{i \in {\rm openers}}  \!\! \sfa_{\crr(i,\pi),\, \nee(i,\pi)}
   \prod\limits_{i \in {\rm closers}}  \! \sfb_{\qne(i,\pi)}
   \prod\limits_{i \in {\rm insiders}} \!\!  \sfd_{\crr(i,\pi),\, \nee(i,\pi)}
   \prod\limits_{i \in {\rm singletons}}  \!\!\! \sfe_{\qne(i,\pi)}
   \;.
        \nonumber \\[-2mm]
 \label{def.Bnhat}
\end{eqnarray}
These polynomials then have a beautiful J-fraction:

\begin{theorem}[First master J-fraction for set partitions]
   \label{thm.setpartitions.Jtype.final1}
The ordinary generating function of the polynomials
$B_n(\bsfa,\bsfb,\bsfd,\bsfe)$
has the J-type continued fraction
\begin{eqnarray}
   & & \hspace*{-12mm}
   \sum_{n=0}^\infty B_n(\bsfa,\bsfb,\bsfd,\bsfe) \: t^n
   \;=\;
       \nonumber \\
   & &
   \cfrac{1}{1 - \sfe_0 t - \cfrac{\sfa_{00} \sfb_{0} t^2}{1 -  (\sfd_{00} + \sfe_1) t - \cfrac{(\sfa_{01} + \sfa_{10}) \sfb_{1} t^2}{1 - (\sfd_{01} + \sfd_{10} + \sfe_2)t - \cfrac{(\sfa_{02} + \sfa_{11} + \sfa_{20}) \sfb_{2} t^2}{1 - \cdots}}}}
       \nonumber \\[1mm]
   \label{eq.thm.setpartitions.Jtype.final1}
\end{eqnarray}
with coefficients
\begin{subeqnarray}
   \gamma_n  & = &   \sfd^\star_{n-1} \,+\, \sfe_n
         \\[1mm]
   \beta_n   & = &   \sfa^\star_{n-1} \, \sfb_{n-1}
 \label{def.weights.setpartitions.Jtype.final1}
\end{subeqnarray}
where
\be
   \sfa^\star_{n-1}  \;\eqdef\;  \sum_{\ell=0}^{n-1} \sfa_{\ell,n-1-\ell}
 \label{def.astar.setpartitions}
\ee
and likewise for $\sfd$.
\end{theorem}

\noindent
We will prove this theorem in Section~\ref{subsec.setpartitions.J}.


\medskip

{\bf Remarks.}
1.  It seems far from obvious (at least to us)
why $B_n(\bsfa,\bsfb,\bsfd,\bsfe)$ depends on
$\bsfa,\bsfb,\bsfd,\bsfe$ only via the
combinations (\ref{def.weights.setpartitions.Jtype.final1}a,b).
Of course, this behavior is a consequence of
the bijection onto labeled Motzkin paths
that we will use in Section~\ref{subsec.setpartitions.J}
to prove Theorem~\ref{thm.setpartitions.Jtype.final1}.
But it would be interesting to understand it combinatorially,
directly at the level of set partitions.

2. It is unfortunate that the polynomial \reff{def.Bnhat}
treats openers and closers asymmetrically,
but we do not see any way to avoid this.
One can, of course, interchange the roles of openers and closers
by passing to the reversed partition;
indeed, this reversal will be employed, for technical reasons,
in our proof in Section~\ref{subsec.setpartitions.J}.
But, whichever way one does it,
one is left with a polynomial that uses
doubly-indexed indeterminates $\sfa_{\ell,\ell'}$ for one class
and singly-indexed indeterminates $\sfb_{\ell}$ for the other.
We do not see any way to obtain a continued fraction
for a polynomial with two doubly-indexed indeterminates.
(But perhaps we are missing something.)
\myendremark

\medskip

Let us now show how to recover
$B_n(x_1,x_2,y_1,y_2,v_1,v_2,p_1,p_2,q_1,q_2,r)$
as a specialization of $B_n(\bsfa,\bsfb,\bsfd,\bsfe)$,
and thus obtain Theorem~\ref{thm.setpartitions.Jtype.pq.refined}
(and hence also Theorem~\ref{thm.setpartitions.Jtype})
as a special case of Theorem~\ref{thm.setpartitions.Jtype.final1}.
The needed specialization is
\begin{subeqnarray}
   \sfa_{\ell,\ell'}
   & = &
   p_2^\ell \:\times\:
   \begin{cases}
       y_2             &  \textrm{if $\ell' = 0$}  \\[1mm]
       q_2^{\ell'} v_2   &  \textrm{if $\ell' \ge 1$}
   \end{cases}
       \\[1mm]
   \sfb_{\ell}   & = &  x_2
       \\[1mm]
   \sfd_{\ell,\ell'}
   & = &
   p_1^\ell \:\times\:
   \begin{cases}
       y_1             &  \textrm{if $\ell' = 0$}  \\[1mm]
       q_1^{\ell'} v_1   &  \textrm{if $\ell' \ge 1$}
   \end{cases}
       \\[1mm]
   \sfe_\ell   & = &  r^\ell \, x_1
 \label{eq.specialization.thm.setpartitions.Jtype.pq}
\end{subeqnarray}
Most of this is obvious:  singletons get a weight $\sfe_\ell = r^\ell x_1$;
we count blocks of size~$\ge 2$ at their closers,
hence $\sfb_{\ell} = x_2$;
and by \reff{eq.crne.total.refined} we count crossings and nestings
at openers and insiders, which explains
the factors $p_2^\ell$ and $q_2^{\ell'}$ in $\sfa_{\ell,\ell'}$,
and $p_1^\ell$~and $q_1^{\ell'}$ in $\sfd_{\ell,\ell'}$.
Finally, recall that $j$ is an exclusive record of $\pi$
if it is not the largest element of its block
(i.e.\ it is either an opener or an insider)
and the next element in its block sticks out farther to the right
than any right neighbor (within its block) of a vertex $< j$.
In other words, $j$ is an exclusive record if and only if
it is an opener or insider
and is not the second element in a nesting, i.e.\ $\nee(j,\pi) = 0$.
This explains the factor $y_2$ in $\sfa_{\ell,\ell'}$ when $\ell' = 0$,
and the factor $v_2$ when $\ell' \ge 1$;
and likewise the factors $y_1$ and $v_1$ in $\sfd_{\ell,\ell'}$.
%
This completes the proof that
$B_n(x_1,x_2,y_1,y_2,v_1,v_2,p_1,p_2,q_1,q_2,r)$
is obtained from the specialization
\reff{eq.specialization.thm.setpartitions.Jtype.pq}
of $B_n(\bsfa,\bsfb,\bsfd,\bsfe)$.
Inserting this specialization into
\reff{def.weights.setpartitions.Jtype.final1}/\reff{def.astar.setpartitions}
yields the weights \reff{def.weights.setpartitions.Jtype.pq.refined}.


\bigskip

{\bf Remarks.}
1.  The definitions \reff{def.cr.ne.k.pi}--\reff{def.Bnhat}
can be motivated by following closely the analogy
with the first master J-fraction for permutations
(Section~\ref{subsec.intro.permutations.firstmaster}).
We begin by mapping set partitions into permutations as follows:
Given a set partition $\pi \in \Pi_n$,
we define the permutation $\sigma \in \Sym_n$
such that the disjoint cycles of~$\sigma$
are the blocks of $\pi$, each traversed in increasing order
(with the largest element of~course followed by the smallest element).
The mapping $\pi \mapsto \sigma$
is clearly a bijection of $\Pi_n$ onto $\Sym_n^\star$,
where $\Sym_n^\star$ denotes the subset of $\Sym_n$ consisting of permutations
in which each cycle of length $\ell \ge 2$ contains
precisely one cycle peak (namely, the cycle maximum),
one cycle valley (namely, the cycle minimum),
$\ell-2$ cycle double rises, and no cycle double falls.
In particular, openers, closers, insiders and singletons of $\pi \in \Pi_n$
map, respectively, into
cycle valleys, cycle peaks, cycle double rises and fixed points
of $\sigma \in \Sym_n^\star$.

We now define set-partition statistics
that are simply the images of the permutation statistics
$\ucross$, $\unest$, $\lcross$, $\lnest$ and $\lev$
[cf.\ \reff{def.ucrossnestjk}/\reff{def.level.bis}]
under this mapping.

The images of $\ucross(j,\sigma)$ and $\unest(j,\sigma)$
are precisely $\crr(j,\pi)$ and $\nee(j,\pi)$
as defined in \reff{def.cr.ne.k.pi}.
This explains why in the definitions \reff{def.cr.ne.k.pi}
we have put the distinguished index in {\em second}\/ position.
More specifically,
the images of $\ucrosscval$, $\ucrosscdrise$, $\unestcval$ and $\unestcdrise$
are $\crop$, $\crin$, $\neop$ and $\nein$, respectively.

As for $\lcross$ and $\lnest$,
we observe that the only lower arcs in the permutation~$\sigma$
are those that map the closer of a non-singleton block to that block's opener;
therefore, the images of $\lcross(k,\sigma)$ and $\lnest(k,\sigma)$ are
\begin{subeqnarray}
    \ovbar(k,\pi)
    & = &
    \#\{ (B,B') \colon\:  \min B < \min B' < k = \max B < \max B'  \}
        \qquad\quad \\[2mm]
    \covbar(k,\pi)
    & = &
    \#\{ (B,B') \colon\:  \min B' < \min B < k = \max B < \max B'  \}
        \qquad\quad
  \label{def.ov.cov.kpi.0}
 \end{subeqnarray}
which can be nonzero only when $k$ is a closer
(of a block $B$ of size~$\ge 2$).
And then
\begin{subeqnarray}
   \ov(\pi)
   & = &
   \sum\limits_{k \in {\rm closers}}  \!\! \ovbar(k,\pi)
        \\[2mm]
   \cov(\pi)
   & = &
   \sum\limits_{k \in {\rm closers}}  \!\! \covbar(k,\pi)
 \label{eq.ovbar.covbar.closers}
\end{subeqnarray}
since we can count the block $B$ at its closer.
So the analogues of $\lcross$ and $\lnest$ are overlaps and coverings.
(Of course, we here have only $\lcrosscpeak$ and $\lnestcpeak$,
since $\lcrosscdfall = \lnestcdfall = 0$ for $\sigma \in \Sym_n^\star$.)

Finally, the image of $\lev(j,\sigma)$
can be defined in two equivalent ways,
corresponding to upper pseudo-nestings and lower pseudo-nestings.
If $j$ is a singleton, we define
\be
   \psne(j,\pi)  \;=\;  \#\{ i<j<k \colon\: (i,k) \in \scrg_\pi  \}
   \;.
 \label{eq.psne}
\ee
Summing over singletons $j$ gives the total number of pseudo-nestings:
\be
   \psne(\pi)   \;=\;  \sum\limits_{j \in {\rm singletons}}  \psne(j,\pi)
   \;.
\ee
On the other hand, if $j$ is a singleton we also define
\be
   \pscov(j,\pi)  \;=\;  \#\{ B \colon\: \min B < j < \max B \}
   \;.
 \label{eq.pscov}
\ee
Summing over singletons $j$ gives the total number of pseudo-coverings:
\be
   \pscov(\pi)   \;=\;  \sum\limits_{j \in {\rm singletons}}  \pscov(j,\pi)
   \;.
\ee
But it is easy to see that $\psne(j,\pi) = \pscov(j,\pi)$,
since for each block $B$ satisfying $\min B < j < \max B$,
there is precisely one edge $(i,k) \in \scrg_\pi$ with $i,k \in B$
that satisfies $i<j<k$.
So these are simply two different names for the same object,
which we have here called $\qne(j,\pi)$ [cf.\ \reff{eq.qne}],
specialized now to the case in which $j$ is a singleton.

The polynomial \reff{def.Bnhat} is then obtained from the
corresponding permutation polynomial \reff{def.Qn.firstmaster}
by specializing to $\bsfc = \bzero$
(corresponding to the absence of cycle double falls
 for $\sigma \in \Sym_n^\star$)
and then further specializing $\bsfb$
to be independent of $\lcross$ and $\lnest$.

Unfortunately, we are unable to employ this mapping
to permutations to {\em prove}\/ Theorem~\ref{thm.setpartitions.Jtype.final1}.
The trouble is that it does not seem possible to encode
the other property defining the subset $\Sym_n^\star \subseteq \Sym_n$
--- namely, that each cycle of length $\ell \ge 2$ contains
precisely one cycle peak and one cycle valley  ---
in the polynomial \reff{def.Qn.firstmaster}.
Perhaps this problem could be alleviated
by using instead the second master J-fraction
for permutations [cf.\ \reff{def.Qn.secondmaster}],
which includes the statistic $\cyc$:
for we could substitute $\sfa \to \lambda^{-1} \sfa$
and $\sfe \to \lambda^{-1} \sfe$ and take $\lambda \to\infty$,
thereby forcing $\sigma \in \Sym_n^\star$.
But then we would only be able to handle $\ucross$ and $\unest$
via their sum, which would amount to specializing (for~instance)
$v_1 = y_1$ and $v_2 = y_2$ in Theorem~\ref{thm.setpartitions.Jtype},
which would be a severe limitation.
Or perhaps Theorem~\ref{thm.setpartitions.Jtype.final1}
could be proven by finding a {\em different}\/ bijection
of set~partitions onto a subclass of permutations.
But we have been unable (thus far) to find a suitable bijection,
so we are instead obliged to prove Theorem~\ref{thm.setpartitions.Jtype.final1}
by a direct argument on set partitions
(see Section~\ref{subsec.setpartitions.J}).

2.  Instead of (\ref{eq.specialization.thm.setpartitions.Jtype.pq}b)
we could take, more generally, $\sfb_\ell = (q')^\ell x_2$.
But it is easy to see that if $k$ is a closer, then
\be
   \qne(k,\pi)  \;=\;  \ovbar(k,\pi) + \covbar(k,\pi)
\ee
[cf.\ \reff{eq.qne} and \reff{def.ov.cov.kpi.0}], so that
\be
   \sum\limits_{k \in {\rm closers}}  \!\! \qne(k,\pi)
   \;=\;
   \ov(\pi) + \cov(\pi)
   \;=\;
   \crop(\pi) + \neop(\pi)
\ee
[cf.\ (\ref{eq.lemma.cross+nest.ov+cov}a) and \reff{eq.ovbar.covbar.closers}].
So taking $\sfb_\ell = (q')^\ell x_2$ is equivalent to
taking $\sfb_\ell = x_2$ and then multiplying $p_2$ and $q_2$ by $q'$.
This can alternatively be seen by observing that,
according to (\ref{def.weights.setpartitions.Jtype.final1}b),
multiplying $\sfb_\ell$ by $(q')^\ell$ has the same effect on the J-fraction
as multiplying $\sfa_{\ell,\ell'}$ by $(q')^{\ell+\ell'}$.
\myendremark

\subsection{First master S-fraction}
   \label{subsec.intro.setpartitions.firstmaster.S}

We can also obtain a master S-fraction by specializing the parameters
in Theorem~\ref{thm.setpartitions.Jtype.final1}.
Indeed, the J-fraction
\reff{eq.thm.setpartitions.Jtype.final1}/\reff{def.weights.setpartitions.Jtype.final1}
is the contraction \reff{eq.contraction_even.coeffs} of the S-fraction
\be
   \cfrac{1}{1 - \cfrac{\sfb_0 t}{1 - \cfrac{\sfa_{00} t}{1 - \cfrac{\sfb_1 t}{1 - \cfrac{(\sfa_{01} + \sfa_{10}) t}{1 - \cdots}}}}}
\ee
with coefficients $\alpha_{2k-1} = \sfb_{k-1}$
and $\alpha_{2k} = \sfa^\star_{k-1}$
if we choose $\bsfd,\bsfe$ so that
\be
   \sfd^\star_{n-1} \,+\, \sfe_n
   \;=\;
   \sfa^\star_{n-1} \,+\, \sfb_n
   \quad\hbox{for all $n \ge 0$}
   \;.
\ee
There are many ways of doing this;
the simplest is to set $\sfe_n = \sfb_n$ for all $n \ge 0$
and then choose $\bsfa,\bsfd$ in any way such that
\be
   \sfd^\star_{n-1}
   \;=\;
   \sfa^\star_{n-1}
   \quad\hbox{for all $n \ge 1$}
   \;.
\ee
Even this latter choice can be done in many ways;
the simplest is to choose $\bsfa$ freely
and then set $\bsfd = \bsfa$.
These choices lead to the following result:

\begin{theorem}[Master S-fraction for set partitions]
   \label{thm.setpartitions.Stype.final1}
The ordinary generating function of the polynomials
$B_n(\bsfa,\bsfb,\bsfa,\bsfb)$
has the S-type continued fraction
\be
   \sum_{n=0}^\infty B_n(\bsfa,\bsfb,\bsfa,\bsfb) \: t^n
   \;\:=\;\:
   \cfrac{1}{1 - \cfrac{\sfb_0 t}{1 - \cfrac{\sfa_{00} t}{1 - \cfrac{\sfb_1 t}{1 - \cfrac{(\sfa_{01} + \sfa_{10}) t}{1 - \cdots}}}}}
   \label{eq.thm.setpartitions.Stype.final1}
\ee
with coefficients
\begin{subeqnarray}
   \alpha_{2k-1}  & = &  \sfb_{k-1}  \\
   \alpha_{2k}    & = &  \sfa^\star_{k-1}
 \label{def.weights.setpartitions.Stype.final1}
\end{subeqnarray}
where
$\displaystyle
   \sfa^\star_{n-1}  \;\eqdef\;  \sum_{\ell=0}^{n-1} \sfa_{\ell,n-1-\ell}
$.
\end{theorem}

To obtain the S-fraction \reff{eq.cor.setpartitions.Jtype.pq}
from Theorem~\ref{thm.setpartitions.Stype.final1},
we make the specializations
\begin{subeqnarray}
   \sfa_{\ell,\ell'}
   & = &
   p^\ell \:\times\:
   \begin{cases}
       y             &  \textrm{if $\ell' = 0$}  \\[1mm]
       q^{\ell'} v   &  \textrm{if $\ell' \ge 1$}
   \end{cases}
       \\[1mm]
   \sfb_\ell   & = &  x
 \label{eq.specialization.thm.setpartitions.Stype.pq}
\end{subeqnarray}

\subsection{Second master J-fraction}
   \label{subsec.intro.setpartitions.secondmaster}

Let us now define a second master J-fraction,
following the same scheme as in
Section~\ref{subsec.intro.setpartitions.firstmaster}
but now using overlaps and coverings
(as defined in Section~\ref{subsec.intro.setpartitions.pq.2.0})
in place of crossings and nestings.

We begin by defining
\begin{subeqnarray}
   \ov(j,\pi)
   & = &
   \#\{ (B_1,B_2) \colon\:  j \in B_2 \,\hbox{ and }\,
                          \min B_1 < j < \max B_1 < \max B_2  \}
         \nonumber \\ \\[2mm]
   \cov(j,\pi)
   & = &
   \#\{ (B_1,B_2) \colon\:  j \in B_2 \,\hbox{ and }\,
                          \min B_1 < j < \max B_2 < \max B_1  \}
         \nonumber \\
 \label{def.ov.cov.kpi}
\end{subeqnarray}
Note that $\ov(j,\pi)$ and $\cov(j,\pi)$ can be nonzero
only when $j$ is either an opener or an insider
(since $j \in B_2$ and $j < \max B_2$).
If we sum over openers, then each block~$B_2$ gets counted once
(namely, with $j = \min B_2$),
and we obtain the total numbers of overlaps and coverings:
\begin{subeqnarray}
   \ov(\pi)
   & = &
   \!\!
   \sum\limits_{j \in {\rm openers}}  \!\! \ov(j,\pi)
         \\[2mm]
   \cov(\pi)
   & = &
   \!\!
   \sum\limits_{j \in {\rm openers}}  \!\! \cov(j,\pi)
 \label{eq.ovcov.total}
\end{subeqnarray}
On the other hand, if we sum over insiders,
then we obtain the quantities $\ovin$ and $\covin$
defined in \reff{def.ovin.covin}:
\begin{subeqnarray}
   \ovin(\pi)
   & = &
   \!\!
   \sum\limits_{j \in {\rm insiders}}  \!\! \ov(j,\pi)
         \\[2mm]
   \covin(\pi)
   & = &
   \!\!
   \sum\limits_{j \in {\rm insiders}}  \!\! \cov(j,\pi)
 \label{eq.ovin.covin.total}
\end{subeqnarray}

\medskip

{\bf Remark.}
Compare \reff{eq.crne.total.refined}
with \reff{eq.ovcov.total}/\reff{eq.ovin.covin.total}:
we see that the total numbers of overlaps and coverings
are {\em not}\/ analogous
to the total numbers of crossings and nestings;
rather, they are analogous 
to the total numbers of crossings and nestings {\em of~opener type}\/.
\myendremark

\medskip

In addition, let us define
\be
   \qcov(j,\pi)  \;=\;  \#\{ B \colon\:  B \not\ni j \hbox{ and }
                                         \min B < j < \max B \}
   \;.
 \label{eq.qcov}
\ee
We call this a {\em quasi-covering}\/ of the vertex $j$;
please note that here $j$ can be a vertex of any type.
When $j$ is a singleton, this gives the pseudo-coverings:
\be
   \sum\limits_{j \in {\rm singletons}}  \qcov(j,\pi)
   \;=\;
   \pscov(\pi)
   \;.
 \label{eq.qcov.singletons}
\ee
When $j$ is an opener or an insider, we have simply
\be
   \qcov(j,\pi)  \;=\;  \ov(j,\pi)  \,+\, \cov(j,\pi)
   \;,
 \label{eq.qcov.ov+cov}
\ee
so that in particular
\begin{subeqnarray}
   \sum\limits_{j \in {\rm openers}}  \qcov(j,\pi)
   & = &
   \ov(\pi)  \,+\, \cov(\pi)
       \\[2mm]
   \sum\limits_{j \in {\rm insiders}}  \qcov(j,\pi)
   & = &
   \ovin(\pi)  \,+\, \covin(\pi)
       \slabel{eq.ovin+covin}
\end{subeqnarray}
And finally, when $j$ is a closer,
$\qcov(j,\pi)$ counts the number of blocks $B_1$
that either overlap or cover the block $B_2$ whose closer is $j$;
in particular, we have
\be
   \sum\limits_{j \in {\rm closers}}  \qcov(j,\pi)
   \;=\;
   \ov(\pi)  \,+\, \cov(\pi)
   \;.
 \label{eq.qcovclosers.ov+cov}
\ee

But the quantity $\qcov(j,\pi)$ has already been introduced
under a different name: we have in fact
\be
   \qcov(j,\pi)  \;=\;  \qne(j,\pi)
 \label{eq.qcov=qne}
\ee
[cf.\ \reff{eq.qne}],
no matter what is the type of the vertex $j$.
And this is easy to see:
if a vertex $j$ and a block $B$ satisfy
$j \notin B$ and $\min B < j < \max B$,
then there is precisely one edge $(i,l) \in \scrg_\pi$ with $i,l \in B$
that satisfies $i<j<l$.
So quasi-nestings and quasi-coverings are just two different names
for the same quantity.
For future reference, let us record the relevant facts about this quantity:

\begin{lemma}[Quasi-nestings of openers and insiders]
   \label{lemma.qne}
\nopagebreak
\quad\hfill
\vspace*{-1mm}
\begin{itemize}
   \item[(a)]  If $j$ is an opener or an insider, we have
\be
   \qne(j,\pi)
   \;=\;
   \qcov(j,\pi)
   \;=\; 
   \crr(j,\pi)  + \nee(j,\pi)
   \;=\;
   \ov(j,\pi)  + \cov(j,\pi)
   \;.
   \quad
\ee
   \item[(b)]  We have
\begin{subeqnarray}
   & &  \hspace*{-15mm}
   \sum\limits_{j \in {\rm openers}}  \!\! \qne(j,\pi)
   \;=\;
   \!\!
   \sum\limits_{j \in {\rm openers}}  \!\! \qcov(j,\pi)
   \;=\;
   \crop(\pi) + \neop(\pi)
   \;=\;
   \ov(\pi) + \cov(\pi)
         \nonumber \\[-4mm] \\[2mm]
   & &  \hspace*{-15mm}
   \sum\limits_{j \in {\rm insiders}}  \!\! \qne(j,\pi)
   \;=\;
   \!\!
   \sum\limits_{j \in {\rm insiders}}  \!\! \qcov(j,\pi)
   \;=\;
   \crin(\pi) + \nein(\pi)
   \;=\;
   \ovin(\pi) + \covin(\pi)
         \nonumber \\[-4mm] \\[2mm]
   & &  \hspace*{-15mm}
   \sum\limits_{j \in {\rm closers}}  \!\! \qne(j,\pi)
   \;=\;
   \!\!
   \sum\limits_{j \in {\rm closers}}  \!\! \qcov(j,\pi)
   \;=\;
   \crop(\pi) + \neop(\pi)
   \;=\;
   \ov(\pi) + \cov(\pi)
         \nonumber \\[-4mm] \\[2mm]
   & &  \hspace*{-15mm}
   \sum\limits_{j \in {\rm singletons}}  \!\! \qne(j,\pi)
   \;=\;
   \!\!
   \sum\limits_{j \in {\rm singletons}}  \!\! \qcov(j,\pi)
   \;=\;
   \psne(\pi) + \pscov(\pi)
         \nonumber \\[-4mm]
\end{subeqnarray}
\end{itemize}
\end{lemma}

\noindent
Here we have simply recalled \reff{eq.crne.total.refined},
\reff{eq.qne.crr+nee}, \reff{eq.ovcov.total}, \reff{eq.ovin.covin.total},
\reff{eq.qcov.ov+cov}, \reff{eq.qcovclosers.ov+cov} and \reff{eq.qcov=qne}.
Note that Lemma~\ref{lemma.qne}(b) refines Lemma~\ref{lemma.cross+nest.ov+cov},
while Lemma~\ref{lemma.qne}(a) is a further refinement.

Let us now introduce four infinite families of indeterminates
$\bsfa = (\sfa_{\ell,\ell'})_{\ell,\ell' \ge 0}$,
$\bsfb = (\sfb_{\ell})_{\ell \ge 0}$,
$\bsfd = (\sfd_{\ell,\ell'})_{\ell,\ell' \ge 0}$,
$\bsfe = (\sfe_\ell)_{\ell \ge 0}$
and define the polynomials
$B^{(2)}_n(\bsfa,\bsfb,\bsfd,\bsfe)$ by
\begin{eqnarray}
   & &  \hspace*{-8mm}
   B^{(2)}_n(\bsfa,\bsfb,\bsfd,\bsfe)
   \;=\;
       \nonumber \\[4mm]
   & & 
   \sum_{\pi \in \Pi_n}
   \;
   \prod\limits_{i \in {\rm openers}}  \!\! \sfa_{\ov(i,\pi),\, \cov(i,\pi)}
   \prod\limits_{i \in {\rm closers}}  \! \sfb_{\qne(i,\pi)}
   \prod\limits_{i \in {\rm insiders}} \!\!  \sfd_{\ov(i,\pi),\, \cov(i,\pi)}
   \prod\limits_{i \in {\rm singletons}}  \!\!\! \sfe_{\qne(i,\pi)}
   \;.
        \nonumber \\[-2mm]
 \label{def.Bnhat2}
\end{eqnarray}
But it turns out that these polynomials
$B^{(2)}_n(\bsfa,\bsfb,\bsfd,\bsfe)$
are not simply {\em analogues}\/
of the polynomials $B_n(\bsfa,\bsfb,\bsfd,\bsfe)$
defined in Section~\ref{subsec.intro.setpartitions.firstmaster};
they {\em are}\/ the polynomials $B_n(\bsfa,\bsfb,\bsfd,\bsfe)$
in disguise:

\begin{theorem}[Second master J-fraction for set partitions]
   \label{thm.setpartitions.Jtype.final2}
The polynomials
\linebreak
$B^{(2)}_n(\bsfa,\bsfb,\bsfd,\bsfe)$
defined in \reff{def.Bnhat2}
are identical to the polynomials $B_n(\bsfa,\bsfb,\bsfd,\bsfe)$
defined in \reff{def.Bnhat}.
In particular, their ordinary generating function
has the same J-fraction
\reff{eq.thm.setpartitions.Jtype.final1}/%
\reff{def.weights.setpartitions.Jtype.final1}.
\end{theorem}

\smallskip
\noindent
Indeed, we will prove this theorem in Section~\ref{subsec.setpartitions.J.2}
by showing that the polynomials $B^{(2)}_n(\bsfa,\bsfb,\bsfd,\bsfe)$
have the J-fraction
\reff{eq.thm.setpartitions.Jtype.final1}/%
\reff{def.weights.setpartitions.Jtype.final1}.
It is an interesting open problem to find a natural bijection $\Pi_n \to \Pi_n$
that proves $B_n(\bsfa,\bsfb,\bsfd,\bsfe) = B^{(2)}_n(\bsfa,\bsfb,\bsfd,\bsfe)$
by mapping directly the relevant statistics.
(One possibility, of course, is to use the bijection obtained
 by composing the two bijections to labeled Motzkin paths
 that will be constructed in
 Sections~\ref{subsec.setpartitions.J} and \ref{subsec.setpartitions.J.2}.
 But we do not know, at present, how to give a
 {\em simple}\/ and {\em explicit}\/ definition of this bijection.)

\bigskip

Let us now show how to recover
$B^{(2)}_n(x_1,x_2,y_1,y_2,v_1,v_2,p_1,p_2,q_1,q_2,r)$
defined in \reff{def.Bn.ovcov.brec}
as a specialization
of $B^{(2)}_n(\bsfa,\bsfb,\bsfd,\bsfe)$,
and thus obtain the J-fraction
corresponding to Theorem~\ref{thm.setpartitions.pq.ovcov.brec}
as a special case of Theorem~\ref{thm.setpartitions.Jtype.final2}.
The needed specialization is
precisely \reff{eq.specialization.thm.setpartitions.Jtype.pq},
and the reasoning is very similar to that used in deriving
\reff{eq.specialization.thm.setpartitions.Jtype.pq}.
The only difference is that we now use overlaps and coverings
in place of crossings and nestings,
and block-records in place of exclusive records.
We have defined block-records in such a way that
$j$ is a block-record if and only if it is an opener or insider
and $\cov(j,\pi) = 0$;
so the reasoning used in deriving
\reff{eq.specialization.thm.setpartitions.Jtype.pq}
applies verbatim, with nestings replaced by coverings.

\subsection{Third and fourth master J-fractions}
   \label{subsec.intro.setpartitions.thirdmaster}

We now introduce some polynomials that mix the statistics
that were used in the first and second master J-fractions
(Sections~\ref{subsec.intro.setpartitions.firstmaster}
 and \ref{subsec.intro.setpartitions.secondmaster}).
So introduce indeterminates
$\bsfa = (\sfa_{\ell,\ell'})_{\ell,\ell' \ge 0}$,
$\bsfb = (\sfb_{\ell})_{\ell \ge 0}$,
$\bsfd = (\sfd_{\ell,\ell'})_{\ell,\ell' \ge 0}$,
$\bsfe = (\sfe_\ell)_{\ell \ge 0}$
as before,
and define the polynomials
$B^{(3)}_n(\bsfa,\bsfb,\bsfd,\bsfe)$
and $B^{(4)}_n(\bsfa,\bsfb,\bsfd,\bsfe)$ by
\begin{eqnarray}
   & &  \hspace*{-8mm}
   B^{(3)}_n(\bsfa,\bsfb,\bsfd,\bsfe)
   \;=\;
       \nonumber \\[4mm]
   & & 
   \sum_{\pi \in \Pi_n}
   \;
   \prod\limits_{i \in {\rm openers}}  \!\! \sfa_{\ov(i,\pi),\, \cov(i,\pi)}
   \prod\limits_{i \in {\rm closers}}  \! \sfb_{\qne(i,\pi)}
   \prod\limits_{i \in {\rm insiders}} \!\!  \sfd_{\crr(i,\pi),\, \nee(i,\pi)}
   \prod\limits_{i \in {\rm singletons}}  \!\!\! \sfe_{\qne(i,\pi)}
        \nonumber \\[-2mm]
 \label{def.Bnhat3}
 \\[2mm]
   & &  \hspace*{-8mm}
   B^{(4)}_n(\bsfa,\bsfb,\bsfd,\bsfe)
   \;=\;
       \nonumber \\[4mm]
   & & 
   \sum_{\pi \in \Pi_n}
   \;
   \prod\limits_{i \in {\rm openers}}  \!\! \sfa_{\crr(i,\pi),\, \nee(i,\pi)}
   \prod\limits_{i \in {\rm closers}}  \! \sfb_{\qne(i,\pi)}
   \prod\limits_{i \in {\rm insiders}} \!\!  \sfd_{\ov(i,\pi),\, \cov(i,\pi)}
   \prod\limits_{i \in {\rm singletons}}  \!\!\! \sfe_{\qne(i,\pi)}
        \nonumber \\[-2mm]
 \label{def.Bnhat4}
\end{eqnarray}
So $B^{(3)}$ employs overlaps and coverings for openers,
and crossings and nestings for insiders, while $B^{(4)}$ does the reverse.

It turns out that the polynomials $B^{(3)}$ and $B^{(4)}$ are,
like $B^{(2)}$, identical to the polynomials $B_n$
defined in Section~\ref{subsec.intro.setpartitions.firstmaster}:

\begin{theorem}[Third and fourth master J-fractions for set partitions]
   \label{thm.setpartitions.Jtype.final34}
The polynomials
$B^{(3)}_n(\bsfa,\bsfb,\bsfd,\bsfe)$ and $B^{(4)}_n(\bsfa,\bsfb,\bsfd,\bsfe)$
defined in \reff{def.Bnhat3}/\reff{def.Bnhat4}
are identical to the polynomials $B_n(\bsfa,\bsfb,\bsfd,\bsfe)$
defined in \reff{def.Bnhat}.
In particular, their ordinary generating function
has the same J-fraction
\reff{eq.thm.setpartitions.Jtype.final1}/%
\reff{def.weights.setpartitions.Jtype.final1}.
\end{theorem}

\smallskip
\noindent
Indeed, we will prove this theorem in Section~\ref{subsec.setpartitions.J.3}
by showing that the polynomials
$B^{(3)}_n(\bsfa,\bsfb,\bsfd,\bsfe)$ and $B^{(4)}_n(\bsfa,\bsfb,\bsfd,\bsfe)$
have the J-fraction
\reff{eq.thm.setpartitions.Jtype.final1}/%
\reff{def.weights.setpartitions.Jtype.final1}.

The polynomials defined in Theorem~\ref{thm.Bnhat.third.equivalences}(iii,iv)
are then obtained from $B^{(3,4)}_n(\bsfa,\bsfb,\bsfd,\bsfe)$
by the same specialization \reff{eq.specialization.thm.setpartitions.Jtype.pq}
that was used for $B_n$ and $B_n^{(2)}$
in Sections~\ref{subsec.intro.setpartitions.firstmaster}
and \ref{subsec.intro.setpartitions.secondmaster}, respectively;
the reasoning is identical to that used there.

\subsection{A remark on the Wachs--White statistics and inversion statistics}
   \label{subsec.setpartitions.inv}

Wachs and White \cite{Wachs_91} have introduced four statistics
on set partitions, which can be defined as follows:\footnote{
   Wachs and White \cite{Wachs_91} actually defined their statistics
   on words ${\bf w} = w_1 \cdots w_n \in [k]^n$.
   Now, words ${\bf w} \in [k]^n$ satisfying
   $w_i \le \max\limits_{1 \le j < i} w_j + 1$
   and $\max\limits_{1 \le i \le n} w_i = k$
   --- termed {\em restricted growth functions}\/
   of length $n$ and maximum $k$ ---
   are in bijection with partitions of $[n]$ with $k$ blocks:
   we write $\pi = \{B_1,\ldots,B_k\}$
   where $\min B_1 < \min B_2 < \ldots < \min B_k$,
   and set $w_i = r$ if $i \in B_r$.
   The statistics (\ref{def.wachs-white}a--d) then arise
   by restricting the Wachs--White statistics to restricted growth functions
   and mapping them to set partitions via the bijection.
}
\begin{subeqnarray}
   \lb(\pi)
   & = &
   \#\{ (B_1,B_2,k) \colon\:  \min B_1 < \min B_2 < k \in B_1  \}
         \\[2mm]
   \ls(\pi)
   & = &
   \#\{ (B_1,B_2,k) \colon\:  \min B_1 < \min B_2 \le k \in B_2  \}
         \\[2mm]
   \rb(\pi)
   & = &
   \#\{ (B_1,B_2,k) \colon\:  \min B_1 < \min B_2 
            \,\hbox{ with }\, k \in B_1
            \,\hbox{ and }\, k < \max B_2 \}
         \nonumber \\ \\
   \rs(\pi)
   & = &
   \#\{ (B_1,B_2,k) \colon\:  \min B_1 < \min B_2 \le k < \max B_1
            \,\hbox{ with }\, k \in B_2 \}
   \qquad\quad
 \label{def.wachs-white}
\end{subeqnarray}
We can also modify the definition of $\ls$
to make it more closely analogous to that of $\lb$,
by removing the case $\min B_2 = k$:
\begin{subeqnarray}
   \ls'(\pi)
   & \eqdef &
   \#\{ (B_1,B_2,k) \colon\:  \min B_1 < \min B_2 < k \in B_2  \}
        \\[2mm]
   & = &
   \ls(\pi)  \,-\, \binom{|\pi|}{2}
   \;.
 \label{def.lsprime}
\end{subeqnarray}
The statistics $\lb$ and $\ls$ had been introduced earlier by
Milne \cite[Remark~4.13]{Milne_82}:
a~triplet $(B_1,B_2,k)$ contributing to $\lb$ (resp.\ $\ls$)
is called an {\em inversion}\/ (resp.\ {\em dual inversion}\/) of $\pi$,
and we also write $\inv(\pi) = \lb(\pi)$ and $\invtilde(\pi) = \ls(\pi)$.
See e.g.\ \cite[section~3]{Zeng_95}.

Wachs and White \cite{Wachs_91} called $\lb$ and $\ls$ the ``easy'' statistics,
because it is straightforward to show that the coefficient array
\be
   S_{p,q}(n,k)
   \;\eqdef\;
   \sum_{\pi \in \Pi_{n,k}}  q^{\lb(\pi)} p^{\ls(\pi)}
\ee
(where $\Pi_{n,k}$ denotes the partitions of $[n]$ with $k$ blocks)
satisfies the recurrence
\be
   S_{p,q}(n,k)
   \;=\;
   p^{k-1} S_{p,q}(n-1,k-1)  \,+\,  [k]_{p,q} S_{p,q}(n-1,k)
   \;.
 \label{eq.Spq.recurrence}
\ee
Equivalently, the coefficient array
\be
   S'_{p,q}(n,k)
   \;\eqdef\;
   \sum_{\pi \in \Pi_{n,k}}  q^{\lb(\pi)} p^{\ls'(\pi)}
   \;=\;
   p^{-\binom{k}{2}} \, S_{p,q}(n,k)
\ee
satisfies the recurrence
\be
   S'_{p,q}(n,k)
   \;=\;
   S'_{p,q}(n-1,k-1)  \,+\,  [k]_{p,q} S'_{p,q}(n-1,k)
   \;,
\ee
which shows in particular that the pair $(\lb,\ls')$ has a symmetric distribution
on $\Pi_{n,k}$.

By contrast with the ``easy'' statistics $\lb$ and $\ls$,
Wachs and White \cite{Wachs_91} called $\rb$ and $\rs$ the ``hard'' statistics,
because there is no obvious way to prove a recurrence for them.
Wachs and White nevertheless showed, by a nontrivial bijection,
that the pair $(\rs,\rb)$ is equidistributed on $\Pi_{n,k}$ with $(\lb,\ls)$.

Here we will show that, curiously, one of the ``hard'' statistics
--- namely, $\rs$ ---
has a simple interpretation in terms of our overlap and covering statistics:

\begin{proposition}
   \label{prop.rs}
For partitions $\pi \in \Pi_n$, we have
\be
   \rs(\pi)
   \;=\;
   \ov(\pitilde) \,+\, 2\,\cov(\pitilde) \,+\, \covin(\pitilde) \,+\, \pscov(\pitilde)
         \label{eq.rs}
\ee
where $\pitilde$ denotes the reversal of $\pi$,
i.e.\ the image of $\pi$ under the map $i \mapsto \itilde \eqdef n+1-i$.
\end{proposition}

We remark that $\ov$, $\cov$ and $\pscov$ are manifestly reversal-invariant
[i.e.\ $\ov(\pitilde) = \ov(\pi)$, etc.].
By contrast, $\ovin$ and $\covin$ are not reversal-invariant for $n \ge 5$,
but their sum $\ovin + \covin$ is reversal-invariant by \reff{eq.ovin+covin}.

As preparation for the proof of Proposition~\ref{prop.rs},
we define reversals of the statistics \reff{def.ov.cov.kpi}:
\begin{subeqnarray}
   \ovtilde(k,\pi)
   & \eqdef &
   \ov(\ktilde,\pitilde)
        \nonumber \\
   & = &
   \#\{ (B_1,B_2) \colon\:  k \in B_1 \,\hbox{ and }\,
                          \min B_1 < \min B_2 < k < \max B_2  \}
       \qquad\quad \\[5mm]
   \covtilde(k,\pi)
   & \eqdef &
   \cov(\ktilde,\pitilde)
        \nonumber \\
   & = &
   \#\{ (B_1,B_2) \colon\:  k \in B_2 \,\hbox{ and }\,
                          \min B_1 < \min B_2 < k < \max B_1  \}
       \qquad\quad \\[5mm]
   \qcovtilde(k,\pi)
   & \eqdef &
   \qcov(\ktilde,\pitilde)
        \nonumber \\
   & = &
   \#\{ B \colon\:  B \not\ni k \,\hbox{ and }\, \min B < k < \max B \}
 \label{def.ov.cov.kpi.reversed.0}
\end{subeqnarray}
[The definition of $\qcov$ is in fact reversal-invariant,
so that $\qcovtilde(k,\pi) = \qcov(k,\pi)$.]
 
\proofof{Proposition~\ref{prop.rs}}
In the definition of $\rs(\pi)$,
we separate the cases with $\min B_2 < k$ from those with $\min B_2 = k$:
\begin{subeqnarray}
   \rs(\pi)
   & = &
   \sum_{k \in {\rm insiders} \,\cup\, {\rm closers}}  \! \covtilde(k,\pi)
   \;\;+
   \sum_{k \in {\rm openers} \,\cup\, {\rm singletons}}  \!\! \qcovtilde(k,\pi)
   \qquad
        \\[2mm]
   & = &
   \sum_{k \in {\rm insiders} \,\cup\, {\rm openers}}  \!\!\! \cov(k,\pitilde)
   \;\;+
   \sum_{k \in {\rm closers} \,\cup\, {\rm singletons}}  \qcov(k,\pitilde)
        \\[2mm]
   & = &
   \covin(\pitilde) \,+\, \cov(\pitilde) \,+\,
     [\ov(\pitilde) \,+\, \cov(\pitilde)]  \,+\, \pscov(\pitilde)
\end{subeqnarray}
by (\ref{def.ov.cov.kpi.reversed.0}b,c),
(\ref{eq.ovcov.total}b), (\ref{eq.ovin.covin.total}b),
\reff{eq.qcov.singletons} and \reff{eq.qcovclosers.ov+cov}.
%
\qed

{\bf Remark.}
Straightforward computation shows that $\lb$, $\ls$ and $\rb$
--- in contrast to $\rs$ ---
{\em cannot}\/ be written as a linear combination of
$\ov$, $\cov$, $\ovin$, $\covin$, $\ovinrev$, $\covinrev$,
$\binom{|\pi|}{2}$, $n |\pi|$ and $\binom{n}{2}$
[where $\ovinrev(\pi) = \ovin(\pitilde)$ and $\covinrev(\pi) = \covin(\pitilde)$].
In fact, a linear combination $a_1 \lb + a_2 \ls + a_3 \rb + a_4 \rs$
can be written in this way {\em only}\/ if $a_1 = a_2 = a_3 = 0$.
To see this, for $\ls$ and $\rb$ it suffices to consider $n=3$,
while for $\lb$ and general linear combinations it suffices to consider $n=4$.
\myendremark

\bigskip

We can now rederive an S-fraction
for the generating polynomials associated to the $q$-Stirling numbers
\be
   S_q(n,k)
   \;\eqdef\;
   S_{1,q}(n,k)
   \;=\;
   S'_{1,q}(n,k)
   \;=\;
   \sum_{\pi \in \Pi_{n,k}}  q^{\lb(\pi)}
   \;,
\ee
which was obtained some years ago by one of us \cite[eq.~(2.1)]{Zeng_95}.
Indeed, using the Wachs--White equidistribution result $\lb \sim \rs$
together with the involution $\pi \mapsto \pitilde$
and the identity \reff{eq.rs},
and then applying the definition \reff{def.Bn.ovcov.brec}
and Theorem~\ref{thm.setpartitions.pq.ovcov.brec},
we see that
\begin{eqnarray}
   & &
   \sum_{k=0}^n S_q(n,k) \, x^k
   \;\eqdef\;
   \sum_{\pi \in \Pi_n} x^{|\pi|} q^{\lb(\pi)}
   \;=\;
   \sum_{\pi \in \Pi_n} x^{|\pi|} q^{\rs(\pi)}
   \;=\;
   \sum_{\pi \in \Pi_n} x^{|\pi|} q^{\rs(\pitilde)}
          \nonumber \\[1mm]
   & & \quad
   \;=\;
   B_n^{(2)}(x,x,1,1,1,1,1,q,q,q^2,q)
   \;=\;
   B_n(x,x,1,1,1,1,1,q,q,q^2,q)
   \;.
   \qquad\quad
 \label{eq.zeng.inv}
\end{eqnarray}
Specializing Corollary~\ref{cor.setpartitions.Jtype.pq}
to $y = v = 1$, $p = 1$ and $r = q$,
we recover the S-fraction of \cite[eq.~(2.1)]{Zeng_95}:
\be
   \sum_{n=0}^\infty \sum_{k=0}^n S_q(n,k) \, x^k \: t^n
   \;=\;
   \cfrac{1}{1 - \cfrac{xt}{1 - \cfrac{t}{1 - \cfrac{qxt}{1- \cfrac{(1+q)t}{1 - \cfrac{q^2 xt}{1 - \cfrac{(1+q+q^2)t}{1-\cdots}}}}}}}
   \label{eq.Sfrac.zeng1}
\ee
with coefficients
\begin{subeqnarray}
   \alpha_{2k-1}  & = &  q^{k-1} x   \\[1mm]
   \alpha_{2k}    & = &  [k]_q
 \label{eq.Sfrac.zeng1.coeffs}
\end{subeqnarray}

  
On the other hand, in the same paper
Zeng also obtained \cite[eqn.~(2.2)]{Zeng_95}
an S-fraction for the generating polynomials associated to
the modified $q$-Stirling numbers
\be
   \widetilde{S}_q(n,k)
   \;\eqdef\;
   q^{\binom{k}{2}}
   S_{1,q}(n,k)
   \;=\;
   S_{q,1}(n,k)
   \;=\;
   \sum_{\pi \in \Pi_{n,k}}  q^{\ls(\pi)}
   \;.
\ee
Namely,\footnote{
   The formula following \cite[eqn.~(2.2)]{Zeng_95} has a typographical error:
   it should read $\lambda_{2n-1} = a q^{2n-2}$, not $a q^{2n}$.
   Note that the correct formula is given in \cite[eqn.~(2.11)]{Zeng_95}.

   Also, the definition \cite[eqn.~(1.2)]{Zeng_95} has a typographical error:
   the coefficient should be $q^k$, not $q^{k-1}$.
}
\be
   \sum_{n=0}^\infty \sum_{k=0}^n \widetilde{S}_q(n,k) \, x^k \: t^n
   \;=\;
   \cfrac{1}{1 - \cfrac{xt}{1 - \cfrac{(1+(q-1)x)t}{1 - \cfrac{q^2 xt}{1- \cfrac{(1+q(q-1)x)t}{1 - \cfrac{q^4 xt}{1 - \cfrac{(1+q^2(q-1)x)t}{1-\cdots}}}}}}}
   \label{eq.Sfrac.zeng2}
\ee
with coefficients
\begin{subeqnarray}
   \alpha_{2k-1}  & = &  q^{2k-2} x   \\[1mm]
   \alpha_{2k}    & = &  (1 + q^{k-1}(q-1)x) \, [k]_q
 \slabel{eq.Sfrac.zeng2.coeffs.b}
 \label{eq.Sfrac.zeng2.coeffs}
\end{subeqnarray}
Unfortunately, we do {\em not}\/ know how to obtain this S-fraction
as a special case of our results here;
its combinatorial meaning remains quite mysterious (at least to us).
We leave it as an open problem to understand
\reff{eq.Sfrac.zeng2}/\reff{eq.Sfrac.zeng2.coeffs}
as a special case of some more general result.
Please note that \reff{eq.Sfrac.zeng2}/\reff{eq.Sfrac.zeng2.coeffs}
differs from all of the other continued fractions in this paper
in that the coefficient \reff{eq.Sfrac.zeng2.coeffs.b}
contains a term with a minus sign;
this may be an indication of its combinatorial subtlety.

\subsection{A remark on the Ehrenborg--Readdy intertwining statistic}
   \label{subsec.setpartitions.intertwining}

\begin{figure}[t]
\centering
\begin{tikzpicture}
\draw (0,0) node {$\bullet$};
\draw (0,-.2) node [below] {$1$};
\draw (2,0) node {$\bullet$};
\draw (2,-.2) node[below] {$2$};\draw (4,0) node{$\bullet$};
\draw (4,-.2) node[below] {$3$};
\draw (6,0) node{$\bullet$};
\draw (6,-.2) node[below] {$4$};
\draw (8,0) node{$\bullet$};
\draw (8,-.2) node[below] {$5$};
\draw (10,0) node{$\bullet$};
\draw (10,-.2) node[below] {$6$};
\draw [very thick] (0,0) arc (160:20:2.1);
\draw [very thick] (4,0) arc (150:30:3.45);
\draw [very thick] (0,0) arc (10:80:2);
\draw [very thick] (12,2) arc (90:180:2);
\draw [very thick, dotted]  (6,0) arc (0:180:2);
\draw[very thick, dotted] (2,0) arc (0:80:2);
\draw[very thick, dotted] (8,0) arc (0:180:1);
\draw [very thick,dotted] (10,2) arc (90:180:2);
\end{tikzpicture}
\caption{
    Computation of the interwining number of the partition
    $\pi=\{\{1, 3, 6\}, \{2,4,5\}\}$:
    we have $\iota(\pi)=4$.
}
   \label{fig.ehrenborg}
\end{figure}
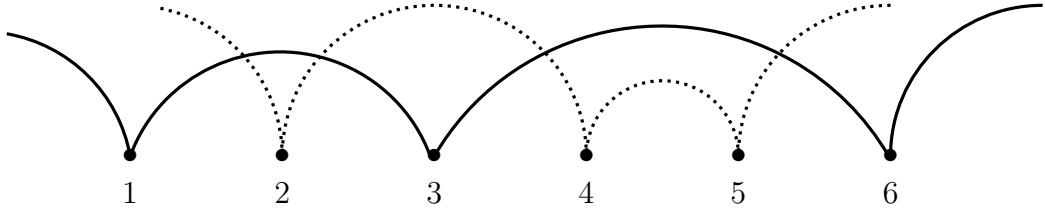

Ehrenborg and Readdy \cite[section~6]{Ehrenborg_96}
have introduced a statistic on set partitions that can be defined as follows:
For $i,j \in \Z$, let $\interval(i,j)$ denote the open interval
\be
   \interval(i,j)
   \;=\;
   \{ m \in \Z \colon\:  \min(i,j) < m < \max(i,j) \}
   \;.
\ee
By definition $\interval(i,j) = \interval(j,i)$.
Then, for two disjoint nonempty finite subsets $B,C \subset \Z$,
define the {\em intertwining number}\/
\be
   \iota(B,C)
   \;=\;
   \# \{(b,c) \in B \times C \colon\:  \interval(b,c) \cap (B \cup C) = \emptyset \}
   \;.
\ee
Of course $\iota(B,C) = \iota(C,B)$.
This intertwining number can be interpreted graphically as follows:
Draw solid (resp.\ dashed) arcs between consecutive elements of $B$ (resp.~$C$)
as in the usual graphical representation of a set partition;
but now also draw a solid (resp.\ dashed) arc
from the smallest element of $B$ (resp.~$C$) to $-\infty$,
and a solid (resp.\ dashed) arc
from the largest element of $B$ (resp.~$C$) to $+\infty$.
Then $\iota(B,C)$ is the total number of crossings between
solid and dashed arcs,
with the understanding that the arcs are drawn so that
two arcs to $-\infty$, or two arcs to $+\infty$, never cross
(see Figure~\ref{fig.ehrenborg}).
Indeed, if $b \in B$ and $c \in C$
(let's say for concreteness that $b < c$)
and there is no point of $B$ or $C$ between $b$ and $c$,
then the arc upwards from $b$
(whether to the next element of $B$ or to $+\infty$)
necessarily intersects the arc downwards from $c$
(whether to the previous element of $C$ or to $-\infty$);
but if there is an element of $B$ and/or $C$ between $b$ and $c$,
then these arcs will not intersect.
It is now easy to see that $\iota(B,C) \ge 1$.

Now, for a partition $\pi = \{B_1,\ldots,B_k\}$ of $[n]$,
define the {\em intertwining number}\/
\be
   \iota(\pi)
   \;=\;
   \sum_{1 \le i < j \le k} \! \iota(B_i,B_j)
 \label{def.iota.pi}
\ee
(of course this quantity does not depend on
 how the blocks of $\pi$ are ordered).
Since $\iota(B_i,B_j) \ge 1$ for all $i \neq j$,
we have $\iota(\pi) \ge \binom{|\pi|}{2}$;
we therefore define the {\em reduced intertwining number}\/
\be
   \iota'(\pi)  \;\eqdef\;  \iota(\pi) \,-\, \binom{|\pi|}{2}  \;\ge\; 0
   \;.
 \label{def.iotaprime.pi}
\ee

Ehrenborg and Readdy showed \cite[Proposition~6.3]{Ehrenborg_96} that
\begin{subeqnarray}
    S_q(n,k)  & = &  \sum_{\pi \in \Pi_{n,k}} q^{\iota'(\pi)}
  \slabel{eq.ehrenborg.iotaprime}
        \\[2mm]
    \widetilde{S}_q(n,k)  & = &  \sum_{\pi \in \Pi_{n,k}} q^{\iota(\pi)}
\end{subeqnarray}
(Note that their $q$-Stirling numbers $S[n,k]$
 correspond to our $\widetilde{S}_q(n,k)$.)
Here we would like to observe that their intertwining number
can be written as a combination of our crossing, nesting, overlap and covering
statistics:

\begin{proposition}
   \label{prop.intertwining}
For partitions $\pi \in \Pi_n$, we have
\begin{subeqnarray}
   \iota'(\pi)
   & = &
   \crr(\pi) \,+\, \ov(\pi) \,+\, \cov(\pi) \,+\, \pscov(\pi)
 \slabel{eq.prop.intertwining.a}
      \\[2mm]
   & = &
   \crin(\pi) \,+\, 2\,\crop(\pi) \,+\, \neop(\pi) \,+\, \psne(\pi)
   \;.
 \slabel{eq.prop.intertwining.b}
 \label{eq.prop.intertwining}
\end{subeqnarray}
\end{proposition}

\proof
Consider a pair of blocks $B,C$ contributing to the sum \reff{def.iota.pi},
and let us fix the order by assuming that $\min B < \min C$.
For each block, there are three types of arcs:
internal arcs, the arc to $-\infty$, and the arc to $+\infty$.
So there are nine types of possible crossings:

1,2,3) The arc from $-\infty$ to $\min B$ cannot intersect any arc of $C$.

4,5) The arc from $-\infty$ to $\min C$ will always intersect an arc of $B$
(namely, an internal arc of $B$ in case $\min B < \min C < \max B$,
 and the arc from $\max B$ to $+\infty$ in case
 $\min B \le \max B < \min C$).
So there are $\binom{|\pi|}{2}$ intersections of this type,
which compensates the term subtracted in \reff{def.iotaprime.pi}.

6) The arc from $\max B$ to $+\infty$ will intersect an internal arc of $C$
in case $\min C < \max B < \max C$,
or in other words $\min B < \min C < \max B < \max C$,
i.e.\ blocks $B$ and $C$ form an overlap.

7) The arc from $\max B$ to $+\infty$ cannot intersect
     the arc from $\max C$ to $+\infty$.

8) The arc from $\max C$ to $+\infty$ will intersect an internal arc of $B$
in case $\min B < \max C < \max B$,
or in other words $\min B < \min C \le \max C < \max B$
i.e.\ blocks $B$ and $C$ form a covering or a pseudo-covering.

9) Finally, the intersections between internal arcs of different blocks give
(when summed over pairs $B,C$) the contribution $\crr(\pi)$.

Putting this all together proves \reff{eq.prop.intertwining.a}.

Then \reff{eq.prop.intertwining.b}
follows by using the trivial identity $\crr = \crop + \crin$
along with the identities $\ov + \cov = \crop + \neop$
[cf. (\ref{eq.lemma.cross+nest.ov+cov}a)]
and $\pscov = \psne$ [cf.\ before \reff{def.ovin.covin}].
\qed

Applying now \reff{eq.ehrenborg.iotaprime} and \reff{eq.prop.intertwining.b}
together with the definition \reff{def.Bn.crossnest.refined}, we see that
\be
   \sum_{k=0}^n S_q(n,k) \, x^k
   \;=\;
   \sum_{\pi \in \Pi_n} x^{|\pi|} q^{\iota'(\pi)}
   \;=\;
   B_n(x,x,1,1,1,1,1,q,q^2,1,q,q)
   \;,
\ee
which differs from \reff{eq.zeng.inv} by interchanging
$(p_1,p_2) \leftrightarrow (q_1,q_2)$.
Since the S-fraction of Corollary~\ref{cor.setpartitions.Jtype.pq}
is invariant under $p \leftrightarrow q$ when $y = v = 1$,
we re-obtain the S-fraction \reff{eq.Sfrac.zeng1}/\reff{eq.Sfrac.zeng1.coeffs}.

\bigskip

{\bf Remarks.}
1.  It would be interesting to try to find a statistic
``dual'' to $\iota'$ that would give the two-variable polynomials
$S_{p,q}(n,k)$ or $S'_{p,q}(n,k)$.

2.  Sagan \cite{Sagan_91} has defined statistics $\maj$ and $\majhat$
for set partitions, which have the distributions
$S_q(n,k)$ and $\widetilde{S}_q(n,k)$, respectively.
He also discussed the joint distributions
$(\inv,\maj)$, $(\inv,\majhat)$ and $(\maj,\majhat)$,
each of which satisfies a different $(p,q)$-generalization
of the Stirling recurrence:
$(\inv,\majhat)$ corresponds to \reff{eq.Spq.recurrence},
while the other two are new.
It does not seem that any of the three have nice continued fractions,
except when $p=1$ or $q=1$.
\myendremark

\subsection{Counting connected components; indecomposable set partitions}
   \label{subsec.setpartitions.connected}

Let us now show how to extend our set-partition polynomials
to count also the connected components of a set partition.
As a corollary, we will obtain continued fractions for
indecomposable set partitions.
The method is identical, {\em mutatis mutandis}\/,
to the one used in Section~\ref{subsec.permutations.connected}
to count connected components in permutations.

A {\em divider}\/ of a set partition $\pi \in \Pi_n$
is an index $i \in [n]$ such that the interval $[1,i]$
is a union of blocks of $\pi$;
equivalently, the complementary interval $[i+1,n]$
is a union of blocks of $\pi$.
Clearly, when $n = 0$ (hence $\pi = \emptyset$) there are no dividers;
when $n \ge 1$, the index $n$ is always a divider,
and there may or may not be others.
A {\em connected component}\/ (or {\em atomic part}\/) of $\pi \in \Pi_n$
\cite[A127743]{OEIS}
is a minimal nonempty interval $[i,j] \subseteq [n]$
such that the intervals $[1,i-1]$, $[i,j]$ and $[j+1,n]$
are all unions of blocks of $\pi$.
If $1 \le i_1 < i_2 < \ldots < i_k = n$ are the dividers of $\pi$,
then $[1,i_1]$, $[i_1 + 1, i_2]$, \ldots, $[i_{k-1} + 1, i_k]$
are its connected components.
So the number of connected components equals the number of dividers;
we write it as $\ccc(\pi)$.
Thus $\ccc(\emptyset) = 0$;
for $n \ge 1$ we have $1 \le \ccc(\pi) \le n$,
with $\ccc(\pi) = n$ if and only if $\pi$ is the partition into singletons.
A set partition $\pi$ is called {\em indecomposable}\/
(or {\em atomic}\/) if $\ccc(\pi) = 1$
\cite[A087903/A074664]{OEIS}.\footnote{
   {\bf Warning:}
   We were tempted to use here the term ``irreducible'',
   but we felt obliged to avoid it because this term has been used previously
   for {\em at least two distinct}\/ other classes of partitions!
   Several authors \cite{Bender_85,Beissinger_85} \cite[p.~16]{Mansour_13}
   call a set partition ``irreducible'' if no proper subinterval of $[1,n]$
   is a union of blocks \cite[A099947]{OEIS};
   this is more restrictive than our condition
   that no {\em initial}\/ subinterval $[1,i]$ with $i \le n-1$
   is a union of blocks.
   On the other hand,
   the term ``irreducible'' is sometimes used \cite[A055105]{OEIS}
   to denote yet another class of partitions
   (also called ``unsplittable'').
   This latter class is equinumerous with the indecomposable partitions
   as defined here \cite{Chen_11e}
   but gives rise to a different triangular array
   when refined according to the number of blocks.
}

In any of the set-partition polynomials studied thus far,
we can insert an additional factor $\zeta^{\ccc(\pi)}$.
This affects the continued fractions as follows:

\begin{theorem}[Counting connected components in set partitions]
   \label{thm.connected.setpartitions}
Consider any of the polynomials
\reff{eq.thm.setpartitions},
\reff{def.Bn.x1x2y1y2},
\reff{def.Bn.crossnest.refined},
\reff{def.Bn.ovcov.brec},
\reff{def.Bnhat},
\reff{def.Bnhat2},
\reff{def.Bnhat3},
or \reff{def.Bnhat4},
and insert an additional factor $\zeta^{\ccc(\pi)}$.
Then the continued fractions associated to the ordinary generating functions
are modified as follows:
in each S-fraction, multiply $\alpha_1$ by $\zeta$;
in each J-fraction, multiply $\gamma_0$ and $\beta_1$ by $\zeta$.
\end{theorem}

This result has an easy proof in our labeled-Motzkin-paths formalism,
as we shall remark in
Section~\ref{subsec.setpartitions.J}.
But it also has a simple ``renewal theory'' explanation, as follows:
Given set partitions $\pi = \{B_1,\ldots,B_k\} \in \Pi_m$
and $\pi' = \{B'_1,\ldots,B'_l\} \in \Pi_n$,
let us define their {\em concatenation}\/ $\pi|\pi' \in \Pi_{m+n}$
as $\{B_1,\ldots,B_k, B'_1 + m,\ldots, B'_l+m\}$.
Multiple concatenations are defined in the obvious way.
Then every set partition can be written uniquely as a concatenation
of (zero or more) indecomposable set partitions
(namely, $\pi$ restricted to its connected components,
 with indices relabeled to start at 1).
Now let $P_n$ be any set-partition polynomial based on statistics that are
additive under concatenation, and include also a factor $\zeta^{\ccc(\pi)}$;
and let $P_n^{\rm ind}$ be the corresponding polynomial with the sum
restricted to indecomposable set partitions (without the factor $\zeta$).
Now define the ordinary generating functions
\begin{subeqnarray}
   f(t)  & = &  \sum_{n=0}^\infty P_n t^n  \\[2mm]
   g(t)  & = &  \sum_{n=1}^\infty P_n^{\rm ind} t^n
\end{subeqnarray}
Then it is immediate from the foregoing that
\be
   f(t)  \;=\;  {1 \over 1 \,-\, \zeta g(t)}
   \;.
 \label{eq.set_partitions.renewal}
\ee
Moreover, all of the statistics that have been considered here
are indeed additive under concatenation:
it is easy to see that this holds for statistics
based on counting blocks by size,
on classifying elements as opener/closer/insider/singleton,
on the record or block-record structure,
on crossings and nestings,
or on overlaps and coverings.
Theorem~\ref{thm.connected.setpartitions} is an immediate consequence.

These considerations also allow us to deduce continued fractions
for the ordinary generating functions of the polynomials $P_n^{\rm ind}$
associated to indecomposable set partitions.
Indeed, it follows immediately from \reff{eq.set_partitions.renewal}
that if $f(t)$ is the ordinary generating function associated to all
set partitions (without factors $\zeta^{\ccc(\pi)}$),
then $g(t) = 1 - 1/f(t)$ is the ordinary generating function associated to
indecomposable set partitions.
The continued fractions transform according to the same formulae
\reff{eq.countconn.Pn.S}--\reff{eq.countconn.Pnind.J}
as for indecomposable permutations.

\section{Perfect matchings}
  \label{sec.intro.perfectmatchings}

\subsection{S-fraction}

Euler showed \cite[section~29]{Euler_1760} that the
generating function of the odd semifactorials can be represented
as an S-type continued fraction
\be
   \sum_{n=0}^\infty (2n-1)!! \: t^n
   \;=\;
   \cfrac{1}{1 - \cfrac{1t}{1 - \cfrac{2t}{1 - \cfrac{3t}{1- \cdots}}}}
 \label{eq.2n-1semifact.contfrac}
\ee
with coefficients $\alpha_n = n$.\footnote{
   See also \cite[Section~2.6]{Callan_09} for a combinatorial proof
   of \reff{eq.2n-1semifact.contfrac}
   based on a counting of height-labeled Dyck paths.
}
Since $(2n-1)!!$ enumerates perfect matchings of a $2n$-element set,
it is natural to seek polynomial refinements of this sequence
that enumerate perfect matchings of $[2n]$ according to some
natural statistic(s).
Note that we can regard a perfect matching either as a special type
of set partition (namely, one in which all blocks are of size~2)
or as a special type of permutation
(namely, one in which all cycles are of length~2,
 i.e.\ a fixed-point-free involution).
We write these two interpretations
as $\pi \in \scrm_{2n} \subseteq \Pi_{2n}$
and $\sigma \in \scri_{2n} \subseteq \Sym_{2n}$, respectively.


Inspired by \reff{eq.2n-1semifact.contfrac},
let us introduce the polynomials $M_n(x,y,u,v)$ defined
by the continued fraction
\be
   \sum_{n=0}^\infty M_n(x,y,u,v) \: t^n
   \;=\;
   \cfrac{1}{1 - \cfrac{xt}{1 - \cfrac{(y+v)t}{1 - \cfrac{(x+2u)t}{1- \cfrac{(y+3v)t}{1 - \cfrac{(x+4u)t}{1 - \cfrac{(y+5v)t}{1-\cdots}}}}}}}
 \label{eq.matching.fourvar.contfrac}
\ee
with coefficients
\begin{subeqnarray}
   \alpha_{2k-1}  & = &  x + (2k-2) u \\
   \alpha_{2k}    & = &  y + (2k-1) v
 \label{def.weights.matching.fourvar}
\end{subeqnarray}
Clearly $M_n(x,y,u,v)$ is a homogeneous polynomial of degree $n$.
Since $M_n(1,1,1,1) = (2n-1)!!$,
it is plausible to expect that $M_n(x,y,u,v)$ enumerates
perfect matchings of $[2n]$ according to some natural trivariate statistic.

To show this, let us adopt the interpretation of
perfect matchings as fixed-point-free-involutions.
We recall the classification of indices $i$ of a permutation $\sigma$
into cycle peaks, cycle valleys, cycle double rises, cycle double falls,
and fixed points.
Note that if $\sigma$ is an involution,
then it has no cycle double rises or cycle double falls;
moreover, $i$ is a cycle peak (resp.\ cycle valley)
if and only if it is the largest (resp.\ smallest) element
of a 2-element cycle.

Now let $\sigma$ be a fixed-point-free involution on $[2n]$,
so that it consists of $n$ 2-element cycles.
For each cycle, we look at its largest element (i.e.\ the cycle peak)
and classify it into four types:
\begin{itemize}
   \item {\em even cycle-peak antirecord}\/ (ecpar)\ \ 
      [i.e.\ $i$ is even and is an antirecord];
   \item {\em odd cycle-peak antirecord}\/ (ocpar)\ \ 
      [i.e.\ $i$ is odd and is an antirecord];
   \item {\em even cycle-peak non-antirecord}\/ (ecpnar)\ \ 
      [i.e.\ $i$ is even and is not an antirecord];
   \item {\em odd cycle-peak non-antirecord}\/ (ocpnar)\ \ 
      [i.e.\ $i$ is odd and is not an antirecord].
\end{itemize}
(Note that a cycle peak cannot be a record, but that it can be an antirecord.)
Similarly, we classify the smallest element of each cycle
(i.e.\ the cycle valley) into four types:
\begin{itemize}
   \item {\em even cycle-valley record}\/ (ecvr)\ \ 
      [i.e.\ $i$ is even and is a record];
   \item {\em odd cycle-valley record}\/ (ocvr)\ \ 
      [i.e.\ $i$ is odd and is a record];
   \item {\em even cycle-valley non-record}\/ (ecvnr)\ \ 
      [i.e.\ $i$ is even and is not a record];
   \item {\em odd cycle-valley non-record}\/ (ocvnr)\ \ 
      [i.e.\ $i$ is odd and is not a record].
\end{itemize}
(Note that a cycle valley cannot be an antirecord, but that it can be a record.)

This classification of indices also has an easy translation
into the interpretation of perfect matchings as set partitions
in which every block has size 2.
Obviously ``cycle valley'' is equivalent to ``opener'',
and ``cycle peak'' to ``closer''.
Moreover, an opener is a record
if and only if it is an exclusive record
in the sense defined in Section~\ref{subsec.intro.setpartitions.S}:
that is, an opener $j$ that is paired with a closer $k$ ($> j$)
is a record if and only if there does not exist an opener $i < j$
that is paired with a closer $l > k$.
In terms of the nesting statistic \reff{def.ne.k.pi},
an opener $j$ is a record if and only if
$\nee(j,\pi) = 0$.
Similarly, a closer $k$ that is paired with an opener $j$ ($< k$)
is an antirecord if and only if there does not exist a closer $l > k$
that is paired with an opener $i < j$;
or in other words $\nee(j,\pi) = 0$.
(Note that this latter equation involves $j$, not $k$.)

With these preliminaries, we can now state our result:

\begin{theorem}[S-fraction for perfect matchings]
   \label{thm.matchings}
The polynomials $M_n(x,y,u,v)$ defined by
\reff{eq.matching.fourvar.contfrac}/\reff{def.weights.matching.fourvar}
have the combinatorial interpretation
\begin{subeqnarray}
   M_n(x,y,u,v)
   & = &
   \sum_{\sigma \in \scri_{2n}}
      x^{\ecpar(\sigma)} y^{\ocpar(\sigma)}
         u^{\ecpnar(\sigma)}  v^{\ocpnar(\sigma)}
                        \slabel{eq.matching.fourvar.a}  \\[2mm]
   & = &
   \sum_{\sigma \in \scri_{2n}}
      x^{\ocvr(\sigma)} y^{\ecvr(\sigma)}
         u^{\ocvnr(\sigma)}  v^{\ecvnr(\sigma)}
                        \slabel{eq.matching.fourvar.b}
   \;,
 \label{eq.matching.fourvar}
\end{subeqnarray}
where the sums run over fixed-point-free involutions of $[2n]$.
\end{theorem}

The interpretations \reff{eq.matching.fourvar.a} and
\reff{eq.matching.fourvar.b} are of course trivially equivalent
under the bijection $\sigma \mapsto R \circ \sigma \circ R$
with $R(i) = 2n+1-i$, which preserves the cycle structure of a
permutation but interchanges even with odd,
peak with valley, and record with antirecord.

Comparing the definitions
\reff{eq.eulerian.fourvar.contfrac}/\reff{def.weights.eulerian.fourvar}
and
\reff{eq.matching.fourvar.contfrac}/\reff{def.weights.matching.fourvar},
we see immediately that
\be
   M_n(x,y,u,v)  \;=\;  P_n(x,y+v,2u,2v)
   \;.
 \label{eq.Mn.Pn.identity}
\ee
We leave it as an open problem to give a bijective proof of this identity
based on the combinatorial interpretations
\reff{eq.eulerian.fourvar.arec} [or \reff{eq.eulerian.fourvar.cyc}]
and \reff{eq.matching.fourvar}.

Some special cases of Theorem~\ref{thm.matchings} were previously known,
notably:

\begin{itemize}
   \item The polynomials
       \cite[Proposition~7]{Dumont_86}
       \cite[Corollaire~15]{Randrianarivony_94}
       \cite{Savage_12,Ma_13,Ma_13b,Ma_15,Ma_17}
       \cite[Project~6.6.1]{Boros_04}
       \cite[A185411/A185410/A156919]{OEIS}
\be
   M_n(x,y,x,y)  \;=\;
   \sum_{\sigma \in \scri_{2n}} x^{\ecp(\sigma)} y^{\ocp(\sigma)}
   \;,
\ee
which count perfect matchings of $[2n]$ according to the number of pairs
that have even or odd largest entries.
   \item The polynomials \cite[A127160]{OEIS}
\be
   M_n(x,x,u,u)  \;=\;
   \sum_{\sigma \in \scri_{2n}} x^{\arec(\sigma)} u^{n-\arec(\sigma)}
   \;,
\ee
which count fixed-point-free involutions of $[2n]$ according to the number of
antirecords (or records).
These polynomials arise in several contexts:
\begin{itemize}
   \item $M_n(c,c,1,1)$ are the even moments $\mu_{2n}$
of the so-called Askey--Wimp--Kerov distribution \cite{AW84,K98},
which is the orthogonality measure
for the associated Hermite polynomials $H_n(x;c)$ \cite{Drake_09}.
   \item $M_n(\gamma+1,\gamma+1,1,1)$ are the even moments $\mu_{2n}$
of the limiting distribution for the Gaussian $\beta$ ensemble
when $N \to \infty$, $\beta \to 0$ with $\beta N \to 2\gamma$:
see \cite{ABG12,BP15,DS15,BCG21}.\footnote{
      The statistic $\textrm{roof}(\pi)$ defined in \cite[section~4.1]{BCG21}
      equals the number of openers $j$ such that $\crr(j,\pi) = 0$.
      By Theorem~\ref{thm.matchings.Stype.final1} below,
      the weight $x^{\textrm{roof}(\pi)}$
      gives rise to the S-fraction \reff{eq.matching.fourvar.contfrac}
      for $M_n(x,x,1,1)$.
      The roof statistic is in some sense ``dual'' to the record statistic:
      compare Lemma~\ref{lemma.matchings.pi}(b) below.

      Indeed, we can say more.
      It can be checked that the Kasraoui--Zeng \cite{Kasraoui_06} involution
      on set partitions, $f_{\rm KZ} \colon \Pi_N \to \Pi_N$
      --- which interchanges crossings and nestings ---
      in fact interchanges the numbers of crossings and nestings
      that use any particular index $k$ in {\em third}\/ position
      (that is, $i < j < k < l$);  here $k$ must be an insider or a closer.
      On the other hand, in \reff{def.cr.ne.k.pi} we defined
      $\crr(j,\pi)$ and $\nee(j,\pi)$ using a specified index $j$
      in {\em second}\/ position (so that $j$ must be an insider or an opener).
      It follows that the map $f'_{\rm KZ} \eqdef R \circ f_{\rm KZ} \circ R$,
      where $R(i) = N+1-i$, is an involution on $\Pi_N$
      that interchanges $\crr(j,\pi)$ and $\nee(j,\pi)$
      for each index $j \in [N]$.
      Restricting to perfect matchings of $[2n]$,
      we see that $f'_{\rm KZ}$ interchanges
      $\crr(j,\pi)$ and $\nee(j,\pi)$ for each $j \in [2n]$;
      in particular, it interchanges roofs with records.
      For example, $f'_{\rm KZ}$ maps
      $\pi = \{ \{1,8\}, \{2,4\}, \{3,5\}, \{6,7\} \}$ onto
      $\pi' = \{ \{1,4\}, \{2,7\}, \{3,5\}, \{6,8\} \}$, with
      \begin{eqnarray*}
         & & \crr(2,\pi) = 0 = \nee(2,\pi')  \\
         & & \crr(6,\pi) = 0 = \nee(6,\pi')  \\
         & & \nee(2,\pi) = 1 = \crr(2,\pi')  \\
         & & \nee(6,\pi) = 1 = \crr(6,\pi')
      \end{eqnarray*}
      so that the indices 2 and 6 are roofs (but not records) in $\pi$,
      and records (but not roofs) in $\pi'$.

      We thank Bishal Deb and Vadim Gorin for discussions on this.
}
   \item $x M_n(x+1,x+1,1,1)$ counts rooted maps embeddable on
      an orientable surface (of arbitrary genus), with $n$ edges,
      with respect to the number of vertices
      \cite[Theorem~3]{Arques_00} \cite[A053979]{OEIS}.
   \item Similarly, $2^n x M_n({x \over 2}+1,{x \over 2}+1,1,1)$
      counts rooted maps embeddable on
      an arbitrary (orientable or non-orientable) surface
      (of arbitrary genus), with $n$ edges, 
      with respect to the number of vertices
      \cite[Theorem~2]{Liu_14b}.
   \item $x M_n(x+1,x+1,1,1)$ also arises in a problem concerning
      extreme-value distributions in probability theory \cite{Albrecher_09}.
\end{itemize}
Note that in this case the continued-fraction coefficients simplify to
\be
   \alpha_n  \;=\;  x + (n-1)u
   \;,
 \label{eq.alpha.Mnxxuu}
\ee
as we are no longer distinguishing between even and odd.
This special case of Theorem~\ref{thm.matchings}
can also be deduced from Theorem~\ref{thm.perm.cycle-alternating.second}:
since a fixed-point-free involution is simply a cycle-alternating
permutation in which each cycle has exactly one cycle valley,
we set $\lambda = 1/y_1$ and take $y_1 \to 0$;
then \reff{def.weights.perm.Stype.cycle-alternating.second}
becomes \reff{eq.alpha.Mnxxuu}.
\end{itemize}

We can prove Theorem~\ref{thm.matchings}
as a corollary of our master J-fraction for set partitions
(Theorem~\ref{thm.setpartitions.Jtype.final1})
by specializing variables;
the reasoning is similar to our treatment of
cycle-alternating permutations in
Section~\ref{subsec.permutations.cycle-alternating}.
We need a simple combinatorial lemma:


\begin{lemma}[Openers in perfect matchings]
   \label{lemma.matchings.pi}
Let $\pi \in \scrm_{2n} \subseteq \Pi_{2n}$ be a perfect matching of $[2n]$,
and let $j \in [2n]$ be an opener of $\pi$.  Then:
\begin{itemize}
  \item[(a)]  $j$ has opposite parity to $\crr(j,\pi) + \nee(j,\pi)$.
  \item[(b)]  $j$ is a record if and only if $\nee(j,\pi) = 0$.
\end{itemize}
\end{lemma}

\proof
(a) By \reff{eq.qne.crr+nee} and \reff{eq.qne} we have
\be
   \crr(j,\pi)  \,+\, \nee(j,\pi)
   \;=\;
   \qne(j,\pi)  \;=\;  \#\{ i<j<l \colon\: (i,l) \in \scrg_\pi  \}
   \;.
\ee
Let $m \ge 0$ be the number of pairs $(i,i') \in \scrg_\pi$
with both $i,i' < j$.
Since $j$ is an opener, it is paired with some element $j' > j$.
Therefore, every element $i < j$ is either paired with another
element $i' < j$ or else with some element $l > j$.
Hence
\be
   j-1  \;=\; \#\{i \colon i<j\}  \;=\;  \crr(j,\pi)  \,+\, \nee(j,\pi) \,+\, 2m
   \;,
\ee
which proves (a).

(b) was already observed in the paragraph preceding
Theorem~\ref{thm.matchings}.
\qed

\proofof{Theorem~\ref{thm.matchings}}
In \reff{def.Bnhat} we set $\bsfd = \bsfe = \bzero$
to force $\pi$ to be a perfect matching.
We also set $\sfb_\ell = 1$ for all $\ell \ge 0$
(i.e.\ we do not weight closers).
So the weight is simply
$\!\!\prod\limits_{i \in {\rm openers}}  \!\! \sfa_{\crr(i,\pi),\, \nee(i,\pi)}$.
By Lemma~\ref{lemma.matchings.pi}(a,b) we obtain the polynomial
\reff{eq.matching.fourvar.b} if we set
\be
   \sfa_{\ell,\ell'}
   \;=\;
   \begin{cases}
       x     &  \textrm{if $\ell' = 0$ and $\ell$ is even}  \\
       y     &  \textrm{if $\ell' = 0$ and $\ell$ is odd}   \\
       u     &  \textrm{if $\ell' \ge 1$ and $\ell+\ell'$ is even}
                                                            \\
       v     &  \textrm{if $\ell' \ge 1$ and $\ell+\ell'$ is odd}
   \end{cases}
\ee
Then
\be
   \sfa^\star_{n-1}  \;\eqdef\;  \sum_{\ell=0}^{n-1} \sfa_{\ell,n-1-\ell}
                     \;=\;
   \begin{cases}
       x + (n-1)u     &  \textrm{if $n$ is odd}  \\
       y + (n-1)v     &  \textrm{if $n$ is even}
   \end{cases}
\ee
With these specializations,
the J-fraction \reff{eq.thm.setpartitions.Jtype.final1}
becomes the S-fraction \reff{eq.matching.fourvar.contfrac}
if we identify $B_{2n} = M_n$ and replace $t^2$ by $t$.
This proves \reff{eq.matching.fourvar.b};
and \reff{eq.matching.fourvar.a} then follows
by the trivial bijection noted earlier.
\qed

We can alternatively prove Theorem~\ref{thm.matchings} as a corollary of our
second master S-fraction for cycle-alternating permutations
(Theorem~\ref{thm.perm.cycle-alternating.3.second}).
Here we use the interpretation of perfect matchings as
fixed-point-free involutions
(i.e.\ permutations in which every cycle is of length~2);
and we observe that a fixed-point-free involution of $[2n]$
is simply a cycle-alternating permutation of $[2n]$
in which the number of cycles is maximal (namely, $n$).
So we can obtain perfect matchings from
Theorem~\ref{thm.perm.cycle-alternating.3.second}
by replacing $t \to t/\lambda$ and then taking $\lambda \to\infty$.
The details are as follows.
We begin with a simple combinatorial lemma,
which is a close analogue of Lemma~\ref{lemma.matchings.pi}:

\begin{lemma}
   \label{lemma.matchings}
Let $\sigma \in \scri_{2n} \subseteq \Sym_{2n}$
be a fixed-point-free involution of $[2n]$.
If $i \in [2n]$ is a cycle peak of $\sigma$, then:
\begin{itemize}
  \item[(a)]  $i$ has the same parity as $\lcross(i,\sigma) + \lnest(i,\sigma)$.
  \item[(b)]  $i$ is an antirecord if and only if $\lnest(i,\sigma) = 0$.
\end{itemize}
If $i \in [2n]$ is a cycle valley of $\sigma$, then:
\begin{itemize}
  \item[(c)]  $i$ has the opposite parity from
      $\ucross(i,\sigma) + \unest(i,\sigma)$.
  \item[(d)]  $i$ is a record if and only if $\unest(i,\sigma) = 0$.
\end{itemize}
\end{lemma}

\proof
(a)  Let $k \in [2n]$ be a cycle peak of $\sigma$, so that $\sigma(k) < k$.
Then the set $\{j \colon\; j < k \}$, which has cardinality $k-1$,
can be partitioned as
\be
   \{\sigma(k)\}
   \;\cup\;
   \{j < k \colon\: \sigma(j) < k \}
   \;\cup\;
   \{j < \sigma(k) \colon\: \sigma(j) > k \}
   \;\cup\;
   \{\sigma(k) < j < k \colon\: \sigma(j) > k \}
   \;.
\ee
The first of these sets has cardinality 1;
the second has even cardinality;
the third has cardinality $\lnest(k,\sigma)$;
and the fourth has cardinality $\lcross(k,\sigma)$.

(b) Again let $k \in [2n]$ be a cycle peak, so that $\sigma(k) < k$.
Then $k$ fails to be an antirecord in case there exists an index $l > k$
such that $\sigma(l) < \sigma(k)$.
But this is precisely the statement that $\lnest(k,\sigma) > 0$.

%

(c,d)  The proofs are similar.
\qed

\secondproofof{Theorem~\ref{thm.matchings}}
We start from the polynomial
\reff{def.Qn.firstmaster.cycle-alternating.second}
for cycle-alternating permutations,
multiply by $\lambda^{-n}$, and take $\lambda \to \infty$:
this restricts the sum to fixed-point-free involutions.
We also set $\sfa_\ell = 1$ for all $\ell \ge 0$
(i.e.\ we do not weight cycle valleys).
The result is
\be
   \lim_{\lambda \to\infty}
   \lambda^{-n} \, \widehat{Q}_{2n}(\bone,\bsfb,\bzero,\bzero,\bzero,\lambda)
   \;=\;
   \sum_{\sigma \in \scri_{2n}}
   \;\:
   \prod\limits_{i \in {\rm cpeak}} \!\! \sfb_{\lcross(i,\sigma),\,\lnest(i,\sigma)}
   \;.
\ee
By Lemma~\ref{lemma.matchings}(a,b) we obtain the polynomial
\reff{eq.matching.fourvar.a} if we set
\be
   \sfb_{\ell,\ell'}
   \;=\;
   \begin{cases}
       x     &  \textrm{if $\ell' = 0$ and $\ell$ is even}  \\
       y     &  \textrm{if $\ell' = 0$ and $\ell$ is odd}   \\
       u     &  \textrm{if $\ell' \ge 1$ and $\ell+\ell'$ is even}
                                                            \\
       v     &  \textrm{if $\ell' \ge 1$ and $\ell+\ell'$ is odd}
   \end{cases}
\ee
Then
\be
   \sfb^\star_{n-1}  \;=\;  \sum_{\ell=0}^{n-1} \sfb_{\ell,n-1-\ell}
                     \;=\;
   \begin{cases}
       x + (n-1)u     &  \textrm{if $n$ is odd}  \\
       y + (n-1)v     &  \textrm{if $n$ is even}
   \end{cases}
\ee
{}From \reff{def.weights.permutations.Stype.final1.cycle-alternating.second}
we have
\be
   \lim_{\lambda \to\infty}
   \lambda^{-1} \, \alpha_n
   \;=\;
   \sfb^\star_{n-1}
   \;,
\ee
which completes the proof.
\qed

\medskip

{\bf Remarks.}
1.  It is natural to try amalgamating (\ref{eq.matching.fourvar}a,b)
into an eight-variable polynomial
\begin{eqnarray}
   \widehat{M}_n(x,y,u,v,\bar{x}, \bar{y}, \bar{u}, \bar{v})
   & = &
   \sum_{\sigma \in \scri_{2n}}
      x^{\ecpar(\sigma)} y^{\ocpar(\sigma)}
         u^{\ecpnar(\sigma)}  v^{\ocpnar(\sigma)}  \:\times
       \qquad\qquad
       \nonumber \\[-1mm]
   & & \qquad\;
      \bar{x}^{\ocvr(\sigma)} \bar{y}^{\ecvr(\sigma)}
         \bar{u}^{\ocvnr(\sigma)}  \bar{v}^{\ecvnr(\sigma)}
      \;,
\end{eqnarray}
so that $\widehat{M}_n(x,y,u,v,1,1,1,1) = M_n(x,y,u,v)$
and $\widehat{M}_n(1,1,1,1,\bar{x}, \bar{y}, \bar{u}, \bar{v}) =
 M_n(\bar{x}, \bar{y}, \bar{u}, \bar{v})$.
But it seems that we can get a J-fraction with polynomial coefficients
only if we specialize to six variables:
either $u=x$ and $v=y$, or $\bar{u}=\bar{x}$ and $\bar{v}=\bar{y}$.
And in these cases we get an S-fraction:  for instance,
\be
   \sum_{n=0}^\infty \widehat{M}_n(x,y,u,v,\bar{x},\bar{y},\bar{x},\bar{y}) \: t^n
   \;=\;
   \cfrac{1}{1 - \cfrac{x \bar{x}t}{1 - \cfrac{(y+v) \bar{y}t}{1 - \cfrac{(x+2u) \bar{x}t}{1- \cfrac{(y+3v) \bar{y}t}{1 - \cdots}}}}}
 \label{eq.matching.sixvar.contfrac}
\ee
with coefficients
\begin{subeqnarray}
   \alpha_{2k-1}  & = &  [x + (2k-2) u] \, \bar{x} \\
   \alpha_{2k}    & = &  [y + (2k-1) v] \, \bar{y}
 \label{def.weights.matching.sixvar}
\end{subeqnarray}
But this is actually an immediate consequence of Theorem~\ref{thm.matchings},
once we realize that the number of odd (resp.\ even) cycle valleys
is equal to the number of even (resp.\ odd) cycle peaks;
so if we count cycle valleys only by parity and not by record status,
then $\bar{x}$ simply multiplies $x$ and $u$,
and $\bar{y}$ multiplies $y$ and $v$.

2. One could try to go farther,
by introducing the \emph{disjoint} classification
of each 2-element cycle into 16 categories according to the
status of its cycle peak (ecpar, ocpar, ecpnar, ocpnar)
\emph{and} its cycle valley (ecvr, ocvr, ecvnr, ocvnr),
and then defining a \emph{homogeneous} 16-variable polynomial.
We do not know whether any interesting continued fractions
can be obtained from this polynomial.
\myendremark

\subsection[$p,q$-generalizations]{$\bm{p,q}$-generalizations}

Let us now generalize the four-variable polynomial $M_n(x,y,u,v)$
by adding weights for crossings and nestings
as in Section~\ref{subsec.intro.permutations.pq}.
Note that if $\sigma \in \Sym_n$ is an involution
(not necessarily fixed-point-free),
we trivially have $\ucross(\sigma) = \lcross(\sigma)$
and $\unest(\sigma) = \lnest(\sigma)$;
so we write them simply as $\crr(\sigma)$ and $\nee(\sigma)$, respectively.
These quantities of course coincide with $\crr(\pi)$ and $\nee(\pi)$
as defined in Section~\ref{subsec.intro.setpartitions.pq}
for the matching (not necessarily perfect) $\pi \in \Pi_n$
that corresponds to the involution $\sigma \in \Sym_n$;
and we here have $\crr(\pi) = \crop(\pi)$ and $\nee(\pi) = \neop(\pi)$
since a matching has no insiders.
We now define
\begin{subeqnarray}
   M_n(x,y,u,v,p,q)
   & = &
   \sum_{\sigma \in \scri_{2n}}
      x^{\ecpar(\sigma)} y^{\ocpar(\sigma)}
         u^{\ecpnar(\sigma)}  v^{\ocpnar(\sigma)}
         p^{\crr(\sigma)} q^{\nee(\sigma)}
     \qquad
      \\[2mm]
   & = &
   \sum_{\sigma \in \scri_{2n}}
      x^{\ocvr(\sigma)} y^{\ecvr(\sigma)}
         u^{\ocvnr(\sigma)}  v^{\ecvnr(\sigma)}
         p^{\crr(\sigma)} q^{\nee(\sigma)}
   \;,
 \label{def.matching.fourvar.pq}
\end{subeqnarray}
where the sums run over fixed-point-free involutions of $[2n]$,
and the equality of (\ref{def.matching.fourvar.pq}a)
and (\ref{def.matching.fourvar.pq}b) again follows from
the bijection $\sigma \mapsto R \circ \sigma \circ R$
(which preserves the numbers of crossings and nestings).

In fact, we can go farther, by distinguishing crossings and nestings
according to whether the element in second position is even or odd.
That is, let us say that a crossing or nesting $i < j < k < l$
is {\em even}\/ (resp.\ {\em odd}\/) if $j$ is even (resp.\ odd).
We denote by $\ecr(\sigma)$, $\ocr(\sigma)$, $\ene(\sigma)$, $\oone(\sigma)$
the numbers of even crossings, odd crossings, even nestings and odd nestings,
respectively.
We then define
\begin{subeqnarray}
   & &
   \!\!
   M_n(x,y,u,v,p_+,p_-,q_+,q_-)
         \nonumber \\[2mm]
   & & \qquad =\;
   \sum_{\sigma \in \scri_{2n}}
      x^{\ecpar(\sigma)} y^{\ocpar(\sigma)}
         u^{\ecpnar(\sigma)}  v^{\ocpnar(\sigma)}
         p_+^{\ocr(\sigma)} p_-^{\ecr(\sigma)}
         q_+^{\oone(\sigma)} q_-^{\ene(\sigma)}
     \qquad\qquad
      \\[2mm]
   & & \qquad =\;
   \sum_{\sigma \in \scri_{2n}}
      x^{\ocvr(\sigma)} y^{\ecvr(\sigma)}
         u^{\ocvnr(\sigma)}  v^{\ecvnr(\sigma)}
         p_+^{\ecr(\sigma)} p_-^{\ocr(\sigma)}
         q_+^{\ene(\sigma)} q_-^{\oone(\sigma)}
   \;,
 \slabel{def.matching.fourvar.pq.plusminus.b}
 \label{def.matching.fourvar.pq.plusminus}
\end{subeqnarray}
where the two formulae are again related by
$\sigma \mapsto R \circ \sigma \circ R$.
We find:

\begin{theorem}[S-fraction for perfect matchings, $p,q$-generalization]
   \label{thm.matching.fourvar.pq.Stype}
The ordinary generating function of the polynomials
$M_n(x,y,u,v,p_+,p_-,q_+,q_-)$ has the S-type continued fraction
\be
   \sum_{n=0}^\infty M_n(x,y,u,v,p_+,p_-,q_+,q_-) \: t^n
   \;\,=\;\,
   \cfrac{1}{1 - \cfrac{xt}{1 - \cfrac{(p_+ y + q_+ v)t}{1 - \cfrac{(p_-^2 x+q_- \, [2]_{p_-,q_-} u)t}{1- \cfrac{(p_+^3 y+ q_+ \, [3]_{p_+,q_+} v)t}{1 - \cdots}}}}}
   \label{eq.thm.matching.fourvar.pq.Stype}
\ee
with coefficients
\begin{subeqnarray}
   \alpha_{2k-1}  & = &  p_-^{2k-2} x + q_- \, [2k-2]_{p_-,q_-} u \\[1mm]
   \alpha_{2k}    & = &  p_+^{2k-1} y + q_+ \, [2k-1]_{p_+,q_+} v
 \label{def.weights.matching.fourvar.pq.Stype}
\end{subeqnarray}
\end{theorem}


Note that if $u=x$ and/or $v=y$,
then the weights \reff{def.weights.matching.fourvar.pq.Stype} simplify to
$\alpha_{2k-1} = [2k-1]_{p_-,q_-} \, x$ and $\alpha_{2k} = [2k]_{p_+,q_+} \, y$,
respectively.
For the special case $x=y=u=v$, $p_+ = p_-$ and $q_+ = q_-$,
the S-fraction \reff{eq.thm.matching.fourvar.pq.Stype}
was obtained previously by Kasraoui and Zeng \cite{Kasraoui_06}
(see also \cite[p.~3280]{Blitvic_12}).

The proof of Theorem~\ref{thm.matching.fourvar.pq.Stype}
is a straightforward extension of the method used for
Theorem~\ref{thm.matchings}:

\proofof{Theorem~\ref{thm.matching.fourvar.pq.Stype}}
In \reff{def.Bnhat} we set $\bsfd = \bsfe = \bzero$
to force $\pi$ to be a perfect matching.
We also set $\sfb_\ell = 1$ for all $\ell \ge 0$
(i.e.\ we do not weight closers).
So the weight is simply
$\!\!\prod\limits_{i \in {\rm openers}}  \!\! \sfa_{\crr(i,\pi),\, \nee(i,\pi)}$.
By Lemma~\ref{lemma.matchings.pi}(a,b)
we obtain the polynomial \reff{def.matching.fourvar.pq.plusminus.b} if we set
\be
   \sfa_{\ell,\ell'}
   \;=\;
   \begin{cases}
       p_-^\ell q_-^{\ell'} x
           &  \textrm{if $\ell' = 0$ and $\ell$ is even}  \\
       p_+^\ell q_+^{\ell'} y
           &  \textrm{if $\ell' = 0$ and $\ell$ is odd}   \\
       p_-^\ell q_-^{\ell'} u
           &  \textrm{if $\ell' \ge 1$ and $\ell+\ell'$ is even} \\
       p_+^\ell q_+^{\ell'} v
           &  \textrm{if $\ell' \ge 1$ and $\ell+\ell'$ is odd}
   \end{cases}
\ee
Then
\be
   \sfa^\star_{n-1}  \;\eqdef\;  \sum_{\ell=0}^{n-1} \sfa_{\ell,n-1-\ell}
                     \;=\;
   \begin{cases}
       p_-^{n-1} x + q_- [n-1]_{p_-,q_-} u     &  \textrm{if $n$ is odd}  \\[1mm]
       p_+^{n-1} y + q_+ [n-1]_{p_+,q_+} v     &  \textrm{if $n$ is even}
   \end{cases}
\ee
With these specializations,
the J-fraction \reff{eq.thm.setpartitions.Jtype.final1}
becomes the S-fraction \reff{eq.thm.matching.fourvar.pq.Stype}
if we identify $B_{2n} = M_n$ and replace $t^2$ by $t$.
\qed

{\bf Remarks.}
1.  Comparing the continued fractions
\reff{eq.thm.perm.pq.Stype.BIG.1}/\reff{def.weights.thm.perm.pq.Stype.BIG.1}
and
\reff{eq.thm.matching.fourvar.pq.Stype}/\reff{def.weights.matching.fourvar.pq.Stype},
we see that
\begin{eqnarray}
   & &
   \!\!
   M_n(x,y,u,v,p_+,p_-,q_+,q_-)
        \nonumber \\[1mm]
   & & \qquad  =\;
   P_n(x,p_+ y+q_+ v,(p_- + q_-)u/q_-,(p_+ +q_+)v,p_+^2,p_-^2,q_+^2,q_-^2)
   \;,
   \qquad\qquad
 \label{eq.Mn.Pn.identity.pq}
\end{eqnarray}
which generalizes \reff{eq.Mn.Pn.identity}.
We leave it as an open problem to find a bijective proof
of \reff{eq.Mn.Pn.identity.pq}.

Let us observe a curious fact about \reff{eq.Mn.Pn.identity.pq}.
One might think that the appearance of the {\em squares}\/ of $p_\pm,q_\pm$
on the right-hand side of \reff{eq.Mn.Pn.identity.pq}
comes from the fact that each crossing or nesting in a perfect matching
corresponds to {\em two}\/ crossings or nestings
--- one upper and one lower --- in the corresponding permutation.
But this does {\em not}\/ seem to be the correct explanation,
since the meaning of the subscripts $+$ and $-$ is different
on the two sides of \reff{eq.Mn.Pn.identity.pq}:
in $M_n$ it distinguishes even from odd,
while in $P_n$ it distinguishes upper from lower.
So we really do not understand the combinatorial meaning of
\reff{eq.Mn.Pn.identity.pq}.

2.  An explicit expression for $M_n(1,1,1,1,p,1)$
--- which counts perfect matchings of $[2n]$
by the number of crossings (or nestings) ---
was found implicitly by Touchard \cite{Touchard_52}
and explicitly by Riordan \cite{Riordan_75}
(see also \cite{Read_79,Penaud_95,Prodinger_12,Blitvic_12,Josuat-Verges_13}):
\be
   M_n(1,1,1,1,p,1)  \;=\;
   (1-p)^{-n}  \sum_{k=0}^n (-1)^k \, t_{n,k} \, p^{k(k+1)/2}
 \label{eq.Mn1q}
\ee
where
\be
   t_{n,k}  \;=\;  \binom{2n}{n+k} - \binom{2n}{n+k+1}
            \;=\;  {2k+1 \over n+k+1} \binom{2n}{n+k}
            \;=\;  {2k+1 \over 2n+1} \binom{2n+1}{n+k+1}
 \label{eq.Mn1q.bis}
\ee
are a variant of the ballot numbers.
No explicit expression for $M_n(1,1,1,1,p,q)$ seems to be known.
\myendremark

%

\subsection{Master S-fraction}

Finally, we can get a master S-fraction for perfect matchings
by specializing the first master J-fraction for set partitions
(Theorem~\ref{thm.setpartitions.Jtype.final1}).
We introduce two infinite families of indeterminates
$\bsfa = (\sfa_{\ell,\ell'})_{\ell,\ell' \ge 0}$ and
$\bsfb = (\sfb_{\ell})_{\ell \ge 0}$,
and define the polynomials $M_n(\bsfa,\bsfb)$ by
\be
   M_n(\bsfa,\bsfb)
   \;=\;
   \sum_{\pi \in \scrm_{2n}}
   \;
   \prod\limits_{i \in {\rm openers}}  \!\! \sfa_{\crr(i,\pi),\, \nee(i,\pi)}
   \prod\limits_{i \in {\rm closers}}  \! \sfb_{\qne(i,\pi)}
   \;,
 \label{def.Mn.master}
\ee
where $\crr(i,\pi)$, $\nee(i,\pi)$ and $\qne(i,\pi)$
are as defined in (\ref{def.cr.ne.k.pi}a,b) and \reff{eq.qne}.
Of course, $M_n(\bsfa,\bsfb) = B_{2n}(\bsfa,\bsfb,\bzero,\bzero)$
since setting $\bsfd = \bsfe = \bzero$ in \reff{def.Bnhat}
is precisely what is needed to restrict the summation to perfect matchings.
{}From Theorem~\ref{thm.setpartitions.Jtype.final1} we immediately deduce:

\begin{theorem}[Master S-fraction for perfect matchings]
   \label{thm.matchings.Stype.final1}
The ordinary generating function of the polynomials
$M_n(\bsfa,\bsfb)$ has the S-type continued fraction
\be
   \sum_{n=0}^\infty M_n(\bsfa,\bsfb) \: t^n
   \;=\;
   \cfrac{1}{1 - \cfrac{\sfa_{00} \sfb_{0} t}{1  - \cfrac{(\sfa_{01} + \sfa_{10}) \sfb_{1} t}{1  - \cfrac{(\sfa_{02} + \sfa_{11} + \sfa_{20}) \sfb_{2} t}{1 - \cdots}}}}
   \label{eq.thm.matchings.Stype.final1}
\ee
with coefficients
\be
   \alpha_n   \;=\;   \sfa^\star_{n-1} \, \sfb_{n-1}
 \label{def.weights.matchings.Stype.final1}
\ee
where
$\displaystyle
   \sfa^\star_{n-1}  \;\eqdef\;  \sum_{\ell=0}^{n-1} \sfa_{\ell,n-1-\ell}
$.
\end{theorem}

Alternatively, we can prove Theorem~\ref{thm.matchings.Stype.final1}
as a corollary of our
second master S-fraction for cycle-alternating permutations
(Theorem~\ref{thm.perm.cycle-alternating.3.second}),
using the interpretation of perfect matchings as fixed-point-free involutions:

\secondproofof{Theorem~\ref{thm.matchings.Stype.final1}}
We start by applying the bijection $\sigma \mapsto R \circ \sigma \circ R$
to \reff{def.Qn.firstmaster.cycle-alternating.second}:
this interchanges valleys with peaks, and upper with lower, yielding
\be
   \widehat{Q}_{2n}(\bsfa,\bsfb,\bzero,\bzero,\bzero,\lambda)
   \;=\;
   \sum_{\sigma \in \Sym^{\rm ca}_{2n}}
   \;\:
   \lambda^{\cyc(\sigma)} \;
   \prod\limits_{i \in {\rm cpeak}}  \! \sfa_{\lcross(i,\sigma)+\lnest(i,\sigma)}
   \prod\limits_{i \in {\rm cval}} \!\! \sfb_{\ucross(i,\sigma),\,\unest(i,\sigma)}
   \;.
 \label{def.Qn.firstmaster.cycle-alternating.second.REVERSED}
\ee
We then multiply by $\lambda^{-n}$ and take $\lambda \to \infty$:
this restricts the sum to fixed-point-free involutions.
We have $\crr(j,\pi) = \ucross(j,\sigma)$
and $\nee(j,\pi) = \unest(j,\sigma)$
since in both cases the distinguished index is in second position
[compare (\ref{def.cr.ne.k.pi}a,b) to (\ref{def.ucrossnestjk}a,b)].
And we have $\qne(k,\pi) = \lcross(k,\sigma)+\lnest(k,\sigma)$
since in both cases the distinguished index is in third position
[compare \reff{eq.qne} to (\ref{def.ucrossnestjk}c,d)].
So \reff{def.Qn.firstmaster.cycle-alternating.second.REVERSED}
corresponds to \reff{def.Mn.master}
with the interchange of letters $\sfa \leftrightarrow \sfb$.
Then from \reff{def.weights.permutations.Stype.final1.cycle-alternating.second}
we have
\be
   \lim_{\lambda \to\infty}
   \lambda^{-1} \, \alpha_n
   \;=\;
   \sfa_{n-1} \, \sfb^\star_{n-1}
   \;,
\ee
which corresponds to \reff{def.weights.matchings.Stype.final1}
with $\sfa \leftrightarrow \sfb$.
\qed

\bigskip

{\bf Remark.}
It is unfortunate that Theorem~\ref{thm.matchings.Stype.final1}
treats openers and closers asymmetrically.
But we do not know how to avoid this.
If we consider perfect matchings as a special case of set partitions,
this asymmetry is imposed by \reff{def.Bnhat}
and Theorem~\ref{thm.setpartitions.Jtype.final1},
as remarked already in Section~\ref{subsec.intro.setpartitions.firstmaster}.
If we consider perfect matchings as a limiting case of
cycle-alternating permutations,
we are obliged to use the {\em second}\/ master S-fraction
--- that is, Theorem~\ref{thm.perm.cycle-alternating.3.second}
instead of Theorem~\ref{thm.perm.cycle-alternating.3} ---
in order to get access to the statistic $\cyc(\sigma)$;
and Theorem~\ref{thm.perm.cycle-alternating.3.second}
treats cycle valleys differently from cycle peaks.
\myendremark

\subsection{Counting connected components; indecomposable perfect matchings}  
   \label{subsec.perfectmatchings.connected}

We can extend our polynomials to count also the connected components
of a perfect matching.
Since the method is essentially identical to the one used previously for
permutations (Section~\ref{subsec.permutations.connected})
and set partitions (Section~\ref{subsec.setpartitions.connected}),
we will be brief.

We can consider a perfect matching
either as a permutation in which all cycles are of length~2,
or as a set partition in which all blocks are of size~2.
We then specialize the definitions of ``divider'', ``connected component''
and ``indecomposable'' from
either permutations (Section~\ref{subsec.permutations.connected})
or set partitions (Section~\ref{subsec.setpartitions.connected});
both methods give the same notions for perfect matchings.
We denote the number of connected components in a perfect matching
by $\ccc(\sigma) = \ccc(\pi)$.
The enumeration of indecomposable perfect matchings
can be found in \cite[A000698]{OEIS}.
The enumeration of perfect matchings by number of connected components
is (to our surprise) not in \cite{OEIS} at present;  it begins as
%
%
\vspace*{-5mm}
\begin{table}[H]
\centering
\small
\begin{tabular}{c|rrrrrrrrr|r}
$n \setminus k$ & 0 & 1 & 2 & 3 & 4 & 5 & 6 & 7 & 8 & Row sums \\
\hline
0 & 1 &  &  &  &  &  &  &  &  & 1  \\
1 & 0 & 1 &  &  &  &  &  &  &  & 1  \\
2 & 0 & 2 & 1 &  &  &  &  &  &  & 3  \\
3 & 0 & 10 & 4 & 1 &  &  &  &  &  & 15  \\
4 & 0 & 74 & 24 & 6 & 1 &  &  &  &  & 105  \\
5 & 0 & 706 & 188 & 42 & 8 & 1 &  &  &  & 945  \\
6 & 0 & 8162 & 1808 & 350 & 64 & 10 & 1 &  &  & 10395  \\
7 & 0 & 110410 & 20628 & 3426 & 568 & 90 & 12 & 1 &  & 135135  \\
8 & 0 & 1708394 & 273064 & 38886 & 5696 & 850 & 120 & 14 & 1 & 2027025  \\
\end{tabular}
\end{table}
\vspace*{-5mm}

In any of the set-partition polynomials studied thus far,
we can insert an additional factor $\zeta^{\ccc(\sigma)}$.
We then have:

\begin{theorem}[Counting connected components in perfect matchings]
   \label{thm.connected.perfectmatchings}
Consider any of the polynomials
\reff{eq.matching.fourvar},
\reff{def.matching.fourvar.pq}
or \reff{def.Mn.master},
and insert an additional factor $\zeta^{\ccc(\sigma)} = \zeta^{\ccc(\pi)}$.
Then the S-fractions associated to the ordinary generating functions
are modified by multiplying $\alpha_1$ by $\zeta$.
\end{theorem}

The reasoning is identical to that in
Sections~\ref{subsec.permutations.connected}
and \ref{subsec.setpartitions.connected}.

\section{Preliminaries for the proofs}   \label{sec.prelim}

Our proofs are based on Flajolet's \cite{Flajolet_80}
combinatorial interpretation of continued fractions
in terms of Dyck and Motzkin paths,
together with some bijections mapping combinatorial objects
(permutations, set partitions or perfect matchings)
to labeled Dyck or Motzkin paths.
We begin by reviewing briefly these two ingredients.

\subsection{Combinatorial interpretation of continued fractions}
   \label{subsec.prelim.1}

Recall that a {\em Motzkin path}\/ of length $n \ge 0$
is a path $\omega = (\omega_0,\ldots,\omega_n)$
in the right quadrant $\N \times \N$,
starting at $\omega_0 = (0,0)$ and ending at $\omega_n = (n,0)$,
whose steps $s_j = \omega_j - \omega_{j-1}$
are $(1,1)$ [``rise''], $(1,-1)$ [``fall''] or $(1,0)$ [``level''].
We write $\scrm_n$ for the set of Motzkin paths of length~$n$,
and $\scrm = \bigcup_{n=0}^\infty \scrm_n$.
A Motzkin path is called a {\em Dyck path}\/ if it has no level steps.
A Dyck path always has even length;
we write $\scrd_{2n}$ for the set of Dyck paths of length~$2n$,
and $\scrd = \bigcup_{n=0}^\infty \scrd_{2n}$.

Let ${\bf a} = (a_k)_{k \ge 0}$, ${\bf b} = (b_k)_{k \ge 1}$
and ${\bf c} = (c_k)_{k \ge 0}$ be indeterminates;
we will work in the ring $\Z[[{\bf a},{\bf b},{\bf c}]]$
of formal power series in these indeterminates.
To each Motzkin path $\omega$ we assign a weight
$W(\omega) \in \Z[{\bf a},{\bf b},{\bf c}]$
that is the product of the weights for the individual steps,
where a rise starting at height~$k$ gets weight~$a_k$,
a~fall starting at height~$k$ gets weight~$b_k$,
and a level step at height~$k$ gets weight~$c_k$.
Flajolet \cite{Flajolet_80} showed that
the generating function of Motzkin paths
can be expressed as a continued fraction:

\begin{theorem}[Flajolet's master theorem]
   \label{thm.flajolet}
We have
\be
   \sum_{\omega \in \scrm}  W(\omega)
   \;=\;
   \cfrac{1}{1 - c_0 - \cfrac{a_0 b_1}{1 - c_1 - \cfrac{a_1 b_2}{1- c_2 - \cfrac{a_2 b_3}{1- \cdots}}}}
 \label{eq.thm.flajolet}
\ee
as an identity in $\Z[[{\bf a},{\bf b},{\bf c}]]$.
Equivalently, we have
\be
   \sum_{n=0}^\infty t^n \sum_{\omega \in \scrm_n}  W(\omega)
   \;=\;
   \cfrac{1}{1 - c_0 t - \cfrac{a_0 b_1 t^2}{1 - c_1 t - \cfrac{a_1 b_2 t^2}{1 - \cdots}}}
 \label{eq.flajolet.motzkin}
\ee
as an identity in $\Z[{\bf a},{\bf b},{\bf c}][[t]]$,
i.e.\ a J-fraction \reff{def.Jtype} with coefficients
$\gamma_n = c_n$ and $\beta_n = a_{n-1} b_n$.
\end{theorem}

Specializing \reff{eq.flajolet.motzkin} to ${\bf c} = \bzero$
and replacing $t^2$ by $t$, we obtain:

\begin{corollary}[Flajolet's master theorem for Dyck paths]
   \label{cor.flajolet.dyck}
We have
\be
   \sum_{n=0}^\infty t^n \sum_{\omega \in \scrd_{2n}}  W(\omega)
   \;=\;
   \cfrac{1}{1 - \cfrac{a_0 b_1 t}{1 - \cfrac{a_1 b_2 t}{1 - \cdots}}}
 \label{eq.flajolet.dyck}
\ee
as an identity in $\Z[{\bf a},{\bf b}][[t]]$,
i.e.\ an S-fraction \reff{def.Stype} with coefficients $\alpha_n = a_{n-1} b_n$.
\end{corollary}

\subsection{Labeled Dyck and Motzkin paths}

Let ${\bf A} = (A_k)_{k \ge 0}$, ${\bf B} = (B_k)_{k \ge 1}$
and ${\bf C} = (C_k)_{k \ge 0}$ be sequences of nonnegative integers.
An {\em $({\bf A},{\bf B},{\bf C})$-labeled Motzkin path of length $n$}\/
is a pair $(\omega,\xi)$
where $\omega = (\omega_0,\ldots,\omega_n)$
is a Motzkin path of length $n$,
and $\xi = (\xi_1,\ldots,\xi_n)$ is a sequence of integers satisfying
\be
   1  \:\le\: \xi_i  \:\le\:
   \begin{cases}
       A(h_{i-1})  & \textrm{if $h_i = h_{i-1} + 1$ (i.e.\ step $i$ is a rise)}
              \\[1mm]
       B(h_{i-1})  & \textrm{if $h_i = h_{i-1} - 1$ (i.e.\ step $i$ is a fall)}
              \\[1mm]
       C(h_{i-1})  & \textrm{if $h_i = h_{i-1}$ (i.e.\ step $i$ is a level step)}
   \end{cases}
 \label{eq.xi.ineq}
\ee
where $h_i$ is the height of the Motzkin path after step $i$,
i.e.\ $\omega_i = (i, h_i)$.
[For typographical clarity
 we have here written $A(k)$ as a synonym for $A_k$, etc.]
We call the pair $(\omega,\xi)$
an {\em $({\bf A},{\bf B})$-labeled Dyck path}\/
if $\omega$ is a Dyck path (in this case ${\bf C}$ plays no role).
We denote by $\scrm_n({\bf A},{\bf B},{\bf C})$
the set of $({\bf A},{\bf B},{\bf C})$-labeled Motzkin paths of length $n$,
and by $\scrd_{2n}({\bf A},{\bf B})$
the set of $({\bf A},{\bf B})$-labeled Dyck paths of length $2n$.

Let us stress that the numbers $A_k$, $B_k$ and $C_k$ are allowed
to take the value 0.
Whenever this happens, the path $\omega$ is forbidden to take a step
of the specified kind at the specified height.

We shall also make use of multicolored Motzkin paths.
An {\em $\ell$-colored Motzkin path}\/ is simply a Motzkin path
in which each level step has been given a ``color''
from the set $\{1,2,\ldots,\ell\}$.
In other words, we distinguish $\ell$ different types of level steps.
An {\em $({\bf A},{\bf B},{\bf C}^{(1)},\ldots,{\bf C}^{(\ell)})$-labeled
 $\ell$-colored Motzkin path of length $n$}\/
is then defined in the obvious way,
where we use the sequence ${\bf C}^{(j)}$ to bound
the label $\xi_i$ when step $i$ is a level step of type $j$.
We denote by $\scrm_n({\bf A},{\bf B},{\bf C}^{(1)},\ldots,{\bf C}^{(\ell)})$
the set of $({\bf A},{\bf B},{\bf C}^{(1)},\ldots,{\bf C}^{(\ell)})$-labeled
$\ell$-colored Motzkin paths of length $n$.

\bigskip

{\bf Remark.}  What we have called an
$({\bf A},{\bf B},{\bf C})$-labeled Motzkin path
is (up to small changes in notation)
called a {\em path diagramme}\/ by Flajolet \cite[p.~136]{Flajolet_80}
and a {\em history}\/ by Viennot \cite[p.~II-9]{Viennot_83}.
The triplet $({\bf A},{\bf B},{\bf C})$ is called a
{\em possibility function}\/.
\myendremark

\bigskip

Following Flajolet \cite[Proposition~7A]{Flajolet_80},
we can state a ``master J-fraction'' for
$({\bf A},{\bf B},{\bf C})$-labeled Motzkin paths.
Let ${\bf a} = (a_{k,\xi})_{k \ge 0 ,\; 1 \le \xi \le A(k)}$,
${\bf b} = (b_{k,\xi})_{k \ge 1 ,\; 1 \le \xi \le B(k)}$
and ${\bf c} = (c_{k,\xi})_{k \ge 0 ,\; 1 \le \xi \le C(k)}$
be indeterminates;
we give an $({\bf A},{\bf B},{\bf C})$-labeled Motzkin path $(\omega,\xi)$
a weight $W(\omega,\xi)$
that is the product of the weights for the individual steps,
where a rise starting at height~$k$ with label $\xi$ gets weight~$a_{k,\xi}$,
a~fall starting at height~$k$ with label $\xi$ gets weight~$b_{k,\xi}$,
and a level step at height~$k$ with label $\xi$ gets weight~$c_{k,\xi}$.
Then:

\begin{theorem}[Flajolet's master theorem for labeled Motzkin paths]
   \label{thm.flajolet_master_labeled_Motzkin}
We have
\be
   \sum_{n=0}^\infty t^n
   \sum_{(\omega,\xi) \in \scrm_n({\bf A},{\bf B},{\bf C})} \!\!\!  W(\omega)
   \;=\;
   \cfrac{1}{1 - c_0 t - \cfrac{a_0 b_1 t^2}{1 - c_1 t - \cfrac{a_1 b_2 t^2}{1- c_2 t - \cfrac{a_2 b_3 t^2}{1- \cdots}}}}
\ee
as an identity in $\Z[{\bf a},{\bf b},{\bf c}][[t]]$, where
\be
   a_k  \;=\;  \sum_{\xi=1}^{A(k)} a_{k,\xi}
   \;,\qquad
   b_k  \;=\;  \sum_{\xi=1}^{B(k)} b_{k,\xi}
   \;,\qquad
   c_k  \;=\;  \sum_{\xi=1}^{C(k)} c_{k,\xi}
   \;.
 \label{def.weights.akbkck}
\ee
\end{theorem}

\noindent
This is an immediate consequence of Theorem~\ref{thm.flajolet}
together with the definitions.
There is obviously also a similar theorem for
$({\bf A},{\bf B},{\bf C}^{(1)},\ldots,{\bf C}^{(\ell)})$-labeled
$\ell$-colored Motzkin paths,
in which $c_k$ involves a sum over the colors of the level steps.

By specializing to ${\bf c} = \bzero$ and replacing $t^2$ by $t$,
we obtain the corresponding theorem
for $({\bf A},{\bf B})$-labeled Dyck paths:

\begin{corollary}[Flajolet's master theorem for labeled Dyck paths]
   \label{cor.flajolet_master_labeled_Dyck}
We have
\be
   \sum_{n=0}^\infty t^n
   \sum_{(\omega,\xi) \in \scrd_{2n}({\bf A},{\bf B})} \!\!\!  W(\omega)
   \;=\;
   \cfrac{1}{1 - \cfrac{a_0 b_1 t}{1 - \cfrac{a_1 b_2 t}{1- \cfrac{a_2 b_3 t}{1- \cdots}}}}
\ee
as an identity in $\Z[{\bf a},{\bf b}][[t]]$, where
$a_k$ and $b_k$ are defined by \reff{def.weights.akbkck}.
\end{corollary}


We will also use (following Biane \cite{Biane_93})
doubly labeled Motzkin paths.
Let ${\bf A'} = (A'_k)_{k \ge 0}$, ${\bf A''} = (A''_k)_{k \ge 0}$,
${\bf B'} = (B'_k)_{k \ge 1}$ ${\bf B''} = (B''_k)_{k \ge 1}$,
${\bf C'} = (C'_k)_{k \ge 0}$, ${\bf C''} = (C''_k)_{k \ge 0}$
be sequences of nonnegative integers.
An {\em $({\bf A'},{\bf A''},{\bf B'},{\bf B''},{\bf C'},{\bf C''})$-doubly
labeled Motzkin path of length $n$}\/
is a pair $(\omega,\xi)$
where $\omega = (\omega_0,\ldots,\omega_n)$
is a Motzkin path of length~$n$,
and $\xi = \big( (\xi'_1,\xi''_1),\ldots, (\xi'_n,\xi''_n) \big)$
is a sequence of pairs of integers satisfying
\be
   1  \:\le\: \xi'_i  \:\le\:
   \begin{cases}
       A'(h_{i-1})  & \textrm{if $h_i = h_{i-1} + 1$ (i.e.\ step $i$ is a rise)}
              \\[1mm]
       B'(h_{i-1})  & \textrm{if $h_i = h_{i-1} - 1$ (i.e.\ step $i$ is a fall)}
              \\[1mm]
       C'(h_{i-1})  & \textrm{if $h_i = h_{i-1}$ (i.e.\ step $i$ is a level step)}
   \end{cases}
 \label{eq.xi.ineq.prime}
\ee
and likewise for $\xi''_i$.
Of course, doubly labeled paths can be mapped bijectively
onto singly labeled paths with $A_k = A'_k A''_k$ etc.;
but this bijection is in most cases unnatural, 
so we prefer to work with doubly labeled paths
whenever they express a combinatorially natural construction.
We also define {\em $\ell$-colored doubly labeled Motzkin paths}\/
in the obvious way.
Theorem~\ref{thm.flajolet_master_labeled_Motzkin}
has an obvious extension to doubly labeled Motzkin paths
(and to $\ell$-colored doubly labeled Motzkin paths),
which we refrain from writing out.

\section{Permutations: Proofs}   \label{sec.proofs.permutations}

%
%
%
%
%

\subsection{First master J-fraction: \\
      Proof of Theorems~\ref{thm.perms.S}(a),
      \ref{thm.perm.Jtype}, \ref{thm.perm.pq.Jtype},
      \ref{thm.perm.pq.Jtype.BIG} and \ref{thm.permutations.Jtype.final1}}
  \label{subsec.permutations.J}

In this section we will prove the first master J-fraction
for permutations (Theorem~\ref{thm.permutations.Jtype.final1}).
As a corollary we will also obtain Theorem~\ref{thm.perm.pq.Jtype.BIG},
which is obtained from Theorem~\ref{thm.permutations.Jtype.final1}
by the specialization~\reff{eq.Qn.BIG.specializations};
and from this we will in turn obtain
Theorems \ref{thm.perm.Jtype} and \ref{thm.perm.pq.Jtype},
which are special cases of Theorem~\ref{thm.perm.pq.Jtype.BIG},
as well as Theorem~\ref{thm.perms.S}(a), which is linked by
contraction \reff{eq.contraction_even.coeffs}
to the specialization \reff{def.specialization.xyuv}
of Theorem \ref{thm.perm.Jtype}.

To prove Theorem~\ref{thm.permutations.Jtype.final1},
we will employ a variant of the Foata--Zeilberger \cite{Foata_90} bijection.
More precisely, we will construct a bijection from $\Sym_n$ to the set of
$({\bf A},{\bf B},{\bf C}^{(1)}, \linebreak
 {\bf C}^{(2)},{\bf C}^{(3)})$-labeled
3-colored Motzkin paths of length $n$, where
\begin{subeqnarray}
   A_k        & = &  k+1       \quad\hbox{for $k \ge 0$}  \\
   B_k        & = &  k         \qquad\;\;\,\hbox{for $k \ge 1$}  \\
   C_k^{(1)}  & = &  k         \qquad\;\;\,\hbox{for $k \ge 0$}  \\
   C_k^{(2)}  & = &  k         \qquad\;\;\,\hbox{for $k \ge 0$}  \\
   C_k^{(3)}  & = &  1         \qquad\;\;\:\hbox{for $k \ge 0$}
 \label{def.abc}
\end{subeqnarray}
We will begin by explaining how the Motzkin path $\omega$ is defined;
then we will explain how the labels $\xi$ are defined;
next we will prove that the mapping is indeed a bijection;
next we will translate the various statistics from
$\Sym_n$ to our labeled Motzkin paths;
and finally we will sum over labels $\xi$ to obtain the weight $W(\omega)$
associated to a Motzkin path $\omega$,
which upon applying \reff{eq.flajolet.motzkin}
will yield Theorem~\ref{thm.permutations.Jtype.final1}.

\bigskip

{\bf Step 1: Definition of the Motzkin path.}
Given a permutation $\sigma \in \Sym_n$,
we classify the indices $i \in [n]$
in the usual way as cycle peak, cycle valley, cycle double rise,
cycle double fall or fixed point.
We then define a path $\omega = (\omega_0,\ldots,\omega_n)$
starting at $\omega_0 = (0,0)$ and ending at $\omega_n = (n,0)$,
with steps $s_1,\ldots,s_n$, as follows:
\begin{itemize}
   \item If $i$ is a cycle valley, then $s_i$ is a rise.
   \item If $i$ is a cycle peak, then $s_i$ is a fall.
   \item If $i$ is a cycle double fall, then $s_i$ is a level step of type 1.
   \item If $i$ is a cycle double rise, then $s_i$ is a level step of type 2.
   \item If $i$ is a fixed point, then $s_i$ is a level step of type 3.
\end{itemize}
Of course we need to prove that this is indeed a Motzkin path,
i.e.\ that all the heights $h_i$ are nonnegative and that $h_n = 0$.
We do this by obtaining a precise interpretation of the height $h_i$:

\begin{lemma}
   \label{lemma.heights}
For $i \in [n+1]$ we have
\begin{subeqnarray}
   h_{i-1}  & = &  \# \{j < i \colon\:  \sigma(j) \ge i \}
       \slabel{eq.lemma.heights.a}  \\[2mm]
            & = &  \# \{j < i \colon\:  \sigma^{-1}(j) \ge i \}
       \slabel{eq.lemma.heights.b}
 \label{eq.lemma.heights}
\end{subeqnarray}
\end{lemma}

\noindent
In particular, if $i$ is a fixed point,
then by comparing \reff{eq.lemma.heights.a} with \reff{def.level}
we see that the height of the Motzkin path after (or before) step $i$
equals the level of the fixed point:
\be
   h_{i-1}  \;=\;  h_i  \;=\;  \lev(i,\sigma)
   \;.
 \label{eq.height.fix}
\ee

\proofof{Lemma~\ref{lemma.heights}}
We shall prove \reff{eq.lemma.heights} in the equivalent form
\begin{subeqnarray}
   h_i      & = &  \# \{j \le i \colon\:  \sigma(j) > i \}
       \slabel{eq.lemma.heights2.a}  \\[2mm]
            & = &  \# \{j \le i \colon\:  \sigma^{-1}(j) > i \}
       \slabel{eq.lemma.heights2.b}  
 \label{eq.lemma.heights2}
\end{subeqnarray}
for $0 \le i \le n$ (which implies in particular that $h_0 = h_n = 0$).
To prove \reff{eq.lemma.heights2},
we represent a permutation $\sigma \in \Sym_n$ by a bipartite digraph
in which the top row of vertices is labeled $1,\ldots,n$
and the bottom row $1',\ldots,n'$,
and we draw an arrow from $i$ to $j'$ in case $\sigma(i) = j$
(see Figure~\ref{fig.FZ}a).
\begin{figure}[p]
\vspace*{4cm}
\centering
\hspace*{-8cm}
\begin{picture}(8,3)(0,0)
\setlength{\unitlength}{10mm}
\linethickness{.5mm}
\put(-0.1,2.4){$1$}\put(0,2){\circle*{0.2}}
\put(1.4,2.4){$2$}\put(1.5,2){\circle*{0.2}}
\put(2.9,2.4){$3$}\put(3,2){\circle*{0.2}}
\put(4.4,2.4){$4$}\put(4.5,2){\circle*{0.2}}
\put(5.9,2.4){$5$}\put(6,2){\circle*{0.2}}
\put(7.4,2.4){$6$}\put(7.5,2){\circle*{0.2}}
\put(8.9,2.4){$7$}\put(9,2){\circle*{0.2}}
\put(-0.1,-0.7){$1'$}\put(0,0){\circle*{0.2}}
\put(1.4,-0.7){$2'$}\put(1.5,0){\circle*{0.2}}
\put(2.9,-0.7){$3'$}\put(3,0){\circle*{0.2}}
\put(4.4,-0.7){$4'$}\put(4.5,0){\circle*{0.2}}
\put(5.9,-0.7){$5'$}\put(6,0){\circle*{0.2}}
\put(7.4,-0.7){$6'$}\put(7.5,0){\circle*{0.2}}
\put(8.9,-0.7){$7'$}\put(9,0){\circle*{0.2}}
\put(4.25,-1.7){(a)}
\qbezier(0,2)(3,1)(6,0)
\qbezier(1.5,2)(4.5,1)(7.5,0)
\qbezier(3,2)(1.5,1)(0,0)
\qbezier(3,0)(6,1)(9,2)
\qbezier(1.5,0)(3.75,1)(6,2)
\qbezier(7.5,2)(8.25,1)(9,0)
\qbezier(4.5,0)(4.5,1)(4.5,2)
\end{picture}
\\[7.5cm]
\centering
\hspace*{-9.6cm}
\begin{picture}(8,2.5)(0,0)
\setlength{\unitlength}{1.5cm}
\linethickness{.3mm}
\put(0,0){\vector(1,0){7.8}}\put(0,0){\vector(0,1){3}}
\put(-0.05,-0.5){$0$}
\put(0.95,-0.5){$1$}\put(1,0){\circle*{.15}}
\put(1.95,-0.5){$2$}\put(2,0){\circle*{.15}}
\put(2.95,-0.5){$3$}\put(3,0){\circle*{.15}}
\put(3.95,-0.5){$4$}\put(4,0){\circle*{.15}}
\put(4.95,-0.5){$5$}\put(5,0){\circle*{.15}}
\put(5.95,-0.5){$6$}\put(6,0){\circle*{.15}}
\put(6.95,-0.5){$7$}\put(7,0){\circle*{.15}}
\put(3.35,-1.2){(b)}
\put(-0.4,-0.1){$0$}
\put(-0.4,0.9){$1$}\put(0,1){\circle*{.15}}
\put(-0.4,1.9){$2$}\put(0,2){\circle*{.15}}
\put(0,0){\red{\line(1,1){1}}}
\put(1,1){\red{\line(1,1){1}}}
\put(2,2){\red{\line(1,0){1}}} \put(2.4,2.1){\small 1}
\put(3,2){\red{\line(1,0){1}}}\put(3.4,2.1){\small 3}
\put(4,2){\red{\line(1,-1){1}}}
\put(5,1){\red{\line(1,0){1}}}\put(5.4,1.1){\small 2}
\put(6,1){\red{\line(1,-1){1}}}
\put(0,0){\red{\circle*{.17}}}
\put(1,1){\red{\circle*{.17}}}
\put(2,2){\red{\circle*{.17}}}
\put(3,2){\red{\circle*{.17}}}
\put(4,2){\red{\circle*{.17}}}
\put(5,1){\red{\circle*{.17}}}
\put(6,1){\red{\circle*{.17}}}
\put(7,0){\red{\circle*{.17}}}
\end{picture}
\\[5cm]
\centering
\hspace*{-14cm}
\begin{picture}(10,3)(8,0)
\linethickness{.4mm}
\setlength{\unitlength}{0.55cm}
\put(-0.1,2.5){$i$}\put(0,2){\circle*{.2}}
\put(-0.1,-1){$i'$}\put(0,0){\circle*{.2}}
\qbezier(0,2)(0,1)(0,0)
\put(0,-1.5){\makebox(0,-1)[c]{fix}}
\put(4.9,2.5){$i$}\put(5,2){\circle*{.2}}
\put(4.9,-1){$i'$}\put(5,0){\circle*{.2}}
\put(6.9,2.5){$k$}\put(7,2){\circle*{.2}}
\put(2.8,-1){$j'$}\put(3,0){\circle*{.2}}
\qbezier(3,0)(4,1)(5,2)
\qbezier[15](5,0)(6,1)(7,2)
\put(5,-1.5){\makebox(0,-1)[c]{cdfall}}
\put(11.9,2.5){$i$}\put(12,2){\circle*{.2}}
\put(12,-1){$i'$}\put(12,0){\circle*{.2}}
\put(14,-1){$k'$}\put(14,0){\circle*{.2}}
\put(9.9,2.5){$j$}\put(10,2){\circle*{.2}}
\qbezier(10,2)(11,1)(12,0)
\qbezier[15](12,2)(13,1)(14,0)
\put(12,-1.5){\makebox(0,-1)[c]{cdrise}}
\put(17.4,2.5){$i$}\put(17.5,2){\circle*{.2}}
\put(17.4,-1){$i'$}\put(17.5,0){\circle*{.2}}
\put(19.4,2.5){$k$}\put(19.5,2){\circle*{.2}}
\put(20.4,-1){$j'$}\put(20.5,0){\circle*{.2}}
\qbezier[15](17.5,2)(19,1)(20.5,0)
\qbezier[15](17.5,0)(18.5,1)(19.5,2)
\put(18.5,-1.5){\makebox(0,-1)[c]{cval}}
\put(26.9,2.5){$i$}\put(27,2){\circle*{.2}}
\put(26.9,-1){$i'$}\put(27,0){\circle*{.2}}
\put(24.4,-1){$j'$}\put(24.5,0){\circle*{.2}}
\put(24.9,2.5){$k$}\put(25,2){\circle*{.2}}
\qbezier(25,2)(26,1)(27,0)
\qbezier(24.5,0)(25.75,1)(27,2)
\put(26,-1.5){\makebox(0,-1)[c]{cpeak}}
\put(13.8,-4){(c)}
\end{picture}
\vspace*{3cm}
%
\caption{
   (a) Bipartite digraph representing the permutation
       $\sigma = 5614273 = (152673)\,(4) \in \Sym_7$.
       Arrows run from the top row to the bottom row
       and are suppressed for clarity.
   (b) Motzkin path corresponding to the same permutation $\sigma$,
       with the types (1,2,3) of the level steps indicated.
   (c) The five possibilities for a column $ii'$.
       The dotted edges are those which are not yet ``seen'' at step $i$.
 \label{fig.FZ}
}
\end{figure}
We then ``read'' this diagram from left to right,
adding one column at each step
and taking account of all arrows that have been ``completely seen''
at the given stage:
namely, after $i$ steps we include all arrows $a \to b'$
for which both $a$ and $b$ are $\le i$.
We now claim that the height of the Motzkin path
(Figure~\ref{fig.FZ}b)
after $i$ steps
equals the number of unconnected dots in the top row after $i$ steps,
and also equals the number of unconnected dots in the bottom row
after $i$ steps;  these claims are, respectively,
\reff{eq.lemma.heights2.a} and \reff{eq.lemma.heights2.b}.
We prove these claims by induction on $i$,
by considering the five possibilities for what can happen
when we include the new column $ii'$ (see Figure~\ref{fig.FZ}c):
\begin{itemize}
   \item If $i$ is a fixed point, we have added an arrow $i \to i'$,
      and have thus added one connected dot to each row;
      the number of unconnected dots in each row remains unchanged
      from the previous step.  This agrees with $h_i = h_{i-1}$,
      since $s_i$ is a level step.
   \item If $i$ is a cycle double fall, we have added an arrow $i \to j'$
      with $j < i$.  This means that we have added the connected dot $i$
      to the top row;  in the bottom row we have added the unconnected dot $i'$
      but have also changed $j'$ from unconnected to connected.
      Therefore the number of unconnected dots in each row
      remains unchanged from the previous step;
      and this agrees with $h_i = h_{i-1}$, since $s_i$ is a level step.
   \item If $i$ is a cycle double rise, we have added an arrow $j \to i'$
      with $j < i$.  This means that we have added the connected dot $i'$
      to the bottom row;  in the top row we have added the unconnected dot $i$
      but have also changed $j$ from unconnected to connected.
      Once again $h_i = h_{i-1}$, with $s_i$ being a level step.
   \item If $i$ is a cycle valley, then no new arrows are added;
      we have therefore added the unconnected dot $i$ to the top row
      and the unconnected dot $i'$ to the bottom row.
      This agrees with $h_i = h_{i-1} + 1$, since $s_i$ is a rise.
   \item If $i$ is a cycle peak, then two new arrows are added:
      $i \to j'$ with $j = \sigma(i) < i$,
      and $k \to i'$ with $k = \sigma^{-1}(i) < i$.
      Therefore, in each row we have added one connected dot ($i$ or $i'$)
      and also changed one dot from unconnected to connected;
      therefore, the number of unconnected dots in each row decreases by 1.
      This agrees with $h_i = h_{i-1} - 1$, since $s_i$ is a fall.
\end{itemize}
\qed

{\bf Remark.}
It follows from the definition of the mapping $\sigma \mapsto \omega$
that the permutations $\sigma$ and $\sigma^{-1}$
map onto the same Motzkin path,
except that the level steps of types 1 and 2 are interchanged.
This explains why \reff{eq.lemma.heights.a} equals \reff{eq.lemma.heights.b}.
\myendremark

\bigskip

{\bf Step 2: Definition of the labels $\bm{\xi_i}$.}
We now define
\be
   \xi_i
   \;=\;
   \begin{cases}
       1 \,+\, \#\{j \colon\: j < i \hbox{ and } \sigma(j) > \sigma(i) \}
           & \textrm{if $\sigma(i) > i$ \ (cval, cdrise)}
              \\[1mm]
       1 \,+\, \#\{j \colon\: j > i \hbox{ and }  \sigma(j) < \sigma(i) \}
           & \textrm{if $\sigma(i) < i$ \ (cpeak, cdfall)}
              \\[1mm]
       1   & \textrm{if $\sigma(i) = i$ \ (fix)}
   \end{cases}
 \label{def.xi}
\ee
These definitions have a simple interpretation in terms of the
nesting statistics defined in (\ref{def.ucrossnestjk}b,d):
\be
   \xi_i - 1
   \;=\;
   \begin{cases}
       \unest(i,\sigma) & \textrm{if $\sigma(i) > i$ \ (cval, cdrise)}
              \\[1mm]
       \lnest(i,\sigma) & \textrm{if $\sigma(i) < i$ \ (cpeak, cdfall)}
              \\[1mm]
       0   & \textrm{if $\sigma(i) = i$ \ (fix)}
   \end{cases}
 \label{def.xi.bis}
\ee
Of course, we must verify that the inequalities \reff{eq.xi.ineq}/\reff{def.abc}
are satisfied; to do this, we interpret $h_{i-1} - \xi_i$
in terms of the crossing statistics defined in (\ref{def.ucrossnestjk}a,c):

\begin{lemma}[Crossing statistics]
   \label{lemma.crossing}
We have
\begin{eqnarray}
   h_{i-1} + 1 - \xi_i   & = &  \ucross(i,\sigma)
         \quad\textrm{\rm if $i \in \cval$}
     \label{eq.lemma.ucrosscval} \\[1mm]
   h_{i-1} - \xi_i   & = &  \ucross(i,\sigma)
         \quad\textrm{\rm if $i \in \cdrise$}
     \label{eq.lemma.ucrosscdrise} \\[1mm]
   h_{i-1} - \xi_i   & = &  \lcross(i,\sigma)
         \quad\:\textrm{\rm if $i \in \cpeak \,\cup\, \cdfall$}
     \label{eq.lemma.lcrosscpeak}
\end{eqnarray}
\end{lemma}

\proof
(a) If $i$ is a cycle valley
    (so that $\sigma(i) > i$ and $\sigma^{-1}(i) > i$), then
\begin{subeqnarray}
   h_{i-1} + 1 - \xi_i
   & = &
   \# \{j < i \colon\:  \sigma(j) \ge i \}
     \:-\:
   \# \{j < i \colon\:  \sigma(j) > \sigma(i) \}
          \qquad \\[2mm]
   & = &
   \# \{j < i \colon\:  \sigma(i) \ge \sigma(j) \ge i \}
          \\[2mm]
   & = &
   \# \{j < i \colon\:  \sigma(i) > \sigma(j) > i \}
          \\[2mm]
   & = &
   \ucross(i,\sigma)
   \;.
 \label{eq.Axi.cval}
\end{subeqnarray}
(b) If $i$ is a cycle double rise
      (so that $\sigma(i) > i$ and $\sigma^{-1}(i) < i$), then
\begin{subeqnarray}
   h_{i-1} - \xi_i
   & = &
   \# \{j < i \colon\,  \sigma(j) \ge i \}
     \:-\: 1 \:-\:
   \# \{j < i \colon\,  \sigma(j) > \sigma(i) \}
          \qquad \\[2mm]
   & = &
   \# \{j < i \colon\,  \sigma(j) > i \}
     \:-\: 
   \# \{j < i \colon\,  \sigma(j) > \sigma(i) \}
          \\[2mm]
   & = &
   \# \{j < i \colon\,  \sigma(i) \ge \sigma(j) > i \}
          \\[2mm]
   & = &
   \# \{j < i \colon\,  \sigma(i) > \sigma(j) > i \}
          \\[2mm]
   & = &
   \ucross(i,\sigma)
   \;.
 \label{eq.Cxi.cdrise}
\end{subeqnarray}
(c) If $i$ is a cycle peak or a cycle double fall
      (so that $\sigma(i) < i$), then
\begin{subeqnarray}
   h_{i-1} - \xi_i
   & = &
   \# \{j \colon\,  j < i \hbox{ and } \sigma^{-1}(j) \ge i \}
     \:-\: 1 \:-\:
   \# \{j \colon\,  \sigma(j) < \sigma(i) \hbox{ and } j > i \}
          \nonumber \\ \\
   & = &
   \# \{j \colon\,  \sigma(j) < i \hbox{ and } j \ge i \}
     \:-\: 1 \:-\:
   \# \{j \colon\,  \sigma(j) < \sigma(i) \hbox{ and } j > i \}
          \nonumber \\ \\
   & = &
   \# \{j \colon\,  \sigma(j) < i \hbox{ and } j > i \}
     \:-\:
   \# \{j \colon\,  \sigma(j) < \sigma(i) \hbox{ and } j > i \}
          \\[2mm]
   & = &
   \# \{j > i \colon\:  \sigma(i) \le \sigma(j) < i \}
          \\[2mm]
   & = &
   \# \{j > i \colon\:  \sigma(i) < \sigma(j) < i \}
          \\[2mm]
   & = &
   \lcross(i,\sigma)
   \;.
 \label{eq.Bxi.cpeak.cdfall}
\end{subeqnarray}
\qed

Since the quantities \reff{eq.lemma.ucrosscval}--\reff{eq.lemma.lcrosscpeak}
are manifestly nonnegative,
it follows immediately 
that the inequalities \reff{eq.xi.ineq}/\reff{def.abc} are satisfied.

For future use, let us also (partially) interpret the labels $\xi_i$
in terms of the bipartite digraph employed in the proof of
Lemma~\ref{lemma.heights} (Figure~\ref{fig.FZ}).
First recall that $h_{i-1}$ equals
the number of unconnected dots in the top row after $i-1$ steps,
and also equals the number of unconnected dots in the bottom row
after $i-1$ steps.
Now, if $i$ is a cycle double fall, then at stage $i$ we add an arrow
from $i$ on the top row to an unconnected dot $j'$ on the bottom row,
where $j = \sigma(i) < i$;
and if $i$ is a cycle peak, then we add the just-mentioned arrow
and also add an arrow from an unconnected dot $k$ on the top row to $i'$,
where $k = \sigma^{-1}(i) < i$.
We now claim that, in these two cases, $\xi_i$ is the index
of the unconnected dot $j'$ among all the unconnected dots on the bottom row:
that is, $\xi_i = r$ if and only if $j'$
is the $r$th unconnected dot on the bottom row, reading from left to right.
Indeed, by definition $\xi_i - 1$
equals $\#\{k \colon\: k > i \hbox{ and }  \sigma(k) < \sigma(i) \}$,
which is precisely the number of unconnected dots on the bottom row
to the left of $j' = \sigma(i)'$.

For cycle double rises and cycle valleys, by contrast,
the labels $\xi_i$ do not have any simple interpretation
in terms of the bipartite digraph {\em when read from left to right}\/,
as they depend on the value of $\sigma(i)$,
which is $> i$ and hence unknown at time $i$.
(See also the Remark after Step~2 in Section~\ref{subsec.permutations.J.v2}.)

\bigskip

{\bf Step 3: Proof of bijection.}
We prove that the mapping $\sigma \mapsto (\omega,\xi)$ is a bijection by
explicitly describing the inverse map.

First, some preliminaries:
Given a permutation $\sigma \in \Sym_n$,
we define five subsets of $[n]$:
\begin{subeqnarray}
   F & = & \{i \colon\: \sigma(i) > i\}  \;=\; \hbox{positions of excedances}
       \\[1mm]
   F' & = & \{i \colon\: i > \sigma^{-1}(i)\}  \;=\; \hbox{values of excedances}
       \\[1mm]
   G & = & \{i \colon\: \sigma(i) < i\}  \;=\; \hbox{positions of anti-excedances}
       \\[1mm]
   G' & = & \{i \colon\: i < \sigma^{-1}(i)\}  \;=\; \hbox{values of anti-excedances}
       \\[1mm]
   H  & = & \{i \colon\: \sigma(i) = i\}  \;=\;  \hbox{fixed points}
\end{subeqnarray}
Let us observe that
\begin{subeqnarray}
   F \cap F'  & = &  \hbox{cycle double rises}  \\[1mm]
   G \cap G'  & = &  \hbox{cycle double falls}  \\[1mm]
   F \cap G'  & = &  \hbox{cycle valleys}  \\[1mm]
   F' \cap G  & = &  \hbox{cycle peaks}  \\[1mm]
   F \cap G   & = &  \emptyset \\[1mm]
   F' \cap G'   & = &  \emptyset
\end{subeqnarray}
and of course $H$ is disjoint from $F,F',G,G'$.

Let us also recall the notion of an {\em inversion table}\/:
Let $S$ be a totally ordered set of cardinality $k$,
and let $\bm{x} = (x_1,\ldots,x_k)$ be an enumeration of $S$;
then the (left-to-right) inversion table corresponding to $\bm{x}$
is the sequence $\bm{p} = (p_1,\ldots,p_k)$ of nonnegative integers
defined by $p_\alpha = \#\{\beta < \alpha \colon\: x_\beta > x_\alpha \}$.
Note that $0 \le p_\alpha \le \alpha-1$ for all $\alpha \in [k]$,
so there are exactly $k!$ possible inversion tables.
Given the inversion table $\bm{p}$,
we can reconstruct the sequence $\bm{x}$ 
by working from right to left, as follows:
There are $p_k$ elements of $S$ larger than $x_k$,
so $x_k$ must be the $(p_k+1)$th largest element of $S$.
Then there are $p_{k-1}$ elements of $S \setminus \{x_k\}$
larger than $x_{k-1}$,
so $x_{k-1}$ must be the $(p_{k-1}+1)$th largest element
of $S \setminus \{x_k\}$.
And so forth.
[Analogously, the right-to-left inversion table corresponding to $\bm{x}$
is the sequence $\bm{p} = (p_1,\ldots,p_k)$ of nonnegative integers
defined by $p_\alpha = \#\{\beta > \alpha \colon\: x_\beta < x_\alpha \}$.]

With these preliminaries out of the way,
we can now describe the map $(\omega,\xi) \mapsto \sigma$.
Given the 3-colored Motzkin path $\omega$,
we read off which indices $i$ correspond to cycle valleys,
cycle peaks, cycle double falls, cycle double rises, and fixed points;
this allows us to reconstruct the sets $F,F',G,G',H$.
We now use the labels $\xi$ to reconstruct the maps
$\sigma \restrict F \colon\: F \to F'$ and
$\sigma \restrict G \colon\: G \to G'$, as follows:
Let $i_1,\ldots,i_k$ be the elements of $F$ written in increasing order;
then the sequence $j_1,\ldots,j_k$ defined by $j_\alpha = \sigma(i_\alpha)$
is a listing of $F'$ whose inversion table is given by
$p_\alpha = \xi_{i_\alpha} - 1$:
this is the content of \reff{def.xi} in the case $\sigma(i) > i$.
So we can use $\xi \restrict F$ to reconstruct $\sigma \restrict F$.
In a similar way we can use $\xi \restrict G$
to reconstruct $\sigma \restrict G$,
but now we must use the right-to-left inversion table
because of how \reff{def.xi} is written in the case $\sigma(i) < i$.

\bigskip

{\bf Step 4: Translation of the statistics.}
We have already translated the crossing and nesting statistics
\reff{def.ucrossnestjk}
in terms of the heights $h_{i-1}$ and labels $\xi_i$:
see \reff{def.xi.bis} and Lemma~\ref{lemma.crossing}.
And we have also translated the level of a fixed point
in terms of the height $h_{i-1} = h_i$:
see \reff{eq.height.fix}.
These are all the statistics arising in
Theorem~\ref{thm.permutations.Jtype.final1}.

\bigskip

{\bf Step 5: Computation of the weights \reff{def.weights.permutations.Jtype.final1}.}
Using the bijection, we transfer the weights \reff{def.Qn.firstmaster}
from $\sigma$ to $(\omega,\xi)$
and then sum over $\xi$ to obtain the weight $W(\omega)$.
This weight is factorized over the individual steps $s_i$, as follows:
\begin{itemize}
   \item If $s_i$ is a rise starting at height $h_{i-1} = k$
      (so that $i$ is a cycle valley),
      then from \reff{eq.lemma.ucrosscval} and \reff{def.xi.bis} the weight is
\be
   a_k
   \;=\;
   \sum_{\xi_i = 1}^{k+1}  \sfa_{k+1-\xi_i,\,\xi_i-1}
   \;=\;
   \sfa^\star_k
\ee
where $\sfa^\star_k$ was defined in \reff{def.astar}.
   \item If $s_i$ is a fall starting at height $h_{i-1} = k$
      (so that $i$ is a cycle peak and $k \ge 1$),
      then from \reff{eq.lemma.lcrosscpeak} and \reff{def.xi.bis} the weight is
\be
   b_k
   \;=\;
   \sum_{\xi_i = 1}^{k}  \sfb_{k-\xi_i,\,\xi_i-1}
   \;=\;
   \sfb^\star_{k-1}
   \;.
\ee
   \item If $s_i$ is a level step of type 1 at height $h_{i-1} = k$
      (so that $i$ is a cycle double fall and $k \ge 1$),
      then from \reff{eq.lemma.lcrosscpeak} and \reff{def.xi.bis} the weight is
\be
   c^{(1)}_k
   \;=\;
   \sum_{\xi_i = 1}^{k}  \sfc_{k-\xi_i,\,\xi_i-1}
   \;=\;
   \sfc^\star_{k-1}
   \;.
\ee
   \item If $s_i$ is a level step of type 2 at height $h_{i-1} = k$
      (so that $i$ is a cycle double rise and $k \ge 1$),
      then from \reff{eq.lemma.ucrosscdrise} and \reff{def.xi.bis} the weight is
\be
   c^{(2)}_k
   \;=\;
   \sum_{\xi_i = 1}^{k}  \sfd_{k-\xi_i,\,\xi_i-1}
   \;=\;
   \sfd^\star_{k-1}
   \;.
\ee
   \item If $s_i$ is a level step of type 3 at height $h_{i-1} = k$
      (so that $i$ is a fixed point),
      then from \reff{eq.height.fix} the weight is
\be
   c^{(3)}_k
   \;=\;
   \sfe_k
   \;.
\ee
\end{itemize}
Setting $\gamma_k = c^{(1)}_k + c^{(2)}_k + c^{(3)}_k$
and $\beta_k = a_{k-1} b_k$
as instructed in \reff{eq.flajolet.motzkin},
we obtain the weights \reff{def.weights.permutations.Jtype.final1}.
This completes the proof of Theorem~\ref{thm.permutations.Jtype.final1}.
\qed

\bigskip

{\bf Remark.}
Theorem~\ref{thm.connected} on counting connected components in permutations,
as applied to Theorem~\ref{thm.permutations.Jtype.final1},
has an easy proof in our labeled-Motzkin-paths formalism.
{}From \reff{eq.lemma.heights2} we see that
$i$ is a divider (see Section~\ref{subsec.permutations.connected})
if and only if $h_i = 0$.
And this happens if and only if step $s_i$
is either a fall starting at height $h_{i-1} = 1$
or a level step of type~3 at height $h_{i-1} = 0$.
So, giving each connected component a weight $\zeta$
amounts to multiplying $\sfb^\star_0$ and $\sfe_0$ by $\zeta$.
In the J-fraction coefficients \reff{def.weights.permutations.Jtype.final1}
this multiplies $\gamma_0$ and $\beta_1$ by $\zeta$,
exactly as asserted in Theorem~\ref{thm.connected}.
\myendremark

\bigskip

Let us conclude by giving a formula
for the inversion statistic \reff{def.inv}
in terms of the Foata--Zeilberger heights and labels:

\begin{lemma}[Inversion statistic]
   \label{lemma.FZ.inversion}
We have
\be
   \inv(\sigma)
   \;=\;
   \sum\limits_{i} (h_{i-1} + \xi_i - 1)
        \:+\:
   \sum\limits_{i \in \fix} h_{i-1}
   \;.
 \label{eq.lemma.FZ.inversion}
\ee
\end{lemma}

\proof
We use Proposition~\ref{prop.inv} to write $\inv(\sigma)$
in terms of our crossing and nesting statistics,
then \reff{eq.height.fix}, \reff{def.xi.bis} and
\reff{eq.lemma.ucrosscval}--\reff{eq.lemma.lcrosscpeak}
to translate the crossing and nesting statistics to heights and labels.
We get
\begin{subeqnarray}
   \cval + \ucrosscval + 2\,\unestcval
   & = &
   \sum_{i \in \cval} (h_{i-1} + \xi_i)
       \\[2mm]
   \cdrise + \ucrosscdrise + 2\,\unestcdrise
   & = &
   \sum_{i \in \cdrise} (h_{i-1} + \xi_i - 1)
       \\[2mm]
   \cdfall + \lcrosscdfall + 2\,\lnestcdfall
   & = &
   \sum_{i \in \cdfall} (h_{i-1} + \xi_i - 1)
       \\[2mm]
   \lcrosscpeak + 2\,\lnestcpeak
   & = &
   \sum_{i \in \cpeak} (h_{i-1} + \xi_i - 2)
       \\[2mm]
   2\,\psnest
   & = &
   2 \sum_{i \in \fix} h_{i-1}
      \:=\: 2 \sum_{i \in \fix} (h_{i-1} + \xi_i - 1)
   \qquad\qquad
\end{subeqnarray}
and hence
\be
   \inv
   \;=\;
   (\cval - \cpeak) \:+\: \sum\limits_{i} (h_{i-1} + \xi_i - 1)
        \:+\:
   \sum\limits_{i \in \fix} h_{i-1}
   \;.
\ee
Since $\cval = \cpeak$, this proves \reff{eq.lemma.FZ.inversion}.
\qed

A formula equivalent to \reff{eq.lemma.FZ.inversion} was given,
in a different notation, in \cite[eqn.~(8)]{Clarke_97}.

\bigskip

{\bf Final remarks.}
1.  Our definition of the Motzkin path $\omega$ is
essentially the same as that used by
Foata and Zeilberger \cite{Foata_90},
Randrianarivony \cite[Section~2]{Randrianarivony_98b},
Corteel \cite[Section~3.1]{Corteel_07},
and Shin and Zeng \cite[Section~4]{Shin_10} \cite[Section~5]{Shin_12};
the only difference is that we have used three rather than two types
of level steps, for conceptual clarity.
Our definition of the labels $\xi$ is different from the ones in these papers,
but similar in spirit.

2.  Other variants (and other presentations) of the Foata--Zeilberger bijection
can be found in
\cite{Foata_90,Biane_93,DeMedicis_94,Randrianarivony_95,Clarke_97,%
Randrianarivony_98b,Shin_16,Elizalde_18}.
\myendremark

\subsection{Second master J-fraction: \\
     Proof of Theorems~\ref{thm.perms.S}(b),
     \ref{thm.perm.Jtype.v2.weaker0}, \ref{thm.perm.Jtype.v2.weaker0.pq}
     and \ref{thm.permutations.Jtype.final2}}
  \label{subsec.permutations.J.v2}


In this section we will prove the second master J-fraction for permutations
(Theorem~\ref{thm.permutations.Jtype.final2}).
As a corollary we will also obtain Theorem~\ref{thm.perm.Jtype.v2.weaker0.pq},
which is obtained from Theorem~\ref{thm.permutations.Jtype.final2}
by the specialization~\reff{eq.Qn.BIG.specializations.v2};
and from this we will in turn obtain Theorem~\ref{thm.perm.Jtype.v2.weaker0},
as well as Theorem~\ref{thm.perms.S}(b),
which is linked by contraction \reff{eq.contraction_even.coeffs}
to the specialization \reff{def.specialization.xyuv.cyc}
of Theorem \ref{thm.perm.Jtype.v2.weaker0}.

%
%
%
%

Here we need to construct a bijection that will allow us to count the
number of cycles (cyc), which is a global variable.
To do this, we will employ (a slight variant of)
the Biane \cite{Biane_93} bijection;
it is very similar in spirit to the Foata--Zeilberger bijection
used in Section~\ref{subsec.permutations.J},
but organized slightly differently.
The Biane bijection (in our version) maps $\Sym_n$ to the set of
$({\bf A'},{\bf A''},{\bf B'},{\bf B''},
{\bf C^{(1)\prime}},{\bf C^{(1)\prime\prime}},
{\bf C^{(2)\prime}},{\bf C^{(2)\prime\prime}},
{\bf C^{(3)\prime}},{\bf C^{(3)\prime\prime}})$-doubly labeled
3-colored Motzkin paths of length $n$, where
\begin{subeqnarray}
   (A'_k,A''_k)    & = &  (1,1)       \qquad\;\;\,\hbox{for $k \ge 0$}  \\
   (B'_k,B''_k)    & = &  (k,k)       \qquad\;\;\hbox{for $k \ge 1$}  \\
   (C^{(1)\prime}_k, C^{(1)\prime\prime}_k)  & = &  (1,k)
                                      \qquad\;\;\,\hbox{for $k \ge 0$}  \\
   (C^{(2)\prime}_k, C^{(2)\prime\prime}_k)  & = &  (k,1)
                                      \qquad\;\;\,\hbox{for $k \ge 0$}  \\
   (C^{(3)\prime}_k, C^{(3)\prime\prime}_k)  & = &  (1,1)
                                      \qquad\;\;\,\hbox{for $k \ge 0$}
 \label{def.abc.biane.0}
\end{subeqnarray}
Our presentation of this bijection will follow the same steps
as in Section~\ref{subsec.permutations.J}.

\bigskip

{\bf Step 1: Definition of the Motzkin path.}
The Motzkin path $\omega$ associated to a permutation $\sigma \in \Sym_n$
is identical to the one employed in the Foata--Zeilberger bijection.
That is:
\begin{itemize}
   \item If $i$ is a cycle valley, then $s_i$ is a rise.
   \item If $i$ is a cycle peak, then $s_i$ is a fall.
   \item If $i$ is a cycle double fall, then $s_i$ is a level step of type 1.
   \item If $i$ is a cycle double rise, then $s_i$ is a level step of type 2.
   \item If $i$ is a fixed point, then $s_i$ is a level step of type 3.
\end{itemize}
The interpretation of the heights $h_i$ is thus exactly as in
Lemma~\ref{lemma.heights}.

\bigskip

{\bf Step 2: Definition of the labels $\bm{\xi_i = (\xi'_i,\xi''_i)}$.}
\begin{itemize}
   \item If $i$ is a cycle valley, then $\xi'_i = \xi''_i = 1$.
   \item If $i$ is a cycle double fall, then $\xi'_i = 1$ and
    $\xi''_i = 1 + \#\{k \colon\: k > i \hbox{ and }  \sigma(k) < \sigma(i) \}$.
   \item If $i$ is a cycle double rise, then $\xi''_i = 1$ and
    \hbox{$\xi'_i = 1 + \#\{k \colon\: k < \sigma^{-1}(i) \hbox{ and }  \sigma(k) > i \}$}.
   \item If $i$ is a cycle peak, then
    $\xi'_i = 1 + \#\{k \colon\: k < \sigma^{-1}(i) \hbox{ and }  \sigma(k) > i \}$
    and
    $\xi''_i = 1 + \#\{k \colon\: k > i \hbox{ and }  \sigma(k) < \sigma(i) \}$.
   \item If $i$ is a fixed point, then $\xi'_i = \xi''_i = 1$.
\end{itemize}

These labels have a nice interpretation
in terms of the bipartite digraph employed in the proof of
Lemma~\ref{lemma.heights} (Figure~\ref{fig.FZ}).
First recall that $h_{i-1}$ equals
the number of unconnected dots in the top row after $i-1$ steps,
and also equals the number of unconnected dots in the bottom row
after $i-1$ steps.
We then look at what happens at stage $i$:
\begin{itemize}
   \item If $i$ is a cycle valley, then at stage $i$ we add no arrows.
       Since no choices are being made at this stage, we set
       $\xi'_i = \xi''_i = 1$.
   \item If $i$ is a cycle double fall, then at stage $i$ we add an arrow
       from $i$ on the top row to an unconnected dot $j'$ on the bottom row,
       where $j = \sigma(i) < i$;
       then $\xi''_i$ is the index of the unconnected dot $j'$
       among all the unconnected dots on the bottom row.
       Since no unconnected dot on the top row was touched,
       we set $\xi'_i = 1$.
   \item Similarly, if $i$ is a cycle double rise, then at stage $i$
       we add an arrow from an unconnected dot $j$ on the top row
       to $i'$ on the bottom row, where $j = \sigma^{-1}(i) < i$;
       then $\xi'_i$ is the index of the unconnected dot $j$
       among all the unconnected dots on the top row.
       Since no unconnected dot on the bottom row was touched,
       we set $\xi''_i = 1$.
   \item If $i$ is a cycle peak, then we add two arrows:
       from $i$ on the top row to the unconnected dot $j'$ on the bottom row,
       where $j = \sigma(i) < i$;
       and also from the unconnected dot $k$ on the top row
       to $i'$ on the bottom row, where $k = \sigma^{-1}(i) < i$.
       Then $\xi'_i$ (resp.\ $\xi''_i$) is the index of $k$ (resp.\ $j'$)
       among the unconnected dots on the top (resp.\ bottom) row.
   \item If $i$ is a fixed point, then at stage $i$ we add an arrow
       $i \to i'$.
       Since no choices are being made at this stage, we set
       $\xi'_i = \xi''_i = 1$.
\end{itemize}
This interpretation shows in particular that
the inequalities \reff{eq.xi.ineq}/\reff{def.abc.biane.0} are satisfied.

\medskip

{\bf Remark.}
The Biane labels $(\xi'_i,\xi''_i)$
are related to the Foata--Zeilberger labels $\xi_i$
[defined in \reff{def.xi}] as follows:
\begin{subeqnarray}
   & &
   i \in \cdfall \cup \cpeak\colon\:
   \xi''_i(\textrm{Biane}) = \xi_i(\textrm{FZ})
       \slabel{eq.FZ.Biane.a} \\[2mm]
   & &
   i \in \cdrise \cup \cpeak\colon\:
   \xi'_i(\textrm{Biane}) = \xi_{\sigma^{-1}(i)}(\textrm{FZ})
      \hbox{ or equivalently }
      \xi_{i}(\textrm{FZ}) = \xi'_{\sigma(i)}(\textrm{Biane})
              \nonumber \\[-1mm]
 \slabel{eq.FZ.Biane.b}
 \label{eq.FZ.Biane}
\end{subeqnarray}
since if $i$ is a cycle double fall or cycle peak, then $\sigma(i) < i$,
while if $i$ is a cycle double rise or cycle peak, then $\sigma^{-1}(i) < i$.
\myendremark

{\bf Step 3: Proof of bijection.}
The foregoing interpretation shows how to build the bipartite digraph,
and hence reconstruct the permutation $\sigma$,
by successively reading the steps $s_i$ and labels $\xi_i$.
Specifically, at stage $i$ one proceeds as follows \cite[p.~280]{Biane_93}:
\begin{itemize}
   \item  If $s_i$ is a rise (corresponding to $i$ being a cycle valley),
      then we add no arrows.
      [In this case we necessarily have $\xi_i = (1,1)$.]
   \item  If $s_i$ is a level step of type 1
      (corresponding to $i$ being a cycle double fall),
      then $\xi_i = (1,m)$ for some $m \in [h_{i-1}]$,
      and we add an arrow from $i$ on the top row
      to the $m$th (from left to right) unconnected dot on the bottom row.
   \item  If $s_i$ is a level step of type 2
      (corresponding to $i$ being a cycle double rise),
      then $\xi_i = (l,1)$ for some $l \in [h_{i-1}]$,
      and we add an arrow from
      the $l$th (from left to right) unconnected dot on the top row
      to $i'$ on the bottom row.
   \item  If $s_i$ is a fall (corresponding to $i$ being a cycle peak),
      then $\xi_i = (l,m)$ for some $l,m \in [h_{i-1}]$,
      and we add two arrows:
      one going from $i$ on the top row to the $m$th (from left to right)
      unconnected dot on the bottom row;
      and the other going from the $l$th (from left to right)
      unconnected dot on the top row to $i'$ on the bottom row.
   \item If $s_i$ is a level step of type 3
      (corresponding to $i$ being a fixed point),
      we add an arrow from $i$ on the top row to $i'$ on the bottom row.
      [In this case we necessarily have $\xi_i = (1,1)$.]
\end{itemize}
Clearly, once a dot becomes the source or sink of an arrow,
it plays no further role in the construction
and in particular receives no further arrows.
Moreover, since $h_n = 0$, at the end of the construction
there are no unconnected dots.
The final result of the construction thus corresponds to a bijection
between $\{1,\ldots,n\}$ and $\{1',\ldots,n'\}$,
or in other words to a permutation $\sigma \in \Sym_n$.

\bigskip

{\bf Step 4: Translation of the statistics.}

\begin{lemma}[Crossing and nesting statistics]
   \label{lemma.crossing.2a}
We have
\begin{eqnarray}
   h_{i-1} - \xi''_i  & = & \lcross(i,\sigma)
         \hspace*{2.9cm}\textrm{\rm if $i \in \cpeak \,\cup\, \cdfall$}
     \label{eq.lemma.2a.lcrosscpeak} \\[1mm]
   \xi''_i - 1 & = & \lnest(i,\sigma)
         \hspace*{3.05cm}\textrm{\rm if $i \in \cpeak \,\cup\, \cdfall$}
     \label{eq.lemma.2a.lnestcpeak} \\[1mm]
   h_i - 1  & = &  \ucross(i,\sigma) + \unest(i,\sigma)
         \quad\textrm{\rm if $i \in \cval \,\cup\, \cdrise$}
     \label{eq.lemma.2a.ucrosscval} \\[1mm]
   \xi'_i - 1  & = & \unest(\sigma^{-1}(i),\sigma)
         \hspace*{1.95cm}\textrm{\rm if $i \in \cpeak \,\cup\, \cdrise$}
     \label{eq.lemma.2a.unestcpeak} \\[1mm]
   h_{i-1}  \:=\:  h_i  & = &  \lev(i,\sigma)
         \hspace*{3.4cm}\textrm{\rm if $i \in \fix$}
     \label{eq.lemma.2a.fix}
\end{eqnarray}
\end{lemma}

\noindent
Note that \reff{eq.lemma.2a.ucrosscval} is written in terms of $h_i$,
while \reff{eq.lemma.2a.lcrosscpeak} is written in terms of $h_{i-1}$.

\proof
(a) If $i$ is a cycle peak or a cycle double fall, then
$h_{i-1} - \xi''_i = \lcross(i,\sigma)$
exactly as in \reff{eq.Bxi.cpeak.cdfall},
and $\xi''_i - 1 = \lnest(i,\sigma)$
exactly as in \reff{def.xi.bis}.

(b) If $i$ is a cycle valley or a cycle double rise,
then $h_i - 1 = \ucross(i,\sigma) + \unest(i,\sigma)$
is an immediate consequence of \reff{eq.lemma.heights2.a}
and the definitions (\ref{def.ucrossnestjk}a,b).

(c) If $i$ is a cycle peak or a cycle double rise, then
$\xi'_i - 1 = \unest(\sigma^{-1}(i),\sigma)$
is an immediate consequence of the definition of $\xi'_i$
and the definition (\ref{def.ucrossnestjk}b).

(d) If $i$ is a fixed point, then $h_{i-1} = h_i = \lev(i,\sigma)$
was proven in \reff{eq.height.fix}.
\qed

It is instructive to compare this result with Lemma~\ref{lemma.crossing}.
For cycle peaks and cycle double falls,
the results here for $\lcross$ and $\lnest$
are identical to \reff{eq.lemma.lcrosscpeak} and \reff{def.xi.bis}
but with $\xi_i$ replaced by $\xi''_i$;
this is, of course, an immediate consequence of \reff{eq.FZ.Biane.a}.
For cycle valleys and cycle double rises, by contrast,
here we do not learn about $\ucross$ and $\unest$ individually,
but only about their sum.
And finally, for cycle peaks and cycle double rises,
we obtain $\unest$, but evaluated at $\sigma^{-1}(i)$ rather than at $i$;
this is an immediate consequence of \reff{eq.FZ.Biane.b}.

Finally, we come to the counting of cycles (cyc).
We use the term {\em cycle closer}\/ to denote the largest element
in a non-singleton cycle.
Obviously every non-singleton cycle has precisely one cycle closer.
A cycle closer is always a cycle peak, but not conversely:
for instance, in the cycle $(1324)$, both 3 and 4 are cycle peaks,
but only 4 is a cycle closer.
So we need to know which cycle peaks are cycle closers,
or at least how many of the former are the latter.
The answer is as follows:

\begin{lemma}[Counting of cycles]
   \label{lemma.cycles.2}
Fix $i \in [n]$,
and fix $(s_1,\ldots,s_{i-1})$ and $(\xi_1,\ldots,\xi_{i-1})$.
Consider all permutations $\sigma \in \Sym_n$
that have those given values for the first $i-1$ steps and labels
and for which $i$ is a cycle peak.
Then:
\begin{itemize}
   \item[(a)] The value of $\xi_i = (\xi'_i,\xi''_i)$ completely determines
      whether $i$ is a cycle closer or not.
   \item[(b)] For each value $\xi'_i \in [h_{i-1}]$
      there is precisely one value $\xi''_i \in [h_{i-1}]$
      that makes $i$ a~cycle closer, and conversely.
\end{itemize}
\end{lemma}

\proof
We use once again the bipartite digraph
of Figure~\ref{fig.FZ}(a),
and let us also draw a vertical dotted line
(with an upwards arrow) to connect each pair $j' \to j$.
Now consider the restriction of this digraph to
the vertex set $\{1,\ldots,{i-1},1',\ldots,{(i-1)'}\}$:
as discussed in Step~3, this restriction can be reconstructed
from the steps $(s_1,\ldots,s_{i-1})$
and the labels $(\xi_1,\ldots,\xi_{i-1})$.
The connected components of this restriction are of two types:
complete directed cycles and directed open chains;
they correspond to cycles of $\sigma$
whose cycle closers are, respectively, $\le i-1$ and $> i-1$.
Each directed open chain runs from an unconnected dot on the bottom row
to an unconnected dot on the top row.

Now suppose that $i$ is a cycle peak.
Then at stage $i$ we add two arrows:
from $i$ on the top row to an unconnected dot $j'$ on the bottom row;
and also from an unconnected dot $k$ on the top row
to $i'$ on the bottom row.
Here $\xi'_i$ (resp.\ $\xi''_i$) is the index of $k$ (resp.\ $j'$)
among the unconnected dots on the top (resp.\ bottom) row.

Now the point is simply this: $i$ is a cycle closer if and only if
$j'$ and $k$ belong to the {\em same}\/ directed open chain
(with $j'$ being its starting point and $k$ being its ending point).
So for each value $\xi'_i \in [h_{i-1}]$
there is precisely one value $\xi''_i \in [h_{i-1}]$
that makes $i$ a cycle closer, and conversely.
\qed

\bigskip

{\bf Step 5: Computation of the weights
  \reff{def.weights.permutations.Jtype.final2}.}
Using the bijection, we transfer the weights \reff{def.Qn.secondmaster}
from $\sigma$ to $(\omega,\xi)$
and then sum over $\xi$ to obtain the weight $W(\omega)$.
This weight is factorized over the individual steps $s_i$, as follows:
\begin{itemize}
   \item If $s_i$ is a rise starting at height $h_{i-1} = h_i - 1 = k$
      (so that $i$ is a cycle valley),
      then necessarily $\xi_i = (1,1)$,
      and it follows from \reff{eq.lemma.2a.ucrosscval}
      that the weight is
\be
   a_k
   \;=\;
   \sfa_k
   \;.
\ee
   \item If $s_i$ is a fall starting at height $h_{i-1} = k$
      (so that $i$ is a cycle peak and $k \ge 1$),
      then for each choice of $\xi''_i \in [k]$
      there are $k$ possible choices of $\xi'_i$,
      of which one closes a cycle and the rest don't
      (Lemma~\ref{lemma.cycles.2}).
      Therefore, using \reff{eq.lemma.2a.lcrosscpeak}
      and \reff{eq.lemma.2a.lnestcpeak}, the weight is
\be
   b_k
   \;=\;
   (\lambda+k-1) \sum_{\xi''_i = 1}^{k}  \sfb_{k-\xi''_i,\,\xi''_i-1}
   \;=\;
   (\lambda+k-1) \, \sfb^\star_{k-1}
   \;.
\ee
   \item If $s_i$ is a level step of type 1 at height $h_{i-1} = k$
      (so that $i$ is a cycle double fall and $k \ge 1$),
      then from \reff{eq.lemma.2a.lcrosscpeak}
      and \reff{eq.lemma.2a.lnestcpeak} the weight is
\be
   c^{(1)}_k
   \;=\;
   \sum_{\xi''_i = 1}^{k}  \sfc_{k-\xi''_i,\,\xi''_i-1}
   \;=\;
   \sfc^\star_{k-1}
   \;.
\ee
   \item If $s_i$ is a level step of type 2 at height $h_{i-1} = h_i = k$
      (so that $i$ is a cycle double rise and $k \ge 1$),
      then from \reff{eq.lemma.2a.ucrosscval}
      and \reff{eq.lemma.2a.unestcpeak} the weight is
\be
   c^{(2)}_k
   \;=\;
   \sum_{\xi'_i = 1}^{k}  \sfd_{k-1,\,\xi'_i-1}
   \;=\;
   \sfd^\natural_{k-1}
   \;.
\ee
   \item If $s_i$ is a level step of type 3 at height $h_{i-1} = k$
      (so that $i$ is a fixed point),
      then from \reff{eq.lemma.2a.fix} the weight is
\be
   c^{(3)}_k
   \;=\;
   \lambda \sfe_k
   \;.
\ee
\end{itemize}
Setting $\gamma_k = c^{(1)}_k + c^{(2)}_k + c^{(3)}_k$
and $\beta_k = a_{k-1} b_k$
as instructed in \reff{eq.flajolet.motzkin},
we obtain the weights \reff{def.weights.permutations.Jtype.final2}.
This completes the proof of Theorem~\ref{thm.permutations.Jtype.final2}.
\qed

\bigskip

{\bf Remark.}
Theorem~\ref{thm.connected} on counting connected components in permutations,
as applied to Theorem~\ref{thm.permutations.Jtype.final2},
has an easy proof in our labeled-Motzkin-paths formalism;
the argument is identical to the one presented in
Section~\ref{subsec.permutations.J}
for the first master J-fraction.
\myendremark

\bigskip

Let us also observe that Biane \cite[eqn.~(3.2.6)]{Biane_93}
has given a formula, in terms of his heights and labels,
for the inversion statistic \reff{def.inv};
in our notation it is:

\begin{lemma}[Inversion statistic]
   \label{lemma.biane.inversion}
We have
\begin{subeqnarray}
   \inv(\sigma)
   & = &
   2 \sum\limits_{i \in \fix} h_{i-1}
        \:+\:
   \sum\limits_{i \in \cdrise} (h_{i-1} + \xi'_i - 1)
        \:+\:
   \sum\limits_{i \in \cdfall} (h_{i-1} + \xi''_i - 1)
        \qquad  \nonumber \\[1mm]
   & & \qquad\qquad\:
        \:+\:
   \sum\limits_{i \in \cpeak} (2h_{i-1} + \xi'_i + \xi''_i - 3)
 \slabel{eq.lemma.biane.inversion.a}
         \\[2mm]
   & = &
   \sum\limits_{i} (h_{i-1} + \xi'_i + \xi''_i - 2)
        \:+\:
   \sum\limits_{i \in \fix} h_{i-1}
   \;.
 \slabel{eq.lemma.biane.inversion.b}
\end{subeqnarray}
\end{lemma}

Since Biane's proof of this formula is a bit complicated (by induction on $n$),
let us give a simple proof,
following ideas of Elizalde \cite[p.~6, item~(v)]{Elizalde_18}:

\proof
The number of inversions is the number of crossings of lines
in the bipartite digraph of Figure~\ref{fig.FZ}(a).
We now count the crossings of lines $(L,L')$
according to the stage $i$ at which the {\em first}\/ of the two lines $L,L'$
appears in the ``reading'' of the digraph from left to right,
as shown in Figure~\ref{fig.FZ}(c):

1) A fixed point $i$ corresponds to a vertical line $i \to i'$
in the bipartite digraph;
and since at stage $i$
there are $h_{i-1} = h_i$ unconnected dots $j < i$ in the top row
and also $h_{i-1} = h_i$ unconnected dots $j'$ with $j < i$ in the bottom row,
each of these unconnected dots will later connect
to a vertex $k'$ (resp.\ $k$) with $k > i$,
thereby crossing the  vertical line $i \to i'$.
This gives the term $2 \sum\limits_{i \in \fix} h_{i-1}$.

2) A cycle double rise $i$ corresponds to a line $j \to i'$
where $j < i$ (and there is not yet any line emanating from $i$).
The line $j \to i'$ will later be crossed by lines emanating from
$\xi'_i - 1$ unconnected dots $< j$ on the top row,
and lines arriving at $h_{i-1}$ unconnected dots $< i'$ on the bottom row.
This gives the term $\sum\limits_{i \in \cdrise} (h_{i-1} + \xi'_i - 1)$.

3) Similarly, a cycle double fall $i$ corresponds to a line $i \to j'$
where $j < i$ (and there is not yet any line arriving at $i'$).
The line $i \to j'$ will later be crossed by lines emanating from
$h_{i-1}$ unconnected dots $< i$ on the top row,
and lines arriving at $\xi''_i - 1$ unconnected dots $< j'$ on the bottom row.
This gives the term $\sum\limits_{i \in \cdfall} (h_{i-1} + \xi''_i - 1)$.

4) A cycle peak $i$ corresponds to a pair of lines $i \to j'$
and $k \to i'$, where $j,k < i$.
The line $i \to j'$ will later be crossed by
lines emanating from $h_{i-1}$ unconnected dots $< i$ on the top row,
and lines arriving at $\xi''_i - 1$ unconnected dots $< j'$ on the bottom row.
Similarly, the line $k \to i'$ will later be crossed by
lines emanating from $\xi'_i - 1$ unconnected dots $< k$ on the top row,
and lines arriving at $h_{i-1}$ unconnected dots $< i'$ on the bottom row.
However, this double-counts the crossing between the lines
$i \to j'$ and $k \to i'$, so we need to subtract 1.
This gives the term
$\sum\limits_{i \in \cpeak} (2h_{i-1} + \xi'_i + \xi''_i - 3)$.

5) A cycle valley $i$ does not correspond to any line at stage $i$.

Putting this all together gives \reff{eq.lemma.biane.inversion.a}.
This can be rewritten as \reff{eq.lemma.biane.inversion.b}
once we observe that
\be
   \sum\limits_{i \in \cval} h_{i-1}
   \;=\;
   \sum\limits_{i \in \cpeak} (h_{i-1} - 1)
 \label{eq.cvalcpeak}
\ee
(by pairing rises with falls on the Motzkin path)
and recall that the ``null'' values of $\xi'_i$ or $\xi''_i$ are 1.
\qed

\medskip

{\bf Remarks.}
1.  A formula equivalent to \reff{eq.lemma.biane.inversion.b} was given,
in a different notation, in \cite[after eqn.~(8)]{Clarke_97}.

2.  The result of Lemma~\ref{lemma.biane.inversion} can also be used
to give an alternate proof of Proposition~\ref{prop.inv}, as follows:
Start from \reff{eq.lemma.biane.inversion.a}
and rewrite it slightly by using
\be
   \sum\limits_{i \in \cpeak}  h_{i-1}
   \;=\;
   \sum\limits_{i \in \cval}  h_i
   \;.
\ee
Then translate the heights and labels back to crossing and nesting statistics,
using \reff{eq.lemma.2a.lcrosscpeak}--\reff{eq.lemma.2a.fix};
this yields
\begin{eqnarray}
   \hspace*{-1.5cm}
   \inv
   & \!=\! & \cval + \cdrise + \cdfall + \ucross + \lcross
                    + 2\,\lnest + 2\,\psnest + \unest
          \nonumber \\[1mm]
   & & \qquad
    + \sum_{i \in \cdrise \cup \cpeak} \unest(\sigma^{-1}(i),\sigma)
   \;.
 \label{eq.newproof.prop.inv}
\end{eqnarray}
But $i \in \cdrise \cup \cpeak$ $\iff$
$\sigma^{-1}(i) < i$ $\iff$
$\sigma^{-1}(i) \in \cdrise \cup \cval$,
so
\be
   \sum_{i \in \cdrise \cup \cpeak} \unest(\sigma^{-1}(i),\sigma)
   \;=\;
   \sum_{j \in \cdrise \cup \cval} \unest(j,\sigma)
   \;=\;
   \unest(\sigma)
   \;.
 \label{eq.Biane.sigma-1}
\ee
Combining this with \reff{eq.newproof.prop.inv}
proves \reff{eq.prop.inv}.

3.  Note the close similarity between
the Foata--Zeilberger formula \reff{eq.lemma.FZ.inversion}
and the Biane formula \reff{eq.lemma.biane.inversion.b}.
Indeed, by using \reff{eq.FZ.Biane} and \reff{eq.Biane.sigma-1}
together with \reff{def.xi.bis},
it is straightforward to show that the right-hand sides
of \reff{eq.lemma.FZ.inversion} and \reff{eq.lemma.biane.inversion.b}
are equal.\!
\myendremark


Let us explain, finally, why it is {\em not}\/ possible
to use the Biane bijection to combine the counting of cycles
with the counting of inversions.
The trouble is that the last term in \reff{eq.lemma.biane.inversion.a}
--- that is, the sum over cpeak ---
involves both $\xi'_i$ and $\xi''_i$.
By Lemma~\ref{lemma.cycles.2} we know that for each possible value of $\xi'_i$
there is exactly one value of $\xi''_i$ that makes $i$ a cycle closer ---
{\em but we don't know which one it is}\/.
Therefore, we are unable to evaluate the sum over $\xi'_i,\xi''_i \in [k]$
of (for instance) a weight $q^{\xi'_i+\xi''_i}$
times a weight $\lambda$ for each cycle closer.
For instance, suppose that the cycle closer occurs when $\xi'_i = \xi''_i$:
then the sum is
\be
   \biggl( \sum_{j=1}^k q^j \biggr)^{\! 2}
     \:+\:
   (\lambda-1) \sum_{j=1}^k q^{2j}
   \;.
\ee
But if we suppose, by contrast, that the cycle closer occurs when
$\xi'_i = k+1 - \xi''_i$,  then the sum is
\be
   \biggl( \sum_{j=1}^k q^j \biggr)^{\! 2}
     \:+\:
   (\lambda-1) k q^{k+1}
   \;.
\ee
Of course, as explained in Section~\ref{subsec.permutations.inv},
this inability to combine
the counting of cycles with the counting of inversions
is not merely a limitation of our method of proof,
but is inherent in the problem:
the weight $q^{\inv(\sigma)} \lambda^{\cyc(\sigma)}$
gives rise to a J-fraction with coefficients that are rational functions
rather than polynomials.

\section{Set partitions: Proofs}   \label{sec.proofs.setpartitions}

\subsection{S-fraction:  Proof of Theorem~\ref{thm.setpartitions}}
   \label{subsec.setpartitions.S}

It is immediate from the definitions
\reff{eq.eulerian.fourvar.contfrac} and \reff{eq.setpartitions.contfrac}
that $B_n(x,y,v) = P_n(x,y,0,v)$.
We shall prove Theorem~\ref{thm.setpartitions}
by translating the interpretation \reff{eq.eulerian.fourvar.cyc}
of $P_n(x,y,0,v)$ to set~partitions via a suitable bijection.

We define a mapping from set partitions to permutations as follows:
Given a set partition $\pi \in \Pi_n$,
we define the permutation $\sigma \in \Sym_n$
such that the disjoint cycles of~$\sigma$
are the blocks of $\pi$, each traversed in increasing order
(with the largest element of~course followed by the smallest element).
The mapping $\pi \mapsto \sigma$
is clearly a bijection of $\Pi_n$ onto $\Sym_n^\star$,
where $\Sym_n^\star$ denotes the subset of $\Sym_n$ consisting of permutations
in which each cycle of length $\ell \ge 2$ contains
precisely one cycle peak (namely, the cycle maximum),
one cycle valley (namely, the cycle minimum),
$\ell-2$ cycle double rises, and no cycle double falls.
Observe now that these are precisely the permutations in which
each cycle has exactly one non-excedance,
i.e.\ in which $\exc(\sigma) + \cyc(\sigma) = n$.
So setting $u=0$ in \reff{eq.eulerian.fourvar.cyc}
corresponds to restricting the sum to $\Sym_n^\star$,
which we can then pull back to $\Pi_n$ via the bijection.
We clearly have $\cyc(\sigma) = |\pi|$.
To finish the proof, we need only interpret $\erec(\sigma)$ in terms of $\pi$.

But this is easy.
By definition, $\sigma(i)$ is
the next-larger element of the block $B_\pi(i)$ after $i$,
in case $i$ is not the largest element of $B_\pi(i)$,
and is the smallest element of $B_\pi(i)$,
in case $i$ is the largest element of $B_\pi(i)$.
So $\sigma(i) = i$ if and only if $i$ is a singleton;
$\sigma(i) < i$ if and only if $i$ is the largest element of
a non-singleton block;
and $\sigma(i) > i$ if and only if $i$ is a non-largest element of a
(necessarily non-singleton) block.

When is an index $i$ an exclusive record of $\sigma$?
An exclusive record is simply a record that is not a fixed point,
or equivalently a record that is an excedance.
So we eliminate the non-excedances by defining
\be
   \sigma'(i)
   \;=\;
   \begin{cases}
      \sigma(i)  &  \textrm{if $\sigma(i) > i$} \\[1mm]
      0          &  \textrm{if $\sigma(i) \le i$}
   \end{cases}
\ee
Then the exclusive records of the permutation $\sigma$
are the same as the nonzero records of the word $\sigma'$,
i.e.\ the indices $i$ such that $\sigma'(i) \neq 0$
and $\sigma'(j) < \sigma'(i)$ for all $j < i$.
[Note that the only repeated elements of $\sigma'$ are 0,
 so there is no need to distinguish between records and strict records.]
But this is exactly how we have defined ``exclusive record''
for a set partition $\pi$.
\qed

\medskip

{\bf Remark.}
It would be nicer if we could use the interpretation
\reff{eq.eulerian.fourvar.arec}
instead of \reff{eq.eulerian.fourvar.cyc} on the permutation side,
i.e.\ employing $\arec$ instead of $\cyc$,
since we have available better refinements
for the permutation polynomials that do not include the statistic $\cyc$:
compare Theorem~\ref{thm.perm.Jtype} with \ref{thm.perm.Jtype.v2.weaker0},
or \ref{thm.perm.pq.Jtype.BIG} with \ref{thm.perm.Jtype.v2.weaker0.pq},
or \ref{thm.permutations.Jtype.final1} with \ref{thm.permutations.Jtype.final2}.
So we would like to find an injection of partitions of $[n]$
into permutations of $[n]$ in which blocks map to antirecords,
and in which the image permutations are precisely those in which
every index is either an excedance or an antirecord
(that is, in which there are no $\nrcpeak$, $\nrcdfall$ or $\nrfix$).
But we have been unable to find such a mapping.
%
%
%
%
%
\myendremark


\subsection{First master J-fraction:
   Proof of Theorems~\ref{thm.setpartitions.Jtype},
        \ref{thm.setpartitions.Jtype.pq.refined}
        and \ref{thm.setpartitions.Jtype.final1}}
   \label{subsec.setpartitions.J}

In this section we will prove the first master J-fraction for set partitions
(Theorem~\ref{thm.setpartitions.Jtype.final1}).
As a consequence we will also obtain
Theorem~\ref{thm.setpartitions.Jtype.pq.refined},
which is obtained from Theorem~\ref{thm.setpartitions.Jtype.final1}
by the specialization \reff{eq.specialization.thm.setpartitions.Jtype.pq};
and Theorem~\ref{thm.setpartitions.Jtype},
which is a special case of Theorem~\ref{thm.setpartitions.Jtype.pq.refined}.
We will also obtain a second proof of Theorem~\ref{thm.setpartitions},
which is linked by contraction \reff{eq.contraction_even.coeffs}
to the specialization $x_1 = x_2$, $y_1 = y_2$, $v_1 = v_2$
of Theorem \ref{thm.setpartitions.Jtype}.

To prove Theorem~\ref{thm.setpartitions.Jtype.final1},
we will employ the Kasraoui--Zeng \cite{Kasraoui_06} bijection,
which is a variant of one proposed earlier by Flajolet \cite{Flajolet_80}.
(We will discuss the Flajolet bijection in the next subsection.)
However, before introducing this bijection we need first to reinterpret
the polynomial $B_n(\bsfa,\bsfb,\bsfd,\bsfe)$ defined in \reff{def.Bnhat}
by reversing the order of the vertices $1,\ldots,n$.
(The reason for this somewhat embarrassing reversal
 will be discussed after the proof.)
So, given any set partition $\pi \in \Pi_n$,
we define $\pitilde \in \Pi_n$ to be the reversal of $\pi$,
i.e.\ the image of $\pi$ under the map $i \mapsto \itilde \eqdef n+1-i$.
We then define reversals of the statistics \reff{def.cr.ne.k.pi}
employed in \reff{def.Bnhat}:
\begin{subeqnarray}
   \crrtilde(k,\pi)
   \;\eqdef\;
   \crr(\ktilde,\pitilde)
   & = &
   \#\{ i<j<k<l \colon\: (i,k) \in \scrg_\pi \hbox{ and } (j,l) \in \scrg_\pi \}
         \qquad \\[2mm]
   \neetilde(k,\pi)
   \;\eqdef\;
   \nee(\ktilde,\pitilde)
   & = &
   \#\{ i<j<k<l \colon\: (i,l) \in \scrg_\pi \hbox{ and } (j,k) \in \scrg_\pi \}
         \qquad \\[2mm]
   \qnetilde(k,\pi)
   \;\eqdef\;
   \qne(\ktilde,\pitilde)
   & = &
   \#\{ i<k<l \colon\: (i,l) \in \scrg_\pi  \}
 \label{def.cr.ne.kpi.reversed}
\end{subeqnarray}
Thus, $\crrtilde$ and $\neetilde$ are like $\crr$ and $\nee$
but put the distinguished index in third rather than second position.
On the other hand, the definition of $\qne$ is reversal-invariant,
so that in fact $\qnetilde(k,\pi) = \qne(k,\pi)$.
Now summing over $\pi \in \Pi_n$
is of course equivalent to summing over $\pitilde \in \Pi_n$;
and reversal interchanges openers with closers.
It follows that the polynomial 
$B_n(\bsfa,\bsfb,\bsfd,\bsfe)$ defined in \reff{def.Bnhat}
can equivalently be written as
\begin{eqnarray}
   & &  \hspace*{-8mm}
   B_n(\bsfa,\bsfb,\bsfd,\bsfe)
   \;=\;
       \nonumber \\[4mm]
   & &
   \sum_{\pi \in \Pi_n}
   \;
   \prod\limits_{i \in {\rm closers}}  \! \sfa_{\crrtilde(i,\pi),\, \neetilde(i,\pi)}
   \prod\limits_{i \in {\rm openers}}  \!\! \sfb_{\qne(i,\pi)}
   \prod\limits_{i \in {\rm insiders}} \!\!  \sfd_{\crrtilde(i,\pi),\, \neetilde(i,\pi)}
   \prod\limits_{i \in {\rm singletons}}  \!\!\! \sfe_{\qne(i,\pi)}
   \;.
        \nonumber \\[-2mm]
 \label{def.Bnhat.reversed}
\end{eqnarray}
We will employ the reinterpretation \reff{def.Bnhat.reversed} in our proof,
because its groupings of closers-and-insiders and openers-and-singletons
are better adapted to the Kasraoui--Zeng bijection
than the groupings of openers-and-insiders and closers-and-singletons
employed in our original (and in our opinion more natural)
definition \reff{def.Bnhat}.

Let us now define the Kasraoui--Zeng \cite{Kasraoui_06} bijection,
which is a bijection from $\Pi_n$ to the set of
$({\bf A},{\bf B},{\bf C}^{(1)},{\bf C}^{(2)})$-labeled
2-colored Motzkin paths of length $n$, where
\begin{subeqnarray}
   A_k        & = &  1         \quad\hbox{for $k \ge 0$}  \\
   B_k        & = &  k         \quad\hbox{for $k \ge 1$}  \\
   C_k^{(1)}  & = &  k         \quad\hbox{for $k \ge 0$}  \\
   C_k^{(2)}  & = &  1         \quad\hbox{for $k \ge 0$}
 \label{def.abc.KZ}
\end{subeqnarray}
As before, we will begin by explaining how the Motzkin path $\omega$ is defined;
then we will explain how the labels $\xi$ are defined;
next we will prove that the mapping is indeed a bijection;
next we will translate the various statistics from
$\Pi_n$ to our labeled Motzkin paths;
and finally we will sum over labels $\xi$ to obtain the weight $W(\omega)$
associated to a Motzkin path $\omega$,
which upon applying \reff{eq.flajolet.motzkin}
will yield Theorem~\ref{thm.setpartitions.Jtype.final1}.

\bigskip

{\bf Step 1: Definition of the Motzkin path.}
Given a set partition $\pi \in \Pi_n$,
we~classify the indices $i \in [n]$
in the usual way as opener, closer, insider or singleton.
We then define a path $\omega = (\omega_0,\ldots,\omega_n)$
starting at $\omega_0 = (0,0)$ and ending at $\omega_n = (n,0)$,
with steps $s_1,\ldots,s_n$, as follows:
\begin{itemize}
   \item If $i$ is an opener, then $s_i$ is a rise.
   \item If $i$ is a closer, then $s_i$ is a fall.
   \item If $i$ is an insider, then $s_i$ is a level step of type 1.
   \item If $i$ is a singleton, then $s_i$ is a level step of type 2.
\end{itemize}
The interpretation of the heights $h_i$ is almost immediate
from this definition:

\begin{lemma}
   \label{lemma.heights.KZ}
For $i \in \{0,\ldots,n\}$,
$h_i$ is the number of blocks that are ``started but unfinished''
after stage $i$, i.e.
\be
   h_i  \;=\;  \#\{ B \in \pi \colon\: \min B \le i < \max B \}
   \;.
\ee
\end{lemma}

\noindent
In particular, it follows that $\omega$ is indeed a Motzkin path,
i.e.\ all the heights $h_i$ are nonnegative and $h_n = 0$.

\bigskip

{\bf Step 2: Definition of the labels $\bm{\xi_i}$.}
If $i$ is an opener or a singleton,
we set $\xi_i = 1$ as required by \reff{def.abc.KZ}.
If $i$ is an insider or a closer,
we look at the $h_{i-1}$ blocks $B_1,\ldots,B_{h_{i-1}}$
that are ``started but unfinished'' after stage $i-1$
(note that we must have $h_{i-1} \ge 1$).
For each $j \in [h_{i-1}]$,
let $y_j$ be the maximal element of $B_j \cap [i-1]$;
it is necessarily an opener or insider in $B_j$,
and its successor in $B_j$ will be $\ge i$.
(Kasraoui and Zeng \cite{Kasraoui_06} call $y_j$ the ``vacant vertex''.)
We order the blocks $B_j$ so that $y_1 < y_2 < \ldots < y_{h_{i-1}}$.
Then the vertex $i$ belongs to precisely one of these blocks $B_j$
(and is thus the successor of $y_j$ within this block);
we set $\xi_i = j$.

\bigskip

{\bf Step 3: Proof of bijection.}
It is easy to describe the inverse map to $\sigma \mapsto (\omega,\xi)$.
Successively for $i=1,\ldots,n$,
we use the 2-colored Motzkin path $\omega$
to read off the type of the vertex $i$ (opener, closer, insider or singleton);
and if $i$ is an insider or closer,
we use the label $\xi_i$ to decide to which ``started but unfinished'' block
the vertex $i$ should be attached.

\bigskip

{\bf Step 4: Translation of the statistics.}
This too is straightforward:

\begin{lemma}
   \label{lemma.statistics.KZ}
\quad\hfill
\vspace*{-1mm}
\begin{itemize}
   \item[(a)]  If $i$ is an insider or a closer, then
\begin{subeqnarray}
   \crrtilde(i,\pi)   & = &   h_{i-1} - \xi_i  \\[2mm]
   \neetilde(i,\pi)   & = &   \xi_i - 1
 \label{eq.lemma.statistics.KZ.a}
\end{subeqnarray}
   \item[(b)]  If $i$ is a singleton or an opener, then
\be
   \qne(i,\pi)  \;=\;  h_{i-1}
   \;.
 \label{eq.lemma.statistics.KZ.b}
\ee
\end{itemize}
\end{lemma}

\proof
(a) This is \cite[Proposition~3.3]{Kasraoui_06},
but for completeness we give the proof.
By definition, $\crrtilde(i,\pi)$ is the number of quadruplets $r < s < i < l$
such that $(r,i) \in \scrg_\pi$ and $(s,l) \in \scrg_\pi$.
But this means that, of the $h_{i-1}$ vacant vertices existing
at the beginning of stage $i$,
$r$ is the $\xi_i$th vacant vertex and $s$ is a later vacant vertex.
The number of such vertices $s$ is therefore
$\crrtilde(i,\pi) = h_{i-1} - \xi_i$.
Similarly, $\neetilde(i,\pi)$ is the number of quadruplets $r < s < i < l$
such that $(r,l) \in \scrg_\pi$ and $(s,i) \in \scrg_\pi$.
But this means that $s$ is the $\xi_i$th vacant vertex
and that $r$ is an earlier vacant vertex.
Therefore the number of such vertices $r$ is
$\neetilde(i,\pi) = \xi_i - 1$.

(b) Let $B_1,\ldots,B_{h_{i-1}}$ be the blocks
that are ``started but unfinished'' after stage $i-1$,
and let $y_j$ be the maximal element of $B_j \cap [i-1]$.
If $i$ is a singleton or an opener, then each of the vertices $y_j$
is the initial point of an arc that ends at a vertex $> i$;
and these are the only vertices that do so.
So $\qne(i,\pi) = h_{i-1}$.
[If, by contrast, $i$ is an insider or a closer, then one of the $y_j$
 is the initial point of an arc that ends at $i$,
 so $\qne(i,\pi) = h_{i-1} - 1$, in agreement with \reff{eq.qne.crr+nee}
 and part (a).]
\qed

\bigskip

{\bf Step 5: Computation of the weights \reff{def.weights.setpartitions.Jtype.final1}.}
Using the bijection, we transfer the weights \reff{def.Bnhat.reversed}
from $\pi$ to $(\omega,\xi)$
and then sum over $\xi$ to obtain the weight $W(\omega)$.
This weight is factorized over the individual steps $s_i$, as follows:
\begin{itemize}
   \item If $s_i$ is a rise starting at height $h_{i-1} = k$
      (so that $i$ is an opener),
      then from \reff{eq.lemma.statistics.KZ.b} the weight is
\be
   a_k  \;=\;  \sfb_k   \;.
\ee
   \item If $s_i$ is a fall starting at height $h_{i-1} = k$
      (so that $i$ is a closer and $k \ge 1$),
      then from \reff{eq.lemma.statistics.KZ.a} the weight is
\be
   b_k
   \;=\;
   \sum_{\xi_i=1}^k \sfa_{k-\xi_i,\, \xi_i -1}
   \;=\;
   \sfa^\star_{k-1}
\ee
where $\sfa^\star_{k-1}$ was defined in \reff{def.astar.setpartitions}.
   \item If $s_i$ is a level step of type 1 at height $h_{i-1} = k$
      (so that $i$ is an insider and $k \ge 1$),
      then from \reff{eq.lemma.statistics.KZ.a} the weight is
\be
   c^{(1)}_k
   \;=\;
   \sum_{\xi_i=1}^k \sfd_{k-\xi_i,\, \xi_i -1}
   \;=\;
   \sfd^\star_{k-1}
      \;.
\ee
   \item If $s_i$ is a level step of type 2 at height $h_{i-1} = k$
      (so that $i$ is a singleton),
      then from \reff{eq.lemma.statistics.KZ.b} the weight is
\be
   c^{(2)}_k   \;=\;   \sfe_k  \;.
\ee
\end{itemize}
Setting $\gamma_k = c^{(1)}_k + c^{(2)}_k$ and $\beta_k = a_{k-1} b_k$
as instructed in \reff{eq.flajolet.motzkin},
we obtain the weights \reff{def.weights.setpartitions.Jtype.final1}.
This completes the proof of
Theorem~\ref{thm.setpartitions.Jtype.final1}.
\qed

\bigskip

{\bf Remark.}
Theorem~\ref{thm.connected.setpartitions} on counting connected components
in set partitions,
as applied to Theorem~\ref{thm.setpartitions.Jtype.final1},
has an easy proof in our labeled-Motzkin-paths formalism.
{}From Lemma~\ref{lemma.heights.KZ} we see that
$i$ is a divider (see Section~\ref{subsec.setpartitions.connected})
if and only if $h_i = 0$.
And this happens if and only if step $s_i$
is either a fall starting at height $h_{i-1} = 1$
or a level step of type~2 at height $h_{i-1} = 0$.
So, giving each connected component a weight $\zeta$
amounts to multiplying $\sfb_0$ and $\sfe_0$ by $\zeta$.
In the J-fraction coefficients \reff{def.weights.setpartitions.Jtype.final1}
this multiplies $\gamma_0$ and $\beta_1$ by $\zeta$,
exactly as asserted in Theorem~\ref{thm.connected.setpartitions}.
\myendremark

\bigskip

Let us now explain why the reversal $\pi \mapsto \pitilde$
seems to be needed in our proof.
One reason was already explained:
the Kasraoui--Zeng bijection naturally treats
closers and insiders on the same footing,
and openers and singletons on the same footing
(cf.~Lemma~\ref{lemma.statistics.KZ}),
whereas our original definition \reff{def.Bnhat}
interchanged closers with openers in this regard.
So we need to pass to the reversed definition \reff{def.Bnhat.reversed}
in order to apply the bijection.
(Alternatively, we could have applied the Kasraoui--Zeng bijection
 to the reversed partition $\pitilde$,
 but that strikes us as even more unnatural.)

Here is another perspective on the problem:
In Section~\ref{subsec.intro.setpartitions.firstmaster}
we defined the polynomial $B_n(\bsfa,\bsfb,\bsfd,\bsfe)$
[cf.\ \reff{def.Bnhat}]
by close analogy with the permutation polynomial
$Q_n(\bsfa,\bsfb,\bsfc,\bsfd,\bsfe)$
defined in \reff{def.Qn.firstmaster},
when a set partition $\pi \in \Pi_n$ is mapped onto a permutation
$\sigma \in \Sym_n$ by specifying that the
disjoint cycles of~$\sigma$ are the blocks of $\pi$,
each traversed in increasing order.
In particular, openers correspond to cycle valleys,
closers to cycle peaks, insiders to cycle double rises,
and singletons to fixed points;  cycle double falls are forbidden.
But the Foata--Zeilberger bijection employed in our permutation proof
(Section~\ref{subsec.permutations.J}) does {\em not}\/ correspond nicely
to the Kasraoui--Zeng bijection used here in our set-partition proof.
Both bijections ``read'' the input object (permutation or set partition)
from left to right, but they employ very different senses of ``reading''.
In the Foata--Zeilberger bijection,
at stage $i$ we employ the entire permutation $\sigma$
--- not just its restriction $\sigma(1) \cdots \sigma(i)$ ---
in defining both the Motzkin path and the labels.
In the Kasraoui--Zeng bijection, by contrast,
at stage $i$ we know only the restriction of $\pi$ to $[1,i]$,
together with the status of vertex $i$ as
opener, closer, insider or singleton;
but if $i$ is an opener or insider,
we do {\em not}\/ know anything about the part of its block to its right
(except that it is nonempty).
So the definition \reff{def.Bnhat}
of the polynomial $B_n(\bsfa,\bsfb,\bsfd,\bsfe)$
does not correspond nicely to what is needed in the proof;
it turns out that the reformulation \reff{def.Bnhat.reversed}
is more appropriate.

%
%

\subsection{Second master J-fraction:
       Proof of Theorems~\ref{thm.setpartitions.pq.ovcov.brec}
          and \ref{thm.setpartitions.Jtype.final2}}
   \label{subsec.setpartitions.J.2}

In this section we will prove the second master J-fraction for set partitions
(Theorem~\ref{thm.setpartitions.Jtype.final2}).
As a consequence we will also obtain
Theorem~\ref{thm.setpartitions.pq.ovcov.brec},
which is obtained by comparing the specialization
\reff{eq.specialization.thm.setpartitions.Jtype.pq}
of Theorem~\ref{thm.setpartitions.Jtype.final1}
with the same specialization
of Theorem~\ref{thm.setpartitions.Jtype.final2}.

To prove Theorem~\ref{thm.setpartitions.Jtype.final2},
we will employ the Flajolet \cite{Flajolet_80} bijection,
which is a very slight variant of the
Kasraoui--Zeng \cite{Kasraoui_06} bijection
employed in the previous subsection.
We will therefore be brief in our description.

Analogously to what was done in the preceding subsection,
we need to use reversals of the statistics \reff{def.ov.cov.kpi}
employed in \reff{def.Bnhat2}.
These were already defined in \reff{def.ov.cov.kpi.reversed.0}:
\begin{subeqnarray}
   \ovtilde(k,\pi)
   \;\eqdef\;
   \ov(\ktilde,\pitilde)
   & = &
   \#\{ (B_1,B_2) \colon\:  k \in B_1 \,\hbox{ and }\,
                          \min B_1 < \min B_2 < k < \max B_2  \}
     \hspace*{-2cm}
       \nonumber \\ \\
   \covtilde(k,\pi)
   \;\eqdef\;
   \cov(\ktilde,\pitilde)
   & = &
   \#\{ (B_1,B_2) \colon\:  k \in B_2 \,\hbox{ and }\,
                          \min B_1 < \min B_2 < k < \max B_1  \}
     \hspace*{-2cm}
       \nonumber \\ \\
   \qcovtilde(k,\pi)
   \;\eqdef\;
   \qcov(\ktilde,\pitilde)
   & = &
   \#\{ B \colon\:  B \not\ni k \,\hbox{ and }\, \min B < k < \max B \}
 \label{def.ov.cov.kpi.reversed}
\end{subeqnarray}
But the definition of $\qcov$ is reversal-invariant,
so that in fact $\qcovtilde(k,\pi) = \qcov(k,\pi)$.
The polynomial $B^{(2)}_n(\bsfa,\bsfb,\bsfd,\bsfe)$
defined in \reff{def.Bnhat2} can then equivalently be written as
\begin{eqnarray}
   & &  \hspace*{-8mm}
   B^{(2)}_n(\bsfa,\bsfb,\bsfd,\bsfe)
   \;=\;
       \nonumber \\[4mm]
   & & 
   \sum_{\pi \in \Pi_n}
   \;
   \prod\limits_{i \in {\rm closers}}  \!\! \sfa_{\ovtilde(i,\pi),\, \covtilde(i,\pi)}
   \prod\limits_{i \in {\rm openers}}  \! \sfb_{\qcov(i,\pi)}
   \prod\limits_{i \in {\rm insiders}} \!\!  \sfd_{\ovtilde(i,\pi),\, \covtilde(i,\pi)}
   \prod\limits_{i \in {\rm singletons}}  \!\!\! \sfe_{\qcov(i,\pi)}
   \;.
        \nonumber \\[-2mm]
 \label{def.Bnhat2.reversed}
\end{eqnarray}

The Flajolet bijection --- just like the Kasraoui--Zeng bijection ---
takes $\Pi_n$ to the set of
$({\bf A},{\bf B},{\bf C}^{(1)},{\bf C}^{(2)})$-labeled
2-colored Motzkin paths of length $n$, where
${\bf A},{\bf B},{\bf C}^{(1)},{\bf C}^{(2)}$
are given by \reff{def.abc.KZ}.
The details are as follows:

\bigskip

{\bf Step 1: Definition of the Motzkin path.}
This is identical to the Kasraoui--Zeng bijection.

\bigskip

{\bf Step 2: Definition of the labels $\bm{\xi_i}$.}
If $i$ is an opener or a singleton,
we set $\xi_i = 1$ as required by \reff{def.abc.KZ}.
If $i$ is an insider or a closer,
we look at the $h_{i-1}$ blocks $B_1,\ldots,B_{h_{i-1}}$
that are ``started but unfinished'' after stage $i-1$
(note that we must have $h_{i-1} \ge 1$).
For each $j \in [h_{i-1}]$,
let $x_j$ be the {\em minimal}\/ element of $B_j \cap [i-1]$,
or in other words the opener of $B_j$.
(This use of the minimal rather than maximal element
 of $B_j \cap [i-1]$ is the {\em only}\/ change from Kasraoui--Zeng.)
We order the blocks $B_j$ so that $x_1 < x_2 < \ldots < x_{h_{i-1}}$.
Then the vertex $i$ belongs to precisely one of these blocks $B_j$;
we set $\xi_i = j$.

\bigskip

{\bf Step 3: Proof of bijection.}
Exactly as in Kasraoui--Zeng.

\bigskip

{\bf Step 4: Translation of the statistics.}
This too is straightforward,
and is a direct analogue of Lemma~\ref{lemma.statistics.KZ}:

\begin{lemma}
   \label{lemma.statistics.F}
\quad\hfill
\vspace*{-1mm}
\begin{itemize}
   \item[(a)]  If $i$ is an insider or a closer, then
\begin{subeqnarray}
   \ovtilde(i,\pi)   & = &   h_{i-1} - \xi_i  \\[2mm]
   \covtilde(i,\pi)   & = &   \xi_i - 1
 \label{eq.lemma.statistics.F.a}
\end{subeqnarray}
   \item[(b)]  If $i$ is a singleton or an opener, then
\be
   \qcov(i,\pi)  \;=\;  h_{i-1}
   \;.
 \label{eq.lemma.statistics.F.b}
\ee
\end{itemize}
\end{lemma}

\proof
(a) By definition, $\ovtilde(i,\pi)$ is the number of blocks $B'$
such that $\min B < \min B' < i < \max B'$,
where $B$ is the block containing $i$.
But this is exactly the definition of $h_{i-1} - \xi_i$.
Similarly, $\covtilde(i,\pi)$ is the number of blocks $B'$
such that $\min B' < \min B < i < \max B'$,
where $B$ is the block containing $i$.
But this is exactly the definition of $\xi_i - 1$.

(b) By definition, $\qcov(i,\pi)$ is the number of blocks $B \not\ni i$
such that $\min B < i < \max B$.
If $i$ is a singleton or an opener, then this equals $h_{i-1}$.
[If, by contrast, $i$ is an insider or a closer, then one of the blocks
 $B_1,\ldots,B_{h_{i-1}}$ contains $i$,
 so $\qcov(i,\pi) = h_{i-1} - 1$, in agreement with \reff{eq.qcov.ov+cov}
 and part (a).]
\qed

\bigskip

{\bf Step 5: Computation of the weights.}
Identical to Section~\ref{subsec.setpartitions.J},
with the obvious substitutions
$\crrtilde \to \ovtilde$, $\neetilde \to \covtilde$, $\qne \to \qcov$.
This completes the proof of
Theorem~\ref{thm.setpartitions.Jtype.final2}.
\qed

\subsection{Third and fourth master J-fractions:
       Proof of Theorems~\ref{thm.Bnhat.third.equivalences} and
       \ref{thm.setpartitions.Jtype.final34}}
   \label{subsec.setpartitions.J.3}

In this section we will prove the third and fourth
master J-fractions for set partitions
(Theorem~\ref{thm.setpartitions.Jtype.final34}).
As a consequence we will also obtain
Theorem~\ref{thm.Bnhat.third.equivalences},
which is obtained by comparing the specialization
\reff{eq.specialization.thm.setpartitions.Jtype.pq}
of Theorems~\ref{thm.setpartitions.Jtype.final1},
\ref{thm.setpartitions.Jtype.final2} and \ref{thm.setpartitions.Jtype.final34}.

To prove Theorem~\ref{thm.setpartitions.Jtype.final34},
we will employ a bizarre amalgam
of the Flajolet and Kasraoui--Zeng bijections,
in which the labels for insiders are given by Flajolet
and those for closers by Kasraoui--Zeng, or vice versa.
We will again be brief in our description.

Once again we begin by defining the reversed statistics
\reff{def.cr.ne.kpi.reversed} and \reff{def.ov.cov.kpi.reversed}.
We then rewrite the polynomials
$B^{(3)}_n(\bsfa,\bsfb,\bsfd,\bsfe)$
and
$B^{(4)}_n(\bsfa,\bsfb,\bsfd,\bsfe)$
in terms of these reversed statistics:
\begin{eqnarray}
   & &  \hspace*{-8mm}
   B^{(3)}_n(\bsfa,\bsfb,\bsfd,\bsfe)
   \;=\;
       \nonumber \\[1mm]
   & & 
   \sum_{\pi \in \Pi_n}
   \;
   \prod\limits_{i \in {\rm closers}}  \!\! \sfa_{\ovtilde(i,\pi),\, \covtilde(i,\pi)}
   \prod\limits_{i \in {\rm openers}}  \! \sfb_{\qne(i,\pi)}
   \prod\limits_{i \in {\rm insiders}} \!\!  \sfd_{\crrtilde(i,\pi),\, \neetilde(i,\pi)}
   \prod\limits_{i \in {\rm singletons}}  \!\!\! \sfe_{\qne(i,\pi)}
   \;.
        \nonumber \\[-2mm] \\[2mm]
 \label{def.Bnhat3.reversed}
   & &  \hspace*{-8mm}
   B^{(4)}_n(\bsfa,\bsfb,\bsfd,\bsfe)
   \;=\;
       \nonumber \\[1mm]
   & & 
   \sum_{\pi \in \Pi_n}
   \;
   \prod\limits_{i \in {\rm closers}}  \!\! \sfa_{\crrtilde(i,\pi),\, \neetilde(i,\pi)}
   \prod\limits_{i \in {\rm openers}}  \! \sfb_{\qne(i,\pi)}
   \prod\limits_{i \in {\rm insiders}} \!\!  \sfd_{\ovtilde(i,\pi),\, \covtilde(i,\pi)}
   \prod\limits_{i \in {\rm singletons}}  \!\!\! \sfe_{\qne(i,\pi)}
   \;.
        \nonumber \\[-2mm]
 \label{def.Bnhat4.reversed}
\end{eqnarray}
Then the bijections --- call them \#3 and \#4 ---
are defined as follows:

\bigskip

{\bf Step 1: Definition of the Motzkin path.}
Exactly as in the Kasraoui--Zeng and Flajolet bijections.

\bigskip

{\bf Step 2: Definition of the labels $\bm{\xi_i}$.}
In both bijections \#3 and \#4,
if $i$ is an opener or a singleton, we set $\xi_i = 1$.
Then, in bijection \#3:
\begin{itemize}
   \item If $i$ is an insider, we define $\xi_i$ as in the
       Kasraoui--Zeng bijection.
   \item If $i$ is a closer, we define $\xi_i$ as in the
       Flajolet bijection.
\end{itemize}
Bijection \#4 is defined by the reverse scheme:
\begin{itemize}
   \item If $i$ is an insider, we define $\xi_i$ as in the
       Flajolet bijection.
   \item If $i$ is a closer, we define $\xi_i$ as in the
       Kasraoui--Zeng bijection.
\end{itemize}

\bigskip

{\bf Step 3: Proof of bijection.}
Exactly as in Kasraoui--Zeng.

\bigskip

{\bf Step 4: Translation of the statistics.}

\medskip

\begin{lemma}
   \label{lemma.statistics.34}
\quad\hfill \\[1mm]
\noindent
In bijection \#3:
\begin{itemize}
   \item[(a)]  If $i$ is an insider, then
\begin{subeqnarray}
   \crrtilde(i,\pi)   & = &   h_{i-1} - \xi_i  \\[2mm]
   \neetilde(i,\pi)   & = &   \xi_i - 1
 \label{eq.lemma.statistics.3.a}
\end{subeqnarray}
   \item[(b)]  If $i$ is a closer, then
\begin{subeqnarray}
   \ovtilde(i,\pi)   & = &   h_{i-1} - \xi_i  \\[2mm]
   \covtilde(i,\pi)   & = &   \xi_i - 1
 \label{eq.lemma.statistics.3.b}
\end{subeqnarray}
\end{itemize}
\noindent
In bijection \#4:
\begin{itemize}
   \item[(c)]  If $i$ is an insider, then
\begin{subeqnarray}
   \ovtilde(i,\pi)   & = &   h_{i-1} - \xi_i  \\[2mm]
   \covtilde(i,\pi)   & = &   \xi_i - 1
 \label{eq.lemma.statistics.4.c}
\end{subeqnarray}
   \item[(d)]  If $i$ is a closer, then
\begin{subeqnarray}
   \crrtilde(i,\pi)   & = &   h_{i-1} - \xi_i  \\[2mm]
   \neetilde(i,\pi)   & = &   \xi_i - 1
 \label{eq.lemma.statistics.4.d}
\end{subeqnarray}
\end{itemize}
\noindent
In both bijections:
\begin{itemize}
   \item[(e)]  If $i$ is a singleton or an opener, then
\be
   \qne(i,\pi)  \;=\;  \qcov(i,\pi)  \;=\;  h_{i-1}
   \;.
 \label{eq.lemma.statistics.34.e}
\ee
\end{itemize}
\end{lemma}

\proof
This is an immediate consequence of
Lemmas~\ref{lemma.statistics.KZ} and \ref{lemma.statistics.F}
together with \reff{eq.qcov=qne}.
\qed

\bigskip

{\bf Step 5: Computation of the weights.}
Identical to Sections~\ref{subsec.setpartitions.J}
and \ref{subsec.setpartitions.J.2},
with the obvious substitutions.
This completes the proof of
Theorem~\ref{thm.setpartitions.Jtype.final34}.
\qed

\section*{Acknowledgments}
\addcontentsline{toc}{section}{Acknowledgments}

We wish to thank Andrew Elvey Price, Mathias P\'etr\'eolle and Jes\'us Salas
for helpful conversations.
We also thank Natasha Blitvi\'c and Einar Steingr\'{\i}msson
for drawing our attention to their recent work \cite{Blitvic_20}.

This work has benefited greatly from the existence of
the On-Line Encyclopedia of Integer Sequences \cite{OEIS}.
We warmly thank Neil Sloane for founding this indispensable resource,
and the hundreds of volunteers for helping to maintain and expand it.

This work has depended heavily on symbolic computation
(in {\sc Mathematica}) carried out on large-memory computers
funded by the U.S.~taxpayers via
National Science Foundation grant PHY--0424082,
by the U.K.~taxpayers via
Engineering and Physical Sciences Research Council grant EP/N025636/1,
and by a computer donation from the Dell Corporation;
we are extremely grateful to all three.

This research was supported in part by
Engineering and Physical Sciences Research Council grant EP/N025636/1.

\end{document}